\documentclass[12pt]{article}
\usepackage{amssymb}
\usepackage{amsmath}
\usepackage{makeidx}
\usepackage[dvips]{graphicx}
\usepackage{sw20elba}

\newtheorem{theorem}{Theorem}[section]

\newtheorem{condition}[theorem]{Condition}

\newtheorem{corollary}[theorem]{Corollary}

\newtheorem{definition}[theorem]{Definition}
\newtheorem{example}[theorem]{Example}

\newtheorem{lemma}[theorem]{Lemma}

\newtheorem{proposition}[theorem]{Proposition}
\newtheorem{remark}[theorem]{Remark}

\newenvironment{proof}[1][Proof]{\textbf{#1.} }{\ \rule{0.5em}{0.5em}}
\input{tcilatex}
\numberwithin{equation}{section}

\begin{document}

\title{Nonlinear dynamics of a system of particle-like wavepackets }
\author{A. Babin and A. Figotin \\
Department of Mathematics, University of California at Irvine, CA 92697}
\maketitle

\begin{abstract}
This work continues our studies of nonlinear evolution of a system of
wavepackets. We study a wave propagation governed by a nonlinear system of
hyperbolic PDE's with constant coefficients with the initial data being a
multi-wavepacket. By definition a general wavepacket has a well defined
principal wave vector, and, as we proved in previous works, the nonlinear
dynamics preserves systems of wavepackets and their principal wave vectors.
Here we study the nonlinear evolution of a special class of wavepackets,
namely particle-like wavepackets. A particle-like wavepacket is of a dual
nature: on one hand, it is a wave with a well defined principal wave vector,
on the other hand, it is a particle in the sense that it can be assigned a
well defined position in the space. We prove that under the nonlinear
evolution a generic multi-particle wavepacket remains to be a multi-particle
wavepacket with a high accuracy, and every constituting single particle-like
wavepacket not only preserves its principal wave number\ but also it has a
well-defined space position evolving with a constant velocity which is its
group velocity. Remarkably the described properties hold though the involved
single particle-like wavepackets undergo nonlinear interactions and multiple
collisions in the space.\ We also prove that if principal wavevectors of
multi-particle wavepacket are generic, the result of nonlinear interactions
between different wavepackets is small and the approximate linear
superposition principle holds uniformly with respect to the initial spatial
positions of wavepackets.
\end{abstract}

\section{Introduction}

The principal object of our studies here is a general nonlinear evolutionary
system which describes wave propagation in homogeneous media governed by a
hyperbolic PDE's in $\mathbb{R}^{d}$, $d=1,2,3,\ldots $ is the space
dimension, of the form 
\begin{equation}
\partial _{\tau }\mathbf{U}=-\frac{\mathrm{i}}{\varrho }\mathbf{L}\left( -%
\mathrm{i}\nabla \right) \mathbf{U}+\mathbf{F}\left( \mathbf{U}\right) ,\
\left. \mathbf{U}\left( \mathbf{r},\tau \right) \right\vert _{\tau =0}=%
\mathbf{h}\left( \mathbf{r}\right) ,\ \mathbf{r}\in \mathbb{R}^{d},
\label{difeqintr}
\end{equation}%
where (i) $\mathbf{U}=\mathbf{U}\left( \mathbf{r},\tau \right) $, $\mathbf{r}%
\in \mathbb{R}^{d}$, $\mathbf{U}\in \mathbb{C}^{2J}$ is a $2J$ dimensional
vector; (ii) $\mathbf{L}\left( -\mathrm{i}\nabla \right) $ is a linear
self-adjoint differential (pseudodifferential) operator with constant
coefficients with the symbol $\mathbf{L}\left( \mathbf{k}\right) $, which is
a Hermitian $2J\times 2J$ matrix; (iii) $\ \mathbf{F}$ is a general
polynomial nonlinearity; (iv) $\varrho >0$ is a \emph{small parameter}. The
properties of the linear part are described in terms of dispersion relations 
$\omega _{n}\left( \mathbf{k}\right) $ (eigenvalues of the matrix $\mathbf{L}%
\left( \mathbf{k}\right) $).\ The form of the equation suggests that the
processes described by it involve two time scales. Since the nonlinearity $%
\mathbf{F}\left( \mathbf{U}\right) $ is of order one, nonlinear effects
occur at times $\tau $ of order one, whereas the natural time scale of
linear effects, governed by the operator $\mathbf{L}$ with the coefficient $%
1/\varrho $, is of order $\varrho $. Consequently, the small parameter $%
\varrho $ measures the ratio of the slow (nonlinear effects) time scale and
the fast (linear effects) time scale. A typical example of an equation of
the form (\ref{difeqintr}) is the nonlinear Schrodinger equation (NLS) or a
system of NLS's. Many more examples including a general nonlinear wave
equation and the Maxwell equations in periodic media truncated to a finite
number of bands are considered in \cite{BF7}, \cite{BF8}.

As in our previous works \cite{BF7}, \cite{BF8} we consider here the
nonlinear evolutionary system (\ref{difeqintr}) with the initial data $%
\mathbf{h}\left( \mathbf{r}\right) $ being a sum of wavepackets. \emph{The
special focus of this paper is particle-like localized wavepackets which can
be viewed as quasi-particles}. Recall that a general wavepacket is defined
as such a function $\mathbf{h}\left( \mathbf{r}\right) $ that its Fourier
transform $\mathbf{\hat{h}}\left( \mathbf{k}\right) $ is localized in $\beta 
$-neighborhood of a single wavevector $\mathbf{k}_{\ast }$, called \emph{%
principal wavevector}, where $\beta $ is a \emph{small parameter}. The
simplest example of a wavepacket is a function of the form%
\begin{equation}
\mathbf{\hat{h}}\left( \beta ;\mathbf{k}\right) =\beta ^{-d}\mathrm{e}^{-%
\mathrm{i}\mathbf{kr}_{\ast }}\hat{h}\left( \frac{\mathbf{k}-\mathbf{k}%
_{\ast }}{\beta }\right) \mathbf{g}_{n}\left( \mathbf{k}_{\ast }\right) ,\ 
\mathbf{k}\in \mathbb{R}^{d},  \label{hh}
\end{equation}%
where $\mathbf{g}_{n}\left( \mathbf{k}_{\ast }\right) $ is an eigenvector of
the matrix $\mathbf{L}\left( \mathbf{k}_{\ast }\right) $ and $\hat{h}\left( 
\mathbf{k}\right) $ is a scalar Schwartz function (i.e. it is infinitely
smooth and rapidly decaying one). Note that for $\mathbf{\hat{h}}\left(
\beta ,\mathbf{k}\right) $ of the form (\ref{hh}) we have its inverse
Fourier transform 
\begin{equation}
\mathbf{h}\left( \beta ;\mathbf{r}\right) =h\left( \beta \left( \mathbf{r}-%
\mathbf{r}_{\ast }\right) \right) \mathrm{e}^{\mathrm{i}\mathbf{k}_{\ast
}\left( \mathbf{r-r}_{\ast }\right) }\mathbf{g}_{n}\left( \mathbf{k}_{\ast
}\right) ,\ \mathbf{r}\in \mathbb{R}^{d}.  \label{hr}
\end{equation}%
Evidently, $\mathbf{h}\left( \beta ,\mathbf{r}\right) $ described by the
above formula is a plane wave $\mathrm{e}^{\mathrm{i}\mathbf{k}_{\ast }%
\mathbf{r}}\mathbf{g}_{n}\left( \mathbf{k}_{\ast }\right) $ modulated by a
slowly varying amplitude $h\left( \beta \left( \mathbf{r}-\mathbf{r}_{\ast
}\right) \right) $ obtained from $h\left( \mathbf{z}\right) $ by a spatial
shift along the vector $\mathbf{r}_{\ast }$ with a subsequent dilation with
a large factor $\frac{1}{\beta }$. Clearly, the resulting amplitude has a
typical spatial extension proportional to $\beta ^{-1}$ and the spatial
shift produces a noticeable effect if $\left\vert \mathbf{r}_{\ast
}\right\vert \gg \beta ^{-1}$. The spatial form of the wavepacket (\ref{hr})
naturally allows to interpret $\mathbf{r}_{\ast }\in \mathbb{R}^{d}$ as its
position and, consequently, to consider the wavepacket as a particle-like
one with the position $\mathbf{r}_{\ast }\in \mathbb{R}^{d}$. But how one
can define a position for a general wavepacket? Note that that not every
wavepacket is a particle-like one. For example, let, as before, the function 
$h\left( \mathbf{r}\right) $ be a scalar Schwartz function and let us
consider a slightly more general than (\ref{hr}) function%
\begin{equation}
\mathbf{h}\left( \beta ;\mathbf{r}\right) =\left[ h\left( \beta \left( 
\mathbf{r}-\mathbf{r}_{\ast 1}\right) \right) +h\left( \beta \left( \mathbf{r%
}-\mathbf{r}_{\ast 2}\right) \right) \right] \mathrm{e}^{\mathrm{i}\mathbf{k}%
_{\ast }\mathbf{r}}\mathbf{g}_{n}\left( \mathbf{k}_{\ast }\right) ,\ \mathbf{%
r}\in \mathbb{R}^{d},  \label{hr1}
\end{equation}%
where $\mathbf{r}_{\ast 1}$ and $\mathbf{r}_{\ast 2}$ are two arbitrary,
independent vector variables. The wave $\mathbf{h}\left( \beta ,\mathbf{r}%
\right) $ defined by (\ref{hr1}) is a wavepacket with the wave number $%
\mathbf{k}_{\ast }$ for any choice of vectors $\mathbf{r}_{\ast 1}$ and $%
\mathbf{r}_{\ast 2}$, but it is not a particle-like wavepacket, since it
does not have a single position $\mathbf{r}_{\ast }$, but rather it is a sum
of two particle-like wavepackets with two positions $\mathbf{r}_{\ast 1}$
and $\mathbf{r}_{\ast 2}$.

Our way to introduce a general particle-like wavepacket $\mathbf{h}\left(
\beta ,\mathbf{k}_{\ast },\mathbf{r}_{\ast 0};\mathbf{r}\right) $ with a
position $\mathbf{r}_{\ast 0}$ is by treating it as a single element of a
family of wavepackets $\mathbf{h}\left( \beta ,\mathbf{k}_{\ast },\mathbf{r}%
_{\ast };\mathbf{r}\right) $ with $\mathbf{r}_{\ast }\in \mathbb{R}^{d}$
being another \emph{independent parameter}. In fact, we define the entire
family of wavepackets $\mathbf{h}\left( \beta ,\mathbf{k}_{\ast },\mathbf{r}%
_{\ast };\mathbf{r}\right) $, $\mathbf{r}_{\ast }\in \mathbb{R}^{d}$ subject
to certain conditions allowing to interpret any fixed $\mathbf{r}_{\ast }\in 
\mathbb{R}^{d}$ as the position of $\mathbf{h}\left( \beta ,\mathbf{k}_{\ast
},\mathbf{r}_{\ast };\mathbf{r}\right) $. Since we would like, of course, a
wavepacket to maintain its particle-like property under the nonlinear
evolution, it is clear that its definition must be sufficiently flexible to
accommodate the wavepacket evolutionary variations. In light of the above
discussion the definition of the particle-like wavepacket with a transparent
interpretation of its particle properties turns into the key element of the
entire construction. It turns out that there is a precise description of a
particle-like wavepacket, which is rather simple and physically transparent
and such a description is provided in Definition \ref{Definition regwave}
below, see also Remarks \ref{aloc}, \ref{alocgen}. The concept of the
position is applicable to very general functions, it does not require a
parametrization of the whole family of solutions, which was used, for
example in \cite{JohnssonFrGS06}, \cite{GangSigal06a}, \cite{GangSigal06}.

As in our previous works we are interested in nonlinear evolution not only a
single particle-like wavepacket $\mathbf{h}\left( \beta ,\mathbf{k}_{\ast },%
\mathbf{r}_{\ast };\mathbf{r}\right) $ but a system $\left\{ \mathbf{h}%
\left( \beta ,\mathbf{k}_{\ast l},\mathbf{r}_{\ast l};\mathbf{r}\right)
\right\} $ of particle-like wavepackets which we call \emph{multi-particle
wavepacket}. Under certain natural conditions of genericity on $\mathbf{k}%
_{\ast l}$ we prove here that under the nonlinear evolution: (i) the
multi-particle wavepacket remains to be a multi-particle wavepacket; (ii)
the principal wavevectors $\mathbf{k}_{\ast l}$ remain constant; (ii) the
spatial position $\mathbf{r}_{\ast l}$ of the corresponding wavepacket
evolves with the constant velocity which is exactly its group velocity $%
\frac{1}{\varrho }\nabla \omega _{n}\left( \mathbf{k}_{\ast l}\right) $. The
evolution of positions of wavepackets becomes the most simple in the case
when at $\tau =0$ we have $\mathbf{r}_{\ast l}=\frac{1}{\varrho }\mathbf{r}%
_{\ast }^{0}$, that is the case when spatial positions are bounded in the
same spatial scale in which their group velocities are bounded. In this case
the evolution of the positions is described by the formula 
\begin{equation}
\mathbf{r}_{l}\left( \tau \right) =\frac{1}{\varrho }\left[ \mathbf{r}_{\ast
}^{0}+\tau \nabla \omega _{n_{l}}\left( \mathbf{k}_{\ast l}\right) \right]
,\ \tau \geq 0.  \label{hr2}
\end{equation}%
The rectilinear motion of positions of particle-like wavepackets is a direct
consequence of the spatial homogeneity of the master system (\ref{difeqintr}%
). If the system were not spatially homogeneous, the motion of the positions
of particle-like wavepackets would not be uniform, but we don't study that
problem in this paper. In the rescaled coordinates $\mathbf{y}=\varrho 
\mathbf{r}$ the trajectory of every particle is a fixed, uniquely defined
straight line defined uniquely if $\frac{\varrho }{\beta }\rightarrow 0$ as $%
\varrho ,\beta \rightarrow 0$. Notice that under above mentioned genericity
condition the uniform and independent motion (\ref{hr2}) of the positions of
all involved particle-like wavepackets $\left\{ \mathbf{h}\left( \beta ,%
\mathbf{k}_{\ast l},\mathbf{r}_{\ast l};\mathbf{r}\right) \right\} $
persists though they can collide in the space. In the latter case they
simply pass through each other without significant nonlinear interactions,
and the nonlinear evolution with a high accuracy is reduced just to a
nonlinear evolution of shapes of the particle-like wavepackets. In the case
when the set of the principal wavevectors $\left\{ \mathbf{k}_{\ast
l}\right\} $ satisfy certain resonance conditions some components of the
original multi-particle wavepacket can evolve into a more complex structure
which can be only partly localized in the space and, for instance, can be
needle or pancake like. We do not study in detail those more complex
structures here.

Now let us discuss in more detail the superposition principle introduced and
studied for general multi-wavepackets in \cite{BF8} in the particular case
when initially all $\mathbf{r}_{\ast l}=0$. Here we consider multi-particle
wavepackets with arbitrary $\mathbf{r}_{\ast l}$ and develop a new argument
based on the analysis of an averaged wavepacket interaction system
introduced in \cite{BF7}. Assume that the initial data $\mathbf{h}$ for the
evolution equation (\ref{difeqintr}) to be the sum of a finite number of 
\emph{wavepackets (particle-like wavepackets)} $\mathbf{h}_{l}$, $l=1,\ldots
,N$, i.e.%
\begin{equation}
\mathbf{h}=\mathbf{h}_{1}+\ldots +\mathbf{h}_{N}\   \label{hhh1}
\end{equation}%
where the monochromaticity of every wavepacket $\mathbf{h}_{l}$ is
characterized by another small parameter $\beta $. The well known \emph{%
superposition principle} is a fundamental property of every linear
evolutionary system, stating that the solution $\mathbf{U}$ corresponding to
the initial data $\mathbf{h}$ as in (\ref{hhh1}) equals%
\begin{equation}
\mathbf{U}=\mathbf{U}_{1}+\ldots +\mathbf{U}_{N},\text{ for }\mathbf{h}=%
\mathbf{h}_{1}+\ldots +\mathbf{h}_{N},  \label{hhh2}
\end{equation}%
where $\mathbf{U}_{l}$ is the solution to the same linear problem with the
initial data $\mathbf{h}_{l}$.

Evidently the standard superposition principle can not hold exactly as a
general principle for a nonlinear system, and, at the first glance, there is
no expectation for it to hold even approximately. We show though that, in
fact, the superposition principle does hold with a high accuracy for general
dispersive nonlinear wave systems such as (\ref{difeqintr}) provided that
the initial data are a sum of generic particle-like wavepackets, and this
constitutes one of the subjects of this paper. Namely, the superposition
principle for nonlinear wave systems states that the solution $\mathbf{U}$
corresponding to the multi-particle wavepacket initial data $\mathbf{h}$ as
in (\ref{hhh1}) satisfies \emph{\ }%
\begin{equation*}
\mathbf{U}=\mathbf{U}_{1}+\ldots +\mathbf{U}_{N}+\mathbf{D},\text{ for }%
\mathbf{h}=\mathbf{h}_{1}+\ldots +\mathbf{h}_{N},\text{ where }\mathbf{D}%
\text{ is small.}
\end{equation*}
More detailed statement of the superposition principle for nonlinear
evolution of wavepackets is as follows. We study the nonlinear evolution
equation (\ref{difeqintr}) on a finite time interval 
\begin{equation}
0\leq \tau \leq \tau _{\ast },\text{ where }\tau _{\ast }>0\text{ is a fixed
number}  \label{taustar}
\end{equation}%
which may depend on the $L^{\infty }$ norm of the initial data $\mathbf{h}$
but, importantly, $\tau _{\ast }$\emph{\ does not depend on }$\varrho $. 
\emph{We consider classes of initial data such that wave evolution governed
by (\ref{difeqintr}) is significantly nonlinear on time interval }$\left[
0,\tau _{\ast }\right] $\emph{\ and the effect of the nonlinearity }$F\left( 
\mathbf{U}\right) $\emph{\ does not vanish as }$\varrho \rightarrow 0$\emph{%
. }We assume that $\beta ,\varrho $ satisfy%
\begin{equation}
0<\beta \leq 1,\ 0<\varrho \leq 1,\ \frac{\beta ^{2}}{\varrho }\leq C_{1}%
\text{ with some }C_{1}>0.  \label{rbb1}
\end{equation}%
The above condition of boundedness on the dispersion parameter $\frac{\beta
^{2}}{\varrho }$ ensures that the dispersive effects are not dominant and
they do not suppress nonlinear effects, see \cite{BF7}, \cite{BF8} for a
discussion.

Let us introduce the solution operator $\mathcal{S}\left( \mathbf{h}\right)
\left( \tau \right) :\mathbf{h}\rightarrow \mathbf{U}\left( \tau \right) $
relating the initial data $\mathbf{h}$ of the nonlinear evolution equation (%
\ref{difeqintr}) to its solution $\mathbf{U}\left( t\right) $. Suppose that
the initial state is a system of particle-like wavepackets or \emph{%
multi-particle wavepacket}, namely $\mathbf{h}=\dsum \mathbf{h}_{l}$, with $%
\mathbf{h}_{l}$, $l=1,\ldots ,N$ being "generic" wavepackets. Then for all
times\emph{\ }$0\leq \tau \leq \tau _{\ast }$ the following\emph{\
superposition principle} holds 
\begin{gather}
\mathcal{S}\left( \dsum\nolimits_{l=1}^{N}\mathbf{h}_{l}\right) \left( \tau
\right) =\dsum\nolimits_{l=1}^{N}\mathcal{S}\left( \mathbf{h}_{l}\right)
\left( \tau \right) +\mathbf{D}\left( \tau \right) ,  \label{Gsum1} \\
\left\Vert \mathbf{D}\left( \tau \right) \right\Vert _{E}=\sup\limits_{0\leq
\tau \leq \tau _{\ast }}\left\Vert \mathbf{D}\left( \tau \right) \right\Vert
_{L^{\infty }}\leq C_{\delta }\frac{\varrho }{\beta ^{1+\delta }}\text{ for
any small }\delta >0.  \label{Dbetkap}
\end{gather}%
Obviously, the right-hand side of (\ref{Dbetkap}) may be small \ only if $%
\varrho \leq C_{1}\beta $. There are examples (see \cite{BF7}) in which $%
\mathbf{D}\left( \tau \right) $ is not small for $\varrho =C_{1}\beta $. \
In what follows we refer to a linear combination of particle-like
wavepackets as a\emph{\ multi-particle wavepacket}, and to single
particle-like wavepackets which constitutes the multi-particle wavepacket as 
\emph{component particle wavepackets}.

Very often in theoretical studies of equations of the form (\ref{difeqintr})
or ones reducible to it a functional dependence between $\varrho $ and $%
\beta $ is imposed, resulting in a single small parameter. The most common
scaling is $\varrho =\beta ^{2}$. The nonlinear evolution of wavepackets for
a variety of equations which can be reduced to the form (\ref{difeqintr})
was studied in numerous physical and mathematical papers, mostly by
asymptotic expansions of solutions with respect to a single small parameter
similar to $\beta $, see \cite{BenYoussefLannes02}, \cite{BonaCL05}, \cite%
{ColinLannes}, \cite{CraigSulemS92}, \cite{GiaMielke}, \cite{JolyMR98}, \cite%
{KalyakinUMN}, \cite{Maslov83}, \cite{PW}, \cite{Schneider98a}, \cite%
{Schneider05} and references therein. Often the asymptotic expansions are
based on a specific ansatz prescribing a certain form\ to the solution. In
our studies here we do not use asymptotic expansions with respect to a small
parameter and do not prescribe a specific form to the solution, but we
impose conditions on the initial data requiring it to be a wavepacket or a
linear combination of wavepackets. Since we want to establish a general
property of a wide class of systems, we apply a general enough dynamical
approach. There is a number of general approaches developed for the studies
of high-dimensional and infinite-dimensional nonlinear evolutionary systems
\ of hyperbolic type, \cite{Bambusi03}, \cite{BM}, \cite{Craig}, \cite{GW}, 
\cite{Iooss}, \cite{Kuksin}, \cite{MSZ}, \cite{Schneider98a}, \cite{SU}, 
\cite{Weinstein}, \cite{W}) and references therein. \ The approach we
develop here is based on the introduction of a wavepacket interaction
system. We show in \cite{BF8} and here that solutions to this system are in
a close relation to solutions of the original system.

The superposition principle implies, in particular, that in the process of
nonlinear evolution every single wavepacket propagates almost \
independently of other wavepackets (even though they may "collide" in
physical space for a certain period of time) and the exact solution equals
the sum of particular single wavepacket solutions with a high precision. In
particular, the dynamics of a solution with multi-wavepacket initial data is
reduced to dynamics of separate solutions with single wavepacket data. Note
that the nonlinear evolution of a single wavepacket solution for many
problems is studied in detail, namely it is well approximated by its own
nonlinear Schrodinger equation (NLS), see \cite{ColinLannes}, \cite%
{GiaMielke}, \cite{KalyakinUMN}, \cite{Kalyakin2}, \cite{Schneider98a}, \cite%
{Schneider05}, \cite{SU}, \cite{BF7} and references therein.

Let us give now an elementary physical argument justifying the superposition
principle which goes as follows. If there would be no nonlinearity, the
system would be linear and, consequently, the superposition principle would
hold exactly. Hence, any deviation from it is due to the nonlinear
interactions between wavepackets, and one has to estimate their impact.
Suppose that initially at time $\tau =0$ the spatial extension $s$ of every
involved wavepacket is characterized by the parameter $\beta ^{-1}$ as in (%
\ref{hr}). Assume also (and it is quite an assumption) that the involved
wavepackets evolving nonlinearly maintain somehow their wavepacket
identities, including the group velocities and the spatial extensions. Then,
consequently, the spatial extension of every involved wavepacket is
propositional to $\beta ^{-1}$ and its group velocity $v_{l}$ is
proportional to $\varrho ^{-1}$. The difference $\Delta v$ between any two
different group velocities is also proportional to $\varrho ^{-1}$. Then the
time when two different wavepackets overlap in the space is proportional to $%
s/\left\vert \Delta v\right\vert $ and, hence, to $\varrho /\beta $. Since
the nonlinear term is of order one, the \ magnitude of the impact of the
nonlinearity during this time interval should be roughly proportional to $%
\varrho /\beta $, which results in the same order of the magnitude of $%
\mathbf{D}$ in (\ref{Gsum1})-(\ref{Dbetkap}). Observe, that this estimate is
in agreement with our rigorous estimate of the magnitude of $\mathbf{D}$ in (%
\ref{Dbetkap}) if we set there $\delta =0$.

The rigorous proof of the superposition principle presented here is not
directly based on the above argument since it already implicitly relies on
the principle. Though some components of the physical argument can be found
in our rigorous proof. For example, we prove that the involved wavepackets
maintain under the nonlinear evolution constant values of their wavevectors
with well defined group velocities (the wavepacket preservation). The
Theorem \ref{Theorem E1sol} allows to estimate spatial extensions of
particle-like wavepackets under the nonlinear evolution. The proof of the
superposition principle for general wavepackets provided in \cite{BF8} \ is
based on general algebraic-functional considerations \ and on the theory of
analytic operator expansions in Banach spaces. Here we develop an
alternative approach with a proof based on properties of the wavepacket
interaction systems introduced in \cite{BF7}.

To provide a flexibility in formulating more specific statements related to
the spatial localization of wavepackets we introduce a few types of
wavepackets:

\begin{itemize}
\item a single \emph{particle-like wavepacket} $w$ which is characterized by
the following properties: (a) its modal decomposition involves only
wavevectors from $\beta $-vicinity of a single wavevector $\mathbf{k}_{\ast
} $, where $\beta >0$ is a small parameter; (b) it is spatially localized in
all directions and can be assigned its position $\mathbf{r}_{\ast }$;

\item a \emph{multi-particle wavepacket} which is a system $\left\{
w_{l}\right\} $ of particle-like wavepackets with the corresponding sets of
wavevectors $\left\{ \mathbf{k}_{\ast l}\right\} $ and positions $\left\{ 
\mathbf{r}_{\ast l}\right\} $;

\item a \emph{spatially localized multi-wavepacket} which is a system $%
\left\{ w_{l}\right\} $ with $w_{l}$ being either a particle-like wavepacket
or a general wavepacket.
\end{itemize}

We would like to note that a more detailed analysis, which is left for
another paper, indicates that under certain resonance conditions nonlinear
interactions of particle-like wavepackets may produce a \emph{spatially
localized wavepacket} $w$ characterized by the following properties: (i) its
modal decomposition involves only wavevectors from $\beta $-vicinity of a
single wavevector $\mathbf{k}_{\ast }$, where $\beta >0$ is a small
parameter; (ii) it is only partly spatially localized in some, not
necessarily all directions, and, for instance, it can be needle-like or
pancake-like.

We also would like to point out that the particular form (\ref{difeqintr})
of dependence on the small parameter $\varrho $ is chosen so that
appreciable nonlinear effects occur at times of order one. In fact, many
important classes of problems involving small or parameters \ can be readily
reduced to the framework of (\ref{difeqintr}) by a simple rescaling. It can
be seen from the following examples. The first example is a system with a 
\emph{small nonlinearity} 
\begin{equation}
\partial _{t}\mathbf{v}=-\mathrm{i}\mathbf{Lv}+\alpha \mathbf{f}\left( 
\mathbf{v}\right) ,\ \left. \mathbf{v}\right\vert _{t=0}=\mathbf{h},\
0<\alpha \ll 1,  \label{smallnon}
\end{equation}%
where the initial data is bounded uniformly in $\alpha $. Such problems are
reduced to (\ref{difeqintr}) by the time rescaling $\tau =t\alpha $.$\ $Note
that here $\varrho =\alpha $ \ and the finite time interval $0\leq \tau \leq
\tau _{\ast }$ corresponds to the \emph{long time interval} $0\leq t\leq
\tau _{\ast }/\alpha $.

The second example is a system with \emph{small initial data} considered on
long time intervals. The system itself has no small parameters but the
initial data are small, namely 
\begin{gather}
\partial _{t}\mathbf{v}=-\mathrm{i}\mathbf{Lv}+\mathbf{f}_{0}\left( \mathbf{v%
}\right) ,\ \left. \mathbf{v}\right\vert _{t=0}=\alpha _{0}\mathbf{h},\
0<\alpha _{0}\ll 1,\text{ where}  \label{smallinit} \\
\mathbf{f}_{0}\left( \mathbf{v}\right) =\mathbf{f}_{0}^{\left( m\right)
}\left( \mathbf{v}\right) +\mathbf{f}_{0}^{\left( m+1\right) }\left( \mathbf{%
v}\right) +\ldots ,  \notag
\end{gather}%
where $\alpha _{0}$ is a small parameter and $\mathbf{f}^{\left( m\right)
}\left( \mathbf{v}\right) $ is a \ homogeneous polynomial of degree $m\geq 2$%
. After rescaling $\mathbf{v}=\alpha _{0}\mathbf{V}$ we obtain the following
equation \ with a small nonlinearity 
\begin{equation}
\partial _{t}\mathbf{V}=-\mathrm{i}\mathbf{LV}+\alpha _{0}^{m-1}\left[ 
\mathbf{f}_{0}^{\left( m\right) }\left( \mathbf{V}\right) +\alpha _{0}%
\mathbf{f}^{0\left( m+1\right) }\left( \mathbf{V}\right) +\ldots \right] ,\
\left. \mathbf{V}\right\vert _{t=0}=\mathbf{h},  \label{smin1}
\end{equation}%
which is of the form of (\ref{smallnon}) with $\alpha =\alpha _{0}^{m-1}$.
Introducing the slow time variable $\tau =t\alpha _{0}^{m-1}$ we get from
the above an equation of the form (\ref{difeqintr}), namely%
\begin{equation}
\partial _{\tau }\mathbf{V}=-\frac{\mathrm{i}}{\alpha _{0}^{m-1}}\mathbf{LV}+%
\left[ \mathbf{f}^{\left( m\right) }\left( \mathbf{V}\right) +\alpha _{0}%
\mathbf{f}^{\left( m+1\right) }\left( \mathbf{V}\right) +\ldots \right] ,\
\left. \mathbf{V}\right\vert _{t=0}=\mathbf{h},  \label{smin2}
\end{equation}%
where the nonlinearity does not vanish as $\alpha _{0}\rightarrow 0$. In
this case $\varrho =\alpha _{0}^{m-1}$ and the finite time interval $0\leq
\tau \leq \tau _{\ast }$ corresponds to the long time interval $0\leq t\leq 
\frac{\tau _{\ast }}{\alpha _{0}^{m-1}}$ \ with small $\alpha _{0}\ll 1$.

The third example is related to a \emph{high-frequency \ carrier wave} \ in
the initial data. To be concrete, we consider the Nonlinear Schrodinger
equation 
\begin{equation}
\partial _{\tau }U=-\mathrm{i}\partial _{x}^{2}U+\mathrm{i}\alpha \left\vert
U\right\vert ^{2}U,\ \left. U\right\vert _{\tau =0}=h_{1}\left( M\beta
x\right) e^{\mathrm{i}Mk_{\ast 1}x}+h_{2}\left( M\beta x\right) e^{\mathrm{i}%
Mk_{\ast 2}x}+c.c.,  \label{NLShM}
\end{equation}%
where $c.c.$ stands for complex conjugate of the prior term, and $M\gg 1$ is
a large parameter. The equation (\ref{NLShM}) can be readily recast into the
form (\ref{difeqintr}) by change of variables $y=Mr$ yielding 
\begin{gather}
\partial _{\tau }U=-\mathrm{i}\frac{1}{\varrho }\partial _{r}^{2}U+\mathrm{i}%
\alpha \left\vert U\right\vert ^{2}U,\ \left. U\right\vert _{\tau
=0}=h_{1}\left( \beta r\right) e^{\mathrm{i}k_{\ast 1}r}+h_{2}\left( \beta
r\right) e^{\mathrm{i}k_{\ast 2}r}+c.c.,  \label{NLShM1} \\
\text{where }\varrho =\frac{1}{M^{2}}\ll 1.  \notag
\end{gather}

Summarizing the above analysis we list below important ingredients of our
approach.

\begin{itemize}
\item The wave nonlinear evolution is analyzed based on the modal
decomposition with respect to the linear part of the system. The
significance of the modal decomposition to the nonlinear analysis is based
on the following properties: (i) the wave modal amplitudes do not evolve
under the linear evolution; (ii) the same amplitudes evolve slowly under the
nonlinear evolution; (iii) modal decomposition is instrumental to the
wavepacket definition including its spatial extension and the group velocity.

\item Components of multi-particle wavepacket are characterized by their
wavevectors $\mathbf{k}_{\ast l}$, band numbers $n_{l}$ and spatial
positions $\mathbf{r}_{\ast l}$. The nonlinear evolution preserves $\mathbf{k%
}_{\ast l}$ and $n_{l}$ whereas the spatial positions evolve uniformly with
the velocities $\frac{1}{\varrho }\nabla \omega _{n_{l}}\left( \mathbf{k}%
_{\ast l}\right) $.

\item The problem involves two small parameters $\beta $ and $\varrho $
respectively in the initial data and coefficients of the master equation (%
\ref{difeqintr}). These parameters scale respectively (i) the range of
wavevectors involved in its modal composition, with $\beta ^{-1}$ scaling
its spatial extension, and (ii) $\varrho $ scaling the ratio of the slow and
the fast time scales. We make no assumption on the functional dependence
between $\beta $ and $\varrho $, which are essentially independent and are
subject only to inequalities.

\item The nonlinear \emph{evolution is studied for a finite time} $\tau
_{\ast }$ which may depend on, say, the amplitude of the initial excitation,
and, importantly, $\tau _{\ast }$ is long enough to observe appreciable
nonlinear phenomena which are not vanishingly small. The superposition
principle can be extended to longer time intervals up to blow-up time or
even infinity if relevant uniform in $\beta $ and $\varrho $ estimates of
solutions in appropriate norms are available.

\item In the chosen slow time scale there are two fast wave processes with
typical time scale of order $\varrho $ which can be attributed to the linear
operator $\mathbf{L}$: (i) fast time oscillations resulting in time
averaging and consequent suppression of many nonlinear interactions; (ii)
fast wavepacket propagation with large group velocities resulting in
effective weakening of nonlinear interactions which are not time-averaged
because of resonances. It is these two processes provide mechanisms leading
to the superposition principle.
\end{itemize}

The rest of the paper is organized as follows. In the following Subsection
2.1 we introduce definitions of wavepackets, multiwavepackets and particle
wavepackets. In subsection 2.1 we also formulate and briefly discuss some
important results of \cite{BF7} which are used in this paper, and in
Subsection 2.2 we formulate new results. In Section 3 we formulate
conditions imposed on the linear and the nonlinear parts of the evolution
equation (\ref{difeqintr}), and also introduce relevant concepts describing
resonance interactions inside of wavepackets. In Section 4 we introduce an
integral form of the basic evolution equation and study basic properties of
involved operators. \ In Section 5 we introduce wavepacket interaction
system describing the dynamics of wavepackets. In Section 6 we, first,
define averaged wavepacket interaction system which plays a fundamental role
in the analysis of the dynamics of multiwavepackets and then prove that
solutions to this system approximate solutions to the original equation with
high accuracy. We also discuss there properties of averaged nonlinearities,
in particular, for universally and conditionally universal invariant
wavepackets, and prove the fundamental theorems on preservation of
multi-particle wavepackets, namely Theorems \ref{Theorem sumwavereg} and
Theorem \ref{Theorem regwave}. In Section 7 we prove the superposition
principle using an approximate decoupling of the averaged wavepacket
interaction system. \ In the last subsection of this section we prove some
generalizations to the cases involving non-generic resonance interactions
such as the second-harmonic and the third-harmonic generations.

\section{Statement of results}

This section consists of two subsections. In the first one we introduce
basic concepts and terminology and formulate relevant results from \cite{BF7}
which are used latter on, and in the second one we formulate new results of
this paper.

\subsection{Wavepackets and their basic properties}

Since the both linear operator $\mathbf{L}\left( -\mathrm{i}\nabla \right) $
and the nonlinearity $\mathbf{F}\left( \mathbf{U}\right) $ are translation
invariant, it is natural and convenient to recast the evolution equation (%
\ref{difeqintr}) by applying to it the Fourier transform with respect to the
space variables $\mathbf{r}$, namely 
\begin{equation}
\partial _{\tau }\mathbf{\hat{U}}\left( \mathbf{k}\right) =-\frac{\mathrm{i}%
}{\varrho }\mathbf{L}\left( \mathbf{k}\right) \mathbf{\hat{U}}\left( \mathbf{%
k}\right) +\hat{F}\left( \mathbf{\hat{U}}\right) \left( \mathbf{k}\right) ,\
\left. \mathbf{\hat{U}}\left( \mathbf{k}\right) \right\vert _{\tau =0}=%
\mathbf{\hat{h}}\left( \mathbf{k}\right) ,  \label{difeqfou}
\end{equation}%
where $\mathbf{\hat{U}}\left( \mathbf{k}\right) $ is the Fourier transform
of $\mathbf{U}\left( \mathbf{r}\right) $, i.e.%
\begin{equation}
\mathbf{\hat{U}}\left( \mathbf{k}\right) =\int_{\mathbb{R}^{d}}\mathbf{U}%
\left( \mathbf{r}\right) \mathrm{e}^{-\mathrm{i}\mathbf{r}\cdot \mathbf{k}}\,%
\mathrm{d}\mathbf{r},\text{ }\mathbf{U}\left( \mathbf{r}\right) =\left( 2\pi
\right) ^{-d}\int_{\mathbb{R}^{d}}\mathbf{\hat{U}}\left( \mathbf{k}\right) 
\mathrm{e}^{\mathrm{i}\mathbf{r}\cdot \mathbf{k}}\,\mathrm{d}\mathbf{r},%
\text{\ where }\mathbf{r},\mathbf{k}\in \mathbb{R}^{d},  \label{Ftrans}
\end{equation}%
and $\hat{F}$ is the Fourier form of the nonlinear operator $\mathbf{F}%
\left( \mathbf{U}\right) $ involving convolutions, see (\ref{Fmintr}) \ for
details. The equation (\ref{difeqfou}) is written in terms of Fourier modes,
and we call it the modal form of the original equation (\ref{difeqintr}).
The most of our studies are conducted first for the modal form (\ref%
{difeqfou}) of the evolution equation and carried over then to the original
equation (\ref{difeqintr}).

The nonlinear evolution equations (\ref{difeqintr}), (\ref{difeqfou}) are
commonly interpreted as describing wave propagation in a nonlinear medium.
We assume that the linear part $\mathbf{L}\left( \mathbf{\mathbf{k}}\right) $
is a $2J\times 2J$ Hermitian matrix with eigenvalues $\omega _{n,\zeta
}\left( \mathbf{k}\right) $ and eigenvectors $\mathbf{g}_{n,\zeta }\left( 
\mathbf{k}\right) $ satisfying%
\begin{equation}
\mathbf{L}\left( \mathbf{\mathbf{k}}\right) \mathbf{g}_{n,\zeta }\left( 
\mathbf{k}\right) =\omega _{n,\zeta }\left( \mathbf{k}\right) \mathbf{g}%
_{n,\zeta }\left( \mathbf{k}\right) ,\ \zeta =\pm ,\ \omega _{n,+}\left( 
\mathbf{k}\right) \geq 0,\ \omega _{n,-}\left( \mathbf{k}\right) \leq 0,\
n=1,\ldots ,J,  \label{OmomL}
\end{equation}%
where $\omega _{n,\zeta }\left( \mathbf{k}\right) $ are real-valued,
continuous for all non-singular $\mathbf{k}$ \ functions, and vectors $%
\mathbf{g}_{n,\zeta }\left( \mathbf{k}\right) \in \mathbb{C}^{2J}$ have unit
length in the standard Euclidean norm. The functions $\omega _{n,\zeta
}\left( \mathbf{k}\right) $, $n=1,\ldots ,J$, are called \emph{dispersion
relations} between the frequency $\omega $ and the \emph{wavevector} $%
\mathbf{k}$ with $n$ being the \emph{band number}. We assume that the
eigenvalues are naturally ordered by 
\begin{equation}
\omega _{J,+}\left( \mathbf{k}\right) \geq \ldots \geq \omega _{1,+}\left( 
\mathbf{k}\right) \geq 0\geq \omega _{1,-}\left( \mathbf{k}\right) \geq
\ldots \geq \omega _{J,-}\left( \mathbf{k}\right) ,  \label{omgr0}
\end{equation}%
and for \emph{almost every} $\mathbf{k}$ (with respect to the standard
Lebesgue measure) the eigenvalues are distinct and, consequently, the above
inequalities become strict. Importantly, we also assume the following \emph{%
diagonal symmetry} condition%
\begin{equation}
\omega _{n,-\zeta }\left( -\mathbf{k}\right) =-\omega _{n,\zeta }\left( 
\mathbf{k}\right) ,\ \zeta =\pm ,\ n=1,\ldots ,J,  \label{invsym}
\end{equation}%
which is naturally present in many physical problems (see also Remark \ref%
{Remark symmetry} below), and is a fundamental condition imposed on the
matrix $\mathbf{L}\left( \mathbf{\mathbf{k}}\right) $. Very often we use the
abbreviation 
\begin{equation}
\omega _{n,+}\left( \mathbf{k}\right) =\omega _{n}\left( \mathbf{k}\right) .
\label{ompl}
\end{equation}%
In particular we obtain from (\ref{invsym})%
\begin{equation}
\omega _{n,-}\left( \mathbf{k}\right) =-\omega _{n}\left( -\mathbf{k}\right)
,\ \omega _{n,\zeta }\left( \mathbf{k}\right) =\zeta \omega _{n}\left( \zeta 
\mathbf{k}\right) ,\ \zeta =\pm .  \label{omz}
\end{equation}%
In addition to that in many examples we also have%
\begin{equation}
\mathbf{g}_{n,\zeta }\left( \mathbf{k}\right) =\mathbf{g}_{n,-\zeta }^{\ast
}\left( -\mathbf{k}\right) ,\text{ where }z^{\ast }\text{ is complex
conjugate to }z.  \label{omgr0a}
\end{equation}%
We also use rather often the orthogonal projection $\Pi _{n,\zeta }\left( 
\mathbf{\mathbf{k}}\right) $ in $\mathbb{C}^{2J}$ onto the complex line
defined by the eigenvector $\mathbf{g}_{n,\zeta }\left( \mathbf{k}\right) $,
namely 
\begin{equation}
\Pi _{n,\zeta }\left( \mathbf{\mathbf{k}}\right) \mathbf{\hat{u}}\left( 
\mathbf{k}\right) =\tilde{u}_{n,\zeta }\left( \mathbf{k}\right) \mathbf{g}%
_{n,\zeta }\left( \mathbf{k}\right) =\mathbf{\hat{u}}_{n,\zeta }\left( 
\mathbf{k}\right) ,\ n=1,\ldots ,J,\ \zeta =\pm .  \label{Pin}
\end{equation}

As it is indicated by the title of this paper we study the nonlinear problem
(\ref{difeqintr}) for initial data $\mathbf{\hat{h}}$ in the form of a
properly defined particle-like \emph{wavepackets} or, more generally, a sum
of such wavepackets to which we refer as \emph{multi-particle wavepacket}.
The simplest example of a wavepacket $\mathbf{w}$ is provided by the
following formula 
\begin{equation}
\mathbf{w}\left( \beta ;\mathbf{r}\right) =\Phi _{+}\left( \beta \left( 
\mathbf{r}-\mathbf{r}_{\ast }\right) \right) \mathrm{e}^{\mathrm{i}\mathbf{k}%
_{\ast }\cdot \left( \mathbf{r-r}_{\ast }\right) }\mathbf{g}_{n,+}\left( 
\mathbf{k}_{\ast }\right) ,\ \mathbf{r}\in \mathbb{R}^{d},  \label{wpint}
\end{equation}%
where $\mathbf{k}_{\ast }\in \mathbb{R}^{d}$ is \emph{wavepacket principal
wavevector}, $n$ is \emph{band number}, and $\beta >0$ is a small parameter.
We refer to the pair $\left( n,\mathbf{k}_{\ast }\right) $ in (\ref{wpint})
as \emph{wavepacket }$nk$\emph{-pair} and $\mathbf{r}_{\ast }$ as \emph{%
wavepacket position. }Observe that the space extension of the wavepacket $%
\mathbf{w}\left( \beta ;\mathbf{r}\right) $ is proportional to $\beta ^{-1}$
and it is large for small $\beta $. Notice also that if $\beta \rightarrow 0$
the wavepacket $\mathbf{w}\left( \beta ;\mathbf{r}\right) $ as in (\ref%
{wpint}) tends, up to a constant factor, to the elementary eigenmode $%
\mathrm{e}^{\mathrm{i}\mathbf{k}_{\ast }\cdot \mathbf{r}}\mathbf{g}_{n,\zeta
}\left( \mathbf{k}_{\ast }\right) $ of the operator $\mathbf{L}\left( -%
\mathrm{i}\nabla \right) $ with the corresponding eigenvalue $\omega
_{n,\zeta }\left( \mathbf{k}_{\ast }\right) $. We refer to wavepackets \ of
the simple form (\ref{wpint}) as \emph{simple wavepackets} to underline the
very special way the parameter $\beta $ enters its representation. The
function $\Phi _{\zeta }\left( \mathbf{r}\right) $, which we call \emph{%
wavepacket envelope}, describes its shape and it can be any scalar
complex-valued regular enough function, for example a function from Schwartz
space. Importantly, as $\beta \rightarrow 0$ the $L^{\infty }$ \emph{norm of
a wavepacket \ (\ref{wpint}) remains constant, hence nonlinear effects in (%
\ref{difeqintr}) remain strong}.

Evolution of wavepackets in problems which can be reduced to the form (\ref%
{difeqintr}) were studied for a variety of equations in numerous physical
and mathematical papers, mostly by asymptotic expansions with respect to a
single small parameter similar to $\beta $, see \cite{BenYoussefLannes02}, 
\cite{BonaCL05}, \cite{ColinLannes}, \cite{CraigSulemS92}, \cite{GiaMielke}, 
\cite{JolyMR98}, \cite{KalyakinUMN}, \cite{Maslov83}, \cite{PW}, \cite%
{Schneider98a}, \cite{Schneider05} and references therein. We are interested
in general properties of evolutionary systems of the form (\ref{difeqintr})
with wavepacket initial data which hold for a wide class of nonlinearities
and all values of the space dimensions $d$ and the number $2J$ of the system
components. Our approach is not based on asymptotic expansions but involves
the two small parameters $\beta $ and $\varrho $ with mild constraints (\ref%
{rbb1}) on their relative smallness. The constraints can be expressed either
in the form of certain inequalities or equalities, and a possible simple
form of such a constraint can be a power law 
\begin{equation}
\beta =C\varrho ^{\varkappa }\text{ where }C>0\text{ and }\varkappa >0\text{
are arbitrary constants.}  \label{powerk}
\end{equation}%
Of course, general features of wavepacket evolution are independent of
particular values of\ the constant $C$. In addition to that, some
fundamental properties such as wavepacket preservation are also totally
independent on particular choice of the values of $\varkappa $ in (\ref%
{powerk}), whereas other properties are independent of $\varkappa $ as it
varies in certain intervals. For for instance, dispersion effects are
dominant for $\varkappa <1/2$, whereas the\ wavepacket superposition
principle of \cite{BF7} holds for $\varkappa <1$.

To eliminate unbounded (as $\varrho \rightarrow 0$) linear term in (\ref%
{difeqfou}) by replacing it with\ a highly oscillatory factor we introduce
the \emph{slow variable} $\mathbf{\hat{u}}\left( \mathbf{k},\tau \right) $
by the formula 
\begin{equation}
\mathbf{\hat{U}}\left( \mathbf{k},\tau \right) =\mathrm{e}^{-\frac{\mathrm{i}%
\tau }{\varrho }\mathbf{L}\left( \mathbf{k}\right) }\mathbf{\hat{u}}\left( 
\mathbf{k},\tau \right) ,  \label{Uu0}
\end{equation}%
and get the following equation for $\mathbf{\hat{u}}\left( \mathbf{k},\tau
\right) $ 
\begin{equation}
\partial _{\tau }\mathbf{\hat{u}}=\mathrm{e}^{\frac{\mathrm{i}\tau }{\varrho 
}\mathbf{L}}\mathbf{\hat{F}}\left( \mathrm{e}^{\frac{-\mathrm{i}\tau }{%
\varrho }\mathbf{L}}\mathbf{\hat{u}}\right) ,\ \left. \mathbf{\hat{u}}%
\right\vert _{\tau =0}=\mathbf{\hat{h}},  \label{dfsF}
\end{equation}%
which, in turn, can be transformed by time integration into the integral
form 
\begin{equation}
\mathbf{\hat{u}}=\mathcal{F}\left( \mathbf{\hat{u}}\right) +\mathbf{\hat{h}}%
,\ \mathcal{F}\left( \mathbf{\hat{u}}\right) =\int_{0}^{\tau }\mathrm{e}^{%
\frac{\mathrm{i}\tau ^{\prime }}{\varrho }\mathbf{L}}\mathbf{\hat{F}}\left( 
\mathrm{e}^{\frac{-\mathrm{i}\tau ^{\prime }}{\varrho }\mathbf{L}}\mathbf{%
\hat{u}}\left( \tau ^{\prime }\right) \right) \mathrm{d}\tau ^{\prime }
\label{ubaseq}
\end{equation}%
with explicitly defined nonlinear polynomial integral operator $\mathcal{F}=%
\mathcal{F}\left( \varrho \right) $. This operator is bounded uniformly with
respect to $\varrho $ in the Banach space $E=C\left( \left[ 0,\tau _{\ast }%
\right] ,L^{1}\right) .$This space has functions $\mathbf{\hat{v}}\left( 
\mathbf{k},\tau \right) $, $0\leq \tau \leq \tau _{\ast }$ as elements \ and
has \ the norm 
\begin{equation}
\left\Vert \mathbf{\hat{v}}\left( \mathbf{k},\tau \right) \right\Vert
_{E}=\left\Vert \mathbf{\hat{v}}\left( \mathbf{k},\tau \right) \right\Vert
_{C\left( \left[ 0,\tau _{\ast }\right] ,L^{1}\right) }=\sup_{0\leq \tau
\leq \tau _{\ast }}\int_{\mathbb{R}^{d}}\left\vert \mathbf{\hat{v}}\left( 
\mathbf{k},\tau \right) \right\vert \,\mathrm{d}\mathbf{k},  \label{Elat}
\end{equation}%
where $L^{1}$ is the Lebesgue space of functions $\mathbf{\hat{v}}\left( 
\mathbf{k}\right) $ with the standard norm 
\begin{equation}
\left\Vert \mathbf{\hat{v}}\left( \mathbf{\cdot }\right) \right\Vert
_{L^{1}}=\int_{\mathbb{R}^{d}}\left\vert \mathbf{\hat{v}}\left( \mathbf{k}%
\right) \right\vert \,\mathrm{d}\mathbf{k}.  \label{L1}
\end{equation}%
Sometimes we use more general weighted spaces $L^{1,a}$ with the norm 
\begin{equation}
\left\Vert \mathbf{\hat{v}}\right\Vert _{L^{1,a}}=\int_{\mathbb{R}%
^{d}}\left( 1+\left\vert \mathbf{k}\right\vert \right) ^{a}\left\vert 
\mathbf{\hat{v}}\left( \mathbf{k}\right) \right\vert \,\mathrm{d}\mathbf{k}%
,\ a\geq 0.  \label{L1a}
\end{equation}%
The space $C\left( \left[ 0,\tau _{\ast }\right] ,L^{1,a}\right) $ with the
norm \ 
\begin{equation}
\left\Vert \mathbf{\hat{v}}\left( \mathbf{k},\tau \right) \right\Vert
_{E_{a}}=\sup_{0\leq \tau \leq \tau _{\ast }}\int_{\mathbb{R}^{d}}\left(
1+\left\vert \mathbf{k}\right\vert \right) ^{a}\left\vert \mathbf{\hat{v}}%
\left( \mathbf{k},\tau \right) \right\vert \mathrm{d}\mathbf{k}.  \label{Ea}
\end{equation}%
is denoted by $E_{a}$, and, obviously, $E_{0}=E$.

A rather elementary existence and uniqueness theorem (Theorem \ref{Theorem
Existence1}) implies that if $\mathbf{\hat{h}}\in L^{1,a}$ then for a small
and, importantly, independent of $\varrho $ constant $\tau _{\ast }>0$ this
equation has a unique solution 
\begin{equation}
\mathbf{\hat{u}}\left( \tau \right) =\mathcal{G}\left( \mathcal{F}\left(
\varrho \right) ,\mathbf{\hat{h}}\right) \left( \tau \right) ,\ \tau \in %
\left[ 0,\tau _{\ast }\right] ,\ \mathbf{\hat{u}}\in C^{1}\left( \left[
0,\tau _{\ast }\right] ,L^{1,a}\right) ,  \label{ugh1}
\end{equation}%
where $\mathcal{G}$ denotes the \emph{solution operator} for the equation (%
\ref{ubaseq}). If $\mathbf{\hat{u}}\left( \mathbf{k},\tau \right) $ is a
solution to the equation (\ref{ubaseq}) we call the function $\mathbf{U}%
\left( \mathbf{r},\tau \right) $ defined by (\ref{Uu0}), (\ref{Ftrans}) an $%
\emph{F}$\emph{-solution} to the equation (\ref{difeqintr}). We denote by $%
\hat{L}^{1}$ the space of functions $\mathbf{V}\left( \mathbf{r}\right) $
such that their Fourier transform $\mathbf{\hat{V}}\left( \mathbf{k}\right) $
belongs to $L^{1}$, and define $\left\Vert \mathbf{V}\right\Vert _{\hat{L}%
^{1}}=\left\Vert \mathbf{\hat{V}}\right\Vert _{L^{1}}$. Since 
\begin{equation}
\left\Vert \mathbf{V}\right\Vert _{L^{\infty }}\leq \left( 2\pi \right)
^{-d}\left\Vert \mathbf{\hat{V}}\right\Vert _{L^{1}}\text{ and }\hat{L}%
^{1}\subset L^{\infty },  \label{Linf}
\end{equation}%
$F$-solutions to (\ref{difeqintr}) belong to $C^{1}\left( \left[ 0,\tau
_{\ast }\right] ,\hat{L}^{1}\right) \subset C^{1}\left( \left[ 0,\tau _{\ast
}\right] ,L^{\infty }\right) $.

We would like to define wavepackets in a form which explicitly allows them
to be real valued. This is accomplished based on the symmetry (\ref{invsym})
of the dispersion relations, which allows to introduce a \emph{doublet
wavepacket} 
\begin{equation}
\mathbf{w}\left( \beta ;\mathbf{r}\right) =\Phi _{+}\left( \beta \left( 
\mathbf{r}-\mathbf{r}_{\ast }\right) \right) \mathrm{e}^{\mathrm{i}\mathbf{k}%
_{\ast }\cdot \left( \mathbf{r-r}_{\ast }\right) }\mathbf{g}_{n,+}\left( 
\mathbf{k}_{\ast }\right) +\Phi _{-}\left( \beta \left( \mathbf{r}-\mathbf{r}%
_{\ast }\right) \right) \mathrm{e}^{-\mathrm{i}\mathbf{k}_{\ast }\cdot
\left( \mathbf{r-r}_{\ast }\right) }\mathbf{g}_{n,-}\left( -\mathbf{k}_{\ast
}\right) .  \label{wpn1}
\end{equation}%
Such a wavepacket is real if $\Phi _{-}\left( \mathbf{r}\right) $, $\mathbf{g%
}_{n,-}\left( -\mathbf{k}_{\ast }\right) $ are complex conjugate
respectively to $\Phi _{+}\left( \mathbf{r}\right) $, $\mathbf{g}%
_{n,+}\left( \mathbf{k}_{\ast }\right) $, i.e. if%
\begin{equation}
\Phi _{-}\left( \mathbf{r}\right) =\Phi _{+}^{\ast }\left( \mathbf{r}\right)
,\ \mathbf{g}_{n,+}\left( \mathbf{k}_{\ast }\right) =\mathbf{g}_{n,-}\left( -%
\mathbf{k}_{\ast }\right) ^{\ast }.  \label{wpn1a}
\end{equation}%
Usually when considering wavepackets with $\ nk$-pair $\left( n,\mathbf{k}%
_{\ast }\right) $ we mean doublet ones as in (\ref{wpn1}), but sometimes we
use the term wavepacket also for an elementary one as defined by (\ref{wpint}%
). Note that the latter use is consistent with the former one since it is
possible to take one of two terms in (\ref{wpn1}) to be zero.

Below we give a precise definition of a wavepacket. To identify
characteristic properties of a wavepacket suitable for our needs, let us
look at the Fourier transform $\mathbf{\hat{w}}\left( \beta ;\mathbf{k}%
\right) $ of an elementary wavepacket $\mathbf{w}\left( \beta ;\mathbf{r}%
\right) $ defined by (\ref{wpint}), that is 
\begin{equation}
\mathbf{\hat{w}}\left( \beta ;\mathbf{k}\right) =\beta ^{-d}\mathrm{e}^{-%
\mathrm{i}\mathbf{k}\cdot \mathbf{r}_{\ast }}\hat{\Phi}\left( \beta
^{-1}\left( \mathbf{k}-\mathbf{k}_{\ast }\right) \right) \mathbf{g}_{n,\zeta
}\left( \mathbf{k}_{\ast }\right) .  \label{wpint1}
\end{equation}%
We call such $\mathbf{\hat{w}}\left( \beta ;\mathbf{k}\right) $ wavepacket
too and it possesses the following properties: (i) its $L^{1}$ norm is
bounded (in fact, constant), uniformly in $\beta \rightarrow 0$; (ii) for
every $\epsilon >0$ the value $\mathbf{\hat{w}}\left( \beta ;\mathbf{k}%
\right) \rightarrow 0$ for every $\mathbf{k}$ outside a $\beta ^{1-\epsilon
} $-neighborhood of $\mathbf{k}_{\ast }$, and the convergence is faster than
any power of $\beta $ if $\Phi $ is a Schwartz function. To explicitly
interpret the last property we introduce a \emph{cutoff function} $\Psi
\left( \mathbf{\eta }\right) $ which is infinitely smooth and such that 
\begin{equation}
\Psi \left( \mathbf{\eta }\right) \geq 0,\ \Psi \left( \mathbf{\eta }\right)
=1\text{ for }\left\vert \mathbf{\eta }\right\vert \leq 1/2,\ \Psi \left( 
\mathbf{\eta }\right) =0\text{ for }\left\vert \mathbf{\eta }\right\vert
\geq 1,  \label{j0}
\end{equation}%
and its shifted/rescaled modification 
\begin{equation}
\Psi \left( \beta ^{1-\epsilon },\mathbf{k}_{\ast };\mathbf{k}\right) =\Psi
\left( \beta ^{-\left( 1-\epsilon \right) }\left( \mathbf{k}-\mathbf{k}%
_{\ast }\right) \right) .  \label{Psik}
\end{equation}%
If an elementary wavepacket $\mathbf{w}\left( \beta ;\mathbf{r}\right) $ is
defined by (\ref{wpint1}) with $\Phi \left( \mathbf{r}\right) $ being a
Schwartz function then%
\begin{equation}
\left\Vert \left( 1-\Psi \left( \beta ^{1-\epsilon },\mathbf{k}_{\ast
};\cdot \right) \right) \mathbf{\hat{w}}\left( \beta ;\cdot \right)
\right\Vert \leq C_{\epsilon ,s}\beta ^{s},\ 0<\beta \leq 1,  \label{halfH0}
\end{equation}%
and the inequality holds for arbitrarily small $\epsilon >0$ and arbitrarily
large $s>0$. Based on the above discussion we give the following definition
of a wavepacket which is a minor variation of \cite[Definiton 8]{BF8}.

\begin{definition}[single-band wavepacket]
\label{dwavepack} Let $\epsilon $ \ be a fixed number, $0<\epsilon <1$. For
a given band number $n\in \left\{ 1,\ldots ,J\right\} $ and a principal
wavevector $\mathbf{k}_{\ast }\in \mathbb{R}^{d}$ a function $\mathbf{\hat{h}%
}\left( \beta ;\mathbf{k}\right) $ is called a \emph{\ wavepacket with }$nk$%
\emph{-pair }$\left( n,\mathbf{k}_{\ast }\right) $ and the degree of
regularity $s>0$ if for small $\beta <\beta _{0}$ with some $\beta _{0}>0$
it satisfies the following conditions: (i) $\mathbf{\hat{h}}\left( \beta ;%
\mathbf{k}\right) $ is $L^{1}$-bounded uniformly in $\beta $, i.e.%
\begin{equation}
\left\Vert \mathbf{\hat{h}}\left( \beta ;\mathbf{\cdot }\right) \right\Vert
_{L^{1}}\leq C,\ 0<\beta <\beta _{0}\text{ for some }C>0;  \label{L1b}
\end{equation}%
(ii) $\mathbf{\hat{h}}\left( \beta ;\mathbf{k}\right) $\ is composed
essentially of two functions $\mathbf{\hat{h}}_{\zeta }\left( \beta ;\mathbf{%
k}\right) $, $\zeta =\pm $, which take values in $n$-th band eigenspace of $%
\mathbf{L}\left( \mathbf{k}\right) $ and are localized near $\zeta \mathbf{k}%
_{\ast }$, namely\ 
\begin{equation}
\mathbf{\hat{h}}\left( \beta ;\mathbf{k}\right) =\mathbf{\hat{h}}_{-}\left(
\beta ;\mathbf{k}\right) +\mathbf{\hat{h}}_{+}\left( \beta ;\mathbf{k}%
\right) +D_{h},\ 0<\beta <\beta _{0},  \label{hbold}
\end{equation}%
where\ the components $\mathbf{\hat{h}}_{\pm }\left( \beta ;\mathbf{k}%
\right) $ satisfy the condition%
\begin{equation}
\mathbf{\hat{h}}_{\zeta }\left( \beta ;\mathbf{k}\right) =\Psi \left( \beta
^{1-\epsilon }/2,\zeta \mathbf{k}_{\ast };\mathbf{k}\right) \Pi _{n,\zeta
}\left( \mathbf{k}\right) \mathbf{\hat{h}}_{\zeta }\left( \beta ;\mathbf{k}%
\right) ,\ \zeta =\pm ,  \label{halfwave}
\end{equation}%
where $\Psi \left( \mathbf{\cdot },\zeta \mathbf{k}_{\ast },\beta
^{1-\epsilon }\right) $ is defined by (\ref{Psik}) and $D_{h}$ is small,
namely it satisfies the inequality 
\begin{equation}
\left\Vert D_{h}\right\Vert _{L^{1}}\leq C^{\prime }\beta ^{s},\ 0<\beta
<\beta _{0},\text{ for some }C^{\prime }>0.  \label{sourloc}
\end{equation}%
The inverse Fourier transform $\mathbf{h}\left( \beta ;\mathbf{r}\right) $
of a wavepacket $\mathbf{\hat{h}}\left( \beta ;\mathbf{k}\right) $ is also
called a wavepacket.
\end{definition}

Evidently, if a wavepacket has the degree of regularity $s$, it also has a
smaller degree of regularity $s^{\prime }\leq s$ with the same $\epsilon $.
Observe that the degree of regularity $s$ is related to the smoothness of $%
\Phi _{\zeta }\left( \mathbf{r}\right) $ as in (\ref{wpint}) so that the
higher is the smoothness the higher $\frac{s}{\epsilon }$ can be taken.
Namely, if $\hat{\Phi}_{\zeta }\in L^{1,a}$ then one can take in (\ref%
{sourloc}) any $s<a\epsilon $ according to the following inequality: 
\begin{equation}
\int \left\vert \left( 1-\Psi \left( \beta ^{\epsilon }\mathbf{\eta }\right)
\right) \hat{\Phi}_{\zeta }\left( \mathbf{\eta }\right) \right\vert d\mathbf{%
\eta }\leq \beta ^{a\epsilon }\left\Vert \hat{\Phi}_{\zeta }\right\Vert
_{L^{1,a}}\leq C\beta ^{s}.  \label{halfH}
\end{equation}

For example, if \ we define $\mathbf{\hat{h}}_{\zeta }$ similarly to (\ref%
{halfwave}) and (\ref{wpint1})\ by the formula 
\begin{equation}
\mathbf{\hat{h}}_{\zeta }\left( \beta ;\mathbf{k}\right) =\Psi \left( \beta
^{-\left( 1-\epsilon \right) }\left( \mathbf{k}-\mathbf{k}_{\ast }\right)
\right) \beta ^{-d}\hat{\Phi}_{\zeta }\left( \beta ^{-1}\left( \mathbf{k}-%
\mathbf{k}_{\ast }\right) \right) \Pi _{n,\zeta }\left( \mathbf{k}\right) 
\mathbf{g}  \label{h2bet}
\end{equation}%
where $\hat{\Phi}_{\zeta }\left( \mathbf{k}\right) $ is a scalar Schwartz
function and $\mathbf{g}$ is a vector, then according to(\ref{halfH})
estimate (\ref{sourloc}) \ holds and \ $\mathbf{\hat{h}}_{\zeta }\left(
\beta ;\mathbf{k}\right) $ is a wavepacket with arbitrarily large degree of
regularity $s$ for any given $\epsilon $ such that $0<\epsilon <1$.

Now let us define a particle-like wavepacket following to the ideas
indicated in the Introduction.

\begin{definition}[ single-band particle-like wavepacket]
\label{Definition regwave} We call a function $\mathbf{\hat{h}}\left( \beta ;%
\mathbf{k}\right) =\mathbf{\hat{h}}\left( \beta ,\mathbf{r}_{\ast };\mathbf{k%
}\right) $, $\mathbf{r}_{\ast }\in \mathbb{R}^{d}$, a \emph{particle-like
wavepacket with the position} $\mathbf{r}_{\ast }$, $nk$-pair\emph{\ }$%
\left( n,\mathbf{k}_{\ast }\right) $ and the degree of regularity $s>0$ if
(i) for every $\mathbf{r}_{\ast }$ it is a wavepacket with the degree of
regularity $s$\ in the sense of the above Definition \ref{dwavepack} \ with
constants $C,C^{\prime }$ independent of $\mathbf{r}_{\ast }\in \mathbb{R}%
^{d}$; (ii) $\mathbf{\hat{h}}_{\zeta }$ in (\ref{hbold}) satisfy the
inequalities 
\begin{equation}
\int_{\mathbb{R}^{d}}\left\vert \nabla _{\mathbf{k}}\left( e^{i\mathbf{r}%
_{\ast }\mathbf{k}}\mathbf{\hat{h}}_{\zeta }\left( \beta ,\mathbf{r}_{\ast };%
\mathbf{k}\right) \right) \right\vert d\mathbf{k}\leq C_{1}\beta
^{-1-\epsilon },\ \zeta =\pm ,\ \mathbf{r}_{\ast }\in \mathbb{R}^{d},
\label{gradker}
\end{equation}%
where $C_{1}>0$ is an independent of $\beta $ and $\mathbf{r}_{\ast }$
constant, $\epsilon $ is the same as in Definition \ref{dwavepack}. The
inverse Fourier transform $\mathbf{h}\left( \beta ;\mathbf{r}\right) $ of a
wavepacket $\mathbf{\hat{h}}\left( \beta ;\mathbf{k}\right) $ is also called
a particle-like wavepacket with the position $\mathbf{r}_{\ast }$. We also
introduce quantity 
\begin{equation}
a\left( \mathbf{r}_{\ast }^{\prime },\mathbf{\hat{h}}_{\zeta }\left( \mathbf{%
r}_{\ast }\right) \right) =\left\Vert \nabla _{\mathbf{k}}\left( e^{i\mathbf{%
r}_{\ast }^{\prime }\mathbf{k}}\mathbf{\hat{h}}_{\zeta }\left( \beta ,%
\mathbf{r}_{\ast };\mathbf{k}\right) \right) \right\Vert _{L^{1}}  \label{ar}
\end{equation}%
which we refer to as the \emph{position detection function }for the
wavepacket $\mathbf{\hat{h}}\left( \beta ,\mathbf{r}_{\ast };\mathbf{k}%
\right) $.
\end{definition}

Note that the left-hand side of (\ref{gradker}) coincides with $a\left( 
\mathbf{r}_{\ast },\mathbf{\hat{h}}_{\zeta }\left( \mathbf{r}_{\ast }\right)
\right) $.

\begin{remark}
\label{Remark hdecay}If $\mathbf{\hat{h}}\left( \beta ;\mathbf{k}\right) =%
\mathbf{\hat{h}}\left( \beta ,\mathbf{r}_{\ast };\mathbf{k}\right) $ is a
particle-like wavepacket with a position $\mathbf{r}_{\ast }$ then, applying
inverse Fourier transform to $\mathbf{\hat{h}}_{\zeta }\left( \beta ,\mathbf{%
r}_{\ast };\mathbf{k}\right) $ and $\nabla _{\mathbf{k}}\mathbf{\hat{h}}%
_{\zeta }\left( \beta ,\mathbf{r}_{\ast };\mathbf{k}\right) $ as in (\ref%
{Ftrans}) we obtain a function $\mathbf{h}\left( \beta ,\mathbf{r}_{\ast };%
\mathbf{r}\right) $ which satisfies 
\begin{equation}
\left\vert \mathbf{r}-\mathbf{r}_{\ast }\right\vert \left\vert \mathbf{h}%
_{\zeta }\left( \beta ,\mathbf{r}_{\ast };\mathbf{r}\right) \right\vert \leq
\left( 2\pi \right) ^{-d}a\left( \mathbf{r}_{\ast },\mathbf{\hat{h}}_{\zeta
}\right)   \label{hdec}
\end{equation}%
implying that $\left\vert \mathbf{h}_{\zeta }\left( \beta ;\mathbf{r}\right)
\right\vert \leq a\left( \mathbf{r}_{\ast },\mathbf{\hat{h}}_{\zeta }\right)
\left\vert \mathbf{r}-\mathbf{r}_{\ast }\right\vert ^{-1}$. This inequality
is useful for large $\left\vert \mathbf{r}-\mathbf{r}_{\ast }\right\vert $,
\ whereas for bounded $\left\vert \mathbf{r}-\mathbf{r}_{\ast }\right\vert $
(\ref{L1b}) implies a simpler inequality 
\begin{equation}
\left\vert \mathbf{h}_{\zeta }\left( \beta ,\mathbf{r}_{\ast };\mathbf{r}%
\right) \right\vert \leq \left( 2\pi \right) ^{-d}\left\Vert \mathbf{\hat{h}}%
\right\Vert _{L^{1}}\leq C.  \label{hdec1}
\end{equation}%
Inequalities (\ref{hdec}) and (\ref{gradker}) suggest that the quantity $%
a\left( \mathbf{r}_{\ast },\mathbf{\hat{h}}_{\zeta }\left( \mathbf{r}_{\ast
}\right) \right) $ can be interpreted as a size of the particle-like
wavepacket $\mathbf{\hat{h}}_{\zeta }\left( \beta ,\mathbf{r}_{\ast };%
\mathbf{k}\right) $.
\end{remark}

Evidently a particle-like wavepacket is a wave and not a point. Hence the
above definition of its position has a degree of uncertainty, allowing, for
example, to replace $\mathbf{r}_{\ast }$ by $\mathbf{r}_{\ast }+\mathbf{a}$
with a \emph{fixed} vector $\mathbf{a}$ (but not allowing unbounded values
of $\mathbf{a}$). The above definition of particle-like wavepacket position
was crafted to meet the following requirements: (i) a system of
particle-like wavepackets remains to be such a system under the nonlinear
evolution; (ii) it is possible (in an appropriate scale) to describe the
trajectories traced out by the positions of a system of particle-wavepackets

\begin{remark}
\label{aloc}Typical dependence of the inverse Fourier transform $\mathbf{h}%
\left( \beta ,\mathbf{r}_{\ast };\mathbf{r}\right) $ of a wavepacket $%
\mathbf{\hat{h}}\left( \beta ,\mathbf{r}_{\ast };\mathbf{k}\right) $ on $%
\mathbf{r}_{\ast }$ is provided by spatial shifts by $\mathbf{r}_{\ast }$ as
in (\ref{wpn1}), namely 
\begin{equation*}
\mathbf{h}\left( \beta ,\mathbf{r}_{\ast };\mathbf{r}\right) =\Phi \left(
\beta \left( \mathbf{r}-\mathbf{r}_{\ast }\right) \right) \mathrm{e}^{%
\mathrm{i}\mathbf{k}_{\ast }\cdot \left( \mathbf{r}-\mathbf{r}_{\ast
}\right) }\mathbf{g}
\end{equation*}
with a constant $\mathbf{g}$. For such a function $\mathbf{h}$ and for any $%
\mathbf{r}_{\ast }^{\prime }\in \mathbb{R}^{d}$ 
\begin{gather*}
a\left( \mathbf{r}_{\ast }^{\prime },\mathbf{\hat{h}}\left( \mathbf{r}_{\ast
}\right) \right) =\left\Vert \nabla _{\mathbf{k}}\left( \beta ^{-d}\mathrm{e}%
^{\mathrm{i}\mathbf{kr}_{\ast }^{\prime }}\mathbf{\hat{h}}\left( \beta ,%
\mathbf{r}_{\ast };\mathbf{k}\right) \right) \right\Vert _{L^{1}}=\left\Vert
\nabla _{\mathbf{k}}\left( \beta ^{-d}\mathrm{e}^{\mathrm{i}\mathbf{kr}%
_{\ast }^{\prime }}\mathrm{e}^{-\mathrm{i}\mathbf{kr}_{\ast }}\hat{\Phi}%
\left( \mathbf{k}\right) \right) \right\Vert _{L^{1}}\left\Vert \mathbf{g}%
\right\Vert  \\
=\left\Vert \mathbf{g}\right\Vert \int \left\vert \mathrm{i}\left( \mathbf{r}%
_{\ast }^{\prime }-\mathbf{r}_{\ast }\right) \hat{\Phi}\left( \mathbf{k}%
^{\prime }\right) +\frac{1}{\beta }\nabla _{\mathbf{k}^{\prime }}\hat{\Phi}%
\left( \mathbf{k}^{\prime }\right) \right\vert d\mathbf{k}^{\prime }.
\end{gather*}%
Hence, taking for simplicity $\left\Vert \mathbf{g}\right\Vert =1$, we obtain%
\begin{equation}
\left\vert \mathbf{r}_{\ast }^{\prime }-\mathbf{r}_{\ast }\right\vert
\left\Vert \hat{\Phi}\right\Vert _{L_{1}}+\frac{1}{\beta }\left\Vert \nabla 
\hat{\Phi}\right\Vert _{L_{1}}\geq a\left( \mathbf{r}_{\ast }^{\prime },%
\mathbf{\hat{h}}\left( \mathbf{r}_{\ast }\right) \right) \geq \left\vert
\left\vert \mathbf{r}_{\ast }^{\prime }-\mathbf{r}_{\ast }\right\vert
\left\Vert \hat{\Phi}\right\Vert _{L_{1}}-\frac{1}{\beta }\left\Vert \nabla 
\hat{\Phi}\right\Vert _{L_{1}}\right\vert .  \label{rstar}
\end{equation}%
For small $\left\vert \mathbf{r}_{\ast }^{\prime }-\mathbf{r}_{\ast
}\right\vert \ll \frac{1}{\beta }$ we see that the position detection
function $a\left( \mathbf{r}_{\ast }^{\prime },\mathbf{\hat{h}}\right) $ is
of order $O\left( \beta ^{-1}\right) $, which is in the agreement with (\ref%
{gradker}). For large $\left\vert \mathbf{r}_{\ast }^{\prime }-\mathbf{r}%
_{\ast }\right\vert \gg \frac{1}{\beta }$ the $a\left( \mathbf{r}_{\ast
}^{\prime },\mathbf{\hat{h}}\right) $ is approximately proportional to $%
\left\vert \mathbf{r}_{\ast }^{\prime }-\mathbf{r}_{\ast }\right\vert $.
Therefore if we know $a\left( \mathbf{r}_{\ast }^{\prime },\mathbf{\hat{h}}%
\left( \mathbf{r}_{\ast }\right) \right) $ as a function of $\mathbf{r}%
_{\ast }^{\prime }$ we can recover the value of $\mathbf{r}_{\ast }$ with
the accuracy of order $O\left( \beta ^{-1-\epsilon }\right) $ with arbitrary
small $\epsilon $. Namely, let us take arbitrary small $\epsilon >0$ and
some $C>0$ and consider the set 
\begin{equation}
B\left( \beta \right) =\left\{ \mathbf{r}_{\ast }^{\prime }:a\left( \mathbf{r%
}_{\ast }^{\prime },\mathbf{\hat{h}}\left( \mathbf{r}_{\ast }\right) \right)
\leq C\beta ^{-1-\epsilon }\right\} \subset \mathbb{R}^{d},  \label{Bbeta}
\end{equation}%
which should provide an approximate location of $\mathbf{r}_{\ast }$.
According to (\ref{rstar}), $\mathbf{r}_{\ast }$ lies in this set for small $%
\beta $. If $\mathbf{r}_{\ast }^{\prime }$ lies in this set then \ 
\begin{equation*}
C\beta ^{-1-\epsilon }\geq a\left( \mathbf{r}_{\ast }^{\prime },\mathbf{\hat{%
h}}\left( \mathbf{r}_{\ast }\right) \right) \geq \left\vert \left\vert 
\mathbf{r}_{\ast }^{\prime }-\mathbf{r}_{\ast }\right\vert \left\Vert \hat{%
\Phi}\right\Vert _{L_{1}}-\frac{1}{\beta }\left\Vert \nabla \hat{\Phi}%
\right\Vert _{L_{1}}\right\vert 
\end{equation*}%
and $\left\vert \mathbf{r}_{\ast }^{\prime }-\mathbf{r}_{\ast }\right\vert
\leq C_{1}\beta ^{-1-\epsilon }+C_{2}\beta ^{-1}$. Hence the diameter of the 
$B\left( \beta \right) $ is of order $O\left( \beta ^{-1-\epsilon }\right) $%
. Observe, taking into account Remark \ref{Remark hdecay}, that the accuracy
of the wavepacket location obviously cannot be better than its size $a\left( 
\mathbf{r}_{\ast },\mathbf{\hat{h}}_{\zeta }\left( \mathbf{r}_{\ast }\right)
\right) \sim \beta ^{-1}$. The above analysis suggest that the function $%
\mathbf{h}\left( \beta ,\mathbf{r}_{\ast };\mathbf{r}\right) $ can be viewed
as \emph{pseudoshifts} of the function $\mathbf{h}\left( \beta ,\mathbf{0};%
\mathbf{r}\right) $ by vectors $\mathbf{r}_{\ast }\in \mathbb{R}^{d}$ in the
sense that the regular spatial shift by $\mathbf{r}_{\ast }$ is combined
with a variation of the shape of $\mathbf{h}\left( \beta ,\mathbf{0};\mathbf{%
r}\right) $ which is limited by the fundamental condition (\ref{gradker}).
In other words, according Definition \ref{Definition regwave} as
wavepacecket moves from $\mathbf{0}$ to $\mathbf{r}_{\ast }$ by a
corresponding spatial shift it is allowed to change its shape subject to the
fundamental condition (\ref{gradker}). The later is instrumental for
capturing nonlinear evolution of particle-like wavepackets governed by an
equation of the form (\ref{difeqintr}).
\end{remark}

\begin{remark}
\label{alocgen} The set $B\left( \beta \right) $ defined by (\ref{Bbeta})
gives an approximate location of the support of the function $\mathbf{\hat{h}%
}\left( \beta ,\mathbf{r}_{\ast };\mathbf{k}\right) $ not only in the
special case considered in Remark \ref{aloc}, but also when $\mathbf{h}%
\left( \beta ,\mathbf{r}_{\ast };\mathbf{r}\right) $ is a general
particle-like wavepacket. One can apply with obvious modifications the above
argument for $\mathrm{e}^{\mathrm{i}\mathbf{kr}_{\ast }}\mathbf{\hat{h}}%
\left( \beta ,\mathbf{r}_{\ast };\mathbf{k}\right) $ in place of $\hat{\Phi}%
\left( \mathbf{k}\right) $ using (\ref{gradker}). Here we give an
alternative argument based on (\ref{hdec}). Notice that condition $a\left( 
\mathbf{r}_{\ast 0},\mathbf{\hat{h}}\left( \mathbf{r}_{\ast }\right) \right)
\leq C\beta ^{-1-\epsilon }$ obviously can be satisfied not only by $\mathbf{%
r}_{\ast 0}=\mathbf{r}_{\ast }$. But one can show that the diameter of the
set of such $\mathbf{r}_{\ast 0}$ is estimated by $O\left( \beta
^{-1-\epsilon }\right) $. Indeed, assume that a given function $\mathbf{h}%
\left( \beta ,\mathbf{r}\right) $ does not vanish at a given point $\mathbf{r%
}_{0}$, that is $\left\vert \mathbf{h}\left( \beta ,\mathbf{r}_{0}\right)
\right\vert \geq c_{0}>0$ \ for all $\beta \leq \beta _{0}$. The fulfillment
of (\ref{gradker}) for the function $\mathbf{h}\left( \beta ,\mathbf{r}%
\right) $ with two different values of $\mathbf{r}_{\ast }$, namely $\mathbf{%
r}_{\ast }=\mathbf{r}_{\ast }^{\prime }$ \ and $\mathbf{r}_{\ast }=\mathbf{r}%
_{\ast }^{\prime \prime }$\ implies that 
\begin{equation*}
a\left( \mathbf{r}_{\ast }^{\prime },\mathbf{\hat{h}}\right) \leq C_{1}\beta
^{-1-\epsilon },\ a\left( \mathbf{r}_{\ast }^{\prime \prime },\mathbf{\hat{h}%
}\right) \leq C_{2}\beta ^{-1-\epsilon },
\end{equation*}%
and according to (\ref{hdec}) for all $\mathbf{r}$\ 
\begin{equation*}
\left\vert \mathbf{r}-\mathbf{r}_{\ast }^{\prime }\right\vert \left\vert 
\mathbf{h}\left( \beta ,\mathbf{r}\right) \right\vert \leq \left( 2\pi
\right) ^{-d}C_{1}\beta ^{-1-\epsilon },\ \left\vert \mathbf{r}-\mathbf{r}%
_{\ast }^{\prime \prime }\right\vert \left\vert \mathbf{h}\left( \beta ,%
\mathbf{r}\right) \right\vert \leq \left( 2\pi \right) ^{-d}C_{2}\beta
^{-1-\epsilon },
\end{equation*}%
Hence,%
\begin{equation*}
\left\vert \mathbf{r}_{0}-\mathbf{r}_{\ast }^{\prime }\right\vert \leq \frac{%
\left( 2\pi \right) ^{-d}C_{1}\beta ^{-1-\epsilon }}{c_{0}},\ \left\vert 
\mathbf{r}_{0}-\mathbf{r}_{\ast }^{\prime \prime }\right\vert \leq \frac{%
\left( 2\pi \right) ^{-d}C_{1}\beta ^{-1-\epsilon }}{c_{0}},
\end{equation*}%
and 
\begin{equation*}
\left\vert \mathbf{r}_{\ast }^{\prime }-\mathbf{r}_{\ast }^{\prime \prime
}\right\vert \leq C_{3}\beta ^{-1-\epsilon }.
\end{equation*}
\end{remark}

Note that if we rescale variables $\mathbf{r}$ and $\mathbf{r}_{\ast }$ as
in Example \ref{Example yscale}, namely $\varrho \mathbf{r}=\mathbf{y}$ and $%
\varrho \mathbf{r}_{\ast }=\mathbf{y}_{\ast }$\ with $\varrho =\beta ^{2}$,
the diameter of the set $B\left( \beta \right) $ in $y$- coordinates is of
order  $\beta ^{1-\epsilon }\ll 1$, and, hence, this set gives a good
approximation for the location of the particle-like wavepacket \ as $\beta
\rightarrow 0$. It is important to notice, that our method to locate the
support of wavepackets is applicable to very general wavepackets and does
not use their specific form. This flexibility allows us to prove that
particle-like wavepackets and their positions are well defined during
nonlinear dynamics \ of generic equations with rather general initial data \
which form infinite-dimensional function spaces. Another approaches to
describe dynamics of waves are applied to situations where solutions under
considerations can be parametrized by a finite number of parameters and the
dynamics of the parameters describes dynamics of the solutions. See for
example \cite{GangSigal06a}, \cite{GangSigal06} where dynamics of centers of
solitions is described.

\begin{remark}
\label{simple}Note that for a single wavepacket initial data $\mathbf{h}%
\left( \beta ,\mathbf{r}-\mathbf{r}_{\ast }^{\prime }\right) $ one can make
a change of variables to a moving frame $\left( \mathbf{x},\tau \right) $,
namely $\left( \mathbf{r},\tau \right) =\left( \mathbf{x}+\mathbf{v}\tau
,\tau \right) $, where $\mathbf{v}=\frac{1}{\varrho }\nabla \omega \left( 
\mathbf{k}_{\ast }\right) $ is the group velocity; this change of variables
makes the group velocity zero. Often it is possible to prove that dynamics
preserves functions which decay at infinity , namely if the initial data $\ 
\mathbf{h}\left( \beta ,\mathbf{x}\right) $ decays at the spatial infinity
then the solution $\mathbf{U}\left( \beta ,\mathbf{x},\tau \right) $ also
decays at infinity (though corresponding proofs can be rather technical). \
This property can be reformulated in rescaled $\mathbf{y}$ variables as
follows: if initial data are localized about zero, then the solution is
localized about zero as well. Then, using the fact that the equation has
constant coefficients, we observe that the solution $\mathbf{U}\left( \beta ,%
\mathbf{y}-\mathbf{y}_{\ast }^{\prime },\tau \right) $, corresponding to $%
\mathbf{h}\left( \beta ,\mathbf{y}-\mathbf{y}_{\ast }^{\prime }\right) $, is
localized about $\mathbf{y}_{\ast }^{\prime }$ provided that $\mathbf{h}%
\left( \beta ,\mathbf{y}\right) $ was localized about the origin. Note that
in this paper we consider much more complicated case of multiple
wavepackets. Even in the simplest case of the initial multiwavepacket which
involves only two components, namely the wavepacket $\mathbf{h}\left( \beta ,%
\mathbf{r}\right) =\mathbf{h}_{1}\left( \beta ,\mathbf{r}-\mathbf{r}_{\ast
}^{\prime }\right) +\mathbf{h}_{2}\left( \beta ,\mathbf{r}-\mathbf{r}_{\ast
}^{\prime \prime }\right) $ with two principal wave vectors $\mathbf{k}%
_{1\ast }\neq \mathbf{k}_{2\ast }$ \ one evidently cannot use the above
considerations based on the change of variables and the translational
invarance. Using other arguments devoloped in this paper we prove that
systems of particle-like wavepackets remain localized in the process of the
nonlinear evolution.
\end{remark}

Note that similarly to (\ref{hh}) and (\ref{hr1}) a function of the form 
\begin{equation*}
\beta ^{-d}\left( \mathrm{e}^{-\mathrm{i}\mathbf{kr}_{\ast 1}}+\mathrm{e}^{-%
\mathrm{i}\mathbf{kr}_{\ast 2}}\right) \left[ \hat{h}\left( \frac{\mathbf{k}-%
\mathbf{k}_{\ast }}{\beta }\right) \right] \mathbf{g}_{n}\left( \mathbf{k}%
_{\ast }\right) ,
\end{equation*}%
defined for any pair of $\mathbf{r}_{\ast 1}$ and $\mathbf{r}_{\ast 2}$
where $\hat{h}$ is a Schwartz function and all constants in Definition \ref%
{dwavepack} are independent of $\mathbf{r}_{\ast 1},\mathbf{r}_{\ast 2}\in 
\mathbb{R}^{d}$, is not a single particle-like wavepacket since it does not
have a single wavepacket position $\mathbf{r}_{\ast }$, but rather it is a
sum of two particle-like wavepackets with two positions $\mathbf{r}_{\ast 1}$
and $\mathbf{r}_{\ast 2}$.

We want to emphasize once more that a particle-like wavepacket is defined as
the family $\mathbf{\hat{h}}\left( \beta ,\mathbf{r}_{\ast };\mathbf{k}%
\right) $ with $\mathbf{r}_{\ast }$ being an independent variable running
the entire space $\mathbb{R}^{d}$, see, for example, (\ref{hh}), (\ref{hr})
and (\ref{wpn1}). In particular, we can choose a dependence of $\mathbf{r}%
_{\ast }$ on $\beta $\ and $\varrho $. An interesting type of such a
dependence is $\mathbf{r}_{\ast }=\mathbf{r}_{\ast }^{0}/\varrho $ where $%
\varrho $ satisfies (\ref{powerk}) as we discuss below in the Example \ref%
{Example yscale}.

Our special interest is in the waves that are finite sums of wavepackets
which we refer to as \emph{multi-wavepackets}.

\begin{definition}[multi-wavepacket]
\label{dmwavepack} Let $S$ be a set of $nk$-pairs: 
\begin{equation}
S=\left\{ \left( n_{l},\mathbf{k}_{\ast l}\right) ,\ l=1,\ldots ,N\right\}
\subset \Sigma =\left\{ 1,\ldots ,J\right\} \times \mathbb{R}^{d},\ \left(
n_{l},\mathbf{k}_{\ast l}\right) \neq \left( n_{l^{\prime }},\mathbf{k}%
_{\ast l^{\prime }}\right) \text{ for }l\neq l^{\prime },  \label{P0}
\end{equation}%
and $N=\left\vert S\right\vert $ be their number. Let $K_{S}$ be a set
consisting of all different wavevectors $\mathbf{k}_{\ast l}$ involved in $S$
with $\left\vert K_{S}\right\vert \leq N$ being the number of its elements. $%
K_{S}$ is called \emph{wavepacket }$k$\emph{-spectrum} and without loss of
genericity we assume the indexing of elements $\ \left( n_{l},\mathbf{k}%
_{\ast l}\right) $ in $S$ to be such that 
\begin{equation}
K_{S}=\left\{ \mathbf{k}_{\ast i},i=1,\ldots ,\left\vert K_{S}\right\vert
\right\} ,\text{ i.e. }l=i\text{ for }1\leq i\leq \left\vert
K_{S}\right\vert \text{.}  \label{K0}
\end{equation}%
A function $\mathbf{\hat{h}}\left( \beta \right) =\mathbf{\hat{h}}\left(
\beta ;\mathbf{k}\right) $ is called a \emph{multi-wavepacket} with $nk$%
\emph{-spectrum} $S$ if it is a finite sum of wavepackets, namely 
\begin{equation}
\mathbf{\hat{h}}\left( \beta ;\mathbf{k}\right) =\sum_{l=1}^{N}\mathbf{\hat{h%
}}_{l}\left( \beta ;\mathbf{k}\right) ,\ 0<\beta <\beta _{0}\text{ for some }%
\beta _{0}>0,  \label{JJ1}
\end{equation}%
where $\mathbf{\hat{h}}_{l}$, $l=1,\ldots ,N$, is a wavepacket with $nk$%
-pair $\left( \mathbf{k}_{\ast l},n_{l}\right) \in S$ as in Definition \ref%
{dwavepack}. If all the wavepackets $\mathbf{\hat{h}}_{l}\left( \beta ;%
\mathbf{k}\right) =\mathbf{\hat{h}}_{l}\left( \beta ,\mathbf{r}_{\ast l};%
\mathbf{k}\right) $ are particle-like ones with respective positions $%
\mathbf{r}_{\ast l}$ then the multi-wavepacket is called \emph{%
multi-particle wavepacket} and we refer to $\left( \mathbf{r}_{\ast
1},\ldots ,\mathbf{r}_{\ast N}\right) $ as its position vector.
\end{definition}

Note that if $\mathbf{\hat{h}}\left( \beta ;\mathbf{k}\right) $ is a
wavepacket then $\mathbf{\hat{h}}\left( \beta ;\mathbf{k}\right) +O\left(
\beta ^{s}\right) $ is also a wavepacket \ with the same $nk$-spectrum, and
the same is true for multi-wavepackets. Hence, we can introduce
multi-wavepackets equivalence relation "$\simeq $" of the degree $s$ by 
\begin{equation}
\mathbf{\hat{h}}_{1}\left( \beta ;\mathbf{k}\right) \simeq \mathbf{\hat{h}}%
_{2}\left( \beta ;\mathbf{k}\right) \text{ if }\left\Vert \mathbf{\hat{h}}%
_{1}\left( \beta ;\mathbf{k}\right) -\mathbf{\hat{h}}_{2}\left( \beta ;%
\mathbf{k}\right) \right\Vert _{L^{1}}\leq C\beta ^{s}\text{ for some
constant \ }C>0.  \label{equiv}
\end{equation}
Note that condition (\ref{gradker}) does not impose restrictions on the term 
$D_{h}$ in (\ref{hbold}), therefore this equivalence can be applied to
particle wavepackets.

Let us turn now to the abstract nonlinear problem (\ref{ubaseq}) where (i) $%
\mathcal{F}=\mathcal{F}\left( \varrho \right) $ depends on $\varrho $ and
(ii) the initial data $\mathbf{\hat{h}}=\mathbf{\hat{h}}\left( \beta \right) 
$\ is a multi-wavepacket depending on $\beta $. We would like to state our
first theorem on multi-wavepacket preservation under the evolution (\ref%
{ubaseq}) for $\beta ,\varrho \rightarrow 0$, which holds provided its $nk$%
\emph{-}spectrum $S$ satisfies a natural condition called \emph{resonance
invariance}. This condition is intimately related to the so-called \emph{%
phase and frequency matching} conditions for stronger nonlinear
interactions, and its concise formulation is as follows. We define for given
dispersion relations $\left\{ \omega _{n}\left( \mathbf{k}\right) \right\} $
and any finite set $S\subset \left\{ 1,\ldots ,J\right\} \times \mathbb{R}%
^{d}$ another finite set $\mathcal{R}\left( S\right) \subset \left\{
1,\ldots ,J\right\} \times \mathbb{R}^{d}$ where $\mathcal{R}$ is a certain
algebraic operation described in Definition \ref{Definition omclos} below.
It turns out that for any $S$ always $S\subseteq \mathcal{R}\left( S\right) $%
, but \emph{if }$\mathcal{R}\left( S\right) =S$\emph{\ we call }$S$\emph{\
resonance invariant}. The condition of resonance invariance is instrumental
for the multi-wavepacket preservation, and there are examples showing that
if it fails, i.e. $\mathcal{R}\left( S\right) \neq S$, the wavepacket
preservation does not hold. \emph{\ }Importantly, the resonance invariance $%
R\left( S\right) =S$\ allows resonances inside the multi-wavepacket, that
includes, in particular, resonances associated with the second and the third
harmonic generations, resonant four-wave interaction etc. \ In this paper we
will use basic results on wavepacket preservation obtained in \cite{BF7} ,
and we formulate theorems from \cite{BF7} we need here. Since we use
constructions from \cite{BF7}, for completeness we provide also their proofs
in the following subsections. The following two theorems are proved in \cite%
{BF7}.

\begin{theorem}[multi-wavepacket preservation]
\label{Theorem invarwave}Suppose that the nonlinear evolution is governed by
(\ref{ubaseq}) and the initial data $\mathbf{\hat{h}}=\mathbf{\hat{h}}\left(
\beta ;\mathbf{k}\right) $ is a multi-wavepacket with $nk$-spectrum $S$ and
the regularity degree $s$, and assume $S$ to be resonance invariant (see
Definition \ref{Definition omclos} below). Let $\rho \left( \beta \right) $
be any function satisfying 
\begin{equation}
0<\rho \left( \beta \right) \leq C\beta ^{s},\text{ for some constant }C>0,
\label{rbr}
\end{equation}%
and let us set $\varrho =\rho \left( \beta \right) $. Then the solution $%
\mathbf{\hat{u}}\left( \tau ,\beta \right) =\mathcal{G}\left( \mathcal{F}%
\left( \rho \left( \beta \right) \right) ,\mathbf{\hat{h}}\left( \beta
\right) \right) \left( \tau \right) $ to (\ref{ubaseq}) for any $\tau \in %
\left[ 0,\tau _{\ast }\right] $ is a multi-wavepacket with $nk$-spectrum $S$
and the regularity degree $s$, i.e.%
\begin{equation}
\mathbf{\hat{u}}\left( \tau ,\beta ;\mathbf{k}\right)
=\dsum\nolimits_{l=1}^{N}\mathbf{\hat{u}}_{l}\left( \tau ,\beta ;\mathbf{k}%
\right) ,\text{ where }\mathbf{\hat{u}}_{l}\text{ is wavepacket with }nk%
\text{-pair }\left( n_{l},\mathbf{k}_{\ast l}\right) \in S.  \label{usumul}
\end{equation}%
The time interval length $\tau _{\ast }>0$ depends only on $L^{1}$-norms of $%
\mathbf{\hat{h}}_{l}\left( \beta ;\mathbf{k}\right) $ and $N$. The
presentation (\ref{usumul}) is unique up to the equivalence (\ref{equiv}) of
degree $s$.
\end{theorem}

The above statement can be interpreted as follows. Modes in $nk$-spectrum $S$
are\ always resonance coupled with modes in $\mathcal{R}\left( S\right) $
through the nonlinear interactions, but if $\mathcal{R}\left( S\right) =S$
then (i) all resonance interactions occur inside $S$ and (ii) only small
vicinity of $S$ is involved in nonlinear interactions leading to the
multi-wavepacket preservation.

The statement of Theorems \ref{Theorem invarwave} directly follows from the
following general theorem proved in \cite{BF7}.

\begin{theorem}[multi-wavepacket approximation]
\label{Theorem sumwave} Let the initial data $\mathbf{\hat{h}}$ in the
integral equation (\ref{ubaseq}) be a multi-wavepacket $\mathbf{\hat{h}}%
\left( \beta ;\mathbf{k}\right) $ with $nk$-spectrum $S$ as in (\ref{P0}),
the regularity degree $s$ and with the parameter $\epsilon >0$ as in
Definition \ref{dwavepack}. \ Assume that $S$ is resonance invariant in the
sense of Definition \ref{Definition omclos} below. Let the cutoff function $%
\Psi \left( \beta ^{1-\epsilon },\mathbf{k}_{\ast };\mathbf{k}\right) $ and
the eigenvector projectors $\Pi _{n,\pm }\left( \mathbf{\mathbf{k}}\right) $
be defined by (\ref{Psik}) and (\ref{Pin}) respectively. For a solution $%
\mathbf{\hat{u}}$ of (\ref{ubaseq}) \ we set%
\begin{equation}
\mathbf{\hat{u}}_{l}\left( \beta ;\tau ,\mathbf{k}\right) =\left[
\dsum\nolimits_{\zeta =\pm }\Psi \left( C\beta ^{1-\epsilon },\zeta \mathbf{k%
}_{\ast l};\mathbf{k}\right) \Pi _{n_{l},\zeta }\left( \mathbf{\mathbf{k}}%
\right) \right] \mathbf{\hat{u}}\left( \beta ;\tau ,\mathbf{k}\right) ,\
l=1,\ldots ,N.  \label{ups1}
\end{equation}%
Then every such $\mathbf{\hat{u}}_{l}\left( \beta ;\tau ,\mathbf{k}\right) $
is a wavepacket and%
\begin{equation}
\sup_{0\leq \tau \leq \tau _{\ast }}\left\Vert \mathbf{\hat{u}}\left( \beta
;\tau ,\mathbf{k}\right) -\dsum\nolimits_{l=1}^{N}\mathbf{\hat{u}}_{l}\left(
\beta ;\tau ,\mathbf{k}\right) \right\Vert _{L^{1}}\leq C_{1}\varrho
+C_{2}\beta ^{s}  \label{uui}
\end{equation}%
where the constants $C,C_{1}$ do not depend on $\epsilon ,s$ and $\beta $,
and the constant $C_{2}$ does not depend on $\beta $and $\varrho $.
\end{theorem}

We would like to point out also that Theorem \ref{Theorem invarwave} allows
to take values $\mathbf{\hat{u}}\left( \tau _{\ast }\right) $ as new
wavepacket initial data for (\ref{difeqintr}) and extend the wavepacket \
invariance of a solution to the next time interval $\tau _{\ast }\leq \tau
\leq \tau _{\ast 1}$. This observation allows to extend the wavepacket \
invariance to larger values of $\tau $ (up to blow-up time or infinity) if
some additional information about solutions with wavepacket initial data is
available, see \cite{BF7}.

Note that the wavepacket form of solutions can be used to obtain long-time
estimates of solutions. Namely, very often behavior of every single
wavepacket is well approximated by its own nonlinear Schrodinger equation
(NLS), see \cite{Colin}, \cite{KSM}, \cite{ColinLannes}, \cite{GiaMielke}, 
\cite{KalyakinUMN}, \cite{Kalyakin2}, \cite{PW}, \cite{Schneider98a}, \cite%
{Schneider05}, \cite{SU} and references therein, see also Section 6. Many
features of the dynamics governed by NLS-type equations are well-understood,
see \cite{Bourgain}, \cite{Caz}, \cite{SchlagK}, \cite{Schlag}, \cite{Sulem}%
, \cite{Weinstein} and references therein. These results can be used to
obtain long-time estimates for every single wavepacket (as, for example, in 
\cite{Kalyakin2}) and, with the help of the superposition principle, for the
multiwavepacket solution.

\subsection{\textbf{Formulation of new results} on particle wavepackets}

In this paper we prove the following refinement of Theorem \ref{Theorem
invarwave} for the case of multi-particle wavepackets.

\begin{theorem}[multi-particle wavepacket preservation]
\label{Theorem regwave}Assume that conditions of Theorem \ref{Theorem
sumwave} hold and, in addition to that, the initial data $\mathbf{\hat{h}}=%
\mathbf{\hat{h}}\left( \beta ;\mathbf{k}\right) $ is a multi-particle
wavepacket of degree $s$ with positions $\mathbf{r}_{\ast 1},\ldots ,\mathbf{%
r}_{\ast N}$ and the multi-particle wavepacket is universally resonance
invariant in the sense of Definition \ref{Definition omclos}. Assume also
that 
\begin{equation}
\rho \left( \beta \right) \leq C\beta ^{s_{0}},s_{0}>0.  \label{ros0}
\end{equation}%
Then the solution $\mathbf{\hat{u}}\left( \beta ;\tau \right) =\mathcal{G}%
\left( \mathcal{F}\left( \rho \left( \beta \right) \right) ,\mathbf{\hat{h}}%
\left( \beta \right) \right) \left( \tau \right) $ to (\ref{ubaseq}) for any 
$\tau \in \left[ 0,\tau _{\ast }\right] $ is a multi-particle wavepacket
with the same $nk$-spectrum $S$ and the same positions $\mathbf{r}_{\ast
1},\ldots ,\mathbf{r}_{\ast N}$. Namely, (\ref{uui}) holds where $\ \mathbf{%
\hat{u}}_{l}$ is wavepacket with $nk$-pair $\left( n_{l},\mathbf{k}_{\ast
l}\right) \in S$ \ defined by (\ref{ups1}), the constants \ $C,C_{1},C_{2}$
do not depend on $\mathbf{r}_{\ast l}$. and every $\ \mathbf{\hat{u}}_{l}$
is equivalent in the sense of the equivalence (\ref{equiv}) of degree $%
s_{1}=\min \left( s,s_{0}\right) $ to a particle wavepacket with the
position $\mathbf{r}_{\ast l}$.
\end{theorem}

\begin{remark}
\label{Pemark pospres}Note that in the statement of the above theorem the
positions $\mathbf{r}_{\ast 1},\ldots ,\mathbf{r}_{\ast N}$ of wavepackets
which compose the solution $\mathbf{\hat{u}}\left( \beta ;\tau ,\mathbf{k}%
\right) $ of (\ref{dfsF}) and (\ref{ubaseq}) do not depend on $\tau $ and,
hence, do not move. Note also that the solution $\mathbf{\hat{U}}\left(
\beta ;\tau ,\mathbf{k}\right) $ of the original equation (\ref{difeqfou}),
related to $\mathbf{\hat{u}}\left( \beta ;\tau ,\mathbf{k}\right) $\ by the
change of variables (\ref{Uu0}), is composed of wavepackets $\mathbf{U}%
_{l}\left( \beta ;\tau ,\mathbf{r}\right) $, corresponding to $\mathbf{u}%
_{l}\left( \beta ;\tau ,\mathbf{r}\right) $, have their positions moving
with respective constant velocities $\nabla _{k}\omega \left( \mathbf{k}%
_{\ast l}\right) $ (see for details Remark \ref{Remark shift}, see also the
following corollary).
\end{remark}

Using Proposition \ref{Proposition Uu} \ we obtain from Theorem \ref{Theorem
regwave} the following corollary.

\begin{corollary}
Let conditions of Theorem \ref{Theorem regwave} hold and $\mathbf{\hat{U}}%
\left( \beta ;\tau ,\mathbf{k}\right) $ be defined by (\ref{Uu0}) in terms
of $\mathbf{\hat{u}}\left( \beta ;\tau ,\mathbf{k}\right) $. Let 
\begin{equation}
\frac{\beta ^{2}}{\varrho }\leq C,\ \text{with \ some }C,\text{\ }0<\beta
\leq \frac{1}{2},\ 0<\varrho \leq \frac{1}{2}.  \label{scale1}
\end{equation}%
Then $\mathbf{\hat{U}}\left( \beta ;\tau ,\mathbf{k}\right) $ is for every $%
\tau \in \left[ 0,\tau _{\ast }\right] $ a particle multi-wavepacket in the
sense of Definition \ref{Definition regwave} with the same $nk$\emph{-}%
spectrum\emph{\ }$S$, regularity $s_{1}$ and with $\tau $-dependent
positions $\mathbf{r}_{\ast l}+\frac{\tau }{\varrho }\nabla _{k}\omega
_{n}\left( \mathbf{k}_{\ast l}\right) $.
\end{corollary}

In the following example we consider the case where spatial positions of
wavepackets have a specific dependence on parameter $\varrho $, namely $%
\mathbf{r}_{\ast }=\mathbf{r}_{\ast }^{0}/\varrho $.

\begin{example}[Wavepacket trajectories and collisions]
\label{Example yscale}Let us rescale the coordinates in the physical space
as follows 
\begin{equation}
\varrho \mathbf{r}=\mathbf{y}  \label{rhoscale}
\end{equation}%
with the consequent rescaling of the wavevector variable (dual with respect
to Fourier transform) $\mathbf{k}=\varrho \mathbf{\eta }$. It follows then
that under the evolution (\ref{difeqintr}) the group velocity of a
wavepacket with a wavevector $\mathbf{k}_{\ast }$ in the new coordinates $%
\mathbf{y}$ becomes $\nabla _{k}\omega \left( \mathbf{k}_{\ast }\right) $
and evidently is of order one. If we set the positions $\mathbf{r}_{\ast l}=%
\mathbf{r}_{\ast l}^{0}/\varrho $ with fixed $\mathbf{r}_{\ast l}^{0}$, then
according to (\ref{hdec}) wavepackets $\left\vert \mathbf{h}\left( \beta ;%
\mathbf{r}\right) \right\vert $ in $\mathbf{y}$-variables have
characteristic spatial scale $\mathbf{y}-\mathbf{r}_{\ast l}^{0}\sim \varrho
a\left( \mathbf{r}_{\ast l},\mathbf{\hat{h}}\right) \sim \varrho \beta ^{-1}$
which is small if $\varrho /\beta \ $is small. The positions of
particle-like wavepackets (quasiparticles) $\mathbf{\hat{U}}\left( \mathbf{y}%
/\varrho ,\tau \right) $ are initially located at $\mathbf{y}_{l}=\mathbf{r}%
_{\ast l}^{0}$ and propagate with the group velocities $\nabla _{k}\omega
\left( \mathbf{k}_{\ast l}\right) $. Their trajectories are straight lines
in the space $\mathbb{R}^{d}$ described by 
\begin{equation*}
\mathbf{y}=\tau \nabla _{k}\omega \left( \mathbf{k}_{\ast l}\right) +\mathbf{%
r}_{\ast l}^{0},\ 0\leq \tau \leq \tau _{\ast }
\end{equation*}%
(compare with (\ref{hr2})). The trajectories may intersect, indicating
"collisions" of quasiparticles. Our results (Theorem \ref{Theorem regwave})
show that if a multi-particle wavepacket initially was universally resonance
invariant, then the involved particle-like wavepackets preserve their
identity in spite of collisions and the fact that the nonlinear interactions
with other wavepackets (quasiparticles) are not small, if fact, they are of
order one. Note that $\mathbf{r}_{\ast l}^{0}$ can be chosen arbitrarily
implying that up to $N\left( N-1\right) $ collisions can occur on the time
interval $\left[ 0,\tau _{\ast }\right] $ on which we study the system
evolution.
\end{example}

To formulate the approximate superposition principle \ for multi-particle
wavepackets, we introduce now the \emph{solution operator} $\mathcal{G}$
mapping the initial data $\mathbf{\hat{h}}$ into the solution $\mathbf{\hat{U%
}}$ $=\mathcal{G}\left( \mathbf{\hat{h}}\right) $ of the modal evolution
equation (\ref{ubaseq}). This operator is defined for $\left\Vert \mathbf{%
\hat{h}}\right\Vert \leq R$ according to the existence and uniqueness
Theorem \ref{Theorem exist}. The main result of this paper is the following
statement.

\begin{theorem}[superposition principle ]
\label{Theorem Superposition}\ Suppose that the initial data $\mathbf{\hat{h}%
}$ of (\ref{ubaseq}) is a multi-particle wavepacket of the form%
\begin{equation}
\mathbf{\hat{h}}=\sum_{l=1}^{N}\mathbf{\hat{h}}_{l},\ N\max_{l}\left\Vert 
\mathbf{\hat{h}}_{l}\right\Vert _{L^{1}}\leq R,  \label{hsumR}
\end{equation}%
satisfying Definition \ref{dmwavepack} and its $nk$-spectrum is universally
resonance invariant in the sense of Definition \ref{Definition omclos}.
Suppose also that that the group velocities of wavepackets are different,
namely 
\begin{equation}
\nabla _{\mathbf{k}}\omega _{n_{l_{1}}}\left( \mathbf{k}_{\ast l_{1}}\right)
\neq \nabla _{\mathbf{k}}\omega _{n_{l_{2}}}\left( \mathbf{k}_{\ast
l_{2}}\right) \text{ if }l_{1}\neq l_{2}  \label{NGVM}
\end{equation}%
and that (\ref{scale1}) holds. Then the solution $\mathbf{\hat{u}}$ $=%
\mathcal{G}\left( \mathbf{\hat{h}}\right) $ to the evolution equation (\ref%
{ubaseq}) satisfies the following approximate superposition principle 
\begin{equation}
\mathcal{G}\left( \sum_{l=1}^{N_{h}}\mathbf{\hat{h}}_{l}\right)
=\sum_{l=1}^{N_{h}}\mathcal{G}\left( \mathbf{\hat{h}}_{l}\right) +\mathbf{%
\tilde{D}},  \label{apprsup}
\end{equation}%
with a small remainder $\mathbf{\tilde{D}}\left( \tau \right) $ such that 
\begin{equation}
\sup_{0\leq \tau \leq \tau _{\ast }}\left\Vert \mathbf{\tilde{D}}\left( \tau
\right) \right\Vert _{L^{1}}\leq C_{\epsilon }\frac{\varrho }{\beta
^{1+\epsilon }}\left\vert \ln \beta \right\vert ,  \label{rem}
\end{equation}%
where (i) $\epsilon $ is the same as in Definition \ref{dwavepack} and can
be arbitrary small; (ii) $\tau _{\ast }$ does not depend on $\beta $, $%
\varrho $, $\mathbf{r}_{\ast l}$\ and $\epsilon $; (iii) $C_{\epsilon }$
does not depend on $\beta ,\varrho $ and positions $\mathbf{r}_{\ast l}$\ .
\end{theorem}

A particular case of the above Theorem in which there was no dependence on $%
\mathbf{r}_{\ast l}$ was proved in \cite{BF8} by a different method based on
the theory of analytic operators in Banach spaces. The condition (\ref{NGVM}%
) can be relaxed if the initial positions of involved particle-like
wavepackets are far apart, and the the corresponding results are formulated
in the theorem below and in Example \ref{Example yscale}.

\begin{theorem}[superposition principle ]
\label{Theorem Superposition1}\ Suppose that the initial data $\mathbf{\hat{h%
}}$ of (\ref{ubaseq}) is a multi-particle wavepacket of the form (\ref{hsumR}%
) with a universally resonance invariant $nk$-spectrum in the sense of
Definition \ref{Definition omclos} and (\ref{scale1}) holds . Suppose also
that either the group velocities of wavepackets are different, namely (\ref%
{NGVM})holds, or the positions $\mathbf{r}_{\ast l}$ \ satisfy the
inequality 
\begin{equation}
\tau _{\ast }\left\vert \mathbf{r}_{\ast l_{1}}-\mathbf{r}_{\ast
l_{2}}\right\vert ^{-1}\leq \frac{\varrho }{2C_{\omega ,2}\beta ^{1-\epsilon
}}\text{ if }\nabla _{\mathbf{k}}\omega _{n_{l_{1}}}\left( \mathbf{k}_{\ast
l_{1}}\right) =\nabla _{\mathbf{k}}\omega _{n_{l_{2}}}\left( \mathbf{k}%
_{\ast l_{2}}\right) \text{,\ }l_{1}\neq l_{2},  \label{NGVMe}
\end{equation}%
where the constant $C_{\omega ,2}$ is the same as in (\ref{Com2}). Then the
solution $\mathbf{\hat{u}}=\mathcal{G}\left( \mathbf{\hat{h}}\right) $ to
the evolution equation (\ref{ubaseq}) satisfies the approximate
superposition principle (\ref{apprsup}), (\ref{rem}).
\end{theorem}

We prove in this paper further generalizations of the particle-like
wavepacket preservation and the superposition principle to the cases where
the $nk$-spectrum of a multi-wavepacket \ is not universal resonance
invariant such as the cases of multi-wavepackets involving the second and
the third harmonic generation. In particular, we prove Theorem \ref{Theorem
regwaveg} showing that that many (but, may be, not all) components of
involved wavepackets remain spatially localized. Another Theorem \ref%
{Theorem Superposition gen} extends the superposition principle to the case
when resonance interactions between components of a multi-wavepackets can
occur.

\section{Conditions and definitions}

In this section we formulate and discuss all definitions and conditions
under which we study the nonlinear evolutionary system (\ref{difeqintr})
through its modal, Fourier form (\ref{difeqfou}). Most of the conditions and
definitions are naturally formulated for the modal form (\ref{difeqfou}),
and this is one of the reasons we use it as the basic one.

\subsection{Linear part}

The basic properties of the linear part $\mathbf{L}\left( \mathbf{\mathbf{k}}%
\right) $ of the system (\ref{difeqfou}), which is a $2J\times 2J$ Hermitian
matrix with eigenvalues $\omega _{n,\zeta }\left( \mathbf{k}\right) $, has
been already discussed in the Introduction. To account for all needed
properties of $\mathbf{L}\left( \mathbf{\mathbf{k}}\right) $ we define the
singular set of points $\mathbf{\mathbf{k}}$.

\begin{definition}[band-crossing points]
\label{Definition band-crossing point} We call $\mathbf{k}_{0}$ a \emph{%
band-crossing point} for $\mathbf{L}\left( \mathbf{\mathbf{k}}\right) $ if $%
\omega _{n+1,\zeta }\left( \mathbf{k}_{0}\right) =\omega _{n,\zeta }\left( 
\mathbf{k}_{0}\right) $ \ for some $n,\zeta $ or $\mathbf{L}\left( \mathbf{%
\mathbf{k}}\right) $ is not continuous at $\mathbf{k}_{0}$ or if $\omega
_{1,\pm }\left( \mathbf{k}_{0}\right) =0$. The set of such points is denoted
by $\sigma _{\mathrm{bc}}$.
\end{definition}

In the next Condition we collect all constraints imposed on the linear
operator $\mathbf{L}\left( \mathbf{\mathbf{k}}\right) $.

\begin{condition}[linear part]
\label{clin}The linear part $\mathbf{L}\left( \mathbf{\mathbf{k}}\right) $
of the system (\ref{difeqfou}) is a $2J\times 2J$ Hermitian matrix with
eigenvalues $\omega _{n,\zeta }\left( \mathbf{k}\right) $ and corresponding
eigenvectors $\mathbf{g}_{n,\zeta }\left( \mathbf{k}\right) $ satisfying for 
$\mathbf{k}\notin \sigma _{\mathrm{bc}}$ the basic relations (\ref{OmomL})-(%
\ref{invsym}). In addition to that we assume:

\begin{enumerate}
\item[(i)] the set of band-crossing points $\sigma _{\mathrm{bc}}$ is a
closed, nowhere dense set in $\mathbb{R}^{d}$ and has zero Lebesgue measure;

\item[(ii)] the entries of the Hermitian matrix $\mathbf{L}\left( \mathbf{%
\mathbf{k}}\right) $ are infinitely differentiable in $\mathbf{k}$ for all $%
\mathbf{k}\notin \sigma _{\mathrm{bc}}$ that readily implies via the
spectral theory, \cite{Kato}, infinite differentiability of all eigenvalues $%
\omega _{n}\left( \mathbf{k}\right) $ in $\mathbf{k}$ for all $\mathbf{k}%
\notin \sigma _{\mathrm{bc}}$;

\item[(iii)] $\mathbf{L}\left( \mathbf{\mathbf{k}}\right) $ satisfies a
polynomial bound 
\begin{equation}
\left\Vert \mathbf{L}\left( \mathbf{\mathbf{k}}\right) \right\Vert \leq
C\left( 1+\left\vert \mathbf{\mathbf{k}}\right\vert ^{p}\right) ,\ \mathbf{k}%
\in \mathbb{R}^{d},\ \text{for some }C>0\text{ and }p>0\text{.}  \label{Lpol}
\end{equation}
\end{enumerate}
\end{condition}

Note that since $\omega _{n,\zeta }\left( \mathbf{k}\right) $ are smooth if $%
\mathbf{k}\notin \sigma _{\mathrm{bc}}$ the following relations hold: 
\begin{equation}
\max_{\left\vert \mathbf{k\pm k}_{\ast l}\right\vert \leq \pi _{0},\
l=1,\ldots ,N,}\left\vert \nabla _{\mathbf{k}}\omega _{n_{l},\zeta
}\right\vert \leq C_{\omega ,1},\ \max_{\left\vert \mathbf{k\pm k}_{\ast
l}\right\vert \leq \pi _{0},\ l=1,\ldots ,N,}\left\vert \nabla _{\mathbf{k}%
}^{2}\omega _{n_{l},\zeta }\right\vert \leq C_{\omega ,2},  \label{Com2}
\end{equation}%
where $C_{\omega ,1}$ and $C_{\omega ,2}$ are positive constants and%
\begin{equation}
\pi _{0}=\frac{1}{2}\min_{l=1,\ldots ,N}\min \left( \limfunc{dist}\left\{
\pm \mathbf{k}_{\ast l},\sigma _{\mathrm{bc}}\right\} ,1\right) .
\label{pi0}
\end{equation}

\begin{remark}[dispersion relations symmetry]
\label{Remark symmetry}The symmetry condition (\ref{invsym}) on the
dispersion relations naturally arise in many physical problems, for example
Maxwell equations in periodic media, see \cite{BF1}-\cite{BF3}, \cite{BF5},
or when $\mathbf{L}\left( \mathbf{k}\right) $ originates from a Hamiltonian.
We would like to stress that this symmetry conditions are not imposed to
simplify studies but rather to take into account fundamental symmetries of
physical media. The symmetry\ causes resonant nonlinear interactions, which
create non-trivial effects. Interestingly, many problems without symmetries
can be put into the framework with the symmetry by a certain extension, \cite%
{BF7}.
\end{remark}

\begin{remark}[band-crossing points]
Band-crossing points are discussed in more details in \cite[Section 5.4]{BF1}%
, \cite[Sections 4.1, 4.2]{BF2}. In particular, generically the set $\sigma
_{\mathrm{bc}}$ of band-crossing point is a manifold of the dimension $d-2$.
Notice also that there is an natural ambiguity in the definition of a
normalized eigenvector $\mathbf{g}_{n,\zeta }\left( \mathbf{k}\right) $ of $%
\mathbf{L}\left( \mathbf{\mathbf{k}}\right) $ which is defined up to a
complex number $\xi $ with $\left\vert \xi \right\vert =1$. This ambiguity
may not allow an eigenvector $\mathbf{g}_{n,\zeta }\left( \mathbf{k}\right) $
which can be a locally smooth function in $\mathbf{k}$\ to be a uniquely
defined continuous function in $\mathbf{k}$ globally for all $\mathbf{k}%
\notin \sigma _{\mathrm{bc}}$ because of a possibility of branching. But,
importantly, the orthogonal projector $\Pi _{n,\zeta }\left( \mathbf{\mathbf{%
k}}\right) $ on $\mathbf{g}_{n,\zeta }\left( \mathbf{k}\right) $ as defined
by (\ref{Pin}) is uniquely defined and, consequently, infinitely
differentiable in $\mathbf{k}$ via the spectral theory, \cite{Kato}, for all 
$\mathbf{k}\notin \sigma _{\mathrm{bc}}$. Since we consider $\mathbf{\hat{U}}%
\left( \mathbf{k}\right) $ as an element of the space $L^{1}$ and $\sigma _{%
\mathrm{bc}}$ is of zero Lebesgue measure considering $\mathbf{k}\notin
\sigma _{\mathrm{bc}}$ is sufficient for us.
\end{remark}

We introduce for vectors $\mathbf{\hat{u}}\in \mathbb{C}^{2J}$ their
expansion with respect to the orthonormal basis$\left\{ \mathbf{g}_{n,\zeta
}\left( \mathbf{k}\right) \right\} $: 
\begin{equation}
\mathbf{\hat{u}}\left( \mathbf{k}\right) =\sum_{n=1}^{J}\sum_{\zeta =\pm }%
\hat{u}_{n,\zeta }\left( \mathbf{k}\right) \mathbf{g}_{n,\zeta }\left( 
\mathbf{k}\right) =\sum_{n=1}^{J}\sum_{\zeta =\pm }\mathbf{\hat{u}}_{n,\zeta
}\left( \mathbf{k}\right) ,\ \mathbf{\hat{u}}_{n,\zeta }\left( \mathbf{k}%
\right) =\Pi _{n,\zeta }\left( \mathbf{\mathbf{k}}\right) \mathbf{\hat{u}}%
\left( \mathbf{k}\right)  \label{Uboldj}
\end{equation}%
and we refer to it as the \emph{modal decomposition} of $\mathbf{\hat{u}}%
\left( \mathbf{k}\right) $ and to $\hat{u}_{n,\zeta }\left( \mathbf{k}%
\right) $ as the \emph{modal coefficients} of $\mathbf{\hat{u}}\left( 
\mathbf{k}\right) $. Evidently 
\begin{equation}
\sum\nolimits_{n=1}^{j}\sum\nolimits_{\zeta =\pm }\Pi _{n,\zeta }\left( 
\mathbf{\mathbf{k}}\right) =I_{2J},\text{ where }I_{2J}\text{ is the }%
2J\times 2J\text{ identity matrix.}  \label{sumPi}
\end{equation}%
Notice that we can define the action of the operator $\mathbf{L}\left( -%
\mathrm{i}\nabla _{\mathbf{r}}\right) $ on any Schwartz function $\mathbf{Y}%
\left( \mathbf{r}\right) $ by the formula%
\begin{equation}
\widehat{\mathbf{L}\left( -\mathrm{i}\nabla _{\mathbf{r}}\right) \mathbf{Y}}%
\left( \mathbf{\mathbf{k}}\right) =\mathbf{L}\left( \mathbf{\mathbf{k}}%
\right) \mathbf{\hat{Y}}\left( \mathbf{\mathbf{k}}\right) ,\text{ }
\label{Ldiff}
\end{equation}%
where in view of the polynomial bound (\ref{Lpol}) the order of $\mathbf{L}$
does not exceed $p.$ In a special case when all the entries of $\mathbf{L}%
\left( \mathbf{\mathbf{k}}\right) $ are polynomials (\ref{Ldiff}) turns into
the action of the differential operator with constant coefficients.

\subsection{Nonlinear part}

The nonlinear term $\hat{F}$ in (\ref{difeqfou}) is assumed to be a general
functional polynomial of the form 
\begin{gather}
\hat{F}\left( \mathbf{\hat{U}}\right) =\sum\nolimits_{m\in \mathfrak{M}_{F}}%
\hat{F}^{\left( m\right) }\left( \mathbf{\hat{U}}^{m}\right) ,\text{ where }%
\hat{F}^{\left( m\right) }\text{ is }m\text{-homogeneous polylinear operator,%
}  \label{Fseries} \\
\mathfrak{M}_{F}=\left\{ m_{1},\ldots ,m_{p}\right\} \subset \left\{
2,3,\ldots \right\} \text{ is a finite set, and }m_{F}=\max \left\{ m:m\in 
\mathfrak{M}_{F}\right\} .  \label{Fseries1}
\end{gather}%
The integer $m_{F}$ in (\ref{Fseries1}) is called the \emph{degree of the
functional polynomial} $\hat{F}$. For instance, if $\mathfrak{M}_{F}=\left\{
2\right\} $ or $\mathfrak{M}_{F}=\left\{ 3\right\} $ the polynomial $\hat{F}%
\,$\ is respectively homogeneous quadratic or cubic. Every $m$-linear
operator $\hat{F}^{\left( m\right) }$ in (\ref{Fseries}) is assumed to be of
the form of a convolution%
\begin{gather}
\hat{F}^{\left( m\right) }\left( \mathbf{\hat{U}}_{1},\ldots ,\mathbf{\hat{U}%
}_{m}\right) \left( \mathbf{k},\tau \right) =\int_{\mathbb{D}_{m}}\chi
^{\left( m\right) }\left( \mathbf{\mathbf{k}},\vec{k}\right) \mathbf{\hat{U}}%
_{1}\left( \mathbf{k}^{\prime }\right) \ldots \mathbf{\hat{U}}_{m}\left( 
\mathbf{k}^{\left( m\right) }\left( \mathbf{k},\vec{k}\right) \right) \,%
\mathrm{\tilde{d}}^{\left( m-1\right) d}\vec{k},  \label{Fmintr} \\
\text{where }\mathbb{D}_{m}=\mathbb{R}^{\left( m-1\right) d},\ \mathrm{%
\tilde{d}}^{\left( m-1\right) d}\vec{k}=\frac{\mathrm{d}\mathbf{k}^{\prime
}\ldots \,\mathrm{d}\mathbf{k}^{\left( m-1\right) }}{\left( 2\pi \right)
^{\left( m-1\right) d}},  \notag \\
\mathbf{k}^{\left( m\right) }\left( \mathbf{k},\vec{k}\right) =\mathbf{k}-%
\mathbf{k}^{\prime }-\ldots -\mathbf{k}^{\left( m-1\right) },\ \vec{k}%
=\left( \mathbf{k}^{\prime },\ldots ,\mathbf{k}^{\left( m\right) }\right) .
\label{conv}
\end{gather}%
indicating that the nonlinear operator $F^{\left( m\right) }\left( \mathbf{U}%
_{1},\ldots ,\mathbf{U}_{m}\right) $ is translation invariant (it may be
local or non-local). The quantities $\chi ^{\left( m\right) }$ in (\ref%
{Fmintr}) are called \emph{susceptibilities}. For numerous examples of
nonlinearities of the form similar to (\ref{Fseries}), (\ref{Fmintr}) see 
\cite{BF1}-\cite{BF7} and references therein. In what follows the nonlinear
term $\hat{F}$ in (\ref{difeqfou}) will satisfy the following conditions.

\begin{condition}[nonlinearity]
\label{cnonbound}The nonlinearity $\hat{F}\left( \mathbf{\hat{U}}\right) $
is assumed to be of the form (\ref{Fseries})-(\ref{Fmintr}). The
susceptibility $\chi ^{\left( m\right) }\left( \mathbf{\mathbf{k}},\mathbf{k}%
^{\prime },\ldots ,\mathbf{k}^{\left( m\right) }\right) $ is infinitely
differentiable for all $\mathbf{\mathbf{k}}$ and $\mathbf{k}^{\left(
j\right) }$ which are not band-crossing points, and is bounded, namely for
some constant $C_{\chi }$ 
\begin{equation}
\left\Vert \chi ^{\left( m\right) }\right\Vert =\left( 2\pi \right)
^{-\left( m-1\right) d}\sup_{\mathbf{\mathbf{k}},\mathbf{k}^{\prime },\ldots
,\mathbf{k}^{\left( m\right) }\in \mathbb{R}^{d}\setminus \sigma _{\mathrm{bc%
}}}\left\vert \chi ^{\left( m\right) }\left( \mathbf{\mathbf{k}},\mathbf{k}%
^{\prime },\ldots ,\mathbf{k}^{\left( m\right) }\right) \right\vert \leq
C_{\chi },\ m\in \mathfrak{M}_{F},  \label{chiCR}
\end{equation}%
where the norm $\left\vert \chi ^{\left( m\right) }\left( \mathbf{k},\vec{k}%
\right) \right\vert $ of $m$-linear tensor $\chi ^{\left( m\right) }:\left( 
\mathbb{C}^{2J}\right) ^{m}\rightarrow \left( \mathbb{C}^{2J}\right) ^{m}$
for fixed $\mathbf{k},\vec{k}$ is defined by 
\begin{equation}
\left\vert \chi ^{\left( m\right) }\left( \mathbf{k},\vec{k}\right)
\right\vert =\sup_{\left\vert \mathbf{x}_{j}\right\vert \leq 1}\left\vert
\chi _{\ }^{\left( m\right) }\left( \mathbf{k},\vec{k}\right) \left( \mathbf{%
x}_{1},\ldots ,\mathbf{x}_{m}\right) \right\vert ,\text{ where }\left\vert 
\mathbf{x}\right\vert \text{ is the Euclidean norm}.  \label{normchi0}
\end{equation}
\end{condition}

Since $\chi _{\zeta ,\vec{\zeta}}^{\left( m\right) }\left( \mathbf{\mathbf{k}%
},\mathbf{k}^{\prime },\ldots ,\mathbf{k}^{\left( m\right) }\right) \ $ are
smooth if $\mathbf{k}\notin \sigma _{\mathrm{bc}}$ the following relation
holds: 
\begin{equation}
\max_{\left\vert \mathbf{k\pm k}_{\ast l}\right\vert \leq \pi _{0},\
l=1,\ldots ,N}\left\vert \nabla \chi _{\zeta ,\vec{\zeta}}^{\left( m\right)
}\left( \mathbf{\mathbf{k}},\mathbf{k}^{\prime },\ldots ,\mathbf{k}^{\left(
m\right) }\right) \right\vert \leq C_{\chi }^{\prime }\text{ }  \label{grchi}
\end{equation}%
if $\mathbf{k}_{\ast l}\notin \sigma _{\mathrm{bc}}$, $\pi _{0}$ is defined
by (\ref{pi0}). The case when $\chi ^{\left( m\right) }\left( \mathbf{k},%
\vec{k}\right) $ depend on small $\varrho $ or, more generally, on $\varrho
^{q}$, $q>0$, \ can be treated similarly, see \cite{BF7}.

\subsection{Resonance invariant $nk$-spectrum}

In this section, being given the dispersion relations $\omega _{n}\left( 
\mathbf{k}\right) \geq 0$, $n\in \left\{ 1,\ldots ,J\right\} $, \ we
consider resonance properties of $nk$-spectra $S$ and the corresponding $k$%
-spectra $K_{S}$ as defined in Definition \ref{dmwavepack}, i.e. 
\begin{equation}
S=\left\{ \left( n_{l},\mathbf{k}_{\ast l}\right) ,\ l=1,\ldots ,N\right\}
\subset \Sigma =\left\{ 1,\ldots ,J\right\} \times \mathbb{R}^{d},\ K_{S}%
\text{ }=\left\{ \mathbf{k}_{\ast _{l}},\ l=1,\ldots ,\left\vert
K_{S}\right\vert \right\} .  \label{ssp1}
\end{equation}%
We precede the formal description of the \emph{resonance invariance} (see
Definition \ref{Definition omclos}) with the following guiding physical
picture. Initially at $\tau =0$ the wave is a multi-wavepacket composed of
modes from a small vicinity of the $nk$-spectrum $S$. As the wave evolves
according to (\ref{difeqfou}) the polynomial nonlinearity inevitably
involves \ a larger set of modes $\left[ S\right] _{\text{out}}\supseteq S$,
but not all modes in $\left[ S\right] _{\text{out}}$\ are "equal" in
developing significant amplitudes. The qualitative picture is that whenever
certain interaction phase function (see (\ref{phim}) below) is not zero, the
fast time oscillations weaken effective nonlinear mode interaction and the
energy transfer from the original modes in $S$ to relevant modes from $\left[
S\right] _{\text{out}}$, keeping their magnitudes vanishingly small as $%
\beta ,\varrho \rightarrow 0$. There is a smaller set of modes $\left[ S%
\right] _{\text{out}}^{\text{res}}$ which can interact with modes from $S$
rather effectively and develop significant amplitudes. Now, 
\begin{equation}
\text{if }\left[ S\right] _{\text{out}}^{\text{res}}\subseteq S\text{ then }S%
\text{ is called resonance invariant.}  \label{ssp2}
\end{equation}%
In simpler situations the resonance invariance conditions turns into the
well-known in nonlinear optics phase and frequency matching conditions. For
instance, if $S$ contains $\left( n_{0},\mathbf{k}_{\ast l_{0}}\right) $ and
the dispersion relations allow for the second harmonic generation in another
band $n_{1}$ so that $2\omega _{n_{0}}\left( \mathbf{k}_{\ast l_{0}}\right)
=\omega _{n_{1}}\left( 2\mathbf{k}_{\ast l_{0}}\right) $, then for $S$ to be
resonance invariant it must contain $\left( n_{1},2\mathbf{k}_{\ast
l_{0}}\right) $ too.

Let us turn now to the rigorous constructions. First we introduce necessary
notations. Let $m\geq 2$ be an integer, $\vec{l}=\left(
l_{1},..,l_{m}\right) $,\ $l_{j}\in \left\{ 1,\ldots ,N\right\} $ be an
integer vector from $\left\{ 1,\ldots ,N\right\} ^{m}$ and $\vec{\zeta}%
=\left( \zeta ^{\left( 1\right) },,..,\zeta ^{\left( m\right) }\right) $, $%
\zeta ^{\left( j\right) }\in \left\{ +1,-1\right\} $ be a binary vector from 
$\left\{ +1,-1\right\} ^{m}$. \emph{\ }Note that a pair $\left( \vec{\zeta},%
\vec{l}\right) $ naturally labels a sample string of the length $m$ composed
of elements $\left( \zeta ^{\left( j\right) },n_{l_{j}},\mathbf{k}_{\ast
l_{j}}\right) $ from the set $\left\{ +1,-1\right\} \times S$. Let us
introduce the sets 
\begin{gather}
\Lambda =\left\{ \left( \zeta ,l\right) :l\in \left\{ 1,\ldots ,N\right\} ,\
\zeta \in \left\{ +1,-1\right\} \right\} ,  \label{setlam} \\
\Lambda ^{m}=\left\{ \vec{\lambda}=\left( \lambda _{1},\ldots ,\lambda
_{m}\right) ,\ \lambda _{j}\in \Lambda ,\ j=1,\ldots ,m\right\} .  \notag
\end{gather}%
There is a natural one-to-one correspondence between $\Lambda ^{m}$ and $%
\left\{ -1,1\right\} ^{m}\times \left\{ 1,\ldots ,N\right\} ^{m}$ \ and we
will write, exploiting this correspondence 
\begin{equation}
\vec{\lambda}=\left( \left( \zeta ^{\prime },l_{1}\right) ,\ldots ,\left(
\zeta ^{\left( m\right) },l_{m}\right) \right) =\left( \vec{\zeta},\vec{l}%
\right) ,\ \vec{\vartheta}\in \left\{ -1,1\right\} ^{m},\ \vec{l}\in \left\{
1,\ldots ,N\right\} ^{m}\text{ for }\vec{\lambda}\in \Lambda ^{m}.
\label{lamprop}
\end{equation}%
Let us introduce the following linear combination \ 
\begin{equation}
\varkappa _{m}\left( \vec{\lambda}\right) =\varkappa _{m}\left( \vec{\zeta},%
\vec{l}\right) =\sum\nolimits_{j=1}^{m}\zeta ^{\left( j\right) }\mathbf{k}%
_{\ast l_{j}}\text{ with }\zeta ^{\left( j\right) }\in \left\{ +1,-1\right\}
,  \label{kapzel}
\end{equation}%
and let $\left[ S\right] _{K,\text{out}}$ be the set of all its values as $%
\mathbf{k}_{\ast l_{j}}\in K_{S}$, $\vec{\lambda}\in \Lambda ^{m}$, namely 
\begin{equation}
\left[ S\right] _{K,\text{out}}=\dbigcup\nolimits_{m\in \mathfrak{M}%
_{F}}\dbigcup\nolimits_{\vec{\lambda}\in \Lambda ^{m}}\left\{ \varkappa
_{m}\left( \vec{\lambda}\right) \right\} .  \label{K1}
\end{equation}%
\ We call $\left[ S\right] _{K,\text{out}}$ \emph{output }$k$\emph{-spectrum}
of $K_{S}$ and assume that%
\begin{equation}
\left[ S\right] _{K,\text{out}}\dbigcap \sigma _{\mathrm{bc}}=\varnothing .
\label{Skempty}
\end{equation}%
We also define the \emph{output }$nk$-\emph{spectrum} of $S$ by 
\begin{equation}
\left[ S\right] _{\text{out}}=\left\{ \left( n,\mathbf{k}\right) \in \left\{
1,\ldots ,J\right\} \times \mathbb{R}^{d}:n\in \left\{ 1,\ldots ,J\right\}
,\ \mathbf{k}\in \left[ S\right] _{K,\text{out}}\right\} .  \label{Sout}
\end{equation}%
We introduce the following functions%
\begin{equation}
\Omega _{1,m}\left( \vec{\lambda}\right) \left( \vec{k}_{\ast }\right)
=\sum\nolimits_{j=1}^{m}\zeta ^{\left( j\right) }\omega _{l_{j}}\left( 
\mathbf{k}_{\ast l_{j}}\right) ,\ \vec{k}_{\ast }=\left( \mathbf{k}_{\ast
1},\ldots ,\mathbf{k}_{\ast \left\vert K_{S}\right\vert }\right) ,\ \text{%
where }\mathbf{k}_{\ast l_{j}}\in K_{S},  \label{Om1}
\end{equation}%
\begin{equation}
\Omega \left( \zeta ,n,\vec{\lambda}\right) \left( \mathbf{k}_{\ast \ast },%
\vec{k}_{\ast }\right) =-\zeta \omega _{n}\left( \mathbf{k}_{\ast \ast
}\right) +\Omega _{1,m}\left( \vec{\lambda}\right) \left( \vec{k}_{\ast
}\right) ,  \label{Omzet}
\end{equation}%
where $\zeta =\pm 1$,$\ m\in \mathfrak{M}_{F}\ $as in (\ref{Fseries}). We
introduce these functions to apply later to phase functions (\ref{phim}).

Now we introduce the \emph{resonance equation} 
\begin{equation}
\Omega \left( \zeta ,n,\vec{\lambda}\right) \left( \zeta \varkappa
_{m}\left( \vec{\lambda}\right) ,\vec{k}_{\ast }\right) =0,\ \vec{l}\in
\left\{ 1,\ldots ,N\right\} ^{m},\ \vec{\zeta}\in \left\{ -1,1\right\} ^{m},
\label{Omeq0}
\end{equation}%
denoting by $P\left( S\right) $ the set of its solutions $\left( m,\zeta ,n,%
\vec{\lambda}\right) $. Such a solution is called $S$\emph{-internal }if%
\begin{equation}
\left( n,\zeta \varkappa _{m}\left( \vec{\lambda}\right) \right) \in S,\text{
that is }n=n_{l_{0}},\ \zeta \varkappa _{m}\left( \vec{\lambda}\right) =%
\mathbf{k}_{\ast l_{0}},\ l_{0}\in \left\{ 1,\ldots ,N\right\} ,
\label{Sint}
\end{equation}%
and we denote the corresponding $l_{0}=I\left( \vec{\lambda}\right) $. We
also denote by $P_{\text{int}}\left( S\right) \subset P\left( S\right) $ the
set of all $S$-internal solutions to (\ref{Omeq0}).

Now we consider the simplest solutions to (\ref{Omeq0}) which play an
important role. Keeping in mind that the string $\vec{l}$ can contain
several copies of a single value $l$, we can recast the sum in (\ref{Om1})
as follows: 
\begin{gather}
\Omega _{1,m}\left( \vec{\lambda}\right) =\Omega _{1,m}\left( \vec{\zeta},%
\vec{l}\right) =\sum\nolimits_{l=1}^{N}\delta _{l}\omega _{l}\left( \mathbf{k%
}_{\ast l}\right) ,\ \text{where }\delta _{l}=\left\{ 
\begin{array}{ccc}
\sum_{j\in \vec{l}^{-1}\left( l\right) }\zeta ^{\left( j\right) } & \text{if}
& \vec{l}^{-1}\left( l\right) \neq \varnothing \\ 
0 & \text{if} & \vec{l}^{-1}\left( l\right) =\varnothing%
\end{array}%
\right. ,  \label{rearr0} \\
\vec{l}^{-1}\left( l\right) =\left\{ j\in \left\{ 1,\ldots ,m\right\}
:l_{j}=l,\right\} ,\ \vec{l}=\left( l_{1},\ldots ,l_{m}\right) ,\ 1\leq
l\leq N.  \notag
\end{gather}

\begin{definition}[universal solutions]
\label{Universal solution}We call a solution $\left( m,\zeta ,n,\vec{\lambda}%
\right) \in P\left( S\right) $ of (\ref{Omeq0}) \emph{universal} if it has
the following properties: (i) only a single coefficient out of all $\delta
_{l}$ in (\ref{rearr0}) is nonzero, namely for some $I_{0}$ we have $\delta
_{I_{0}}=\pm 1$ and $\delta _{l}=0$ for $l\neq I_{0}$; (ii) $n=n_{I_{0}}$ \
and $\zeta =\delta _{I_{0}}$.
\end{definition}

We denote the set of universal solutions to (\ref{Omeq0}) by $P_{\text{univ}%
}\left( S\right) $. \emph{A justification for calling such solution
universal comes from the fact that if a solution is a universal solution
solution for one }$\vec{k}_{\ast }$ \emph{it is a solution for any other }$%
\vec{k}_{\ast }\in \mathbb{R}^{d}$\emph{.} Note that a universal solution is
a $S$-internal solution with $I\left( \vec{\lambda}\right) =I_{0}$ implying 
\begin{equation}
P_{\text{univ}}\left( S\right) \subseteq P_{\text{int}}\left( S\right) .
\label{gint}
\end{equation}%
Indeed, observe that for $\delta _{l}$ as in (\ref{rearr0} 
\begin{equation}
\varkappa _{m}\left( \vec{\lambda}\right) =\varkappa _{m}\left( \vec{\zeta},%
\vec{l}\right) =\sum\nolimits_{j=1}^{m}\zeta ^{\left( j\right) }\mathbf{k}%
_{\ast l_{j}}=\sum\nolimits_{l=1}^{N}\delta _{l}\mathbf{k}_{\ast l}
\label{rearr1}
\end{equation}%
implying $\varkappa _{m}\left( \vec{\lambda}\right) =\delta _{I_{0}}\mathbf{k%
}_{\ast I_{0}}$ and$\ \ \zeta \varkappa _{m}\left( \vec{\lambda}\right)
=\delta _{I_{0}}^{2}\mathbf{k}_{\ast I_{0}}=\mathbf{k}_{\ast I_{0}}$. Then
equation (\ref{Omeq0}) \ is obviously satisfied and $\left( n,\zeta
\varkappa _{m}\left( \vec{\lambda}\right) \right) =\left( n_{I_{0}},\mathbf{k%
}_{\ast I_{0}}\right) \in S$.

\begin{example}[universal solutions]
\label{Example counterprop}Suppose there is just a single band, i.e. $J=1$,
a symmetric dispersion relation $\omega _{1}\left( -\mathbf{k}\right)
=\omega _{1}\left( \mathbf{k}\right) $, a cubic nonlinearity $F$ with $%
\mathfrak{M}_{F}=\left\{ 3\right\} $. We take the $nk$-spectrum\emph{\ }$%
S=\left\{ \left( 1,\mathbf{k}_{\ast }\right) ,\left( 1,-\mathbf{k}_{\ast
}\right) \right\} $, that is $N=2$ and $\mathbf{k}_{\ast 1}=\mathbf{k}_{\ast
},\mathbf{k}_{\ast 2}=-\mathbf{k}_{\ast }$. This example is typical for two
counterpropagating waves. Then $\Omega _{1,3}\left( \vec{\lambda}\right)
\left( \vec{k}_{\ast }\right) =\sum\nolimits_{j=1}^{3}\zeta ^{\left(
j\right) }\omega _{l_{j}}\left( \mathbf{k}_{\ast l_{j}}\right) =\left(
\delta _{1}+\delta _{2}\right) \omega _{1}\left( \mathbf{k}_{\ast }\right) $
\ and $\varkappa _{m}\left( \vec{\lambda}\right)
=\sum\nolimits_{j=1}^{m}\zeta ^{\left( j\right) }\mathbf{k}_{\ast
l_{j}}=\delta _{1}\mathbf{k}_{\ast 1}+\delta _{2}\mathbf{k}_{\ast 2}=\left(
\delta _{1}-\delta _{2}\right) \mathbf{k}_{\ast }$ where we use notation (%
\ref{rearr0}). The universal solution set has the form $P_{\text{univ}%
}\left( S\right) =\left\{ \left( 3,\zeta ,1,\vec{\lambda}\right) :\vec{%
\lambda}\in \Lambda _{\zeta },\zeta =\pm \right\} $ where $\Lambda _{+}$
consists of vectors $\left( \lambda _{1},\lambda _{2},\lambda _{3}\right) $ $%
\ $of the form $\left( \left( +,1\right) ,\left( -,1\right) ,\left(
+,1\right) \right) ,$\ $\left( \left( +,1\right) ,\left( -,1\right) ,\left(
+,2\right) \right) ,\ \left( \left( +,2\right) ,\left( -,2\right) ,\left(
+,1\right) \right) $, $\left( \left( +,2\right) ,\left( -,2\right) ,\left(
+,2\right) \right) $, and vectors obtained from the listed ones by
permutations of coordinates $\lambda _{1},\lambda _{2},\lambda _{3}$. The
solutions from $P_{\text{int}}\left( S\right) $ have to satisfy $\left\vert
\delta _{1}-\delta _{2}\right\vert =1$ and $\left\vert \delta _{1}+\delta
_{2}\right\vert =1$ which is possible only if $\delta _{1}\delta _{2}=0$.
Since $\zeta =\delta _{1}+\delta _{2}$ we have $\zeta \varkappa _{m}\left( 
\vec{\lambda}\right) =\left( \delta _{1}^{2}-\delta _{2}^{2}\right) \mathbf{k%
}_{\ast }$ and $\zeta \varkappa _{m}\left( \vec{\lambda}\right) =\mathbf{k}%
_{\ast 1}$ if $\left\vert \delta _{1}\right\vert =1$ or $\zeta \varkappa
_{m}\left( \vec{\lambda}\right) =\mathbf{k}_{\ast 2}$ if $\left\vert \delta
_{2}\right\vert =1$. Hence $P_{\text{int}}\left( S\right) =P_{\text{univ}%
}\left( S\right) $ in this case. Note that if we set $S_{1}=\left\{ \left( 1,%
\mathbf{k}_{\ast }\right) \right\} $, $S_{2}=\left\{ \left( 1,-\mathbf{k}%
_{\ast }\right) \right\} $ then $S=S_{1}\cup S_{2}$ but $P_{\text{int}%
}\left( S\right) $ is larger than $P_{\text{int}}\left( S_{1}\right) \cup P_{%
\text{int}}\left( S_{2}\right) $. This can be interpreted as follows. When
only \ modes from $S_{1}$ are excited, the modes from $S_{2}$ remain
non-excited. But when the both $S_{1}$ and $S_{2}$ are excited, there is a
resonance effect of $S_{1}$ onto $S_{2}$, represented, for example, by $\vec{%
\lambda}=\left( \left( +,1\right) ,\left( -,1\right) ,\left( +,2\right)
\right) $, which involves the mode $\zeta \varkappa _{m}\left( \vec{\lambda}%
\right) =\mathbf{k}_{\ast 2}$.
\end{example}

Now we are ready to define resonance invariant spectra. First, we introduce
a subset $\left[ S\right] _{\text{out}}^{\text{res}}$ of $\left[ S\right] _{%
\text{out}}$ by the formula%
\begin{eqnarray}
\left[ S\right] _{\text{out}}^{\text{res}} &=&\left\{ \left( n,\mathbf{k}%
_{\ast \ast }\right) \in \left[ S\right] _{\text{out}}:\mathbf{k}_{\ast \ast
}=\zeta ^{\left( 0\right) }\varkappa _{m}\left( \vec{\lambda}\right) ,\ m\in 
\mathfrak{M}_{F},\text{ where}\right.  \label{Sresout} \\
&&\left. \left( m,\zeta ,n,\vec{\lambda}\right) \text{ is a solution of (\ref%
{Omeq0})}\right\} ,  \notag
\end{eqnarray}%
\ calling it \emph{resonant output spectrum} of $S$, and then we define%
\begin{equation}
\text{resonance selection operation }\mathcal{R}\left( S\right) =S\cup \left[
S\right] _{\text{out}}^{\text{res}}.  \label{hull}
\end{equation}

\begin{definition}[resonance invariant $nk$-spectrum]
\label{Definition omclos} The $nk$-spectrum $S$ is called \emph{resonance
invariant} if $\mathcal{R}\left( S\right) =S$ or, equivalently, $\left[ S%
\right] _{\text{out}}^{\text{res}}\subseteq S$. The $nk$-spectrum\ $S$ is
called \emph{universally} \emph{resonance invariant }if $\mathcal{R}\left(
S\right) =S$ and $P_{\text{univ}}\left( S\right) =P_{\text{int}}\left(
S\right) $.
\end{definition}

Obviously, $nk$-spectrum $S$ is resonance invariant \ if and only if all
solutions of (\ref{Omeq0}) are internal, that is $P_{\text{int}}\left(
S\right) =P\left( S\right) $.

It is worth noticing that even when a $nk$-spectrum is not resonance
invariant often it can be easily extended to a resonance invariant one.
Namely, if $\mathcal{R}^{j}\left( S\right) \cap \sigma _{\mathrm{bc}%
}=\varnothing $ for all $j$ then the set 
\begin{equation*}
\mathcal{R}^{\infty }\left( S\right) =\dbigcup\nolimits_{j=1}^{\infty }%
\mathcal{R}^{j}\left( S\right) \subset \Sigma =\left\{ 1,\ldots ,J\right\}
\times \mathbb{R}^{d}\text{ }
\end{equation*}%
is resonance invariant. In addition to that, $\mathcal{R}^{\infty }\left(
S\right) $ is always at most countable. Usually it is finite i.e. $\mathcal{R%
}^{\infty }\left( S\right) =\mathcal{R}^{p}\left( S\right) $ for a finite $p$%
, see examples below; also $\mathcal{R}^{\infty }\left( S\right) =S$ \ for
generic $K_{S}$.

\begin{example}[resonance invariant $nk$-spectra \ for quadratic nonlinearity%
]
\label{Example quad} Suppose there is a single band, i.e. $J=1$, with a
symmetric dispersion relation, and a quadratic nonlinearity $F$, that is $%
\mathfrak{M}_{F}=\left\{ 2\right\} $. Let us assume that $\mathbf{k}_{\ast
}\neq 0$, $\mathbf{k}_{\ast },2\mathbf{k}_{\ast },\mathbf{0}$ are not
band-crossing points and look at two examples. First, suppose that $2\omega
_{1}\left( \mathbf{k}_{\ast }\right) \neq \omega _{1}\left( 2\mathbf{k}%
_{\ast }\right) $ (no second harmonic generation) and $\omega _{1}\left( 
\mathbf{0}\right) \neq 0$. Let us set the $nk$-spectrum to be the set $%
S_{1}=\left\{ \left( 1,\mathbf{k}_{\ast }\right) \right\} $, then $S_{1}$ is
resonance invariant. Indeed, $K_{S_{1}}=\left\{ \mathbf{k}_{\ast }\right\} $%
, $\left[ S_{1}\right] _{K,\text{out}}=\left\{ \mathbf{0},2\mathbf{k}_{\ast
},-2\mathbf{k}_{\ast }\right\} $, $\left[ S_{1}\right] _{\text{out}}=\left\{
\left( 1,\mathbf{0}\right) ,\left( 1,2\mathbf{k}_{\ast }\right) ,\left( 1,-2%
\mathbf{k}_{\ast }\right) \right\} $ and an elementary examination shows
that $\left[ S_{1}\right] _{\text{out}}^{\text{res}}=\varnothing \subset
S_{1}$ implying $\mathcal{R}\left( S_{1}\right) =S_{1}$. For the second
example let us assume $\omega _{1}\left( \mathbf{0}\right) \neq 0$ and $%
2\omega _{1}\left( \mathbf{k}_{\ast }\right) =\omega _{1}\left( 2\mathbf{k}%
_{\ast }\right) $, that is the second harmonic generation is present. Here $%
\left[ S_{1}\right] _{\text{out}}^{\text{res}}=\left\{ \left( 1,2\mathbf{k}%
_{\ast }\right) \right\} $ and $\mathcal{R}\left( S_{1}\right) =\left\{
\left( 1,\mathbf{k}_{\ast }\right) ,\left( 1,2\mathbf{k}_{\ast }\right)
\right\} $ implying $\mathcal{R}\left( S_{1}\right) \neq S_{1}$ and, hence, $%
S_{1}$ is not resonance invariant. Suppose now that $4\mathbf{k}_{\ast },3%
\mathbf{k}_{\ast }\notin \sigma _{\mathrm{bc}}$ and \ $\omega _{1}\left( 
\mathbf{0}\right) \neq 0$, $\omega _{1}\left( 4\mathbf{k}_{\ast }\right)
\neq 2\omega _{1}\left( 2\mathbf{k}_{\ast }\right) $, $\omega _{1}\left( 3%
\mathbf{k}_{\ast }\right) \neq \omega _{1}\left( \mathbf{k}_{\ast }\right)
+\omega _{1}\left( 2\mathbf{k}_{\ast }\right) $ and let us set $%
S_{2}=\left\{ \left( 1,\mathbf{k}_{\ast }\right) ,\left( 1,2\mathbf{k}_{\ast
}\right) \right\} $. An elementary examination shows that $S_{2}$ is
resonance invariant. Note that $S_{2}$ can be obtained by iterating the
resonance selection operator, namely $S_{2}=\mathcal{R}\left( \mathcal{R}%
\left( S_{1}\right) \right) $. Note also that $P_{\text{univ}}\left(
S_{2}\right) \neq P_{\text{int}}\left( S_{2}\right) $. Notice that $\omega
_{1}\left( \mathbf{0}\right) =0$ is a special case since $\mathbf{k}=\mathbf{%
0}$ is a band-crossing point, and it requires a special treatment.
\end{example}

\begin{example}[resonance invariant $nk$-spectra for cubic nonlinearity]
\label{Example cube}Let us consider one-band case with symmetric dispersion
relation and a cubic \ nonlinearity that is $\mathfrak{M}_{F}=\left\{
3\right\} $. First we take $S_{1}=\left\{ \left( 1,\mathbf{k}_{\ast }\right)
\right\} $, we assume that $\mathbf{k}_{\ast },3\mathbf{k}_{\ast }$\ are not
band-crossing points, implying $\left[ S_{1}\right] _{K,\text{out}}=\left\{ 
\mathbf{k}_{\ast },-\mathbf{k}_{\ast },3\mathbf{k}_{\ast },-3\mathbf{k}%
_{\ast }\right\} $. We have $\Omega _{1,3}\left( \vec{\lambda}\right) \left( 
\vec{k}_{\ast }\right) =\sum\nolimits_{j=1}^{3}\zeta ^{\left( j\right)
}\omega _{1}\left( \mathbf{k}_{\ast }\right) =\delta _{1}\omega _{1}\left( 
\mathbf{k}_{\ast }\right) $ \ and $\varkappa _{m}\left( \vec{\lambda}\right)
=\delta _{1}\mathbf{k}_{\ast }$ where we use notation (\ref{rearr0}), $%
\delta _{1}$ takes values $1,-1,3,-3$. If $3\omega _{1}\left( \mathbf{k}%
_{\ast }\right) \neq \omega _{1}\left( 3\mathbf{k}_{\ast }\right) $ then (%
\ref{Omeq0}) has a solution only if $\ \left\vert \delta _{1}\right\vert =1$
and $\delta _{1}=\zeta $, hence $\zeta \varkappa _{m}\left( \vec{\lambda}%
\right) =\mathbf{k}_{\ast }$ and every solution is internal. \ \ Hence, $%
\left[ S_{1}\right] _{\text{out}}^{\text{res}}=\varnothing $ and $\mathcal{R}%
\left( S_{1}\right) =S_{1}$.\ Now consider the case associated with the
third harmonic generation, namely $3\omega _{1}\left( \mathbf{k}_{\ast
}\right) =\omega _{1}\left( 3\mathbf{k}_{\ast }\right) $ and assume that $%
\omega _{1}\left( 3\mathbf{k}_{\ast }\right) +2\omega _{1}\left( \mathbf{k}%
_{\ast }\right) \neq \omega _{1}\left( 5\mathbf{k}_{\ast }\right) $, $%
3\omega _{1}\left( 3\mathbf{k}_{\ast }\right) \neq \omega _{1}\left( 9%
\mathbf{k}_{\ast }\right) $, $2\omega _{1}\left( 3\mathbf{k}_{\ast }\right)
+\omega _{1}\left( \mathbf{k}_{\ast }\right) \neq \omega _{1}\left( 7\mathbf{%
k}_{\ast }\right) $, $2\omega _{1}\left( 3\mathbf{k}_{\ast }\right) -\omega
_{1}\left( \mathbf{k}_{\ast }\right) \neq \omega _{1}\left( 5\mathbf{k}%
_{\ast }\right) $. An\ elementary examination shows that the set $%
S_{4}=\left\{ \left( 1,3\mathbf{k}_{\ast }\right) ,\left( 1,\mathbf{k}_{\ast
}\right) ,\left( 1,-\mathbf{k}_{\ast }\right) \left( 1,-3\mathbf{k}_{\ast
}\right) \right\} $ satisfies $\mathcal{R}\left( S_{4}\right) =S_{4}$.
Consequently, a multiwavepacket having $S_{4}$ as its resonance invariant $%
nk $-spectrum involves the third harmonic generation and, according to
Theorem \ref{Theorem invarwave}, it is preserved under nonlinear evolution.
The above examples indicate that in simple cases the conditions on $\ 
\mathbf{k}_{\ast }$ which can make $S$ non-invariant with respect to $%
\mathcal{R}$ have a form of several algebraic equations, Hence, for almost
all $\mathbf{k}_{\ast }$ such spectra $S$ are resonance invariant. The
examples also show that if we fix $S$\ and dispersion relations then we can
include $S$ in larger spectrum $S^{\prime }=\mathcal{R}^{p}\left( S\right) $
\ using repeated application of the operation $\mathcal{R}$ to $S$, and
often the resulting extended $nk$-spectrum $S^{\prime }$ is resonance
invariant. We show in the following section that $nk$-spectrum $S$ \ with
generic $K_{S}$ is universally resonance invariant.
\end{example}

Note that the concept of resonance invariant $nk$-spectrum gives a
mathematical description of such fundamental concepts of nonlinear optics as
phase matching, frequency matching, four wave interaction in cubic media and
three wave interaction in quadratic media. If a multi-wavepacket has a
resonance invariant spectrum, all these phenomena may take place in the
internal dynamics of \ the multi-wavepacket, but do not lead to resonant \
interactions with continuum of all remaining modes.

\subsection{Genericity of the $nk$-spectrum invariance condition}

In simpler situations, when the number of bands $J$ and wavepackets $N$ are
not too large, the resonance invariance of $nk$- spectrum can be easily
verified as above in Examples \ref{Example quad}, \ref{Example cube}, but
what one can say if $J$ or $N$ are large, or if the dispersion relations are
not explicitly given? We show below that in properly defined non-degenerate
cases a small variation of $K_{S}$ makes $S$ universally resonance
invariant, i.e. the resonance invariance is a generic phenomenon..

Assume that the dispersion relations $\omega _{n}\left( \mathbf{k}\right)
\geq 0$, $n\in \left\{ 1,\ldots ,J\right\} $ are given. Observe then that $%
\Omega _{m}\left( \zeta ,n,\vec{\lambda}\right) =\Omega _{m}\left( \zeta ,n,%
\vec{\lambda}\right) \left( \mathbf{k}_{\ast 1},\ldots ,\mathbf{k}_{\ast
\left\vert K_{S}\right\vert }\right) $ defined by (\ref{Omzet}) is a
continuous function of $\mathbf{k}_{\ast l}\notin \sigma _{\mathrm{bc}}$ for
every $m,\zeta ,n,\vec{\lambda}$.

\begin{definition}[$\protect\omega $-degenerate dispersion relations]
\label{Definition degen}We call dispersion relations $\omega _{n}\left( 
\mathbf{k}\right) $, $n=1,$\ldots $,J$, $\omega $-degenerate if there exists
such a point $\mathbf{k}_{\ast }\in \mathbb{R}^{d}\setminus \sigma _{\mathrm{%
bc}}$ that for all $\mathbf{k}$ \ in a neighborhood of $\mathbf{k}_{\ast }$
at least one of the following four conditions holds: (i) the relations are
linearly dependent, namely $\sum_{n=0}^{J}C_{n}\omega _{n}\left( \mathbf{k}%
\right) =c_{0},$ where \ all $C_{n}$ are integers, one of which is nonzero,
and the $c_{0}$ is a constant; (ii) at least one of $\omega _{n}\left( 
\mathbf{k}\right) $ is a linear function; (iii) at least one of $\omega
_{n}\left( \mathbf{k}\right) $ satisfies equation $C\omega _{n}\left( 
\mathbf{k}\right) =\omega _{n}\left( C\mathbf{k}\right) $ with some $n$ and
integer $C\neq \pm 1$; (iv) at least one of $\omega _{n}\left( \mathbf{k}%
\right) $ satisfies equation $\omega _{n}\left( \mathbf{k}\right) =\omega
_{n^{\prime }}\left( -\mathbf{k}\right) $ where $n^{\prime }\neq n$.
\end{definition}

Note that fulfillment of any of the four conditions in Definition \ref%
{Definition degen} makes impossible turning some non resonance invariant
sets into resonance invariant ones by a variation of $\mathbf{k}_{\ast l}$.
For instance, if $\mathfrak{M}_{F}=\left\{ 2\right\} $ as in Example \ref%
{Example quad} and $2\omega _{1}\left( \mathbf{k}\right) =\omega _{1}\left( 2%
\mathbf{k}\right) $ for all $\mathbf{k}$ in an open set $G$ then the set $%
\left\{ \left( 1,\mathbf{k}_{\ast }\right) \right\} $ with $\mathbf{k}_{\ast
}\in G$ cannot be made resonance invariant by a small variation of $\mathbf{k%
}_{\ast }$. Below we formulate two theorems which show that if dispersion
relations are not $\omega $-degenerate, then a small variation of $\mathbf{k}%
_{\ast l}$ turns non resonance invariant sets into resonance invariant; the
proofs of the theorems are given in \cite{BF7}

\begin{theorem}
\label{Theorem alt}If $\Omega _{m}\left( \zeta ,n_{0},\vec{\lambda}\right)
\left( \mathbf{k}_{\ast 1}^{\prime },\ldots ,\mathbf{k}_{\ast \left\vert
K_{S}\right\vert }^{\prime }\right) =0$ on a cylinder $G$\ in $\left( 
\mathbb{R}^{d}\setminus \sigma _{\mathrm{bc}}\right) ^{\left\vert
K_{S}\right\vert }$ which is a product of small balls $G_{i}\subset \left( 
\mathbb{R}^{d}\setminus \sigma _{\mathrm{bc}}\right) $ then either $\left(
m,\zeta ,n_{0},\vec{\lambda}\right) \in P_{\text{univ}}\left( S\right) $ or
dispersive relations $\omega _{n}\left( \mathbf{k}\right) $ are $\omega $%
-degenerate as in Definition \ref{Definition degen}.
\end{theorem}

\begin{theorem}[genericity of resonance invariance]
\label{Proposition generic} Assume that dispersive relations $\omega
_{n}\left( \mathbf{k}\right) $ are continuous and not $\omega $-degenerate
as in Definition \ref{Definition degen}. Let $\mathcal{K}_{\text{rinv}}$ be
a set of points $\left( \mathbf{k}_{\ast 1},\ldots ,\mathbf{k}_{\ast
\left\vert K_{S}\right\vert }\right) $ such that there exists a universally
resonance invariant $nk$-spectrum $S$ for which its $k$-spectrum $%
K_{S}=\left\{ \mathbf{k}_{\ast 1},\ldots ,\mathbf{k}_{\ast \left\vert
K_{S}\right\vert }\right\} $. Then $\mathcal{K}_{\text{rinv}}$ is open and
everywhere dense set in $\left( \mathbb{R}^{d}\setminus \sigma _{\mathrm{bc}%
}\right) ^{\left\vert K_{S}\right\vert }$.
\end{theorem}

\section{Integrated evolution equation}

Using the variation of constants formula we recast the modal evolution
equation (\ref{difeqfou}) into the following equivalent integral form 
\begin{equation}
\mathbf{\hat{U}}\left( \mathbf{k},\tau \right) =\int_{0}^{\tau }\mathrm{e}^{%
\frac{-\mathrm{i}\left( \tau -\tau ^{\prime }\right) }{\varrho }\mathbf{L}%
\left( \mathbf{k}\right) }\hat{F}\left( \mathbf{\hat{U}}\right) \left( 
\mathbf{k},\tau \right) \,\mathrm{d}\tau ^{\prime }+\mathrm{e}^{\frac{-%
\mathrm{i}\zeta \tau }{\varrho }\mathbf{L}\left( \mathbf{k}\right) }\mathbf{%
\hat{h}}\left( \mathbf{k}\right) ,\ \tau \geq 0.  \label{varc}
\end{equation}%
Then we factor $\mathbf{\hat{U}}\left( \mathbf{k},\tau \right) $ into the
slow variable $\mathbf{\hat{u}}\left( \mathbf{k},\tau \right) $ and the fast
oscillatory term as in (\ref{Uu0}), namely 
\begin{equation}
\mathbf{\hat{U}}\left( \mathbf{k},\tau \right) =\mathrm{e}^{-\frac{\mathrm{i}%
\tau }{\varrho }\mathbf{L}\left( \mathbf{k}\right) }\mathbf{\hat{u}}\left( 
\mathbf{k},\tau \right) ,\ \mathbf{\hat{U}}_{n,\zeta }\left( \mathbf{k},\tau
\right) =\mathbf{\hat{u}}_{n,\zeta }\left( \mathbf{k},\tau \right) \mathrm{e}%
^{-\frac{\mathrm{i}\tau }{\varrho }\zeta \omega _{n}\left( \mathbf{k}\right)
},  \label{Uu}
\end{equation}%
where $\mathbf{\hat{u}}_{n,\zeta }\left( \mathbf{k},\tau \right) $ are the
modal components of $\mathbf{\hat{u}}\left( \mathbf{k},\tau \right) $ as in (%
\ref{Uboldj}). Notice that $\mathbf{\hat{u}}_{n,\zeta }\left( \mathbf{k}%
,\tau \right) $ in (\ref{Uu}) may depend on $\varrho $ and (\ref{Uu}) is
just a change of variables and not an assumption.

\begin{remark}
\label{Remark shift}Note that if $\mathbf{\hat{u}}_{n,\zeta }\left( \mathbf{k%
},\tau \right) $ is a wavepacket, it is localized near its principal
wavevector $\mathbf{k}_{\ast }$. The expansion of $\ \zeta \omega _{n}\left( 
\mathbf{k}\right) $ near the principal wavevector $\zeta \mathbf{k}_{\ast }$
(we take $\zeta =1$ for brevity) takes the form \ 
\begin{equation*}
\omega _{n}\left( \mathbf{k}\right) =\omega \left( \mathbf{k}_{\ast }\right)
+\nabla _{k}\omega _{n}\left( \mathbf{k}_{\ast }\right) \left( \mathbf{k}-%
\mathbf{k}_{\ast }\right) +\frac{1}{2}\nabla _{k}^{2}\omega \left( \mathbf{k}%
_{\ast }\right) \left( \mathbf{k}-\mathbf{k}_{\ast }\right) ^{2}+\ldots
\end{equation*}%
To discuss the impact of the change of variables (\ref{Uu}) we make the
change of variables $\mathbf{k}-\mathbf{k}_{\ast }=\mathbf{\xi }\ $. The
change of variables (\ref{Uu}) 
\begin{gather}
\mathbf{\hat{U}}_{n,+}\left( \mathbf{k},\tau \right) =\mathbf{\hat{u}}%
_{n,+}\left( \mathbf{k},\tau \right) \mathrm{e}^{-\frac{\mathrm{i}\tau }{%
\varrho }\omega _{n}\left( \mathbf{k}\right) }  \label{Usum} \\
=\mathbf{\hat{u}}_{n,\zeta }\left( \mathbf{k},\tau \right) \mathrm{e}^{-%
\frac{\mathrm{i}\tau }{\varrho }\omega _{n}\left( \mathbf{k}_{\ast }\right) }%
\mathrm{e}^{-\frac{\mathrm{i}\tau }{\varrho }\nabla _{k}\omega _{n}\left( 
\mathbf{k}_{\ast }\right) \left( \mathbf{k}-\mathbf{k}_{\ast }\right) }%
\mathrm{e}^{-\frac{\mathrm{i}\tau }{\varrho }\left( \frac{1}{2}\nabla
_{k}^{2}\omega _{n}\left( \mathbf{k}_{\ast }\right) \left( \mathbf{k}-%
\mathbf{k}_{\ast }\right) ^{2}+\ldots \right) }  \notag \\
=\mathbf{\hat{u}}_{n,+}\left( \mathbf{k}_{\ast }+\mathbf{\xi },\tau \right) 
\mathrm{e}^{-\frac{\mathrm{i}\tau }{\varrho }\zeta \omega _{n}\left( \mathbf{%
k}_{\ast }\right) }\mathrm{e}^{-\frac{\mathrm{i}\tau }{\varrho }\nabla
_{k}\omega _{n}\left( \mathbf{k}_{\ast }\right) \mathbf{\xi }}\mathrm{e}^{-%
\frac{\mathrm{i}\tau }{\varrho }R\left( \mathbf{\xi }\right) }  \notag \\
R\left( \mathbf{\xi }\right) =\omega _{n}\left( \mathbf{k}\right) -\omega
_{n}\left( \mathbf{k}_{\ast }\right) -\nabla _{k}\omega _{n}\left( \mathbf{k}%
_{\ast }\right) \left( \mathbf{k}-\mathbf{k}_{\ast }\right) =\frac{1}{2}%
\nabla _{k}^{2}\omega _{n}\left( \mathbf{k}_{\ast }\right) \left( \mathbf{%
\xi }\right) ^{2}+\ldots  \label{Rxi}
\end{gather}%
has the first factor $\mathrm{e}^{-\frac{\mathrm{i}\tau }{\varrho }\omega
_{n}\left( \mathbf{k}_{\ast }\right) }$ responsible for fast time
oscillations \ of $\ \mathbf{\hat{U}}_{n,\zeta }\left( \mathbf{k},\tau
\right) $ and $\mathbf{U}_{n,\zeta }\left( \mathbf{r},\tau \right) $ The
second factor $\mathrm{e}^{-\frac{\mathrm{i}\tau }{\varrho }\nabla
_{k}\omega _{n}\left( \mathbf{k}_{\ast }\right) \mathbf{\xi }}$ is
responsible for the spatial shifts of the inverse Fourier transform by $%
\frac{\tau }{\varrho }\nabla _{k}\omega _{n}\left( \mathbf{k}_{\ast }\right)
,$ since the shifts are time dependent they cause the rectilinear movement
of the wavepacket $\mathbf{U}_{n,\zeta }\left( \mathbf{r},\tau \right) $
with the group velocity $\frac{1}{\varrho }\nabla _{k}\omega _{n}\left( 
\mathbf{k}_{\ast }\right) $, the third factor is responsible for dispersive
effects. Hence the change of variables (\ref{Uu}) effectively introduces the
moving coordinate frame for $\mathbf{\hat{U}}_{n,\zeta }\left( \mathbf{k}%
,\tau \right) $ for every $\mathbf{k}$ and in this coordinate frame $\mathbf{%
\hat{u}}_{n,\zeta }\left( \mathbf{k},\tau \right) $ has zero group velocity
\ and does not have high-frequency time oscillations. The following
proposition shows that if $\mathbf{\hat{u}}_{n,\zeta }\left( \mathbf{k},\tau
\right) $ is a wavepacket with a constant position, $\mathbf{\hat{U}}%
_{n,+}\left( \mathbf{k},\tau \right) $ is a particle wavepacket in the sense
of Definition with position which moves with a constant velocity.
\end{remark}

\begin{proposition}
\label{Proposition Uu}Let $\mathbf{\hat{u}}_{l}\left( \mathbf{k},\tau
\right) $ be for every $\tau \in \left[ 0,\tau _{\ast }\right] $ a particle
wavepacket in the sense of Definition \ref{Definition regwave} with $nk$%
\emph{-}pair\emph{\ }$\left( n,\mathbf{k}_{\ast }\right) $ regularity $s$\
and position $\mathbf{r}_{\ast }\in \mathbb{R}^{d}$ which does not depend on 
$\tau ,$ assume also that and constants $C_{1}$ in (\ref{gradker}) and $%
C,C^{\prime }$ in (\ref{L1b}) and (\ref{sourloc}) do not depend on $\tau $.
Let $\mathbf{\hat{U}}_{l}\left( \mathbf{k},\tau \right) $ be defined in
terms of $\mathbf{\hat{u}}_{l}\left( \mathbf{k},\tau \right) $ by (\ref{Uu}%
). Assume that (\ref{scale1}) holds. Then $\mathbf{\hat{U}}_{l}\left( 
\mathbf{k},\tau \right) $ for every $\tau \in \left[ 0,\tau _{\ast }\right] $
a particle wavepacket in the sense of Definition \ref{Definition regwave}
with $nk$\emph{-}pair\emph{\ }$\left( n,\mathbf{k}_{\ast }\right) $
regularity $s$ and with $\tau $-dependent position $\mathbf{r}_{\ast }+\frac{%
\tau }{\varrho }\nabla _{k}\omega _{n}\left( \mathbf{k}_{\ast }\right) \in 
\mathbb{R}^{d}$.
\end{proposition}

\begin{proof}
The wavepacket $\mathbf{\hat{u}}_{l}\left( \mathbf{k},\tau \right) $ \
involves two components $\mathbf{\hat{u}}_{n,\zeta }\left( \mathbf{k},\tau
\right) $, $\zeta =\pm 1$ for which \ (\ref{halfwave}) holds \ 
\begin{equation}
\mathbf{\hat{u}}_{n,\zeta }\left( \mathbf{k},\tau \right) =\Psi \left( \beta
^{1-\epsilon }/2,\zeta \mathbf{k}_{\ast };\mathbf{k}\right) \Pi _{n,\zeta
}\left( \mathbf{k}\right) \mathbf{\hat{u}}_{n,\zeta }\left( \mathbf{k},\tau
\right) ,\   \label{unz}
\end{equation}%
By (\ref{Uu}) 
\begin{equation*}
\mathbf{\hat{U}}_{n,\zeta }\left( \mathbf{k},\tau \right) =\mathbf{\hat{u}}%
_{n,\zeta }\left( \mathbf{k},\tau \right) \mathrm{e}^{-\frac{\mathrm{i}\tau 
}{\varrho }\zeta \omega _{n}\left( \mathbf{k}\right) }.
\end{equation*}%
According to Definition (\ref{dwavepack}) multiplication by a scalar bounded
continuous function $\mathrm{e}^{-\frac{\mathrm{i}\tau }{\varrho }\zeta
\omega _{n}\left( \mathbf{k}\right) }$ \ \ may only change a constant $%
C^{\prime }$ in (\ref{sourloc}), therefore it transforms wavepackets into
wavepackets. To check that $\mathbf{\hat{U}}_{l}\left( \mathbf{k},\tau
\right) $ \ is a particle-like wavepacket \ we consider (\ref{gradker}) \
with $\mathbf{\hat{h}}_{\zeta }\left( \beta ,\mathbf{r}_{\ast };\mathbf{k}%
\right) $ replaced by $\mathbf{\hat{u}}_{n,\zeta }\left( \mathbf{k},\tau
\right) \mathrm{e}^{-\frac{\mathrm{i}\tau }{\varrho }\zeta \omega _{n}\left( 
\mathbf{k}\right) }$ and $\mathbf{r}_{\ast }$ replaced by $\ \mathbf{r}%
_{\ast }+\frac{\tau }{\varrho }\nabla _{k}\omega _{n}\left( \mathbf{k}_{\ast
}\right) $. \ We consider for brevity $\mathbf{\hat{u}}_{n}\left( \mathbf{k}%
,\tau \right) =\mathbf{\hat{u}}_{n,\zeta }\left( \mathbf{k},\tau \right) $
with $\zeta =1,$the case $\zeta =-1$ is similar. 
\begin{gather*}
\int_{\mathbb{R}^{d}}\left\vert \nabla _{\mathbf{k}}\left( e^{i\left( 
\mathbf{r}_{\ast }+\frac{\tau }{\varrho }\nabla _{k}\omega _{n}\left( 
\mathbf{k}_{\ast }\right) \right) \mathbf{k}}\mathbf{\hat{u}}_{n}\left( 
\mathbf{k},\tau \right) \mathrm{e}^{-\frac{\mathrm{i}\tau }{\varrho }\omega
_{n}\left( \mathbf{k}\right) }\right) \right\vert d\mathbf{k=} \\
\int_{\mathbb{R}^{d}}\left\vert \nabla _{\mathbf{k}}\left( e^{i\left( 
\mathbf{r}_{\ast }+\frac{\tau }{\varrho }\nabla _{k}\omega _{n}\left( 
\mathbf{k}_{\ast }\right) \right) \mathbf{k}}\mathbf{\hat{u}}_{n}\left( 
\mathbf{k},\tau \right) \mathrm{e}^{-\frac{\mathrm{i}\tau }{\varrho }\omega
_{n}\left( \mathbf{k}\right) }\mathrm{e}^{\frac{\mathrm{i}\tau }{\varrho }%
\omega _{n}\left( \mathbf{k}_{\ast }\right) }\right) \right\vert d\mathbf{k}
\\
=\int_{\mathbb{R}^{d}}\left\vert \nabla _{\mathbf{k}}\left( e^{i\mathbf{r}%
_{\ast }\mathbf{k}}\mathbf{\hat{u}}_{n}\left( \mathbf{k},\tau \right) 
\mathrm{e}^{-\frac{\mathrm{i}\tau }{\varrho }R\left( \mathbf{k}-\mathbf{k}%
_{\ast }\right) }\right) \right\vert d\mathbf{k}\leq I_{1}+I_{2}
\end{gather*}%
where where $R\left( \mathbf{\xi }\right) $ is defined by (\ref{Rxi}),%
\begin{eqnarray*}
I_{1} &=&\int_{\mathbb{R}^{d}}\left\vert \mathrm{e}^{-\frac{\mathrm{i}\tau }{%
\varrho }R\left( \mathbf{k}-\mathbf{k}_{\ast }\right) }\nabla _{\mathbf{k}%
}\left( e^{i\mathbf{r}_{\ast }\mathbf{k}}\mathbf{\hat{u}}_{n}\left( \mathbf{k%
},\tau \right) \right) \right\vert d\mathbf{k}, \\
I_{2} &=&\int_{\mathbb{R}^{d}}\left\vert \left( e^{i\mathbf{r}_{\ast }%
\mathbf{k}}\mathbf{\hat{u}}_{n}\left( \mathbf{k},\tau \right) \right) \nabla
_{\mathbf{k}}\mathrm{e}^{-\frac{\mathrm{i}\tau }{\varrho }R\left( \mathbf{k}-%
\mathbf{k}_{\ast }\right) }\right\vert d\mathbf{k}.
\end{eqnarray*}%
The integral $I_{1}$ is bounded uniformly in $\mathbf{r}_{\ast }$ by $%
C_{1}^{\prime }\beta ^{-1-\epsilon }$ since $\mathbf{\hat{u}}_{n,\zeta
}\left( \mathbf{k},\tau \right) $ satisfies (\ref{gradker}). \ Note that 
\begin{eqnarray*}
I_{2} &=&\int_{\mathbb{R}^{d}}\left\vert \left( e^{i\mathbf{r}_{\ast }%
\mathbf{k}}\mathbf{\hat{u}}_{n}\left( \mathbf{k},\tau \right) \right) \nabla
_{\mathbf{k}}\mathrm{e}^{-\frac{\mathrm{i}\tau }{\varrho }R\left( \mathbf{k}-%
\mathbf{k}_{\ast }\right) }\right\vert d\mathbf{k}\leq \\
&&\int_{\mathbb{R}^{d}}\left\vert \mathbf{\hat{u}}_{n}\left( \mathbf{k},\tau
\right) \right\vert \frac{\tau }{\varrho }\left\vert \nabla _{\mathbf{k}%
}R\left( \mathbf{k}-\mathbf{k}_{\ast }\right) \right\vert d\mathbf{k}
\end{eqnarray*}%
Note that according to (\ref{unz}) \ and (\ref{Psik}) $\mathbf{\hat{u}}%
_{n,\zeta }\left( \mathbf{k},\tau \right) \neq 0$ \ only if $\left\vert 
\mathbf{k}-\mathbf{k}_{\ast }\right\vert \leq 2\beta ^{1-\epsilon }$, and
for such $\mathbf{k}-\mathbf{k}_{\ast }$ we have Taylor remainder estimate 
\begin{equation*}
\left\vert \nabla _{\mathbf{k}}R\left( \mathbf{k}-\mathbf{k}_{\ast }\right)
\right\vert \leq C\beta ^{1-\epsilon }.
\end{equation*}%
Therefore $I_{2}\leq C^{\prime }\beta ^{1-\epsilon }/\varrho $ and 
\begin{equation*}
I_{1}+I_{2}\leq C^{\prime }\beta ^{1-\epsilon }/\varrho
\end{equation*}%
Using (\ref{scale1}) we conclude that this inequality implies (\ref{gradker}%
) \ \ for $\mathbf{\hat{U}}_{l}\left( \mathbf{k},\tau \right) $, therefore
it is a particle-like wavepacket.
\end{proof}

From (\ref{varc}) and (\ref{Uu}) we obtain the following \emph{integrated
evolution equation} for $\mathbf{\hat{u}}=\mathbf{\hat{u}}\left( \mathbf{k}%
,\tau \right) $, $\tau \geq 0$, 
\begin{gather}
\mathbf{\hat{u}}\left( \mathbf{k},\tau \right) =\mathcal{F}\left( \mathbf{%
\hat{u}}\right) \left( \mathbf{k},\tau \right) +\mathbf{\hat{h}}\left( 
\mathbf{k}\right) ,\ \mathcal{F}\left( \mathbf{\hat{u}}\right)
=\sum\nolimits_{m\in \mathfrak{M}_{F}}\mathcal{F}^{\left( m\right) }\left( 
\mathbf{\hat{u}}^{m}\left( \mathbf{k},\tau \right) \right) ,  \label{varcu}
\\
\mathcal{F}^{\left( m\right) }\left( \mathbf{\hat{u}}^{m}\right) \left( 
\mathbf{k},\tau \right) =\int_{0}^{\tau }\mathrm{e}^{\frac{\mathrm{i}\tau
^{\prime }}{\varrho }\mathbf{L}\left( \mathbf{k}\right) }\hat{F}_{m}\left(
\left( \mathrm{e}^{\frac{-\mathrm{i}\tau ^{\prime }}{\varrho }\mathbf{L}%
\left( \cdot \right) }\mathbf{\hat{u}}\right) ^{m}\right) \left( \mathbf{k}%
,\tau ^{\prime }\right) \,\mathrm{d}\tau ^{\prime },  \label{Fu}
\end{gather}%
where $\hat{F}_{m}$ are defined by (\ref{Fseries}) and (\ref{Fmintr}) in
terms of the susceptibilities $\chi ^{\left( m\right) }$, and $\mathcal{F}%
^{\left( m\right) }$ are bounded as in the following lemma.

Recall that spaces $L^{1,a}$ are defined by the formula (\ref{L1a}). Below
we formulate basic properties of the spaces. Recall Young's inequality 
\begin{equation}
\left\Vert \mathbf{\hat{u}}\ast \mathbf{\hat{v}}\right\Vert _{L^{1}}\leq
\left\Vert \mathbf{\hat{u}}\right\Vert _{L^{1}}\left\Vert \mathbf{\hat{v}}%
\right\Vert _{L^{1}}.  \label{Yconv}
\end{equation}

This inequality implies boundedness of convolution in $L^{1,a}$, namely the
following Lemma holds.

\begin{lemma}
\label{Lemma convpsi} Let $\hat{H}_{1},\hat{H}_{2}\in L^{1,a}$ be two scalar
functions, $a\geq 0$. Let 
\begin{equation*}
\hat{H}_{3}\left( \mathbf{k}\right) =\int_{\mathbb{R}^{d}}\hat{H}_{1}\left( 
\mathbf{k}-\mathbf{k}^{\prime }\right) \hat{H}_{2}\left( \mathbf{k}^{\prime
}\right) \mathrm{d}\mathbf{k}^{\prime }.
\end{equation*}%
Then%
\begin{equation}
\left\Vert \hat{H}_{3}\left( \mathbf{k}\right) \right\Vert _{L^{1,a}}\leq
\left\Vert \hat{H}_{1}\left( \mathbf{k}\right) \right\Vert
_{L^{1,a}}\left\Vert \hat{H}_{1}\left( \mathbf{k}\right) \right\Vert
_{L^{1,a}}.  \label{L1conpsi}
\end{equation}
\end{lemma}

\begin{proof}
We have 
\begin{gather*}
\left( 1+\left\vert \mathbf{k}\right\vert \right) ^{a}\left\vert \hat{H}%
_{3}\left( \mathbf{k}\right) \right\vert \leq \\
\sup_{\mathbf{k}^{\prime },\mathbf{k}^{\prime \prime }}\frac{\left(
1+\left\vert \mathbf{k}^{\prime }+\mathbf{k}^{\prime \prime }\right\vert
\right) ^{a}}{\left( 1+\left\vert \mathbf{k}^{\prime }\right\vert \right)
^{a}\left( 1+\left\vert \mathbf{k}^{\prime \prime }\right\vert \right) ^{a}}%
\int_{\mathbb{R}^{d}}\left( 1+\left\vert \mathbf{k}-\mathbf{k}^{\prime
}\right\vert \right) ^{a}\left\vert \hat{H}_{1}\left( \mathbf{k}-\mathbf{k}%
^{\prime }\right) \right\vert \left( 1+\left\vert \mathbf{k}^{\prime
}\right\vert \right) ^{a}\left\vert \hat{H}_{2}\left( \mathbf{k}^{\prime
}\right) \right\vert \mathrm{d}\mathbf{k}^{\prime }.
\end{gather*}%
Obviously, 
\begin{equation*}
\frac{1+\left\vert \mathbf{k}^{\prime }+\mathbf{k}^{\prime \prime
}\right\vert }{\left( 1+\left\vert \mathbf{k}^{\prime }\right\vert \right)
\left( 1+\left\vert \mathbf{k}^{\prime \prime }\right\vert \right) }\leq 
\frac{\left( 1+\left\vert \mathbf{k}^{\prime }\right\vert +\left\vert 
\mathbf{k}^{\prime \prime }\right\vert \right) }{\left( 1+\left\vert \mathbf{%
k}^{\prime }\right\vert \right) \left( 1+\left\vert \mathbf{k}^{\prime
\prime }\right\vert \right) }\leq 1.
\end{equation*}%
Applying Young's inequality (\ref{Yconv}) we obtain 
\begin{equation*}
\int_{\mathbb{R}^{d}}\left( 1+\left\vert \mathbf{k}\right\vert \right)
^{a}\left\vert \hat{H}_{3}\left( \mathbf{k}\right) \right\vert \mathrm{d}%
\mathbf{k}\leq \int_{\mathbb{R}^{d}}\left( 1+\left\vert \mathbf{k}^{\prime
}\right\vert \right) \left\vert \hat{H}_{1}\left( \mathbf{k}^{\prime
}\right) \right\vert \mathrm{d}\mathbf{k}^{\prime }\int_{\mathbb{R}%
^{d}}\left( 1+\left\vert \mathbf{k}^{\prime \prime }\right\vert \right)
\left\vert \hat{H}_{2}\left( \mathbf{k}^{\prime \prime }\right) \right\vert 
\mathrm{d}\mathbf{k}^{\prime \prime }
\end{equation*}
Using (\ref{Ea}) we obtain(\ref{L1conpsi}).
\end{proof}

Using Lemma \ref{Lemma convpsi} we derive boundedness of integral operators $%
\mathcal{F}^{\left( m\right) }$.

\begin{lemma}[boundness of multilinear operators]
\label{Lemma bound}Operator $\mathcal{F}^{\left( m\right) }$ defined by (\ref%
{Fmintr}), (\ref{Fu}) is bounded from $E_{a}=C\left( \left[ 0,\tau _{\ast }%
\right] ,L^{1,a}\right) $ into $C^{1}\left( \left[ 0,\tau _{\ast }\right]
,L^{1,a}\right) $, $a\geq 0$ \ and 
\begin{equation}
\left\Vert \mathcal{F}^{\left( m\right) }\left( \mathbf{\hat{u}}_{1}\ldots 
\mathbf{\hat{u}}_{m}\right) \right\Vert _{E_{a}}\leq \tau _{\ast }\left\Vert
\chi ^{\left( m\right) }\right\Vert \dprod\nolimits_{j=1}^{m}\left\Vert 
\mathbf{\hat{u}}_{j}\right\Vert _{E_{a}},  \label{dtf}
\end{equation}%
\begin{equation}
\left\Vert \partial _{\tau }\mathcal{F}^{\left( m\right) }\left( \mathbf{%
\hat{u}}_{1}\ldots \mathbf{\hat{u}}_{m}\right) \right\Vert _{E_{a}}\leq
\left\Vert \chi ^{\left( m\right) }\right\Vert \dprod\nolimits_{j}\left\Vert 
\mathbf{\hat{u}}_{j}\right\Vert _{E_{a}}.  \label{dtf1}
\end{equation}
\end{lemma}

\begin{proof}
Notice that since $\mathbf{L}\left( \mathbf{\mathbf{k}}\right) $ is
Hermitian, $\left\Vert \exp \left\{ -\mathrm{i}\mathbf{L}\left( \mathbf{%
\mathbf{k}}\right) \frac{\tau _{1}}{\varrho }\right\} \right\Vert =1$. Using
the inequality (\ref{L1conpsi}) together with (\ref{Fmintr}), (\ref{Fu}) we
obtain 
\begin{gather*}
\left\Vert \mathcal{F}^{\left( m\right) }\left( \mathbf{\hat{u}}_{1}\ldots 
\mathbf{\hat{u}}_{m}\right) \left( \mathbf{\cdot },\tau \right) \right\Vert
_{L^{1,a}}\leq \sup_{\ \mathbf{\mathbf{k}},\vec{k}}\left\vert \chi _{\
}^{\left( m\right) }\left( \mathbf{\mathbf{k}},\vec{k}\right) \right\vert \\
\int_{\mathbb{R}^{d}}\int_{0}^{\tau }\int_{\mathbb{D}_{m}}\left\vert \left(
1+\left\vert \mathbf{k}^{\prime }\right\vert \right) ^{a}\mathbf{\hat{u}}%
_{1}\left( \mathbf{k}^{\prime }\right) \right\vert \ldots \left\vert \left(
1+\left\vert \mathbf{k}^{\left( m\right) }\right\vert \right) ^{a}\mathbf{%
\hat{u}}_{m}\left( \mathbf{k}^{\left( m\right) }\left( \mathbf{k},\vec{k}%
\right) \right) \right\vert \mathrm{d}\mathbf{k}^{\prime }\ldots \mathrm{d}%
\mathbf{k}^{\left( m-1\right) }\mathrm{d}\tau _{1}\mathrm{d}\mathbf{k}\leq \\
\left\Vert \chi ^{\left( m\right) }\right\Vert \int_{0}^{\tau }\left\Vert 
\mathbf{\hat{u}}_{1}\left( \tau _{1}\right) \right\Vert _{L^{1,a}}\ldots
\left\Vert \mathbf{\hat{u}}_{m}\left( \tau _{1}\right) \right\Vert _{L^{1,a}}%
\mathrm{d}\tau _{1}\leq \tau _{\ast }\left\Vert \chi ^{\left( m\right)
}\right\Vert \left\Vert \mathbf{\hat{u}}_{1}\right\Vert _{E_{a}}\ldots
\left\Vert \mathbf{\hat{u}}_{m}\right\Vert _{E_{a}}.
\end{gather*}%
proving (\ref{dtf}). A similar \ estimate produces we prove (\ref{dtf1}).
\end{proof}

The equation (\ref{varcu}) can be recast as the following abstract equation
in a Banach space 
\begin{equation}
\mathbf{\hat{u}}=\mathcal{F}\left( \mathbf{\hat{u}}\right) +\mathbf{\hat{h}}%
,\ \mathbf{\hat{u}},\mathbf{\hat{h}}\in E_{a},  \label{eqF}
\end{equation}%
and it readily follows from Lemma \ref{Lemma bound} that $\mathcal{F}\left( 
\mathbf{\hat{u}}\right) $ has the following properties.

\begin{lemma}
\label{Lemma Flip} The operator\ $\mathcal{F}\left( \mathbf{\hat{u}}\right) $
defined by (\ref{varcu})-(\ref{Fu}) satisfies the Lipschitz condition 
\begin{equation}
\left\Vert \mathcal{F}\left( \mathbf{\hat{u}}_{1}\right) -\mathcal{F}\left( 
\mathbf{\hat{u}}_{2}\right) \right\Vert _{E_{a}}\leq \tau _{\ast
}C_{F}\left\Vert \mathbf{\hat{u}}_{1}-\mathbf{\hat{u}}_{2}\right\Vert
_{E_{a}}  \label{Flip}
\end{equation}%
where $C_{F}\leq C_{\chi }m_{F}^{2}\left( 4R\right) ^{m_{F}-1}$ if $%
\left\Vert \mathbf{\mathbf{\hat{u}}}_{1}\right\Vert _{E_{a}},\left\Vert 
\mathbf{\mathbf{\hat{u}}}_{2}\right\Vert _{E_{a}}\leq 2R$, with $C_{\chi }$
as in (\ref{chiCR}).
\end{lemma}

We also will use the following form of the contraction principle.

\begin{lemma}[contraction principle]
\label{Lemma contr}Consider equation \ 
\begin{equation}
\mathbf{x}=\mathcal{F}\left( \mathbf{x}\right) +\mathbf{h},\ \mathbf{x},%
\mathbf{h}\in B,  \label{eqFa}
\end{equation}%
where $B$ is a Banach space, $\mathcal{F}$ is an operator in $B$. Suppose
that for some constants $R_{0}>0$ and $0<q<1$ we have%
\begin{eqnarray}
\left\Vert \mathbf{h}\right\Vert &\leq &R_{0},\ \left\Vert \mathcal{F}\left( 
\mathbf{x}\right) \right\Vert \leq R_{0}\text{ if }\left\Vert \mathbf{x}%
\right\Vert \leq 2R_{0},  \label{Fle} \\
\left\Vert \mathcal{F}\left( \mathbf{x}_{1}\right) -\mathcal{F}\left( 
\mathbf{x}_{2}\right) \right\Vert &\leq &q\left\Vert \mathbf{x}_{1}-\mathbf{x%
}_{2}\right\Vert \text{ if }\left\Vert \mathbf{x}_{1}\right\Vert ,\left\Vert 
\mathbf{x}_{2}\right\Vert \leq 2R_{0}.  \label{Flipa}
\end{eqnarray}%
Then there exists a unique solution $\mathbf{x}$ to the equation (\ref{eqFa}%
) such that $\left\Vert \mathbf{x}\right\Vert \leq 2R_{0}$. Let $\left\Vert 
\mathbf{h}_{1}\right\Vert ,\left\Vert \mathbf{h}_{2}\right\Vert \leq R_{0}$
then the two corresponding solutions $\mathbf{x}_{1},\mathbf{x}_{2}$ satisfy 
\begin{equation}
\left\Vert \mathbf{x}_{1}\right\Vert ,\left\Vert \mathbf{x}_{2}\right\Vert
\leq 2R_{0},\ \left\Vert \mathbf{x}_{1}-\mathbf{x}_{2}\right\Vert \leq
\left( 1-q\right) ^{-1}\left\Vert \mathbf{h}_{1}-\mathbf{h}_{2}\right\Vert .
\label{iminu0}
\end{equation}%
Let $\mathbf{x}_{1},\mathbf{x}_{2}$ be the two solutions of correspondingly
two equations of the form (\ref{eqFa}) with $\mathcal{F}_{1}$, $\mathbf{h}%
_{1}\mathbf{\ }$and\ $\mathcal{F}_{2}$, $\mathbf{h}_{2}$. Assume that that $%
\mathcal{F}_{1}\left( \mathbf{u}\right) $ satisfies (\ref{Fle}), (\ref{Flipa}%
) with a Lipschitz constant $q<1$ and that $\left\Vert \mathcal{F}_{1}\left( 
\mathbf{x}\right) -\mathcal{F}_{2}\left( \mathbf{x}\right) \right\Vert \leq
\delta $ for $\left\Vert \mathbf{x}\right\Vert \leq 2R_{0}$. Then 
\begin{equation}
\left\Vert \mathbf{x}_{1}-\mathbf{x}_{2}\right\Vert \leq \left( 1-q\right)
^{-1}\left( \delta +\left\Vert \mathbf{h}_{1}-\mathbf{h}_{2}\right\Vert
\right) .  \label{uminu}
\end{equation}
\end{lemma}

Lemma \ref{Lemma Flip} and the contraction principle as in Lemma \ref{Lemma
contr} imply the following existence and uniqueness theorem.

\begin{theorem}
\label{Theorem exist}Let $\left\Vert \mathbf{\hat{h}}\right\Vert
_{E_{a}}\leq R$, let $\tau _{\ast }<1/C_{F}$ where $C_{F}$ is a constant
from Lemma \ref{Lemma Flip}. Then equation (\ref{varcu}) has a solution $%
\mathbf{\hat{u}}\in E_{a}=C\left( \left[ 0,\tau _{\ast }\right]
,L^{1,a}\right) $ which satisfies $\left\Vert \mathbf{\hat{u}}\right\Vert
_{E_{a}}\leq 2R$, and such a solution is unique. Hence the solution operator 
$\mathbf{\hat{u}}=\mathcal{G}\left( \mathbf{\hat{h}}\right) $ is defined on
\ the ball $\left\Vert \mathbf{\hat{h}}\right\Vert _{E_{a}}\leq R$.
\end{theorem}

The following existence and uniqueness theorem follows from Theorem \ref%
{Theorem exist}.

\begin{theorem}
\label{Theorem Existence1}\ Let $a\geq 0$, (\ref{difeqfou}) satisfy (\ref%
{chiCR}) and $\mathbf{\hat{h}}\in L^{1,a}\left( \mathbb{R}^{d}\right) $,$%
\left\Vert \mathbf{\hat{h}}\right\Vert _{L^{1,a}}\leq R$. Then there exists
a unique solution $\mathbf{\hat{u}}$ to the modal evolution equation (\ref%
{difeqfou}) in the functional space $C^{1}\left( \left[ 0,\tau _{\ast }%
\right] ,L^{1,a}\right) $, $\left\Vert \mathbf{\hat{u}}\right\Vert
_{E_{a}}+\left\Vert \partial \tau \mathbf{\hat{u}}\right\Vert _{E_{a}}\leq
R_{1}\left( R\right) $. The number $\tau _{\ast }$ depends on $R$ and $%
C_{\chi }$.
\end{theorem}

Using the inequality (\ref{Linf}) and applying the inverse Fourier transform
we readily obtain the existence of an $F-$solution of (\ref{difeqintr}) in $%
C^{1}\left( \left[ 0,\tau _{\ast }\right] ,L^{\infty }\left( \mathbb{R}%
^{d}\right) \right) $ from the existence of the solution of equation (\ref%
{difeqfou}) in $C^{1}\left( \left[ 0,\tau _{\ast }\right] ,L^{1}\right) $.
The existence of $F$-solutions with $\left[ a\right] $ bounded spatial
derivatives ($\left[ a\right] $ being an integer part of $a$) follows from
solvability in $C^{1}\left( \left[ 0,\tau _{\ast }\right] ,L^{1,a}\right) $.

Let us recast now the system (\ref{varcu})-(\ref{Fu}) into modal components
using the projections $\Pi _{n,\zeta }\left( \mathbf{\mathbf{k}}\right) $ as
in (\ref{Pin}). The first step to introduce modal susceptibilities $\chi
_{n,\zeta ,\vec{\xi}}^{\left( m\right) }$ having one-dimensional range in $%
\mathbb{C}^{2J}$ and vanishing if one of its arguments $\mathbf{\hat{u}}_{j}$
belongs to a $\left( 2J-1\right) $-dimensional linear subspace in $\mathbb{C}%
^{2J}$ ($j$-th null-space \ of $\chi _{n,\zeta ,\vec{\xi}}^{\left( m\right)
} $ ) as follows.

\begin{definition}[elementary susceptibilities]
\label{Definition elemsus}Let$\ $%
\begin{equation}
\vec{\xi}=\left( \vec{n},\vec{\zeta}\right) \in \left\{ 1,\ldots ,J\right\}
^{m}\times \left\{ -1,1\right\} ^{m}=\Xi ^{m},\left( n,\zeta \right) \in \Xi
\label{Xi}
\end{equation}%
and $\chi _{\ }^{\left( m\right) }\left( \mathbf{\mathbf{k}},\vec{k}\right) %
\left[ \mathbf{\hat{u}}_{1}\left( \mathbf{k}^{\prime }\right) ,\ldots ,%
\mathbf{\hat{u}}_{m}\left( \mathbf{k}^{\left( m\right) }\right) \right] $ be 
$\ m$-linear symmetric\ tensor (susceptibility) as in (\ref{Fmintr}). We
introduce \emph{elementary} \emph{susceptibilities} $\chi _{n,\zeta ,\vec{\xi%
}}^{\left( m\right) }\left( \mathbf{\mathbf{k}},\vec{k}\right) :\left( 
\mathbb{C}^{2J}\right) ^{m}\rightarrow \mathbb{C}^{2J}$ as$\ m$-linear
tensors defined for almost all $\mathbf{\mathbf{k}}$ and $\vec{k}=\left( 
\mathbf{k}^{\prime },\ldots ,\mathbf{k}^{\left( m\right) }\right) $ by the
following formula 
\begin{gather}
\chi _{n,\zeta ,\vec{\xi}}^{\left( m\right) }\left( \mathbf{\mathbf{k}},\vec{%
k}\right) \left[ \mathbf{\hat{u}}_{1}\left( \mathbf{k}^{\prime }\right)
,\ldots ,\mathbf{\hat{u}}_{m}\left( \mathbf{k}^{\left( m\right) }\right) %
\right] =\chi _{n,\zeta ,\vec{n},\vec{\zeta}}^{\left( m\right) }\left( 
\mathbf{\mathbf{k}},\vec{k}\right) \left[ \mathbf{\hat{u}}_{1}\left( \mathbf{%
k}^{\prime }\right) ,\ldots ,\mathbf{\hat{u}}_{m}\left( \mathbf{k}^{\left(
m\right) }\right) \right] =  \label{chim} \\
\Pi _{n,\zeta }\left( \mathbf{\mathbf{k}}\right) \chi ^{\left( m\right)
}\left( \mathbf{\mathbf{k}},\vec{k}\right) \left[ \left( \Pi _{n_{1},\zeta
^{\prime }}\left( \mathbf{k}^{\prime }\right) \mathbf{\hat{u}}_{1}\left( 
\mathbf{k}^{\prime }\right) ,\ldots ,\Pi _{n_{m},\zeta ^{\left( m\right)
}}\left( \mathbf{k}^{\left( m\right) }\left( \mathbf{k},\vec{k}\right)
\right) \mathbf{\hat{u}}_{m}\left( \mathbf{k}^{\left( m\right) }\right)
\right) \right] .  \notag
\end{gather}
\end{definition}

Then using (\ref{sumPi}) and the \emph{elementary} \emph{susceptibilities} (%
\ref{chim}) we get%
\begin{equation}
\chi ^{\left( m\right) }\left( \mathbf{\mathbf{k}},\vec{k}\right) \left[ 
\mathbf{\hat{u}}_{1}\left( \mathbf{k}^{\prime }\right) ,\ldots ,\mathbf{\hat{%
u}}_{m}\left( \mathbf{k}^{\left( m\right) }\right) \right]
=\sum\nolimits_{n,\zeta }\sum\nolimits_{\vec{\xi}}\chi _{n,\zeta ,\vec{\xi}%
}^{\left( m\right) }\left( \mathbf{\mathbf{k}},\vec{k}\right) \left[ \mathbf{%
\hat{u}}_{1}\left( \mathbf{k}^{\prime }\right) ,\ldots ,\mathbf{\hat{u}}%
_{m}\left( \mathbf{k}^{\left( m\right) }\right) \right] .  \label{chisumn}
\end{equation}%
Consequently the modal components $\mathcal{F}_{n,\zeta ,\vec{\xi}}^{\left(
m\right) }$\ of the operators $\mathcal{F}^{\left( m\right) }$ in (\ref{Fu})
are $m$-linear \emph{oscillatory integral operators} defined in terms of the
elementary susceptibilities (\ref{chisumn}) as follows.

\begin{definition}[interaction phase]
Using notations from (\ref{Fmintr}) we introduce for $\vec{\xi}=\left( \vec{n%
},\vec{\zeta}\right) \in \Xi ^{m}$ operator 
\begin{gather}
\mathcal{F}_{n,\zeta ,\vec{\xi}}^{\left( m\right) }\left( \mathbf{\hat{u}}%
_{1}\ldots \mathbf{\hat{u}}_{m}\right) \left( \mathbf{k},\tau \right)
=\int_{0}^{\tau }\int_{\mathbb{D}_{m}}\exp \left\{ \mathrm{i}\phi _{n,\zeta ,%
\vec{\xi}}\left( \mathbf{\mathbf{k}},\vec{k}\right) \frac{\tau _{1}}{\varrho 
}\right\}  \label{Fm} \\
\chi _{n,\zeta ,\vec{\xi}}^{\left( m\right) }\left( \mathbf{\mathbf{k}},\vec{%
k}\right) \left[ \mathbf{\hat{u}}_{1}\left( \mathbf{k}^{\prime },\tau
_{1}\right) ,\ldots ,\mathbf{\hat{u}}_{m}\left( \mathbf{k}^{\left( m\right)
}\left( \mathbf{\mathbf{k}},\vec{k}\right) ,\tau _{1}\right) \right] \mathrm{%
\tilde{d}}^{\left( m-1\right) d}\vec{k}\mathrm{d}\tau _{1},  \notag
\end{gather}%
with the \emph{interaction phase function }$\phi $ defined by 
\begin{gather}
\phi _{n,\zeta ,\vec{\xi}}\left( \mathbf{\mathbf{k}},\vec{k}\right) =\phi
_{n,\zeta ,\vec{n},\vec{\zeta}}\left( \mathbf{\mathbf{k}},\vec{k}\right)
\label{phim} \\
=\zeta \omega _{n}\left( \zeta \mathbf{k}\right) -\zeta ^{\prime }\omega
_{n_{1}}\left( \zeta ^{\prime }\mathbf{k}^{\prime }\right) -\ldots -\zeta
^{\left( m\right) }\omega _{n_{m}}\left( \zeta ^{\left( m\right) }\mathbf{k}%
^{\left( m\right) }\right) ,\ \mathbf{k}^{\left( m\right) }=\mathbf{k}%
^{\left( m\right) }\left( \mathbf{\mathbf{k}},\vec{k}\right)  \notag
\end{gather}%
where $\mathbf{k}^{\left( m\right) }\left( \mathbf{\mathbf{k}},\vec{k}%
\right) $ is \ defined by (\ref{conv}).
\end{definition}

Using $\mathcal{F}_{n,\zeta ,\vec{\xi}}^{\left( m\right) }$ in (\ref{Fm}) we
recast $\mathcal{F}^{\left( m\right) }\left( \mathbf{u}^{m}\right) $ in the
system (\ref{varcu})-(\ref{Fu}) as 
\begin{equation}
\mathcal{F}^{\left( m\right) }\left[ \mathbf{\hat{u}}_{1}\ldots ,\mathbf{%
\hat{u}}_{m}\right] \left( \mathbf{k},\tau \right) =\sum\nolimits_{n,\zeta ,%
\vec{\xi}}\mathcal{F}_{n,\zeta ,\vec{\xi}}^{\left( m\right) }\left[ \mathbf{%
\hat{u}}_{1}\ldots \mathbf{\hat{u}}_{m}\right] \left( \mathbf{k},\tau
\right) ,  \label{Tplusmin}
\end{equation}%
yielding the following system for the modal components $\mathbf{\hat{u}}%
_{n,\zeta }\left( \mathbf{k},\tau \right) $ as in (\ref{Pin}) 
\begin{equation}
\mathbf{\hat{u}}_{n,\zeta }\left( \mathbf{k},\tau \right)
=\sum\nolimits_{m\in \mathfrak{M}_{F}}\sum\nolimits_{\vec{\xi}\in \Xi ^{m}}%
\mathcal{F}_{n,\zeta ,\vec{\xi}}^{\left( m\right) }\left( \mathbf{\hat{u}}%
^{m}\right) \left( \mathbf{k},\tau \right) +\mathbf{\hat{h}}_{n,\zeta
}\left( \mathbf{k}\right) ,\ \left( n,\zeta \right) \in \Xi .  \label{equfa}
\end{equation}

\section{Wavepacket interaction system}

The wavepacket preservation property of the nonlinear evolutionary system in
any of its forms (\ref{difeqintr}), (\ref{difeqfou}), (\ref{varcu}), (\ref%
{eqF}), (\ref{equfa}) is not easy to see directly. It turns out though that
dynamics of wavepackets is well described by a system in a larger space $%
E^{2N}$ based on the original equation (\ref{varcu}) in the space $E$. We
call it \emph{wavepacket interaction system}, which is useful in three ways:
(i) the wavepacket preservation is quite easy to see and verify; (ii) it can
be used to prove the wavepacket preservation for the original nonlinear
problem; (iii) it can be used to study more subtle properties of the
original problem, such as NLS approximation. We start with the system (\ref%
{varcu}) where $\mathbf{\hat{h}}\left( \mathbf{k}\right) $ is a
multiwavepacket with a given $nk$-spectrum $S=\left\{ \left( \mathbf{k}%
_{\ast l},n_{l}\right) ,\ l=1,\ldots ,N\right\} $ as in (\ref{P0}) and $k$%
-spectrum $K_{S}=\left\{ \mathbf{k}_{\ast i},\ i=1,\ldots ,\left\vert
K_{S}\right\vert \right\} $ as in (\ref{K0}). Obviously, for any $\ l$ \ $%
\left( \mathbf{k}_{\ast l},n_{l}\right) =\left( \mathbf{k}_{\ast
i_{l}},n_{l}\right) $ with $i_{l}\leq \left\vert K_{S}\right\vert $ \ and \
indexing $i_{l}=l$ \ for $l\leq \left\vert K_{S}\right\vert $ according to (%
\ref{K0}).

When constructing the wavepacket interaction system it is convenient to have
relevant functions to be explicitly localized about the $k$-spectrum $K_{S}$
of the initial data. We implement that by making up the following cutoff
functions based on (\ref{j0}), (\ref{Psik}) 
\begin{equation}
\Psi _{i,\vartheta }\left( \mathbf{k}\right) =\Psi \left( \mathbf{k}%
,\vartheta \mathbf{k}_{\ast i},\beta ^{1-\epsilon }\right) =\Psi \left(
\beta ^{-\left( 1-\epsilon \right) }\left( \mathbf{k}-\vartheta \mathbf{k}%
_{\ast i}\right) \right) ,\ \mathbf{k}_{\ast i}\in K_{S},\ i=1,\ldots
,\left\vert K_{S}\right\vert ,\ \vartheta =\pm  \label{Psilz}
\end{equation}%
with $\epsilon $ as in Definition \ref{dwavepack} and $\beta >0$ small
enough to satisfy 
\begin{equation}
\beta ^{1/2}\leq \pi _{0},\text{ where }\pi _{0}=\pi _{0}\left( S\right) <%
\frac{1}{2}\min_{\mathbf{k}_{\ast i}\in K_{S}}\limfunc{dist}\left\{ \mathbf{k%
}_{\ast i},\sigma _{\mathrm{bc}}\right\} .  \label{kpi0}
\end{equation}%
In what follows we use notations from (\ref{setlam}) and 
\begin{equation}
\vec{l}=\left( l_{1},\ldots ,l_{m}\right) \in \left\{ 1,\ldots ,N\right\}
^{m},\ \vec{\vartheta}=\left( \vartheta ^{\prime },\ldots ,\vartheta
^{\left( m\right) }\right) \in \left\{ -1,1\right\} ^{m},\ \vec{\lambda}%
=\left( \vec{l},\vec{\vartheta}\right) \in \Lambda ^{m},  \label{narrow}
\end{equation}%
\begin{gather}
\vec{n}=\left( n_{1},\ldots ,n_{m}\right) \in \left\{ 1,\ldots ,J\right\}
^{m},\ \vec{\zeta}\in \left\{ -1,1\right\} ^{m},  \label{zetaar} \\
\vec{\xi}=\left( \vec{n}\ ,\vec{\zeta}\right) \in \Xi ^{m}\ ,\vec{k}=\left( 
\mathbf{\mathbf{k}}^{\prime },\ldots ,\mathbf{\mathbf{k}}^{\left( m\right)
}\right) \in \mathbb{R}^{m},\text{ where }\Xi ^{m}\text{ as in (\ref{Xi})}. 
\notag
\end{gather}%
Based on the above we introduce now the \emph{wavepacket interaction system} 
\begin{gather}
\mathbf{\hat{w}}_{l,\vartheta }\left( \mathbf{\cdot }\right) =\Psi \left( 
\mathbf{\cdot },\vartheta \mathbf{k}_{\ast i_{l}}\right) \Pi
_{n_{l},\vartheta }\left( \mathbf{\cdot }\right) \mathcal{F}\left(
\sum\nolimits_{\left( l^{\prime },\vartheta ^{\prime }\right) \in \Lambda }%
\mathbf{\hat{w}}_{l^{\prime },\vartheta ^{\prime }}\right) +\Psi \left( 
\mathbf{\cdot },\vartheta \mathbf{k}_{\ast i_{l}}\right) \Pi
_{n_{l},\vartheta }\left( \mathbf{\cdot }\right) \mathbf{\hat{h}},\left(
l,\vartheta \right) \in \Lambda ,  \label{sysloc} \\
\mathbf{\vec{w}}=\left( \mathbf{\hat{w}}_{1,+},\mathbf{\hat{w}}_{1,-},...,%
\mathbf{\hat{w}}_{N,+},\mathbf{\hat{w}}_{N,-}\right) \in E^{2N},\ \mathbf{%
\hat{w}}_{l,\vartheta }\in E,  \notag
\end{gather}%
with $\Psi \left( \mathbf{\cdot },\vartheta \mathbf{k}_{\ast i}\right) ,\Pi
_{n,\vartheta }$ being as in (\ref{Psilz}), (\ref{Pin}), $\mathcal{F}$
defined by (\ref{varcu}), and the norm in \ $E^{2N}$ defined based on (\ref%
{Elat}) by the formula 
\begin{equation}
\left\Vert \mathbf{\vec{w}}\right\Vert _{E^{2N}}=\sum\nolimits_{l,\vartheta
}\left\Vert \mathbf{\hat{w}}_{l,\vartheta }\right\Vert _{E},\ E=C\left( %
\left[ 0,\tau _{\ast }\right] ,L^{1}\right) .  \label{E2N}
\end{equation}%
We also use the following concise form of the wave interaction system (\ref%
{sysloc}) 
\begin{gather}
\mathbf{\vec{w}}=\mathcal{F}_{_{\Psi }}\left( \mathbf{\vec{w}}\right) +%
\mathbf{\vec{h}}_{_{\Psi }},\text{ where}  \label{syslocF} \\
\mathbf{\vec{h}}_{_{\Psi }}=\left( \Psi _{i_{1},+}\Pi _{n_{1},+}\mathbf{\hat{%
h}},\Psi _{i_{1},-}\Pi _{n_{1},-}\mathbf{\hat{h}},\ldots ,\Psi _{i_{N},+}\Pi
_{n_{N},+}\mathbf{\hat{h}},\Psi _{i_{N},-}\Pi _{n_{N},-}\mathbf{\hat{h}}%
\right) \in E^{2N}.  \notag
\end{gather}%
The following lemma is analogous to Lemmas \ref{Lemma bound}, \ref{Lemma
Flip}.

\begin{lemma}
\label{Lemma Flippr} Polynomial operator $\mathcal{F}_{\Psi }\left( \mathbf{%
\vec{w}}\right) $ is bounded in $E^{2N}$, $\mathcal{F}_{\Psi }\left( \mathbf{%
0}\right) =\mathbf{0}$, and it satisfies Lipschitz condition 
\begin{equation}
\left\Vert \mathcal{F}_{\Psi }\left( \mathbf{\vec{w}}_{1}\right) -\mathcal{F}%
_{\Psi }\left( \mathbf{\vec{w}}_{2}\right) \right\Vert _{E^{2N}}\leq C\tau
_{\ast }\left\Vert \mathbf{\vec{w}}_{1}-\mathbf{\vec{w}}_{2}\right\Vert
_{E^{2N}},  \label{lipFpsi}
\end{equation}%
where $C$ depends only on $C_{\chi }$ as in (\ref{chiCR}), on the degree of $%
\mathcal{F}$ and on $\left\Vert \mathbf{\vec{w}}_{1}\right\Vert
_{E^{2N}}+\left\Vert \mathbf{\vec{w}}_{2}\right\Vert _{E^{2N}}$, and it does
not depend on $\beta $ and $\varrho $.
\end{lemma}

\begin{proof}
We consider every operator $\mathcal{F}_{n,\zeta ,\vec{\xi}}^{\left(
m\right) }\left( \mathbf{\vec{w}}\right) $ defined by (\ref{Fm}) and prove
its boundedness and the Lipschitz property as in Lemma \ref{Lemma bound}\
using the inequality $\left\vert \exp \left\{ \mathrm{i}\phi _{n,\zeta ,\vec{%
\xi}}\frac{\tau _{1}}{\varrho }\right\} \right\vert \leq 1$ and inequalities
(\ref{j0}), (\ref{chiCR}). Note that the integration in $\tau _{1}$ yields
the factor $\tau _{\ast }$ and consequent summation \ with respect to $%
n,\zeta ,\vec{\xi}$ yields (\ref{lipFpsi}).
\end{proof}

Lemma \ref{Lemma Flippr}, the contraction principle as in Lemma \ref{Lemma
contr} and estimate (\ref{dtf1}) \ for the time derivative yield the
following statement.

\begin{theorem}
\label{Theorem existpr}Let $\left\Vert \mathbf{\vec{h}}_{_{\Psi
}}\right\Vert _{E^{2N}}\leq R.$ Then there exists $\tau _{\ast }>0$ and $\
R_{1}\left( R\right) $ such that equation (\ref{sysloc}) has a solution $%
\mathbf{\vec{w}}\in E^{2N}$ which satisfies 
\begin{equation}
\left\Vert \mathbf{\vec{w}}\right\Vert _{E^{2N}}+\left\Vert \partial _{\tau }%
\mathbf{\vec{w}}\right\Vert _{E^{2N}}\leq R_{1}\left( R\right)  \label{R1R}
\end{equation}%
and such a solution is unique.
\end{theorem}

\begin{lemma}
\label{Lemma wiswave}Every function $\mathbf{\hat{w}}_{l,\zeta }\left( 
\mathbf{k},\tau \right) $ corresponding to the solution of \ (\ref{syslocF}) 
$\ $from $\ E^{2N}$ is a wavepacket with $nk$-pair $\left( \mathbf{k}_{\ast
l},n_{l}\right) $ with the degree of regularity which can be any $s>0$.
\end{lemma}

\begin{proof}
Note that according to (\ref{Psilz}) and (\ref{syslocF}) the function 
\begin{equation*}
\mathbf{\hat{w}}_{l,\vartheta }\left( \mathbf{k},\tau \right) =\Psi \left( 
\mathbf{k},\vartheta \mathbf{k}_{\ast i_{l}},\beta ^{1-\epsilon }\right) \Pi
_{n_{l},\vartheta }\mathcal{F}\left( \mathbf{k},\tau \right) ,\ \left\Vert 
\mathcal{F}\left( \tau \right) \right\Vert _{L^{1}}\leq C,\ 0\leq \tau \leq
\tau _{\ast }
\end{equation*}%
\ involves the factor $\Psi _{l,\vartheta }\left( \mathbf{k}\right) =\Psi
\left( \beta ^{-\left( 1-\epsilon \right) }\left( \mathbf{k}-\vartheta 
\mathbf{k}_{\ast l}\right) \right) $ where $\epsilon $ is as in Definition %
\ref{dwavepack}. Hence,%
\begin{gather}
\Pi _{n,\vartheta ^{\prime }}\mathbf{\hat{w}}_{l,\vartheta }\left( \mathbf{k}%
,\tau \right) =0\text{ if }n\neq n_{l}\text{ or }\vartheta ^{\prime }\neq
\vartheta ,  \label{piweq0} \\
\mathbf{\hat{w}}_{l,\vartheta }\left( \mathbf{k},\tau \right) =\Psi \left( 
\mathbf{k},\vartheta \mathbf{k}_{\ast i_{l}},\beta ^{1-\epsilon }\right) 
\mathbf{\hat{w}}_{l,\vartheta }\left( \mathbf{k},\tau \right) ,\mathbf{\hat{w%
}}_{l,\vartheta }\left( \mathbf{k},\tau \right) =0\text{ if\ }\left\vert 
\mathbf{k}-\vartheta \mathbf{k}_{\ast l}\right\vert \geq \beta ^{1-\epsilon
},  \label{weq0}
\end{gather}%
Since 
\begin{equation}
\Psi \left( \mathbf{k},\vartheta \mathbf{k}_{\ast i_{l}},\beta ^{1-\epsilon
}\right) \Psi \left( \mathbf{k},\vartheta \mathbf{k}_{\ast i_{l}},\beta
^{1-\epsilon }/2\right) =\Psi \left( \mathbf{k},\vartheta \mathbf{k}_{\ast
i_{l}},\beta ^{1-\epsilon }\right) ,  \label{psipsi}
\end{equation}
Definition \ref{dwavepack} for $\mathbf{\hat{w}}_{l,\vartheta }$ is
satisfied with $D_{h}=0$ for any $s>0$ and $C^{\prime }=0$\ in (\ref{sourloc}%
).
\end{proof}

Now we would like to show that if $\mathbf{\hat{h}}$ is a multiwavepacket,
then the function 
\begin{equation}
\mathbf{\hat{w}}\left( \mathbf{k},\tau \right) =\sum\nolimits_{\left(
l,\vartheta \right) \in \Lambda }\mathbf{\hat{w}}_{l,\vartheta }\left( 
\mathbf{k},\tau \right) =\sum\nolimits_{\lambda \in \Lambda }\mathbf{\hat{w}}%
_{\lambda }\left( \mathbf{k},\tau \right)  \label{wkt}
\end{equation}%
constructed based on a solution of (\ref{syslocF}) is an approximate
solution of equation (\ref{eqF}) (see notation (\ref{setlam})). We will
follow the lines of \cite{BF7}. We introduce 
\begin{equation}
\Psi _{\infty }\left( \mathbf{k}\right) =1-\sum\nolimits_{\vartheta =\pm
}\sum\nolimits_{i=1}^{\left\vert K_{S}\right\vert }\Psi \left( \mathbf{k}%
,\vartheta \mathbf{k}_{\ast i}\right) =1-\sum\nolimits_{\vartheta =\pm
}\sum\nolimits_{\mathbf{k}_{\ast i}\in K_{S}}\Psi \left( \frac{\mathbf{k}%
-\vartheta \mathbf{k}_{\ast i}}{\beta ^{1-\epsilon }}\right) .
\label{Psiinf}
\end{equation}%
\ Expanding $m$-linear operator $\mathcal{F}^{\left( m\right) }\left( \left(
\sum_{l,\vartheta }\mathbf{\hat{w}}_{l,\vartheta }\right) ^{m}\right) $ and
using notations (\ref{setlam}), (\ref{lamprop}) we get%
\begin{gather}
\mathcal{F}^{\left( m\right) }\left( \left( \sum\nolimits_{l,\vartheta }%
\mathbf{\hat{w}}_{l,\vartheta }\right) ^{m}\right) =\sum\nolimits_{\vec{%
\lambda}\in \Lambda ^{m}}\mathcal{F}^{\left( m\right) }\left( \mathbf{\vec{w}%
}_{\vec{\lambda}}\right) ,\ \text{where}  \label{binom} \\
\mathbf{\vec{w}}_{\vec{\lambda}}=\mathbf{\hat{w}}_{\lambda _{1}}\ldots 
\mathbf{\hat{w}}_{\lambda _{m}},\ \vec{\lambda}=\left( \lambda _{1},\ldots
,\lambda _{m}\right) \in \Lambda ^{m}.  \label{wlam}
\end{gather}%
The next statement shows that (\ref{wkt}) defines an approximate solution to
integrated evolution equation (\ref{varcu}).

\begin{theorem}
\label{Theorem Dsmall}Let $\mathbf{\hat{h}}$ be a multi-wavepacket with
resonance invariant $nk$-spectrum $S$ with regularity degree $s$, $\mathbf{%
\vec{w}}$ be a solution of (\ref{syslocF}) \ and $\mathbf{\hat{w}}\left( 
\mathbf{k},\tau \right) $ be defined by (\ref{wkt}). Let \ 
\begin{equation}
\mathbf{D}\left( \mathbf{\hat{w}}\right) =\mathbf{\hat{w}}-\mathcal{F}\left( 
\mathbf{\hat{w}}\right) -\mathbf{\hat{h}}.  \label{Dw}
\end{equation}%
Then there exists $\beta _{0}>0$ such that we have the estimate 
\begin{equation}
\left\Vert \mathbf{D}\left( \mathbf{\hat{w}}\right) \right\Vert _{E}\leq
C\varrho +C\beta ^{s},\text{ if }0<\varrho \leq 1,\ \beta \leq \beta _{0}.
\label{Dwb}
\end{equation}
\end{theorem}

\begin{proof}
Let 
\begin{equation}
\mathcal{F}^{-}\left( \mathbf{\hat{w}}\right) =\left(
1-\sum\nolimits_{l,\vartheta }\Psi _{i_{l},\vartheta }\Pi _{n_{l},\vartheta
}\right) \mathcal{F}\left( \mathbf{\hat{w}}\right) ,\ \mathbf{\hat{h}}^{-}=%
\mathbf{\hat{h}}-\sum\nolimits_{l,\vartheta }\Psi _{i_{l},\vartheta }\Pi
_{n_{l},\vartheta }\mathbf{\hat{h}}.  \label{Fhmin}
\end{equation}%
Summation of (\ref{sysloc}) with respect to $l,\vartheta $ yields 
\begin{equation*}
\mathbf{\hat{w}}=\sum\nolimits_{l,\vartheta }\Psi _{i_{l},\vartheta }\Pi
_{n_{l},\vartheta }\mathcal{F}\left( \mathbf{\hat{w}}\right)
+\sum\nolimits_{l,\vartheta }\Psi _{i_{l},\vartheta }\Pi _{n_{l},\vartheta }%
\mathbf{\hat{h}}.
\end{equation*}%
Hence, from (\ref{sysloc}) and (\ref{Dw}) we obtain 
\begin{equation}
\mathbf{D}\left( \mathbf{\hat{w}}\right) =\mathbf{\hat{h}}^{-}-\mathcal{F}%
^{-}\left( \mathbf{\hat{w}}\right) .  \label{Dw1}
\end{equation}%
Using (\ref{hbold}) and (\ref{sourloc}) we consequently obtain 
\begin{equation*}
\left\Vert \Pi _{n_{l},\vartheta }\mathbf{\hat{h}}_{i}\right\Vert
_{L^{1}}\leq C\beta ^{s}\text{ if }n_{l}\neq n_{i};\ \left\Vert \Psi
_{i_{l},\vartheta }\mathbf{\hat{h}}_{i}\right\Vert _{L^{1}}\leq C\beta ^{s}%
\text{ if }\mathbf{k}_{\ast i_{l}}\neq \mathbf{k}_{\ast i},
\end{equation*}%
\begin{equation}
\left\Vert \mathbf{\hat{h}}^{-}\right\Vert _{E}\leq C_{1}\beta ^{s}.
\label{hbet}
\end{equation}%
Now, to show (\ref{Dwb}) it is sufficient to prove that 
\begin{equation}
\left\Vert \mathcal{F}^{-}\left( \mathbf{\hat{w}}\right) \right\Vert
_{E}\leq C_{2}\varrho .  \label{Fbet}
\end{equation}%
Obviously, 
\begin{equation}
\mathcal{F}^{-}\left( \mathbf{\hat{w}}\right) =\left(
1-\sum\nolimits_{l,\vartheta }\Psi _{i_{l},\vartheta }\Pi _{n_{l},\vartheta
}\right) \sum_{m}\mathcal{F}^{\left( m\right) }\left( \mathbf{\hat{w}}%
^{m}\right) .  \label{Fminw}
\end{equation}%
Note that 
\begin{equation}
\sum_{l,\vartheta }\Psi _{i_{l},\vartheta }\Pi _{n_{l},\vartheta
}=\sum\nolimits_{\vartheta =\pm }\sum\nolimits_{\left( n,k_{\ast }\right)
\in S}\Psi \left( \mathbf{\cdot },\vartheta \mathbf{k}_{\ast }\right) \Pi
_{n,\vartheta }.  \label{pspm}
\end{equation}%
Using (\ref{sumPi}) and ( \ref{Psiinf}) we consequently obtain 
\begin{gather}
\sum\nolimits_{\vartheta =\pm }\sum\nolimits_{\left( n,k_{\ast }\right) \in
\Sigma }\Psi \left( \mathbf{\cdot },\vartheta \mathbf{k}_{\ast }\right) \Pi
_{n,\vartheta }+\Psi _{\infty }=1,  \label{sum1} \\
\left( 1-\sum\nolimits_{l,\vartheta }\Psi _{i_{l},\vartheta }\Pi
_{n_{l},\vartheta }\right) =\Psi _{\infty }+\sum\nolimits_{\vartheta =\pm
}\sum\nolimits_{\left( n,k_{\ast }\right) \in \Sigma \setminus S}\Psi \left( 
\mathbf{\cdot },\vartheta \mathbf{k}_{\ast }\right) \Pi _{n,\vartheta }.
\label{psipi}
\end{gather}%
with $\Sigma $ defined in (\ref{ssp1}). Let us expand now $\mathcal{F}%
^{\left( m\right) }\left( \mathbf{\hat{w}}^{m}\right) $ using (\ref{binom}).
According to (\ref{Fminw}) and (\ref{psipi}) to prove (\ref{Fbet}) it is
sufficient to prove that \ for every string $\vec{\lambda}\in \Lambda ^{m}$ $%
\ $the following inequalities hold%
\begin{eqnarray}
\left\Vert \Psi _{\infty }\Pi _{n,\vartheta }\mathcal{F}^{\left( m\right)
}\left( \mathbf{\vec{w}}_{\vec{\lambda}}\right) \right\Vert &\leq
&C_{3}\varrho \text{ for }\left( n,\vartheta \right) \in \Lambda \text{, \
and }  \label{psiinfw} \\
\left\Vert \Psi \left( \mathbf{\cdot },\vartheta \mathbf{k}_{\ast }\right)
\Pi _{n,\vartheta }\mathcal{F}^{\left( m\right) }\left( \mathbf{\vec{w}}_{%
\vec{\lambda}}\right) \right\Vert &\leq &C_{3}\varrho ,\text{ \ if }\left( n,%
\mathbf{k}_{\ast }\right) \in \Sigma \setminus S\text{.}  \label{psiltw}
\end{eqnarray}%
We will use (\ref{piweq0}) and (\ref{weq0}) to obtain the above estimates. \
According to (\ref{Tplusmin})%
\begin{equation}
\mathcal{F}^{\left( m\right) }\left[ \mathbf{\vec{w}}_{\vec{\lambda}}\right]
\left( \mathbf{k},\tau \right) =\sum\nolimits_{n,\zeta }\sum\nolimits_{\vec{%
\xi}}\mathcal{F}_{n,\zeta ,\vec{\xi}}^{\left( m\right) }\left[ \mathbf{\hat{w%
}}_{\lambda _{1}}\ldots \mathbf{\hat{w}}_{\lambda _{m}}\right] \left( 
\mathbf{k},\tau \right) .  \label{Fbin}
\end{equation}%
Note that according to (\ref{piweq0}) if $\lambda _{i}=\left( l,\vartheta
^{\prime }\right) $ 
\begin{equation}
\mathbf{\hat{w}}_{\lambda _{i}}=\Pi _{n,\vartheta }\mathbf{\hat{w}}_{\lambda
_{i}},\ \text{if }n=n_{l}\text{ and }\vartheta ^{\prime }=\vartheta .
\label{wpiw}
\end{equation}%
Let us introduce notation 
\begin{equation}
\vec{n}\left( \vec{l}\right) =\left( n_{l_{1}},\ldots ,n_{l_{m}}\right) ,\ 
\vec{\xi}\left( \vec{\lambda}\right) =\left( \vec{n}\left( \vec{l}\right) ,%
\vec{\vartheta}\right) ,\text{\ for\ }\vec{\lambda}=\left( \vec{l},\vec{%
\vartheta}\right) \in \Lambda ^{m}.  \label{nl}
\end{equation}%
Since 
\begin{equation}
\Pi _{n^{\prime },\vartheta }\Pi _{n,\vartheta ^{\prime }}=0,\text{ if }%
n\neq n^{\prime }\text{ or }\vartheta ^{\prime }\neq \vartheta  \label{PnP'}
\end{equation}%
then (\ref{wpiw}) implies 
\begin{eqnarray}
\mathcal{F}_{n,\zeta ,\vec{\xi}}^{\left( m\right) }\left[ \mathbf{\hat{w}}%
_{\lambda _{1}}\ldots \mathbf{\hat{w}}_{\lambda _{m}}\right] &=&0\text{\ if }%
\vec{\xi}=\left( \vec{n},\vec{\zeta}\right) \neq \vec{\xi}\left( \vec{\lambda%
}\right) ,\text{ and, hence,}  \notag \\
\mathcal{F}^{\left( m\right) }\left[ \mathbf{\vec{w}}_{\vec{\lambda}}\right]
\left( \mathbf{k},\tau \right) &=&\sum\nolimits_{n,\zeta }\mathcal{F}%
_{n,\zeta ,\vec{\xi}\left( \vec{\lambda}\right) }^{\left( m\right) }\left[ 
\mathbf{\hat{w}}_{\lambda _{1}}\ldots \mathbf{\hat{w}}_{\lambda _{m}}\right]
\left( \mathbf{k},\tau \right) ,  \label{nzi}
\end{eqnarray}%
where we use notation (\ref{lamprop}), (\ref{nl}). Note also that 
\begin{equation}
\Pi _{n^{\prime },\vartheta }\mathcal{F}_{n,\zeta ,\vec{\xi}}^{\left(
m\right) }=0\text{\ if }n^{\prime }\neq n\text{\ or\ }\vartheta \neq \zeta ,
\label{Pinth0}
\end{equation}%
and, hence, we have nonzero $\Pi _{n^{\prime },\vartheta }\mathcal{F}%
_{n,\zeta ,\vec{\xi}}^{\left( m\right) }\left( \mathbf{\vec{w}}_{\vec{\lambda%
}}\right) $ only if 
\begin{equation}
\vec{\xi}=\vec{\xi}\left( \vec{\lambda}\right) ,\ n^{\prime }=n\text{, }%
\vartheta =\zeta .  \label{nonzero}
\end{equation}%
By (\ref{Fm})%
\begin{gather}
\mathcal{F}_{n,\zeta ,\vec{\xi}\left( \vec{\lambda}\right) }^{\left(
m\right) }\left( \mathbf{\vec{w}}_{\vec{\lambda}}\right) \left( \mathbf{k}%
,\tau \right) =\int_{0}^{\tau }\int_{\mathbb{D}_{m}}\exp \left\{ \mathrm{i}%
\phi _{n,\zeta ,\vec{\xi}\left( \vec{\lambda}\right) }\left( \mathbf{\mathbf{%
k}},\vec{k}\right) \frac{\tau _{1}}{\varrho }\right\}  \label{Fmw} \\
\chi _{n,\zeta ,\vec{\xi}\left( \vec{\lambda}\right) }^{\left( m\right)
}\left( \mathbf{\mathbf{k}},\vec{k}\right) \left[ \mathbf{\hat{w}}_{\lambda
_{1}}\left( \mathbf{k}^{\prime },\tau _{1}\right) ,\ldots ,\mathbf{\hat{w}}%
_{\lambda _{m}}\left( \mathbf{k}^{\left( m\right) }\left( \mathbf{\mathbf{k}}%
,\vec{k}\right) ,\tau _{1}\right) \right] \mathrm{\tilde{d}}^{\left(
m-1\right) d}\vec{k}\mathrm{d}\tau _{1},  \notag
\end{gather}%
Now we use (\ref{weq0}) and notice that according to the convolution
identity in (\ref{Fmintr}) 
\begin{equation}
\left\vert \mathbf{\hat{w}}_{\lambda _{1}}\left( \mathbf{k}^{\prime },\tau
_{1}\right) \right\vert \cdot \ldots \cdot \left\vert \mathbf{\hat{w}}%
_{\lambda _{m}}\left( \mathbf{k}^{\left( m\right) }\left( \mathbf{\mathbf{k}}%
,\vec{k}\right) ,\tau _{1}\right) \right\vert =0\text{ if }\left\vert 
\mathbf{k}-\sum\nolimits_{i}\vartheta _{i}\mathbf{k}_{\ast l_{i}}\right\vert
\geq m\beta ^{1-\epsilon }.  \label{prweq0}
\end{equation}%
Hence the integral (\ref{Fmw}) is nonzero only if $\left( \mathbf{k},\vec{k}%
\right) $ belongs to the set 
\begin{equation}
B_{\beta }=\left\{ \left( \mathbf{k},\vec{k}\right) :\left\vert \mathbf{k}%
^{\left( i\right) }-\vartheta _{i}\mathbf{k}_{\ast l_{i}}\right\vert \leq
\beta ^{1-\epsilon },\ i=1,\ldots ,m,\ \left\vert \mathbf{k}%
-\sum\nolimits_{i}\vartheta _{i}\mathbf{k}_{\ast l_{i}}\right\vert \leq
m\beta ^{1-\epsilon }\right\} .  \label{konly}
\end{equation}%
We will prove now that if $\left( n,\mathbf{k}_{\ast i}\right) \notin S$
then for small $\beta $ one of the following alternatives holds: \ 
\begin{gather}
\text{either }\Psi \left( \mathbf{\cdot },\vartheta \mathbf{k}_{\ast
i}\right) \Pi _{n^{\prime },\vartheta }\mathcal{F}_{n,\zeta ,\vec{\xi}%
}^{\left( m\right) }\left( \mathbf{\vec{w}}_{\vec{\lambda}}\right) =0
\label{alt1} \\
\text{or (\ref{nonzero}) holds and}\left\vert \phi _{n,\zeta ,\vec{\xi}%
}\left( \mathbf{\mathbf{k}},\vec{k}\right) \right\vert \geq c>0\text{ for }%
\left( \mathbf{k},\vec{k}\right) \in B_{\beta }.  \label{figrc}
\end{gather}%
\ Note then since $\phi _{n,\zeta ,\vec{\xi}}\left( \mathbf{\mathbf{k}},\vec{%
k}\right) $ is smooth then using notation (\ref{kapzel}) we get 
\begin{gather}
\left\vert \phi _{n,\zeta ,\vec{\xi}}\left( \mathbf{\mathbf{k}},\vec{k}%
\right) -\phi _{n^{\prime },\zeta ,\vec{\xi}}\left( \mathbf{\mathbf{k}}%
_{\ast \ast },\vec{k}_{\ast }\right) \right\vert \leq C\beta ^{1-\epsilon }%
\text{ for }\left( \mathbf{k},\vec{k}\right) \in B_{\beta },  \label{fiminfi}
\\
\vec{\vartheta}=\left( \vartheta _{1},\ldots ,\vartheta _{m}\right) ,\ 
\mathbf{\mathbf{k}}_{\ast \ast }=\zeta \sum\nolimits_{i}\vartheta _{i}%
\mathbf{k}_{\ast l_{i}}=\zeta \varkappa _{m}\left( \vec{\vartheta},\vec{l}%
\right) ,  \notag
\end{gather}%
Hence the alternative (\ref{figrc}) holds if%
\begin{equation}
\phi _{n,\zeta ,\vec{\xi}}\left( \mathbf{\mathbf{k}}_{\ast \ast },\vec{k}%
_{\ast }\right) \neq 0,  \label{fineq0}
\end{equation}%
\ and, consequently, it suffices to prove that either (\ref{alt1}) or (\ref%
{fineq0}) holds. Combining (\ref{konly}) with $\Psi \left( \mathbf{k}%
,\vartheta \mathbf{k}_{\ast i}\right) =0$ for $\left\vert \mathbf{k}%
-\vartheta \mathbf{k}_{\ast i}\right\vert \geq \beta ^{1-\epsilon }$ we find
that $\Psi _{i,\vartheta }\mathcal{F}^{\left( m\right) }\left[ \mathbf{\vec{w%
}}_{\vec{\lambda}}\right] $ can be non-zero for small $\beta $ only in a
small neighborhood of a point $\zeta \varkappa _{m}\left( \vec{\vartheta},%
\vec{l}\right) \in \left[ S\right] _{K,\text{out}}$, and that is possible
only if 
\begin{equation}
\mathbf{\mathbf{k}}_{\ast \ast }=\zeta \varkappa _{m}\left( \vec{\vartheta},%
\vec{l}\right) =\vartheta \mathbf{k}_{\ast i},\ \mathbf{k}_{\ast i}\in K_{S}.
\label{kstar}
\end{equation}%
Let us show that the equality 
\begin{equation}
\phi _{n,\zeta ,\vec{\xi}}\left( \mathbf{\mathbf{k}}_{\ast \ast },\vec{k}%
_{\ast }\right) =0  \label{fiopp}
\end{equation}%
is impossible for $\mathbf{\mathbf{k}}_{\ast \ast }$ as in (\ref{kstar}) and 
$n^{\prime }=n$ as in (\ref{Pinth0}), keeping in mind that $\left( n,\mathbf{%
\mathbf{k}}_{\ast i}\right) \notin S$. It follows from (\ref{Omzet}) and (%
\ref{phim}) that the equation (\ref{fiopp}) has the form of the resonance
equation (\ref{Omeq0}). Since $nk$-spectrum $S$ is resonance invariant, in
view of Definition \ref{Definition omclos} the resonance equation (\ref%
{fiopp}) may have a solution only if $\mathbf{\mathbf{k}}_{\ast \ast }=%
\mathbf{\mathbf{k}}_{\ast i}$, $i=i_{l}$, $n=n_{l}$, with $\left( n_{l},%
\mathbf{\mathbf{k}}_{\ast i_{l}}\right) \in S$. Since $\left( n,\mathbf{%
\mathbf{k}}_{\ast i}\right) \notin S$ that implies (\ref{fiopp}) does not
have a solution and, hence, (\ref{fineq0}) holds when $\left( n,\mathbf{%
\mathbf{k}}_{\ast i}\right) \notin S$. Notice that (\ref{R1R}) yields the
following bounds 
\begin{equation}
\left\Vert \mathbf{\hat{w}}_{\lambda _{i}}\right\Vert _{E}\leq R_{1},\
\left\Vert \partial _{\tau }\mathbf{\hat{w}}_{\lambda _{i}}\right\Vert
_{E}\leq C.  \label{R1R1}
\end{equation}%
These bounds combined with Lemma \ref{Lemma intbyparts}, proven below, imply
that if (\ref{fineq0}) holds then (\ref{psiltw}) holds. Now let us turn to (%
\ref{psiinfw}). According to ( \ref{Psiinf}) and (\ref{prweq0}) the term $%
\Psi _{\infty }\Pi _{n^{\prime },\vartheta }\mathcal{F}^{\left( m\right)
}\left( \mathbf{\vec{w}}_{\vec{\lambda}}\right) $ can be non-zero only if $%
\zeta \varkappa _{m}\left( \vec{\lambda}\right) =\mathbf{\mathbf{k}}_{\ast
\ast }\notin K_{S}$. Since $nk$-spectrum $S$ is resonance invariant we
conclude as above that inequality (\ref{fineq0}) holds in this case as well.
The fact that the set of all $\varkappa _{m}\left( \vec{\lambda}\right) $ is
finite, combined with inequality (\ref{fineq0}), imply (\ref{figrc}) for
sufficiently small $\beta $. Using Lemma \ref{Lemma intbyparts} as above we
derive (\ref{psiinfw}). Hence, all terms in the expansion (\ref{Fminw}) are
either zero or satisfy (\ref{psiinfw}) or (\ref{psiltw}) implying
consequently (\ref{Fbet}) and (\ref{Dwb}).
\end{proof}

Here is the lemma used in the above proof.

\begin{lemma}
\label{Lemma intbyparts}Let assume that%
\begin{gather}
\left\vert \Psi _{i,\vartheta ^{\prime }}\Pi _{n^{\prime },\zeta }\chi
_{n,\zeta ,\vec{\xi}}^{\left( m\right) }\left( \mathbf{\mathbf{k}},\vec{k}%
\right) \left[ \mathbf{\hat{w}}_{\lambda _{1}}\left( \mathbf{k}^{\prime
},\tau _{1}\right) ,\ldots ,\mathbf{\hat{w}}_{\lambda _{m}}\left( \mathbf{k}%
^{\left( m\right) }\left( \mathbf{\mathbf{k}},\vec{k}\right) ,\tau
_{1}\right) \right] \right\vert =0\text{ for }\left( \mathbf{\mathbf{k}},%
\vec{k}\right) \in B_{\beta },  \notag \\
\text{and }\left\vert \phi _{n,\zeta ,\vec{\xi}}\left( \mathbf{\mathbf{k}},%
\vec{k}\right) \right\vert \geq \omega _{\ast }>0\text{ for }\left( \mathbf{%
\mathbf{k}},\vec{k}\right) \notin B_{\beta },\text{ where }B_{\beta }\text{
as in (\ref{konly}).}  \label{figrom}
\end{gather}%
Then%
\begin{gather}
\left\Vert \Psi \left( \mathbf{\cdot },\vartheta ^{\prime }\mathbf{k}_{\ast
i}\right) \Pi _{n^{\prime },\zeta }\mathcal{F}_{n,\zeta ,\vec{\xi}}^{\left(
m\right) }\left( \mathbf{\vec{w}}_{\vec{\lambda}}\right) \right\Vert _{E}\leq
\label{estrho} \\
\frac{4\varrho }{\omega _{\ast }}\left\Vert \chi ^{\left( m\right)
}\right\Vert \dprod\nolimits_{j}\left\Vert \mathbf{\hat{w}}_{\lambda
_{j}}\right\Vert _{E}+\frac{2\varrho \tau _{\ast }}{\omega _{\ast }}%
\left\Vert \chi ^{\left( m\right) }\right\Vert \sum\nolimits_{i}\left\Vert
\partial _{\tau }\mathbf{\hat{w}}_{\lambda _{i}}\right\Vert
_{E}\dprod\nolimits_{j\neq i}\left\Vert \mathbf{\hat{w}}_{\lambda
_{j}}\right\Vert _{E}.  \notag
\end{gather}
\end{lemma}

\begin{proof}
Notice that the oscillatory factor in (\ref{Fm}) equals to 
\begin{equation*}
\exp \left\{ \mathrm{i}\phi \left( \mathbf{\mathbf{k}},\vec{k}\right) \frac{%
\tau _{1}}{\varrho }\right\} =\frac{\varrho }{\mathrm{i}\phi \left( \mathbf{%
\mathbf{k}},\vec{k}\right) }\partial _{\tau _{1}}\exp \left\{ \mathrm{i}\phi
\left( \mathbf{\mathbf{k}},\vec{k}\right) \frac{\tau _{1}}{\varrho }\right\}
.
\end{equation*}%
Denoting $\phi _{n,\zeta ,\vec{\xi}}=\phi $, $\Psi _{i,\vartheta ^{\prime
}}\Pi _{n^{\prime },\zeta }\chi _{n,\zeta ,\vec{\xi}}^{\left( m\right)
}=\chi _{\vec{\eta}}^{\left( m\right) }$ and integrating (\ref{Fm})\ by
parts with respect to $\tau _{1}$ we obtain%
\begin{gather}
\Psi \left( \mathbf{k},\vartheta ^{\prime }\mathbf{k}_{\ast i}\right) \Pi
_{n^{\prime },\zeta }\mathcal{F}_{n,\zeta ,\vec{\xi}}^{\left( m\right)
}\left( \mathbf{\vec{w}}_{\vec{\lambda}}\right) \left( \mathbf{k},\tau
\right) =  \label{NFMintbp} \\
\int_{B}\Psi \left( \mathbf{k},\vartheta ^{\prime }\mathbf{k}_{\ast
i}\right) \frac{\varrho \mathrm{e}^{\mathrm{i}\phi \left( \mathbf{\mathbf{k}}%
,\vec{k}\right) \frac{\tau }{\varrho }}}{\mathrm{i}\phi \left( \mathbf{%
\mathbf{k}},\vec{k}\right) }\chi _{\vec{\eta}}^{\left( m\right) }\left( 
\mathbf{\mathbf{k}},\vec{k}\right) \mathbf{\hat{w}}_{\lambda _{1}}\left( 
\mathbf{k}^{\prime },\tau \right) \ldots \mathbf{\hat{w}}_{\lambda
_{m}}\left( \mathbf{k}^{\left( m\right) }\left( \mathbf{k},\vec{k}\right)
,\tau \right) \,\mathrm{\tilde{d}}^{\left( m-1\right) d}\vec{k}  \notag \\
-\int_{B}\Psi \left( \mathbf{k},\vartheta ^{\prime }\mathbf{k}_{\ast
i}\right) \frac{\varrho }{\mathrm{i}\phi \left( \mathbf{\mathbf{k}},\vec{k}%
\right) }\chi _{\vec{\eta}}^{\left( m\right) }\left( \mathbf{\mathbf{k}},%
\vec{k}\right) \mathbf{\hat{w}}_{\lambda _{1}}\left( \mathbf{k}^{\prime
},0\right) \ldots \mathbf{\hat{w}}_{\lambda _{m}}\left( \mathbf{k}^{\left(
m\right) }\left( \mathbf{k},\vec{k}\right) ,0\right) \,\mathrm{\tilde{d}}%
^{\left( m-1\right) d}\vec{k}  \notag \\
-\int_{0}^{\tau }\int_{B}\Psi \left( \mathbf{k},\vartheta ^{\prime }\mathbf{k%
}_{\ast i}\right) \frac{\varrho \mathrm{e}^{\mathrm{i}\phi \left( \mathbf{%
\mathbf{k}},\vec{k}\right) \frac{\tau _{1}}{\varrho }}}{\mathrm{i}\phi
\left( \mathbf{\mathbf{k}},\vec{k}\right) }\chi _{\vec{\eta}}^{\left(
m\right) }\left( \mathbf{\mathbf{k}},\vec{k}\right) \partial _{\tau _{1}}%
\left[ \mathbf{\hat{w}}_{\lambda _{1}}\left( \mathbf{k}^{\prime }\right)
\ldots \mathbf{\hat{w}}_{\lambda _{m}}\left( \mathbf{k}^{\left( m\right)
}\left( \mathbf{k},\vec{k}\right) \right) \right] \,\mathrm{\tilde{d}}%
^{\left( m-1\right) d}\vec{k}d\tau _{1},  \notag
\end{gather}%
where $B$ is the set of $\mathbf{k}^{\left( i\right) }$ for which (\ref%
{konly}) holds. The relations (\ref{chiCR}) and (\ref{j0}) imply $\left\vert
\chi _{\vec{\eta}}^{\left( m\right) }\left( \mathbf{\mathbf{k}},\vec{k}%
\right) \right\vert \leq \left\Vert \chi ^{\left( m\right) }\right\Vert $.
Using then (\ref{figrom}), the Leibnitz formula, (\ref{R1R}) and (\ref{Yconv}%
) we obtain (\ref{estrho}).
\end{proof}

The main result of this subsection is the next theorem which, when combined
with Lemma \ref{Lemma wiswave}, implies the wavepacket preservation, namely
that the solution $\mathbf{\hat{u}}_{n,\vartheta }\left( \mathbf{k},\tau
\right) $ of (\ref{equfa}) is a multi-wavepacket for all $\tau \in \left[
0,\tau _{\ast }\right] $.

\begin{theorem}
\label{Theorem uminw} Assume that conditions of Theorem \ref{Theorem Dsmall}
are fulfilled. Let $\mathbf{\hat{u}}_{n,\vartheta }\left( \mathbf{k},\tau
\right) $ for $n=n_{l}$ and $\mathbf{\hat{w}}_{l,\vartheta }\left( \mathbf{k}%
,\tau \right) \ $be the solutions to respective systems (\ref{equfa}) and (%
\ref{sysloc}), $\mathbf{\hat{w}}$ be defined by (\ref{wkt}). Then for
sufficiently small $\beta _{0}>0$ we have 
\begin{equation}
\left\Vert \mathbf{\hat{u}}_{n_{l},\vartheta }-\Pi _{n_{l},\vartheta }%
\mathbf{\hat{w}}\right\Vert _{E}\leq C\varrho +C^{\prime }\beta ^{s},\
0<\beta \leq \beta _{0},\;l=1,...,N.  \label{uminw}
\end{equation}
\end{theorem}

\begin{proof}
Note that $\mathbf{\hat{u}}_{n,\vartheta }=\Pi _{n,\vartheta }\mathbf{\hat{u}%
}$ where $\mathbf{\hat{u}}$ is a solution of (\ref{varcu}) and, according to
Theorem \ref{Theorem exist}, $\left\Vert \mathbf{\hat{u}}\right\Vert
_{E}\leq 2R$. Comparing the equations (\ref{varcu}) and (\ref{Dw}) , which
are $\mathbf{\hat{u}}=\mathcal{F}\left( \mathbf{\hat{u}}\right) +\mathbf{%
\hat{h}}$ and $\mathbf{\hat{w}}=\mathcal{F}\left( \mathbf{\hat{w}}\right) +%
\mathbf{\hat{h}+D}\left( \mathbf{\hat{w}}\right) $, we find that Lemma \ref%
{Lemma contr} can be applied. Then we notice that by Lemma \ref{Lemma Flip} $%
\mathcal{F}$ has the Lipschitz constant $C_{F}\tau _{\ast }$ for such $%
\mathbf{\hat{u}}$. Taking $C_{F}\tau _{\ast }<1$ as in Theorem \ref{Theorem
exist} we obtain (\ref{uminw}) from (\ref{iminu0}).
\end{proof}

Notice that Theorem \ref{Theorem sumwave} \ is a direct corollary of Theorem %
\ref{Theorem uminw} and Lemma \ref{Lemma wiswave}.

Analogous statement is proven in \cite{BF7} \ for parameter-dependent
equations (\ref{difeqfou}) with $\mathbf{\hat{F}}\left( \mathbf{\hat{U}}%
\right) =\mathbf{\hat{F}}\left( \mathbf{\hat{U}},\varrho \right) $.

The following theorem shows that any multi-wavepacket solution to (\ref%
{varcu}) yields a solution to the wavepacket interaction system (\ref{sysloc}%
).

\begin{theorem}
\label{Theorem necessary}Let $\mathbf{\hat{u}}\left( \mathbf{k},\tau \right) 
$ be a solution of (\ref{varcu}) and assume that $\mathbf{\hat{u}}\left( 
\mathbf{k},\tau \right) $ and $\mathbf{\hat{h}}\left( \mathbf{k}\right) $
are multiwavepackets with $nk$-spectrum $S=\left\{ \left( n_{l},\mathbf{k}%
_{\ast l}\right) \text{, }l=1,\ldots ,N\right\} $ and the regularity degree $%
s$. Let also $\Psi _{i_{l},\vartheta }=\Psi _{i_{l},\vartheta }$ be defined
by (\ref{Psilz}). Then the functions $\mathbf{\hat{w}}_{l,\vartheta
}^{\prime }\left( \mathbf{k},\tau \right) =\Psi _{i_{l},\vartheta }\Pi
_{n_{l},\vartheta }\mathbf{\hat{u}}\left( \mathbf{k},\tau \right) $ are a
solution to the system (\ref{sysloc}) with $\mathbf{\hat{h}}\left( \mathbf{k}%
\right) $ replaced by $\mathbf{\hat{h}}^{\prime }\left( \mathbf{k},\tau
\right) $ satisfying 
\begin{equation}
\left\Vert \mathbf{\hat{h}}\left( \mathbf{k}\right) -\mathbf{\hat{h}}%
^{\prime }\left( \mathbf{k},\tau \right) \right\Vert _{L^{1}}\leq C\beta
^{s},\ 0\leq \tau \leq \tau _{\ast },  \label{hprime}
\end{equation}%
and if $\mathbf{\hat{w}}_{l,\vartheta }$ are solutions of (\ref{sysloc})
with original $\mathbf{\hat{h}}\left( \mathbf{k}\right) $ we have the
inequality 
\begin{equation}
\left\Vert \mathbf{\hat{w}}_{l,\vartheta }^{\prime }\left( \mathbf{k},\tau
\right) -\mathbf{\hat{w}}_{l,\vartheta }\right\Vert _{L^{1}}\leq C\beta
^{s},\ 0\leq \tau \leq \tau _{\ast }.  \label{wwp}
\end{equation}
\end{theorem}

\begin{proof}
Multiplying (\ref{varcu}) by $\Psi _{i_{l},\vartheta }\Pi _{n_{l},\vartheta
} $ we get%
\begin{equation}
\mathbf{\hat{w}}_{l,\vartheta }^{\prime }=\Psi \left( \mathbf{\cdot }%
,\vartheta \mathbf{k}_{\ast i_{l}}\right) \Pi _{n_{l},\vartheta }\mathcal{F}%
\left( \mathbf{\hat{u}}\right) \left( \mathbf{k},\tau \right) +\Psi \left( 
\mathbf{\cdot },\vartheta \mathbf{k}_{\ast i_{l}}\right) \Pi
_{n_{l},\vartheta }\mathbf{\hat{h}}\left( \mathbf{k}\right) ,\ \mathbf{\hat{w%
}}_{l,\vartheta }^{\prime }=\Psi \left( \mathbf{\cdot },\vartheta \mathbf{k}%
_{\ast i_{l}}\right) \Pi _{n_{l},\vartheta }\mathbf{\hat{u}}.  \label{wprime}
\end{equation}%
Since $\mathbf{\hat{u}}\left( \mathbf{k},\tau \right) $ is a multiwavepacket
with regularity $s$ we have 
\begin{equation}
\left\Vert \mathbf{\hat{u}}\left( \mathbf{\cdot },\tau \right) -\mathbf{\hat{%
w}}^{\prime }\left( \mathbf{\cdot },\tau \right) \right\Vert _{L^{1}}\leq
C_{\epsilon }\beta ^{s}\text{ \ where \ }\mathbf{\hat{w}}^{\prime }\left( 
\mathbf{\cdot },\tau \right) =\sum\nolimits_{l,\vartheta }\Psi \left( 
\mathbf{\cdot },\vartheta \mathbf{k}_{\ast i_{l}}\right) \mathbf{\hat{u}}%
\left( \mathbf{\cdot },\tau \right) .  \label{uwprime}
\end{equation}%
Let us recast (\ref{wprime}) in the form 
\begin{gather}
\mathbf{\hat{w}}_{l,\vartheta }^{\prime }=\Psi \left( \mathbf{\cdot }%
,\vartheta \mathbf{k}_{\ast i_{l}}\right) \Pi _{n_{l},\vartheta }\mathcal{F}%
\left( \mathbf{\hat{w}}^{\prime }\right) \left( \mathbf{k},\tau \right)
+\Psi \left( \mathbf{\cdot },\vartheta \mathbf{k}_{\ast i_{l}}\right) \Pi
_{n_{l},\vartheta }\left[ \mathbf{\hat{h}}\left( \mathbf{k}\right) +\mathbf{%
\hat{h}}^{\prime \prime }\left( \mathbf{k},\tau \right) \right] ,
\label{wprime1} \\
\mathbf{\hat{h}}^{\prime \prime }\left( \mathbf{k},\tau \right) =\left[ 
\mathcal{F}\left( \mathbf{\hat{u}}\right) -\mathcal{F}\left( \mathbf{\hat{w}}%
^{\prime }\right) \right] \left( \mathbf{k},\tau \right) .  \notag
\end{gather}%
Denoting $\mathbf{\hat{h}}\left( \mathbf{k}\right) +\mathbf{\hat{h}}^{\prime
\prime }\left( \mathbf{k},\tau \right) =$ $\mathbf{\hat{h}}^{\prime }\left( 
\mathbf{k},\tau \right) $ we observe that (\ref{wprime1}) has the form of (%
\ref{sysloc}) with $\mathbf{\hat{h}}\left( \mathbf{k}\right) $ replaced by $%
\mathbf{\hat{h}}^{\prime }\left( \mathbf{k},\tau \right) $. Inequality (\ref%
{hprime}) follows then from (\ref{uwprime}) and (\ref{Flip}). Using Lemma %
\ref{Lemma contr} we obtain (\ref{wwp}).
\end{proof}

\section{Reduction of wavepacket interaction system to an averaged
interaction system}

\ Our goal in this section is to substitute the wavepacket interaction
system (\ref{sysloc}) with a simpler \emph{averaged interaction system}
which describes the evolution of wavepackets with the same accuracy but has
a simpler nonlinearity, and we follow here the approach developed in \cite%
{BF7}. \ The reduction is a generalization of the classical averaging
principle to the case of continuous spectrum, see \cite{BF7} for a
discussion and further simplification of the averaged interaction system. In
the present paper we do not need the further simplification to a minimal
interaction system leading to a system of NLS-type equations which is done
in \cite{BF7}.

\subsection{Time averaged wavepacket interaction system}

Here we modify the wavepacket interaction system (\ref{sysloc}),
substituting its nonlinearity with another one obtained by the time
averaging, and prove that this substitution produces a small error of order $%
\varrho $. As the first step we recast (\ref{sysloc}) in a slightly
different form by using expansions (\ref{binom}), (\ref{Fbin}) together with
(\ref{nzi})\ and (\ref{Pinth0}) and writing the nonlinearity in the equation
(\ref{sysloc}) in the form 
\begin{eqnarray}
\Psi \left( \mathbf{\cdot },\vartheta \mathbf{k}_{\ast i_{l}}\right) \Pi
_{n_{l},\vartheta }\mathcal{F}\left( \mathbf{\cdot },\tau \right)
&=&\sum\nolimits_{m\in \mathfrak{M}_{F}}\sum\nolimits_{\vec{\lambda}\in
\Lambda ^{m}}\Psi \left( \mathbf{\cdot },\vartheta \mathbf{k}_{\ast
i_{l}}\right) \mathcal{F}_{n_{l},\vartheta ,\vec{\xi}\left( \vec{\lambda}%
\right) }^{\left( m\right) }\left( \mathbf{\vec{w}}_{\vec{\lambda}}\right) ,%
\text{\ }\vec{\lambda}=\left( \vec{l},\vec{\zeta}\right) ,  \label{Fmksi} \\
\mathcal{F}_{n_{l},\vartheta ,\vec{\xi}\left( \vec{\lambda}\right) }^{\left(
m\right) }\left( \mathbf{\vec{w}}_{\vec{\lambda}}\right) \left( \mathbf{k}%
,\tau \right) &=&\left. \mathcal{F}_{n,\zeta ,\vec{n},\vec{\zeta}}^{\left(
m\right) }\left[ \mathbf{\hat{w}}_{\lambda _{1}}\ldots \mathbf{\hat{w}}%
_{\lambda _{m}}\right] \left( \mathbf{k},\tau \right) \right\vert _{\vec{n}=%
\vec{n}\left( \vec{l}\right) ,\text{\ }\left( n,\zeta \right) =\left(
n_{l},\vartheta \right) },  \label{Fmksi1}
\end{eqnarray}%
with $\mathcal{F}_{n,\zeta ,\vec{n},\vec{\zeta}}^{\left( m\right) }$ as in (%
\ref{Fm}) and $\vec{n}\left( \vec{l}\right) \ $as in (\ref{nl}), and we call 
$\mathcal{F}_{n_{l},\vartheta ,\vec{\xi}\left( \vec{\lambda}\right)
}^{\left( m\right) }\left( \mathbf{\vec{w}}_{\vec{\lambda}}\right) $ a \emph{%
decorated monomial} $\mathcal{F}_{n_{l},\vartheta ,\vec{\xi}\left( \vec{%
\lambda}\right) }^{\left( m\right) }$ evaluated at $\mathbf{\vec{w}}_{\vec{%
\lambda}}$. Consequently, the wavepacket interaction system (\ref{sysloc})
can be written in an equivalent form 
\begin{equation}
\mathbf{\hat{w}}_{l,\vartheta }=\sum_{m\in \mathfrak{M}_{F}}\sum_{\vec{%
\lambda}\in \Lambda ^{m}}\Psi \left( \mathbf{\cdot },\vartheta \mathbf{k}%
_{\ast i_{l}}\right) \mathcal{F}_{n_{l},\vartheta ,\vec{\xi}\left( \vec{%
\lambda}\right) }^{\left( m\right) }\left( \mathbf{\vec{w}}_{\vec{\lambda}%
}\right) +\Psi \left( \mathbf{\cdot },\vartheta \mathbf{k}_{\ast
i_{l}}\right) \Pi _{n_{l},\vartheta }\mathbf{\hat{h}},\ l=1,\ldots N,\
\vartheta =\pm .  \label{sysloc1}
\end{equation}%
The construction of the above mentioned time averaged equation reduces to
discarding certain terms in the original system (\ref{sysloc1}). First we
introduce the following sets of indices related to the resonance equation (%
\ref{Omeq0}) and $\Omega _{m}$ defined by (\ref{Omzet}): 
\begin{equation}
\Lambda _{n_{l},\vartheta }^{m}=\left\{ \vec{\lambda}=\left( \vec{l},\vec{%
\zeta}\right) \in \Lambda ^{m}:\Omega _{m}\left( \vartheta ,n_{l},\vec{%
\lambda}\right) =0\right\} ,  \label{resset}
\end{equation}%
and then the \emph{time-averaged nonlinearity }$\mathcal{F}_{\limfunc{av}}$%
\emph{\ by} 
\begin{equation}
\mathcal{F}_{\limfunc{av},n_{l},\vartheta }\left( \mathbf{\vec{w}}\right)
=\sum\nolimits_{m\in \mathfrak{M}_{F}}\mathcal{F}_{n_{l},\vartheta }^{\left(
m\right) },\ \mathcal{F}_{n_{l},\vartheta }^{\left( m\right)
}=\sum\nolimits_{\vec{\lambda}\in \Lambda _{n_{l},\vartheta }^{m}}\mathcal{F}%
_{n_{l},\vartheta ,\vec{\xi}\left( \vec{\lambda}\right) }^{\left( m\right)
}\left( \mathbf{\vec{w}}_{\vec{\lambda}}\right) .  \label{Fav}
\end{equation}%
where $\mathcal{F}_{n_{l},\vartheta ,\vec{\xi}\left( \vec{\lambda}\right)
}^{\left( m\right) }$\ are defined in (\ref{Fmksi1}). \ 

\begin{remark}
\label{Remark averaging}Note that \ the nonlinearity $\mathcal{F}_{\limfunc{%
av},n_{l},\vartheta }^{\left( m\right) }\left( \mathbf{\vec{w}}\right) $ can
be obtained from $\mathcal{F}_{n_{l},\vartheta }^{\left( m\right) }$ by an
averaging \ formula using \emph{averaging operator} $A_{T}$ acting on
polynomial functions $F:\left( \mathbb{C}^{2}\right) ^{N}\rightarrow \left( 
\mathbb{C}^{2}\right) ^{N}$ as follows 
\begin{equation}
\left( A_{T}F\right) _{j,\zeta }=\frac{1}{T}\int_{0}^{T}\mathrm{e}^{-\mathrm{%
i}\zeta \phi _{j}t}F_{j,\zeta }\left( \mathrm{e}^{\mathrm{i}\phi
_{1}t}u_{1,+},\mathrm{e}^{-\mathrm{i}\phi _{1}t}u_{1,-},\ldots ,\mathrm{e}^{%
\mathrm{i}\phi _{N}t}u_{N,+},\mathrm{e}^{-\mathrm{i}\phi
_{N}t}u_{N,-}\right) \mathrm{d}t.  \label{cana1}
\end{equation}%
Using this averaging we define for any polynomial nonlinearity $G:\left( 
\mathbb{C}^{2}\right) ^{N}\rightarrow \left( \mathbb{C}^{2}\right) ^{N}$
averaged polynomial 
\begin{equation}
G_{\limfunc{av},j,\zeta }\left( \vec{u}\right) =\lim_{T\rightarrow \infty
}\left( A_{T}G\right) _{j,\zeta }\left( \vec{u}\right) .  \label{Gav}
\end{equation}%
If frequencies $\phi _{j}$ in (\ref{cana1}) are generic, $G_{\limfunc{av}%
,j,\zeta }\left( \vec{u}\right) $ is always a universal nonlinearity. \ Note
that $\mathcal{F}_{\limfunc{av},n_{l},\vartheta }\left( \mathbf{\vec{w}}%
\right) $ defined by (\ref{Fav}) can be obtained by the formula (\ref{Gav})
where $A_{T}\mathbf{\ }$is defined by formula (\ref{cana1}) with frequencies 
$\phi _{j}=\omega _{n_{j}}\left( \mathbf{k}_{\ast i_{j}}\right) $ (it may be
conditionally universal if the frequencies $\phi _{j}$ are subjected to a
condition of the form (\ref{bnc}), see the following subsection for details,
in particular for definitions of universal and conditionally universal
nonlinearities).
\end{remark}

Finally, we introduce the wave interaction system with time-averaged
nonlinearity as follows: 
\begin{equation}
\mathbf{\hat{v}}_{l,\vartheta }=\Psi \left( \mathbf{\cdot },\vartheta 
\mathbf{k}_{\ast i_{l}}\right) \mathcal{F}_{\limfunc{av},n_{l},\vartheta
}\left( \mathbf{\vec{v}}\right) +\Psi \left( \mathbf{\cdot },\vartheta 
\mathbf{k}_{\ast i_{l}}\right) \Pi _{n_{l},\vartheta }\mathbf{\hat{h}},\
l=1,\ldots N,\vartheta =\pm .  \label{eqav}
\end{equation}%
Similarly to (\ref{syslocF}) we recast this system concisely as 
\begin{equation}
\mathbf{\vec{v}}=\mathcal{F}_{\text{$\limfunc{av}$},\Psi }\left( \mathbf{%
\vec{v}}\right) +\mathbf{\vec{h}}_{_{\Psi }}.  \label{eqavF}
\end{equation}%
The following lemma is analogous to Lemmas \ref{Lemma Flippr}, \ref{Lemma
Flip}.

\begin{lemma}
\label{Lemma Flipprav} Operator\ $\mathcal{F}_{\limfunc{av},\Psi }\left( 
\mathbf{\vec{v}}\right) $ is bounded for bounded $\mathbf{\vec{v}}\in E^{2N}$%
, $\mathcal{F}_{\limfunc{av},\Psi }\left( \mathbf{0}\right) =\mathbf{0}$.
Polynomial operator $\mathcal{F}_{\limfunc{av},\Psi }\left( \mathbf{\vec{v}}%
\right) $ satisfies the Lipschitz condition 
\begin{equation}
\left\Vert \mathcal{F}_{\text{$\limfunc{av}$},\Psi }\left( \mathbf{\vec{v}}%
_{1}\right) -\mathcal{F}_{\text{$\limfunc{av}$},\Psi }\left( \mathbf{\vec{v}}%
_{2}\right) \right\Vert _{E^{2N}}\leq C\tau _{\ast }\left\Vert \mathbf{\vec{v%
}}_{1}-\mathbf{\vec{v}}_{2}\right\Vert _{E^{2N}}
\end{equation}%
where $C$ depends only on $C_{\chi }$ a in (\ref{chiCR}), on the power of $%
\mathcal{F}$ and on $\left\Vert \mathbf{\vec{v}}_{1}\right\Vert
_{E^{2N}}+\left\Vert \mathbf{\vec{v}}_{2}\right\Vert _{E^{2N}}$, and, in
particular, it does not depend on $\beta ,\varrho $.
\end{lemma}

From Lemma \ref{Lemma Flipprav} and the contraction principle we obtain the
following Theorem similarly to Theorem \ref{Theorem existpr}.

\begin{theorem}
\label{Theorem existprav}Let $\left\Vert \mathbf{\vec{h}}_{\Psi }\right\Vert
_{E^{2N}}\leq R$. Then there exists $R_{1}>0$ and $\tau _{\ast }>0$ such
that equation (\ref{eqavF}) has a solution $\mathbf{\vec{v}}\in E^{2N}$
satisfying $\left\Vert \mathbf{\vec{v}}\right\Vert _{E^{2N}}\leq R_{1}$, and
such a solution is unique.
\end{theorem}

The following theorem shows that the averaged interaction system introduced
above provides a good approximation for the wave interaction system.

\begin{theorem}
\label{Theorem uminwav} Let $\mathbf{\hat{v}}_{l,\vartheta }\left( \mathbf{k}%
,\tau \right) \mathbf{\ }$be solution of (\ref{eqav}) and $\mathbf{\hat{w}}%
_{l,\vartheta }\left( \mathbf{k},\tau \right) \ $be the solution of \ (\ref%
{sysloc}). Then for sufficiently small $\beta $ the $\mathbf{\hat{v}}%
_{l,\vartheta }\left( \mathbf{k},\tau \right) $ is a wavepacket satisfying (%
\ref{piweq0}), (\ref{weq0}) with $\mathbf{\hat{w}}$ replaced by $\mathbf{%
\hat{v}}$. In addition to that, there exists $\beta _{0}>0$ such that 
\begin{equation}
\left\Vert \mathbf{\hat{v}}_{l,\vartheta }-\mathbf{\hat{w}}_{l,\vartheta
}\right\Vert _{E}\leq C\varrho ,\ l=1,\ldots ,N;\ \vartheta =\pm ,\text{ for 
}0<\varrho \leq 1,\ 0<\beta \leq \beta _{0}.  \label{vminw}
\end{equation}
\end{theorem}

\begin{proof}
Formula (\ref{piweq0}), (\ref{weq0}) for $\mathbf{\hat{v}}_{l,\vartheta
}\left( \mathbf{k},\tau \right) $ \ follow from (\ref{eqav}). We note that $%
\mathbf{\vec{w}}$ is an approximate solution of (\ref{eqav}), namely we have
an estimate for $\mathbf{D}_{\text{av}}\left( \mathbf{\hat{w}}\right) =%
\mathbf{\hat{w}}-\mathcal{F}_{\limfunc{av},\Psi }-\mathbf{\hat{h}}_{\Psi }$
which is similar to (\ref{Dw}), (\ref{Dwb}):%
\begin{equation}
\left\Vert \mathbf{D}_{\text{av}}\left( \mathbf{\hat{w}}\right) \right\Vert
=\left\Vert \mathbf{\hat{w}}-\mathcal{F}_{\limfunc{av},\Psi }-\mathbf{\hat{h}%
}\right\Vert _{E^{2N}}\leq C\varrho ,\text{ \ if }0<\varrho \leq 1,\beta
\leq \beta _{0}.  \label{wav}
\end{equation}%
The proof of (\ref{wav}) is similar to the proof of (\ref{Fbet}) with minor
simplifications thanks to the absence of terms with $\Psi _{\infty }$. Using
(\ref{wav}) we apply Lemma \ref{Lemma contr} \ and obtain (\ref{vminw}).
\end{proof}

\subsubsection{Properties of averaged nonlinearities}

In this section we discuss elementary properties of nonlinearities obtained
by formula (\ref{Fav}). A key property of such nonlinearities $F_{j,\zeta }$
is the following homogeneity-like property:%
\begin{equation}
F_{j,\zeta }\left( \mathrm{e}^{\mathrm{i}\phi _{1}t}u_{1,+},\mathrm{e}^{-%
\mathrm{i}\phi _{1}t}u_{1,-},\ldots ,\mathrm{e}^{\mathrm{i}\phi
_{N}t}u_{N,+},\mathrm{e}^{-\mathrm{i}\phi _{N}t}u_{N,-}\right) =\mathrm{e}^{%
\mathrm{i}\zeta \phi _{j}t}F_{j,\zeta }\left( u_{1+},u_{1-},\ldots
,u_{N+},u_{N-}\right) .  \label{invfi1}
\end{equation}%
The values of $\phi _{i}$, $i=1,$\ldots $N$ for which this formula holds
depend on the resonance properties of the set $S$ which enters (\ref{Fav})
through the index set $\Lambda _{n_{l},\vartheta }^{m}$. First, let us
consider the simplest case when $\phi _{i}$ are arbitrary. An example of
such a nonlinearity is the function 
\begin{equation*}
F_{2,\zeta }\left( u_{1,+},u_{1,-},u_{2,+},u_{2,-}\right)
=u_{1,+}u_{1,-}u_{2,+}.
\end{equation*}%
We call a nonlinearity which is obtained by the formula (\ref{Fav}) with a
universal resonance invariant set $S$ a \emph{universal nonlinearity. }

\begin{proposition}
\label{Proposition univprop}If $F_{j,\zeta }$ is a universal nonlinearity,
then (\ref{invfi1}) holds for arbitrary set of values $\phi _{i}$, $%
i=1,\ldots ,N$.
\end{proposition}

\begin{proof}
Note that the definition (\ref{Fav}) of the averaged nonlinearity
essentially is based on the selection of vectors $\vec{\lambda}=\left(
\left( \zeta ^{\prime },l^{\prime }\right) ,\ldots ,\left( \zeta ^{\left(
m\right) },l_{m}\right) \right) \in \Lambda _{n_{l},\vartheta }^{m}$ as in (%
\ref{resset}), which is equivalent to the resonance equation (\ref{Omeq0})
with $n=n_{l}$, $\zeta =\vartheta $. This equation has the form 
\begin{equation}
-\zeta \omega _{n}\left( \mathbf{k}_{\ast \ast }\right)
+\sum\nolimits_{l=1}^{N}\delta _{l}\omega _{l}\left( \mathbf{k}_{\ast
l}\right) =0,  \label{Omeq1}
\end{equation}%
with%
\begin{equation}
\mathbf{k}_{\ast \ast }=-\zeta \sum\nolimits_{l=1}^{N}\delta _{l}\mathbf{k}%
_{\ast l},  \label{delk}
\end{equation}%
where $\delta _{l}$ are the same as in (\ref{rearr0}). If $\vec{\lambda}\in
\Lambda _{n_{l},\vartheta }^{m}$ and 
\begin{equation*}
\mathbf{\vec{w}}_{\vec{\lambda}}=\left( \mathbf{\hat{w}}_{\lambda
_{1}}\ldots \mathbf{\hat{w}}_{\lambda _{m}}\right) =\left( \mathbf{\hat{w}}%
_{\zeta ^{\prime },l_{1}}\ldots \mathbf{\hat{w}}_{\zeta ^{\left( m\right)
},l_{m}}\right) =\left( e^{-i\zeta ^{\prime }\phi _{l_{1}}}\mathbf{\hat{v}}%
_{\zeta ^{\prime },l_{1}}\ldots e^{-i\zeta ^{\left( m\right) }\phi _{l_{m}}}%
\mathbf{\hat{v}}_{\zeta ^{\left( m\right) },l_{m}}\right) ,
\end{equation*}%
then, using (\ref{Fmksi1}) and the multilinearity of $\mathcal{F}^{\left(
m\right) }$ we get 
\begin{equation*}
\mathcal{F}_{n_{l},\vartheta ,\vec{\xi}\left( \vec{\lambda}\right) }^{\left(
m\right) }\left( \mathbf{\vec{w}}_{\vec{\lambda}}\right) =\mathrm{e}^{-%
\mathrm{i}\sum \zeta ^{\left( j\right) }\phi _{l_{j}}}\mathcal{F}%
_{n_{l},\vartheta ,\vec{\xi}\left( \vec{\lambda}\right) }^{\left( m\right)
}\left( \mathbf{\vec{v}}_{\vec{\lambda}}\right) ,
\end{equation*}%
and 
\begin{equation}
\sum_{j=1}^{m}\zeta ^{\left( j\right) }\phi
_{l_{j}}=\sum\nolimits_{l=1}^{N}\delta _{l}\phi _{l},  \label{zsumfi}
\end{equation}%
where $\delta _{l}$ are the same as in (\ref{rearr0}). If we have a
universal solution of (\ref{Omeq1}), all coefficients at every $\omega
_{l}\left( \mathbf{k}_{\ast l}\right) $ cancel out ($\omega _{n}\left( 
\mathbf{k}_{\ast \ast }\right) $ also equals one of $\omega _{l}\left( 
\mathbf{k}_{\ast l}\right) $, namely $\omega _{n}\left( \mathbf{k}_{\ast
\ast }\right) =\omega _{n_{I_{0}}}\left( \mathbf{k}_{\ast I_{0}}\right) $).
Using notation (\ref{rearr0}) we see that a universal solution is determined
by the system of equations on binary indices 
\begin{equation}
\delta _{l}=\sum_{j\in \vec{l}^{-1}\left( l\right) }\zeta ^{\left( j\right)
}=0,\ l\neq I_{0},\ \delta _{I_{0}}=\sum_{j\in \vec{l}^{-1}\left(
I_{0}\right) }\zeta ^{\left( j\right) }=\zeta .  \label{univ1}
\end{equation}%
Obviously, the above condition does not involve values of $\omega _{l}$ \
and Hence if $\delta _{l},\zeta $ correspond to a universal solution of (\ref%
{Omeq1}) then we have an identity 
\begin{equation}
-\zeta \phi _{I_{0}}+\sum\nolimits_{l=1}^{N}\delta _{l}\phi _{l}=0,
\label{Omfi}
\end{equation}%
which holds for any $\left( \phi _{1},\ldots ,\phi _{N}\right) \in \mathbb{C}%
^{N}$.
\end{proof}

Consider now the case where the $nk$-spectrum $S$ is resonance invariant but
may be not universal resonance invariant. Definition \ref{Definition omclos}
of resonance invariant $nk$-spectrum \ implies that the set $P\left(
S\right) $ of all solutions of (\ref{Omeq0}) coincides with the set $P_{%
\text{int}}\left( S\right) $ of internal solutions. Hence, all solutions of (%
\ref{Omeq1}), (\ref{delk})) are internal, in particular $\mathbf{k}_{\ast
\ast }=\mathbf{k}_{\ast I_{0}},\omega _{n}\left( \mathbf{k}_{\ast \ast
}\right) =\omega _{n_{I_{0}}}\left( \mathbf{k}_{\ast I_{0}}\right) $ \ with
some $I_{0}$.

If we have a non-universal internal solution of (\ref{Omeq1}), $\omega
_{l}\left( \mathbf{k}_{\ast l}\right) \ $ satisfy the following linear
equation 
\begin{equation}
\zeta \omega _{n_{I_{0}}}\left( \mathbf{k}_{\ast I_{0}}\right)
+\sum\nolimits_{l=1}^{N}\delta _{l}\omega _{l}\left( \mathbf{k}_{\ast
l}\right) =0,\zeta \mathbf{k}_{\ast I_{0}}+\sum\nolimits_{l=1}^{N}\delta _{l}%
\mathbf{k}_{\ast l}=0  \label{bom}
\end{equation}%
where at least one of $b_{j}$ is non-zero. Note that if (\ref{bom}) is
satisfied, we have \ additional (non-universal) solutions of (\ref{Omeq0})
defined by 
\begin{equation}
\sum_{j\in \vec{l}^{-1}\left( l\right) }\zeta ^{\left( j\right) }=\delta
_{l},\ l\neq I_{0},\ \sum_{j\in \vec{l}^{-1}\left( I_{0}\right) }\zeta
^{\left( j\right) }=\zeta +\delta _{I_{0}}.  \label{zeb}
\end{equation}%
Now let us briefly discuss properties of equations (\ref{zeb}). The
right-hand sides of the above system form a vector $\vec{b}=\left(
b_{1},\ldots ,b_{N}\right) $ with $b_{l}=\delta _{l}$, $l\neq I_{0}$, and $%
b_{I_{0}}=\zeta +\delta _{I_{0}}$. Note that $\vec{l}=\left( l_{1},\ldots
,l_{m}\right) $ is uniquely defined by its level sets $\vec{l}^{-1}\left(
l\right) $. For every $l$ the number $\delta _{+l}$ of positive $\zeta
^{\left( j\right) }$ and the number $\delta _{-l}$ of negative $\zeta
^{\left( j\right) }$ with $j\in \vec{l}^{-1}\left( l\right) $ in (\ref{zeb})
satisfy equations 
\begin{equation}
\delta _{+l}-\delta _{-l}=\delta _{l},\ \delta _{+l}+\delta _{-l}=\left\vert 
\vec{l}^{-1}\left( l\right) \right\vert  \label{eqdel}
\end{equation}%
where $\left\vert \vec{l}^{-1}\left( l\right) \right\vert =c_{l}$ is the
cardinality (number of elements) of $\vec{l}^{-1}\left( l\right) $. Hence, $%
\delta _{+l},\delta _{-l}$ are uniquely defined by $\delta _{l},\left\vert 
\vec{l}^{-1}\left( l\right) \right\vert $. Hence, the set of binary
solutions $\vec{\zeta}$ of (\ref{zeb}) with a given $\vec{b}$\ and a given $%
\vec{l}=\left( l_{1},\ldots ,l_{m}\right) $ is determined by subsets of $%
\vec{l}^{-1}\left( l\right) $ with the cardinality $\delta _{+l}$ elements.
Hence, every solution with a given $\vec{b}$\ and a given $\vec{l}$ can be
obtained from one solution by permutations of indices $j$ inside every level
set $\vec{l}^{-1}\left( l\right) $. If $\vec{b}$\ is given and the
cardinalities $\left\vert \vec{l}^{-1}\left( l\right) \right\vert =c_{l}$
are given, we can obtain different $\vec{l}$ which satisfy (\ref{zeb}) by
choosing different decomposition of $\left\{ 1,\ldots ,m\right\} $\ \ into
subsets with given cardinalities $c_{l}$. For given $\vec{b}$ and $\vec{c}%
=\left( c_{1},\ldots ,c_{m}\right) $ we obtain this way the set (may be
empty for some $\vec{b}$, $\vec{c}$) of all solutions of (\ref{zeb}).
Solutions with the same $\vec{b}$ and $\vec{c}$ we call equivalent.

When for a given wavepacket there are several non-equivalent non-universal
solutions, the number of which is denoted by $N_{c}$, we obtain from (\ref%
{bom}) a system of equations with integer coefficients 
\begin{equation}
\sum\nolimits_{l=1}^{N}b_{l,i}\omega _{l}\left( \mathbf{k}_{\ast l}\right)
=0,\ i=1,\ldots ,N_{c}  \label{bnc}
\end{equation}%
and solutions to (\ref{Omeq0}) \ can be found from 
\begin{equation}
\sum_{j\in \vec{l}^{-1}\left( l\right) }\zeta ^{\left( j\right) }=b_{l,i},%
\text{ for some }i,\ 0\leq i\leq N_{c}  \label{zebi}
\end{equation}%
where to include universal solutions we set $b_{l,0}=0$.

Hence, when a wavepacket is universally resonance invariant, we conclude
that all terms in (\ref{Fav}) satisfy (\ref{univ1}). Since (\ref{Omfi})
holds, we get (\ref{invfi1}) for arbitrary $\left( \phi _{1},\ldots ,\phi
_{N}\right) \in \mathbb{C}^{N}$. \ If the wavepacket is conditionally
universal with conditions (\ref{zebi}), then using (\ref{zsumfi}) and (\ref%
{zebi}) we conclude that (\ref{Omfi}) and (\ref{invfi1}) hold if $\left(
\phi _{1},\ldots ,\phi _{N}\right) $ satisfy the system the equations 
\begin{equation}
\sum\nolimits_{l=1}^{N}b_{l,i}\phi _{l}=0,\ i=1,\ldots ,N_{c}.  \label{bfil}
\end{equation}

Now we wold like to describe a special class of solutions of averaged
equations. The evolution equation with an averaged nonlinearity has the form 
\begin{gather}
\partial _{\tau }U_{j,+}=\frac{-\mathrm{i}}{\varrho }\mathcal{L}_{j}\left( -%
\mathrm{i}\nabla \right) U_{j,+}+F_{j,+}\left( U_{1,+},U_{1,-},\ldots
,U_{N,+},U_{N,-}\right) ,  \label{eveqav} \\
\partial _{\tau }U_{j,-}=\frac{\mathrm{i}}{\varrho }\mathcal{L}_{j}\left( 
\mathrm{i}\nabla \right) U_{j,+}+F_{j,-}\left( U_{1,+},U_{1,-},\ldots
,U_{N,+},U_{N,-}\right) ,\ j=1,\ldots ,N,  \notag
\end{gather}%
where $\mathcal{L}\left( -\mathrm{i}\nabla \right) $ \ is a linear scalar
differential operator with constant coefficients. The characteristic
property (\ref{invfi1}) implies that such a system admits special solutions
of the form 
\begin{equation}
U_{j,\zeta }\left( \mathbf{r},\tau \right) =\mathrm{e}^{-\mathrm{i}\phi
_{j}\tau /\varrho }V_{j,\zeta }\left( \mathbf{r}\right)  \label{spec}
\end{equation}%
where $V_{1,\zeta }\left( \mathbf{r}\right) $ solve the time-independent 
\emph{nonlinear eigenvalue problem} 
\begin{gather}
-\mathrm{i}\phi _{j}V_{j,+}=-\mathrm{i}\mathcal{L}_{j}\left( -\mathrm{i}%
\nabla \right) V_{j,+}+\varrho F_{j,+}\left( V_{1,+},V_{1,-},\ldots
,V_{N,+},V_{N,-}\right) ,  \label{nleigen} \\
\mathrm{i}\phi _{j}V_{j,-}=\mathrm{i}\mathcal{L}_{j}\left( \mathrm{i}\nabla
\right) V_{j,+}+\varrho F_{j,-}\left( V_{1,+},V_{1,-},\ldots
,V_{N,+},V_{N,-}\right) ,\ j=1,\ldots ,N.  \notag
\end{gather}

\subsubsection{Examples of universal and conditionally universal
nonlinearities}

Here we give a few examples of equations with averaged nonlinearities. When
the multi-wavepacket is universal resonance invariant, the averaged wave
interaction system involves NLS-type equations.

\begin{example}
\label{Example NLS}The simplest example of (\ref{eveqav}) for one wavepacket
($N=1$) \ and one spatial dimension ($d=1$) is Nonlinear Schrodinger equation%
\begin{gather}
\partial _{\tau }U_{1,+}=-\frac{\mathrm{i}}{\varrho }a_{2}\partial
_{x}^{2}U_{1,+}-\frac{\mathrm{i}}{\varrho }a_{0}U_{1,+}+a_{1}\partial
_{x}U_{1,+}-\mathrm{i}qU_{1,-}U_{1,+}^{2},  \label{exNLS} \\
\partial _{\tau }U_{1,-}=\frac{\mathrm{i}}{\varrho }a_{2}\partial
_{x}^{2}U_{j,-}+\frac{\mathrm{i}}{\varrho }a_{0}U_{1,-}+a_{1}\partial
_{x}U_{1,-}+\mathrm{i}qU_{1,+}U_{1,-}^{2}.  \notag
\end{gather}%
Note that by setting $y=x+a_{1}\tau /\varrho $ we can make $a_{1}=0$.
Obviously, the nonlinearity 
\begin{equation*}
F_{\zeta }\left( U\right) =-\mathrm{i}\zeta qU_{1,-\zeta }U_{1,\zeta }^{2}
\end{equation*}%
satisfies (\ref{invfi1}):%
\begin{equation*}
i\zeta qe^{-i\zeta \phi _{1}}U_{1,-\zeta }\left( e^{i\zeta \phi
_{1}}U_{1,\zeta }\right) ^{2}=e^{i\zeta \phi _{1}}i\zeta qU_{1,-\zeta
}\left( U_{1,\zeta }\right) ^{2}.
\end{equation*}%
The eigenvalue problem in this case takes the form 
\begin{gather}
\mathrm{i}\phi _{1}V_{1,+}=-\mathrm{i}a_{2}\partial _{x}^{2}V_{1,+}-\mathrm{i%
}a_{0}V_{1,+}+a_{1}\partial _{x}V_{1,+}-\mathrm{i}\varrho
qV_{1,-}V_{1,+}^{2},  \label{NLSei} \\
-\mathrm{i}\phi _{1}V_{1,-}=\mathrm{i}a_{2}\partial _{x}^{2}V_{j,-}+\mathrm{i%
}a_{0}V_{1,-}+a_{1}\partial _{x}V_{1,-}+\mathrm{i}\varrho
qV_{1,+}V_{1,-}^{2}.  \notag
\end{gather}%
If $a_{1}=0$ and we consider real-valued $V_{1,+}=V_{1,-}$ we obtain the
equation 
\begin{equation*}
\left( \phi _{1}+a_{0}\right) V_{1,+}=-a_{2}\partial _{x}^{2}V_{1,+}-\varrho
qV_{1,+}^{3}
\end{equation*}%
or, equivalently,%
\begin{equation*}
\frac{\left( \phi _{1}+a_{0}\right) }{\varrho q}V_{1,+}+\frac{a_{2}}{\varrho
q}\partial _{x}^{2}V_{1,+}+V_{1,+}^{3}=0.
\end{equation*}%
If 
\begin{equation}
c^{2}=\frac{a_{2}}{\varrho q}>0,\frac{\left( \phi _{1}+a_{0}\right) }{%
\varrho q}=-b^{2}<0,  \label{cb}
\end{equation}%
the last equation takes the form \ 
\begin{equation*}
-b^{2}V_{1,+}+c^{2}\partial _{x}^{2}V_{1,+}+V_{1,+}^{3}=0.
\end{equation*}%
with a family of classical soliton solutions 
\begin{equation*}
V_{1,+}=2^{1/2}\frac{b}{\cosh \left( b\left( x-x_{0}\right) /c\right) }.
\end{equation*}%
Note that the norm of Fourier transform $\left\Vert \hat{V}_{1,+}\right\Vert
_{L^{1}}=Cb$ where $C$ is an absolute constant, Hence to have $\hat{V}_{1,+}$
bounded in $L^{1}$ uniformly in small $\varrho $ according to (\ref{cb}) we
should take $\phi _{1}=-a_{0}-b^{2}\varrho q$ with a bounded $b$.
\end{example}

If the universal resonance invariant multi-wavepacket involves two
wavepackets ($N=2$) and the nonlinearity $F$ is cubic, that is $\mathfrak{M}%
_{F}=\left\{ 3\right\} $, a semilinear system PDE with averaged nonlinearity
has the form 
\begin{align*}
\partial _{t}U_{2,+}& =-iL_{2}\left( i\nabla \right) U_{2,+}+U_{2,+}\left(
Q_{2,1,+}U_{1,+}U_{1,-}+Q_{2,2,+}U_{2,+}U_{2,-}\right) , \\
\partial _{t}U_{2,-}& =iL_{2}\left( -i\nabla \right) U_{2,-}+U_{2,-}\left(
Q_{2,1,-}U_{1,+}U_{1,-}+Q_{2,2,-}U_{2,+}U_{2,-}\right) , \\
\partial _{t}U_{1,+}& =-iL_{1}\left( i\nabla \right) U_{1,+}+U_{1,+}\left(
Q_{1,1,+}U_{1,+}U_{1,-}+Q_{1,1,+}U_{2,+}U_{2,-}\right) , \\
\partial _{t}U_{1,-}& =iL_{1}\left( -i\nabla \right) U_{1,-}+U_{1,-}\left(
Q_{1,1,-}U_{1,+}U_{1,-}+Q_{1,1,-}U_{2,+}U_{2,-}\right) .
\end{align*}%
Obviously, (\ref{invfi1}) holds with arbitrary $\phi _{1}$, $\phi _{2}$.

Now let us consider quadratic nonlinearities. In particular, let us concider
the one-band symmetric case $\omega _{n}\left( \mathbf{k}\right) =\omega
_{1}\left( \mathbf{k}\right) =\omega _{1}\left( -\mathbf{k}\right) $, i.e. $%
J=1$, $\mathfrak{M}_{F}=\left\{ 2\right\} $ and $m=2$. Suppose that there is
a multi-wavepacket involving two wavepackets with wavevectors $\mathbf{k}%
_{\ast 1},\mathbf{k}_{\ast 2}$ i.e. $N=2$. The resonance equation (\ref%
{Omeq0}) takes now the form 
\begin{equation}
-\zeta \omega _{1}\left( \zeta ^{\prime }\mathbf{k}_{\ast l_{1}}+\zeta
^{\prime \prime }\mathbf{k}_{\ast l_{2}}\right) +\zeta ^{\prime }\omega
_{1}\left( \mathbf{k}_{\ast l_{1}}\right) +\zeta ^{\prime \prime }\omega
_{1}\left( \mathbf{k}_{\ast l_{2}}\right) =0,  \label{ex122}
\end{equation}%
where $l_{1},l_{2}\in \left\{ 1,2\right\} ,$ $\zeta ,\zeta ^{\prime },\zeta
^{\prime \prime }\in \left\{ -1,1\right\} $. All possible cases, and there
are exactly \ four of them, correspond to the four well known effects in the
nonlinear optics: (i) $l_{1}=l_{2}$, $\zeta ^{\prime }=\zeta ^{\prime \prime
}$ and $\zeta ^{\prime }=-\zeta ^{\prime \prime }$ correspond respectively
to second harmonic generation and nonlinear optical rectification; (ii) $%
l_{1}\neq l_{2}$, $\zeta ^{\prime }=\zeta ^{\prime \prime }$and $\zeta
^{\prime }=-\zeta ^{\prime \prime }$ correspond respectively to
sum-frequency and difference-frequency interactions.

Let us suppose now that $\mathbf{k}_{\ast 1},\mathbf{k}_{\ast 2}\neq 0$ and $%
\omega _{1}\left( \mathbf{k}_{\ast 1}\right) \neq 0$, $\omega _{1}\left( 
\mathbf{k}_{\ast 2}\right) \neq 0$, where the last conditions exclude the
optical rectification, and that $\mathbf{k}_{\ast i}\neq 0$ and $\mathbf{k}%
_{\ast i},2\mathbf{k}_{\ast i},\mathbf{0},$ $\pm \mathbf{k}_{\ast 1}\pm 
\mathbf{k}_{\ast 2}$ \ are not band-crossing points. Consider first the case
when the wavepacket is universally resonance invariant.

\begin{example}
\label{Example quad2} Suppose there is a single band, i.e. $J=1$, with a
symmetric dispersion relation, and a quadratic nonlinearity $F$ , that is $%
\mathfrak{M}_{F}=\left\{ 2\right\} $. Let us pick two points $\mathbf{k}%
_{\ast 1}$ and $\mathbf{k}_{\ast 2}\neq \pm \mathbf{k}_{\ast 1}$, and assume
that $\mathbf{k}_{\ast i}\neq 0$ and $\mathbf{k}_{\ast i},2\mathbf{k}_{\ast
i},\mathbf{0},$ $\mathbf{k}_{\ast 1}\pm \mathbf{k}_{\ast 2}$ \ are not
band-crossing points. Assume also that (i) $2\omega _{1}\left( \mathbf{k}%
_{\ast i}\right) \neq \omega _{1}\left( 2\mathbf{k}_{\ast i}\right) ,$ $%
i,j,l=1,2$ , so there is no no second harmonic generation; (ii) $\omega
_{1}\left( \mathbf{k}_{\ast 1}\right) \pm \omega _{1}\left( \mathbf{k}_{\ast
2}\right) \neq \omega _{1}\left( \mathbf{k}_{\ast 1}\pm \mathbf{k}_{\ast
2}\right) ,$ (no sum/difference-frequency interactions); (iii) $\omega
_{1}\left( \mathbf{0}\right) \neq 0$, $\omega _{j}\left( \mathbf{k}_{\ast
1}\right) \pm \omega _{l}\left( \mathbf{k}_{\ast 2}\right) \neq 0$. Let us
set the $nk$-spectrum to be the set $S_{1}=\left\{ \left( 1,\mathbf{k}_{\ast
1}\right) ,\left( 1,\mathbf{k}_{\ast 2}\right) \right\} $. Then $S_{1}$ is
resonance invariant.
\end{example}

In this case (\ref{ex122}) does not have solutions. Hence $\Lambda
_{n_{l},\vartheta }^{m}=\varnothing $ and the averaged nonlinearity equals
zero.

Now let us consider the case where the wavepacket is not universal resonance
invariant, but conditionally universal resonance invariant. In the following
example a conditionally resonance invariant spectrum allows for second
harmonic generation \ in the averaged system.

\begin{example}
\label{Example quad21} Suppose there is a single band, i.e. $J=1$, with a
symmetric dispersion relation, and a quadratic nonlinearity $F$, that is $%
\mathfrak{M}_{F}=\left\{ 2\right\} $. Let us pick two points $\mathbf{k}%
_{\ast 1}$ and $\mathbf{k}_{\ast 2}$ such that $\mathbf{k}_{\ast 2}=2\mathbf{%
k}_{\ast 1}$, and assume that $\mathbf{k}_{\ast i}\neq 0$ and $\mathbf{k}%
_{\ast i},2\mathbf{k}_{\ast i},\mathbf{0},$ $\pm \mathbf{k}_{\ast 1}\pm 
\mathbf{k}_{\ast 2}$ \ are not band-crossing points. Assume also that (i) $%
2\omega _{1}\left( \mathbf{k}_{\ast 1}\right) =\omega _{1}\left( 2\mathbf{k}%
_{\ast 1}\right) $ (second harmonic generation); (ii) $\omega _{i}\left( 
\mathbf{k}_{\ast 1}\right) \pm \omega _{j}\left( \mathbf{k}_{\ast 2}\right)
\neq \omega _{l}\left( \mathbf{k}_{\ast 1}\pm \mathbf{k}_{\ast 2}\right) ,$\ 
$i,j,l=1,2$ (no sum-/difference-frequencies interaction); (iii) $\omega
_{1}\left( \mathbf{0}\right) \neq 0$, $\omega _{j}\left( \mathbf{k}_{\ast
1}\right) \pm \omega _{l}\left( \mathbf{k}_{\ast 2}\right) \neq 0$. Let us
set the $nk$-spectrum to be the set $S=\left\{ \left( 1,\mathbf{k}_{\ast
1}\right) ,\left( 1,\mathbf{k}_{\ast 2}\right) \right\} $. Then $S$ is
resonance invariant. The condition (\ref{bom}) is takes here the form \ 
\begin{equation*}
2\omega _{1}\left( \mathbf{k}_{\ast 1}\right) -\omega _{1}\left( \mathbf{k}%
_{\ast 2}\right) =0,\ 2\mathbf{k}_{\ast 1}-\mathbf{k}_{\ast 2}=0,
\end{equation*}%
and the condition (\ref{bfil}) turns into%
\begin{equation*}
2\omega _{1}\left( \mathbf{k}_{\ast 1}\right) -\omega _{1}\left( \mathbf{k}%
_{\ast 2}\right) =0.
\end{equation*}%
The wavepacket interaction system for such a multiwavepacket has the form 
\begin{align*}
\partial _{t}U_{2,+}& =-iL_{2}\left( i\nabla \right)
U_{2,+}+Q_{2,2,+}U_{1,+}U_{1,+}, \\
\partial _{t}U_{2,-}& =iL_{2}\left( -i\nabla \right)
U_{2,-}+Q_{2,2,-}U_{1,-}U_{1,-}, \\
\partial _{t}U_{1,+}& =-iL_{1}\left( i\nabla \right)
U_{1,+}+Q_{1,2,+}U_{2,+}U_{1,-}, \\
\partial _{t}U_{1,-}& =iL_{1}\left( -i\nabla \right)
U_{1,-}+Q_{1,2,-}U_{2,-}U_{1,+}.
\end{align*}
\end{example}

\subsection{Invariance of multi-particle wavepackets}

The following Lemma shows that particle wavepackets are preserved under
action of certain types of nonlinearities with elementary susceptibilities
as in (\ref{chim}). In the following section we show in particular that
universal nonlinearities are composed of such terms.

\begin{lemma}
\label{Lemma gradF}\ Let \ components $\ \mathbf{\hat{w}}_{l_{i},\zeta }=%
\mathbf{\hat{w}}_{\lambda _{i}}$ of $\mathbf{\vec{w}}_{\vec{\lambda}}=%
\mathbf{\hat{w}}_{\lambda _{1}}\ldots \mathbf{\hat{w}}_{\lambda _{m}}$ be
particle-like wavepackets in the sense of Definition \ref{Definition regwave}
\ and $\mathcal{F}_{n_{l},\vartheta ,\vec{\xi}\left( \vec{\lambda}\right)
}^{\left( m\right) }\left( \mathbf{\vec{w}}_{\vec{\lambda}}\right) $ be as
in (\ref{Fav}). Assume that $\ $%
\begin{equation}
\mathbf{\hat{w}}_{l_{i},\zeta }\left( \mathbf{k},\beta \right) =0\text{ if }%
\left\vert \mathbf{k}-\zeta \mathbf{k}_{\ast li}\right\vert \geq \beta
^{1-\epsilon },\;\zeta =\pm ,\;i=1,\ldots ,m.  \label{hl0}
\end{equation}%
Assume that vector index $\vec{\lambda}\in \Lambda _{n_{l},\vartheta }^{m}$
be such a vector which has at least one \ component $\lambda _{j}=\left(
\zeta _{j},l_{j}\right) $ such that \ 
\begin{equation}
\nabla \omega _{n_{l}}\left( \mathbf{k}_{\ast l}\right) =\nabla \omega
_{n_{l_{j}}}\left( \mathbf{k}_{\ast l_{j}}\right) .  \label{GVM0}
\end{equation}%
\ Then \ for any $\mathbf{r}_{\ast }\in \mathbb{R}^{d}$ 
\begin{gather}
\left\Vert \nabla _{\mathbf{k}}\left( \mathrm{e}^{-\mathrm{i}\mathbf{r}%
_{\ast }\mathbf{k}}\Psi \left( \mathbf{\cdot },\mathbf{k}_{\ast l},\beta
^{1-\epsilon }\right) \mathcal{F}_{n_{l},\vartheta ,\vec{\xi}\left( \vec{%
\lambda}\right) }^{\left( m\right) }\left( \mathbf{\vec{w}}_{\vec{\lambda}%
}\right) \right) \right\Vert _{E}\leq  \label{gradkk} \\
C\tau _{\ast }\left\Vert \nabla _{\mathbf{k}}\mathrm{e}^{-\mathrm{i}\mathbf{r%
}_{\ast }\mathbf{k}^{\left( j\right) }}\mathbf{w}_{l_{j}}\right\Vert
_{E}\dprod\limits_{i\neq j}\left\Vert \mathbf{w}_{l_{j},\zeta
_{j}}\right\Vert _{E}+C\tau _{\ast }\left( \beta ^{-1+\epsilon }+\frac{\beta
^{1-\epsilon }}{\varrho }\right) \dprod\limits_{j=1}^{m}\left\Vert \mathbf{w}%
_{l_{j},\zeta _{j}}\right\Vert _{E}.  \notag
\end{gather}%
where $C$ does not depend on $\mathbf{r}_{\ast }$ and small $\beta ,\varrho $%
.
\end{lemma}

\begin{proof}
Note that \ 
\begin{equation*}
\mathbf{r}_{\ast }\mathbf{k}=\mathbf{r}_{\ast }\left( \mathbf{k}^{\prime
}+\ldots +\mathbf{k}^{\left( m\right) }\right) .
\end{equation*}%
We have by (\ref{Fm}) 
\begin{gather}
\nabla _{\mathbf{k}}\mathrm{e}^{-\mathrm{i}\mathbf{r}_{\ast }\mathbf{k}}%
\mathcal{F}_{n_{l},\vartheta ,\vec{\xi}\left( \vec{\lambda}\right) }^{\left(
m\right) }\left( \mathbf{\vec{w}}_{\vec{\lambda}}\right) \left( \mathbf{k}%
,\tau \right) =\nabla _{\mathbf{k}}\int_{0}^{\tau }\int_{\left[ -\pi ,\pi %
\right] ^{2d}}\exp \left\{ \mathrm{i}\phi _{\ \theta ,\vec{\zeta}}\left( 
\mathbf{\mathbf{k}},\vec{k}\right) \frac{\tau _{1}}{\varrho }\right\}
\label{gradF} \\
\Psi \mathrm{e}^{-\mathrm{i}\mathbf{r}_{\ast }\mathbf{k}}\chi _{\theta ,\vec{%
\zeta}}^{\left( m\right) }\left( \mathbf{\mathbf{k}},\vec{k}\right) \mathbf{w%
}_{l_{1},\zeta ^{\prime }}\left( \mathbf{k}^{\prime }\right) \ldots \mathbf{w%
}_{l_{m},\zeta ^{\left( m\right) }}\left( \mathbf{k}^{\left( m\right)
}\left( \mathbf{k},\vec{k}\right) \right) \,\mathrm{\tilde{d}}^{\left(
m-1\right) d}\vec{k}d\tau _{1}.  \notag
\end{gather}%
Without loss of generality we assume that in (\ref{GVM0}) $l_{j}=l_{m}$ (the
general case is reduced to this one by a renumeration of variables of
integration). By the Leibnitz formula%
\begin{equation}
\nabla _{\mathbf{k}}\left[ \Psi \mathrm{e}^{-\mathrm{i}\mathbf{r}_{\ast }%
\mathbf{k}}\mathcal{F}_{n_{l},\vartheta ,\vec{\xi}\left( \vec{\lambda}%
\right) }^{\left( m\right) }\left( \mathbf{\vec{w}}_{\vec{\lambda}}\right) %
\right] \left( \mathbf{k},\tau \right) =I_{1}+I_{2}+I_{3},  \label{Leib1}
\end{equation}%
where%
\begin{gather*}
I_{1}=\int_{0}^{\tau }\int_{\left[ -\pi ,\pi \right] ^{\left( m-1\right)
d}}\nabla _{\mathbf{k}}\exp \left\{ \mathrm{i}\phi _{\theta ,\vec{\zeta}%
}\left( \mathbf{\mathbf{k}},\vec{k}\right) \frac{\tau _{1}}{\varrho }-%
\mathrm{i}\mathbf{r}_{\ast }\mathbf{k}\right\} \times \\
\Psi \chi _{\theta ,\vec{\zeta}}^{\left( m\right) }\left( \mathbf{\mathbf{k}}%
,\vec{k}\right) \mathrm{e}^{-\mathrm{i}\mathbf{r}_{\ast }\mathbf{k}^{\prime
}}\mathbf{w}_{l_{1},\zeta ^{\prime }}\left( \mathbf{k}^{\prime }\right)
\ldots \mathrm{e}^{-\mathrm{i}\mathbf{r}_{\ast }\mathbf{k}^{\left( m\right)
}}\mathbf{w}_{l_{m},\zeta ^{\left( m\right) }}\left( \mathbf{k}^{\left(
m\right) }\left( \mathbf{k},\vec{k}\right) \right) \mathrm{\tilde{d}}%
^{\left( m-1\right) d}\vec{k}d\tau _{1},
\end{gather*}%
\begin{gather*}
I_{2}=\int_{0}^{\tau }\int_{\left[ -\pi ,\pi \right] ^{\left( m-1\right)
d}}\Psi \exp \left\{ \mathrm{i}\phi _{\theta ,\vec{\zeta}}\left( \mathbf{%
\mathbf{k}},\vec{k}\right) \frac{\tau _{1}}{\varrho }-\mathrm{i}\mathbf{r}%
_{\ast }\mathbf{k}\right\} \times \\
\left[ \nabla _{\mathbf{k}}\left( \Psi \left( \mathbf{\mathbf{k}},\mathbf{k}%
_{\ast l},\beta ^{1-\epsilon }\right) \chi _{\theta ,\vec{\zeta}}^{\left(
m\right) }\left( \mathbf{\mathbf{k}},\vec{k}\right) \right) \right] \mathrm{e%
}^{-\mathrm{i}\mathbf{r}_{\ast }\mathbf{k}^{\prime }}\mathbf{w}_{l_{1},\zeta
^{\prime }}\left( \mathbf{k}^{\prime }\right) \ldots \mathrm{e}^{-\mathrm{i}%
\mathbf{r}_{\ast }\mathbf{k}^{\left( m\right) }}\mathbf{w}_{l_{m},\zeta
^{\left( m\right) }}\left( \mathbf{k}^{\left( m\right) }\left( \mathbf{k},%
\vec{k}\right) \right) \mathrm{\tilde{d}}^{\left( m-1\right) d}\vec{k}d\tau
_{1},
\end{gather*}%
\begin{gather*}
I_{3}=\int_{0}^{\tau }\int_{\left[ -\pi ,\pi \right] ^{\left( m-1\right)
d}}\exp \left\{ \mathrm{i}\phi _{\theta ,\vec{\zeta}}\left( \mathbf{\mathbf{k%
}},\vec{k}\right) \frac{\tau _{1}}{\varrho }-\mathrm{i}\mathbf{r}_{\ast }%
\mathbf{k}\right\} \times \\
\Psi \chi _{\theta ,\vec{\zeta}}^{\left( m\right) }\left( \mathbf{\mathbf{k}}%
,\vec{k}\right) \mathrm{e}^{-\mathrm{i}\mathbf{r}_{\ast }\mathbf{k}^{\prime
}}\mathbf{w}_{l_{1},\zeta ^{\prime }}\left( \mathbf{k}^{\prime }\right)
\ldots \nabla _{\mathbf{k}}\left( \mathrm{e}^{-\mathrm{i}\mathbf{r}_{\ast }%
\mathbf{k}^{\left( m\right) }\left( \mathbf{k},\vec{k}\right) }\mathbf{w}%
_{l_{m},\zeta ^{\left( m\right) }}\left( \mathbf{k}^{\left( m\right) }\left( 
\mathbf{k},\vec{k}\right) \right) \right) \mathrm{\tilde{d}}^{\left(
m-1\right) d}\vec{k}d\tau _{1}.
\end{gather*}%
Since $\mathbf{w}_{j,\zeta }$ are bounded, we have 
\begin{equation}
\left\Vert \mathrm{e}^{-\mathrm{i}\mathbf{r}_{\ast j}\mathbf{k}^{\left(
j\right) }}\mathbf{w}_{l_{j},\zeta ^{\left( j\right) }}\left( \mathbf{k}%
^{\left( j\right) }\right) \right\Vert _{L^{1}}\leq \left\Vert \mathbf{w}%
_{l_{j},\zeta ^{\left( j\right) }}\left( \mathbf{k}^{\left( j\right)
}\right) \right\Vert _{L^{1}}\leq C_{1},\ j=1,\ldots ,m.  \label{ucr}
\end{equation}%
Using (\ref{Yconv}) and (\ref{ucr}) we get%
\begin{equation}
\left\vert I_{3}\right\vert \leq \left\Vert \chi ^{\left( m\right)
}\right\Vert \dprod\limits_{j=1}^{m-1}\left\Vert \mathbf{w}_{l_{j},\zeta
^{\left( j\right) }}\right\Vert _{E}\int_{0}^{\tau }\left\Vert \nabla _{%
\mathbf{k}}\mathrm{e}^{-\mathrm{i}\mathbf{r}_{\ast }\mathbf{k}^{\left(
m\right) }\left( \mathbf{k},\vec{k}\right) }\mathbf{w}_{l_{m},\zeta ^{\left(
m\right) }}\right\Vert _{E}d\tau _{1}.  \label{I3}
\end{equation}%
From (\ref{ucr}), (\ref{Psik}), (\ref{grchi}) and the smoothness of $\Psi
\left( \mathbf{\mathbf{k}},\mathbf{k}_{\ast l},\beta ^{1-\epsilon }\right) $
we get 
\begin{equation}
\left\vert I_{2}\right\vert \leq C_{2}\beta ^{-1+\epsilon
}\dprod\limits_{j=1}^{m}\left\Vert \mathbf{w}_{l_{j},\zeta _{j}}\right\Vert
_{E}.  \label{I2}
\end{equation}%
\bigskip Now let us estimate $I_{1}$. Using (\ref{phim}) we obtain%
\begin{gather}
I_{1}=\int_{0}^{\tau }\int_{\left[ -\pi ,\pi \right] ^{\left( m-1\right) d}}%
\left[ \exp \left\{ \mathrm{i}\phi _{\theta ,\vec{\zeta}}\left( \mathbf{%
\mathbf{k}},\vec{k}\right) \frac{\tau _{1}}{\varrho }\right\} \right]
\label{I1int} \\
\frac{\tau _{1}}{\varrho }\left[ -\theta \nabla _{\mathbf{k}}\omega
_{n_{l}}\left( \mathbf{k}\right) +\zeta ^{\left( m\right) }\nabla _{\mathbf{k%
}}\omega _{n_{l_{m}}}\left( \mathbf{k}^{\left( m\right) }\left( \mathbf{k},%
\vec{k}\right) \right) \right]  \notag \\
\chi _{\theta ,\vec{\zeta}}^{\left( m\right) }\left( \mathbf{\mathbf{k}},%
\vec{k}\right) \mathbf{w}_{l_{1},\zeta ^{\prime }}\left( \mathbf{k}^{\prime
}\right) \ldots \mathbf{w}_{l_{m},\zeta ^{\left( m\right) }}\left( \mathbf{k}%
^{\left( m\right) }\left( \mathbf{k},\vec{k}\right) \right) \mathrm{\tilde{d}%
}^{\left( m-1\right) d}\vec{k}d\tau _{1}.  \notag
\end{gather}%
The difficulty in the estimation of the integral $I_{1}$ comes from the
factor $\frac{\tau _{1}}{\varrho }$ since $\varrho $ is small. \ Since (\ref%
{hl0}) holds, it is sufficient to estimate $I_{1}$ if$.$%
\begin{equation}
\left\vert \mathbf{k}^{\left( j\right) }-\zeta ^{\left( j\right) }\mathbf{k}%
_{\ast n_{j}}\right\vert \leq \beta ^{1-\epsilon }\text{ for all }j.
\label{allj}
\end{equation}%
According to (\ref{kapzel}), since $\vec{\lambda}\in \Lambda
_{n_{l},\vartheta }^{m}$, we have 
\begin{equation*}
\mathbf{k}^{\left( m\right) }\left( \mathbf{\mathbf{k}}_{\ast n_{l}},\vec{k}%
_{\ast }\right) =\mathbf{\mathbf{k}}_{\ast n_{l_{m}}}.
\end{equation*}%
Hence, using (\ref{GVM0}) \ and (\ref{phim}) we obtain $\ $ 
\begin{equation}
\nabla _{\mathbf{\mathbf{k}}}\phi _{\theta ,\vec{\zeta}}\left( \mathbf{%
\mathbf{k}}_{\ast n_{l}},\vec{k}_{\ast }\right) =\left[ -\theta \nabla _{%
\mathbf{k}}\omega _{n_{l}}\left( \mathbf{\mathbf{k}}_{\ast n_{l}}\right)
+\zeta ^{\left( m\right) }\nabla _{\mathbf{k}}\omega _{n_{l_{m}}}\left(
\left( \mathbf{k}^{\left( m\right) }\left( \mathbf{\mathbf{k}}_{\ast n_{l}},%
\vec{k}_{\ast }\right) \right) \right) \right] =0.  \label{grk0}
\end{equation}%
Using (\ref{Com2}) we conclude that in a vicinity of $\vec{k}_{\ast }\mathbf{%
\ }$defined by (\ref{allj}) we have 
\begin{equation*}
\left\vert \left[ -\theta \nabla _{\mathbf{k}}\omega \left( \mathbf{k}%
\right) +\zeta ^{\left( m\right) }\nabla _{\mathbf{k}}\omega \left( \mathbf{k%
}^{\left( m\right) }\left( \mathbf{k},\vec{k}\right) \right) \right]
\right\vert \leq 2\left( m+1\right) C_{\omega ,2}\beta ^{1-\epsilon }.
\end{equation*}%
This yields the estimate%
\begin{equation}
\left\vert I_{1}\right\vert \leq C_{3}\beta ^{1-\epsilon }/\varrho .
\label{I1}
\end{equation}%
Combining (\ref{I1}), (\ref{I2}) and (\ref{I3}) we obtain (\ref{gradk}). \ 
\end{proof}

We introduce a $\beta $-dependent Banach space $E^{1}$ of differentiable
functions of variable $\mathbf{k}$ by the formula%
\begin{equation}
\left\Vert \mathbf{w}\right\Vert _{E^{1}\left( \mathbf{r}_{\ast }\right)
}=\beta ^{1+\epsilon }\left\Vert \nabla _{\mathbf{k}}\left( \mathrm{e}^{-%
\mathrm{i}\mathbf{r}_{\ast }\mathbf{k}}\mathbf{w}\right) \right\Vert
_{E}+\left\Vert \mathbf{w}\right\Vert _{E}.  \label{E1}
\end{equation}%
We use for $2N$- component vectors \ with elements $\mathbf{w}_{i}\left( 
\mathbf{k}\right) \in E^{2}$ the following notation 
\begin{equation}
\mathbf{\vec{w}}\left( \mathbf{k}\right) =\left( \mathbf{w}_{1}\left( 
\mathbf{k}\right) ,\ldots ,\mathbf{w}_{N}\left( \mathbf{k}\right) \right) ,\ 
\mathbf{\vec{r}}_{\ast }=\left( \mathbf{r}_{\ast 1},\ldots ,\mathbf{r}_{\ast
N}\right) ,\ \mathbf{w}_{i}\left( \mathbf{k}\right) =\left( \mathbf{w}%
_{i,+}\left( \mathbf{k}\right) ,\mathbf{w}_{i,-}\left( \mathbf{k}\right)
\right) ,  \label{rarrow}
\end{equation}%
\begin{equation*}
\mathrm{e}^{-\mathrm{i}\mathbf{\vec{r}}_{\ast }\mathbf{k}}\mathbf{\vec{w}}%
\left( \mathbf{k}\right) =\left( \mathrm{e}^{-\mathrm{i}\mathbf{r}_{\ast 1}%
\mathbf{k}}\mathbf{w}_{1}\left( \mathbf{k}\right) ,\ldots ,\mathrm{e}^{-%
\mathrm{i}\mathbf{r}_{\ast N}\mathbf{k}}\mathbf{w}_{N}\left( \mathbf{k}%
\right) \right) ,
\end{equation*}%
Similarly to (\ref{E2N}) we introduce the space $\left( E^{1}\right)
^{2N}\left( \mathbf{\vec{r}}_{\ast }\right) $ with the norm 
\begin{equation}
\left\Vert \mathbf{\vec{w}}\right\Vert _{\left( E^{1}\right) ^{2N}\left( 
\mathbf{\vec{r}}_{\ast }\right) }=\sum\nolimits_{l,\vartheta }\left\Vert 
\mathbf{\hat{w}}_{l,\vartheta }\right\Vert _{E^{1}\left( \mathbf{r}_{\ast
l}\right) }.  \label{E12Nr}
\end{equation}%
The following proposition is obtained by comparing (\ref{E1}) and (\ref%
{gradker}).

\begin{proposition}
\label{Proposition partnorm}A multi-wavepacket $\mathbf{\vec{w}}$ is a
multi-particle one with positions $\mathbf{r}_{\ast 1},$\ldots $,\mathbf{r}%
_{\ast N}$ if and only if 
\begin{equation*}
\left\Vert \mathbf{\vec{w}}\right\Vert _{\left( E^{1}\right) ^{2N}\left( 
\mathbf{\vec{r}}_{\ast }\right) }\leq C
\end{equation*}%
where the constant $C$ does not depend on $\beta ,$ $0<\beta \leq 1/2,$ and $%
\mathbf{\vec{r}}_{\ast }$.
\end{proposition}

In view of the above we will call $E^{1}\left( \mathbf{r}_{\ast }\right) $
and $\left( E^{1}\right) ^{2N}\left( \mathbf{\vec{r}}_{\ast }\right) $
particle spaces. We also use notations 
\begin{equation*}
\Psi _{2}\mathbf{\vec{w}}_{\vec{\lambda}}=\left( \Psi \left( \mathbf{\cdot },%
\mathbf{k}_{\ast l_{1}},\beta ^{1-\epsilon }/2\right) \mathbf{w}_{\lambda
_{1}},\ldots ,\Psi \left( \mathbf{\cdot },\mathbf{k}_{\ast l_{m}},\beta
^{1-\epsilon }/2\right) \mathbf{w}_{\lambda _{m}}\right) ,
\end{equation*}%
\begin{equation*}
\mathcal{F}_{n_{l},\vartheta ,\vec{\lambda},\Psi _{2}}^{\left( m\right)
}\left( \mathbf{\vec{w}}^{m}\right) =\Psi \left( \mathbf{\cdot },\mathbf{k}%
_{\ast l},\beta ^{1-\epsilon }\right) \mathcal{F}_{n_{l},\vartheta ,\vec{\xi}%
\left( \vec{\lambda}\right) }^{\left( m\right) }\left( \Psi _{2}\mathbf{\vec{%
w}}_{\vec{\lambda}}\right) .
\end{equation*}

\begin{lemma}
\label{Lemma E1lip} Let $\mathbf{\vec{w}},\mathbf{\vec{v}}\in \left(
E^{1}\right) ^{2N}\left( \mathbf{\vec{r}}_{\ast }\right) $ \ and $\mathcal{F}%
_{n_{l},\vartheta ,\vec{\xi}\left( \vec{\lambda}\right) }^{\left( m\right)
}\left( \mathbf{\vec{w}}_{\vec{\lambda}}\right) $ be as in (\ref{Fav}).
Assume that vector index $\vec{\lambda}\in \Lambda _{n_{l},\vartheta }^{m}$
be such a vector which has at least one \ component $\lambda _{j}=\left(
\zeta _{j},l_{j}\right) $ \ with $l_{j}=l$. Assume that $\ $ (\ref{rbb1})
holds and $\Psi \left( \mathbf{\cdot },\mathbf{k}_{\ast },\beta ^{1-\epsilon
}\right) $ is defined in (\ref{Psik}). Let $\left\Vert \mathbf{\vec{w}}%
\right\Vert _{\left( E^{1}\right) ^{2N}\left( \mathbf{\vec{r}}_{\ast
}\right) }\leq 2R$. Then 
\begin{equation}
\left\Vert \mathcal{F}_{n_{l},\vartheta ,\vec{\lambda},\Psi _{2}}^{\left(
m\right) }\left( \mathbf{\vec{w}}\right) \right\Vert _{E^{1}\left( \mathbf{r}%
_{\ast l}\right) }\leq C\tau _{\ast }\left\Vert \mathbf{\vec{w}}\right\Vert
_{\left( E\right) ^{2N}}^{m-1}\left\Vert \mathbf{\vec{w}}\right\Vert
_{\left( E^{1}\right) ^{2N}\left( \mathbf{\vec{r}}_{\ast }\right) }
\label{E1bound}
\end{equation}%
where $C$ does not depend on $\beta $, $0<\beta \leq 1/2$, \ and on $\mathbf{%
\vec{r}}_{\ast }$, $\mathbf{r}_{\ast l}$, and $\ \mathbf{\vec{r}}_{\ast }$
is defined by (\ref{rarrow}).\ If $\ $ $\left\Vert \mathbf{\vec{v}}%
\right\Vert _{\left( E^{1}\right) ^{2N}\left( \mathbf{\vec{r}}_{\ast
}\right) }\leq 2R$ the following Lipschitz inequality holds 
\begin{equation}
\left\Vert \mathcal{F}_{n_{l},\vartheta ,\vec{\lambda},\Psi _{2}}^{\left(
m\right) }\left( \mathbf{\vec{w}}\right) -\mathcal{F}_{n_{l},\vartheta ,\vec{%
\lambda},\Psi _{2}}^{\left( m\right) }\left( \mathbf{\vec{v}}\right)
\right\Vert _{E^{1}\left( \mathbf{r}_{\ast l}\right) }\leq C\tau _{\ast
}\left\Vert \mathbf{\vec{w}}-\mathbf{\vec{v}}\right\Vert _{\left(
E^{1}\right) ^{2N}\left( \mathbf{\vec{r}}_{\ast }\right) }  \label{gradk}
\end{equation}%
where $C$ does not depend on $\beta $, $0<\beta \leq 1/2$, \ and on $\mathbf{%
\vec{r}}_{\ast }$, $\mathbf{r}_{\ast l}$.
\end{lemma}

\begin{proof}
Note that $\Psi _{2}\mathbf{\vec{w}}_{\vec{\lambda}}$ and $\Psi _{2}\mathbf{%
\vec{v}}_{\vec{\lambda}}$ are wavepackets in the sense of Definition \ref%
{Definition regwave}. To obtain (\ref{E1bound}) we apply \ inequality (\ref%
{gradkk}) and use (\ref{rbb1}); for the part of $E^{1}$-norm without $%
\mathbf{k}$-derivatives we use (\ref{dtf}). Using multilinearity of $%
\mathcal{F}_{n_{l},\vartheta ,\vec{\lambda},\Psi _{2}}^{\left( m\right) }$
we observe that 
\begin{equation}
\mathcal{F}_{n_{l},\vartheta ,\vec{\lambda},\Psi _{2}}^{\left( m\right)
}\left( \mathbf{\vec{w}}\right) -\mathcal{F}_{n_{l},\vartheta ,\vec{\lambda}%
,\Psi _{2}}^{\left( m\right) }\left( \mathbf{\vec{v}}\right) =\sum_{j=1}^{m}%
\mathcal{F}_{n_{l},\vartheta ,\vec{\lambda},\Psi _{2}}^{\left( m\right)
}\left( \mathbf{w}_{\lambda _{1}},\ldots ,\mathbf{w}_{\lambda _{j}}-\mathbf{v%
}_{\lambda _{j}},\mathbf{v}_{\lambda _{j+1}},\ldots ,\mathbf{v}_{\lambda
_{m}}\right) .  \label{sumdif}
\end{equation}%
We can apply to every term \ inequality (\ref{gradkk}). Multiplying (\ref%
{gradkk}) by $\beta ^{1+\epsilon }$ and using (\ref{rbb1}) \ we deduce (\ref%
{gradk}).
\end{proof}

Now we consider a system similar to (\ref{eqav}), 
\begin{equation}
\mathbf{\hat{v}}_{l,\vartheta }=\mathcal{F}_{\limfunc{av},\Psi
_{2},n_{l},\vartheta }\left( \mathbf{\vec{v}}\right) +\Psi \left( \mathbf{%
\cdot },\vartheta \mathbf{k}_{\ast i_{l}}\right) \Pi _{n_{l},\vartheta }%
\mathbf{\hat{h}},\ l=1,\ldots N,\vartheta =\pm ,  \label{eqav1}
\end{equation}%
where $\mathcal{F}_{\limfunc{av},\Psi ,n_{l},\vartheta }$ is defined by a
formula similar to(\ref{Fav}): 
\begin{equation}
\mathcal{F}_{\limfunc{av},\Psi _{2},n_{l},\vartheta }\left( \mathbf{\vec{v}}%
\right) =\sum\nolimits_{m\in \mathfrak{M}_{F}}\mathcal{F}_{n_{l},\vartheta
}^{\left( m\right) },\ \mathcal{F}_{n_{l},\vartheta }^{\left( m\right)
}=\sum\nolimits_{\vec{\lambda}\in \Lambda _{n_{l},\vartheta }^{m}}\mathcal{F}%
_{n_{l},\vartheta ,\vec{\lambda},\Psi }^{\left( m\right) }\left( \mathbf{%
\vec{v}}\right) .  \label{Fav1}
\end{equation}%
The system (\ref{eqav1}) \ can be written in the form \ similar to (\ref%
{eqavF}) 
\begin{equation}
\mathbf{\vec{v}}=\mathcal{F}_{\limfunc{av},\Psi _{2}}\left( \mathbf{\vec{v}}%
\right) +\mathbf{\vec{h}}_{_{\Psi }}.  \label{eqav1F}
\end{equation}

\begin{theorem}[solvability in particle spaces]
\label{Theorem E1sol} Let the initial data $\mathbf{\vec{h}}$ in the
averaged wavepacket interaction system (\ref{eqav1F}) be a multi-particle
wavepacket $\mathbf{\hat{h}}\left( \beta ,\mathbf{k}\right) $ with $nk$%
-spectrum $S$ as in (\ref{P0}), the regularity degree $s$ and with positions 
$\mathbf{r}_{\ast l}$,$\ l=1,\ldots ,N$. \ Let $\left\Vert \mathbf{\vec{h}}%
\right\Vert _{\left( E^{1}\right) ^{2N}\left( \mathbf{\vec{r}}_{\ast
}\right) }\leq R$. Assume that $S$ is universally resonance invariant in the
sense of Definition \ref{Definition omclos}. Then there exists $\tau _{\ast
\ast }>0$ which does not depend on $\mathbf{\vec{r}}_{\ast }$, $\beta $ and $%
\varrho $ \ such that if $\tau _{\ast }\leq \tau _{\ast \ast }$ \ equation (%
\ref{eqav1F}) has a unique solution $\mathbf{\vec{v}}$ \ in $\left(
E^{1}\right) ^{2N}\left( \mathbf{\vec{r}}_{\ast }\right) $, such that \ 
\begin{equation}
\left\Vert \mathbf{\vec{v}}\right\Vert _{\left( E^{1}\right) ^{2N}\left( 
\mathbf{\vec{r}}_{\ast }\right) }\leq 2R  \label{E1est}
\end{equation}%
where $R$ does not depend on $\varrho ,\beta $ \ and on $\mathbf{\vec{r}}%
_{\ast }$. This solution is is a multi-particle wavepacket with positions $%
\mathbf{r}_{\ast l}$
\end{theorem}

\begin{proof}
Since $S$ is universally resonance invariant \ every vector index $\vec{%
\lambda}\in \Lambda _{n_{l},\vartheta }^{m}$ has at least one component $%
\lambda _{j}=\left( \zeta _{j},l_{j}\right) $ \ with $l_{j}=l$. Hence Lemma %
\ref{Lemma E1lip} is applicable and according to (\ref{gradk}) the operator $%
\mathcal{F}_{\limfunc{av},\Psi _{2}}$ defined by (\ref{Fav1}) is Lipschitz \
in the ball $\left\Vert \mathbf{\vec{v}}\right\Vert _{\left( E^{1}\right)
^{2N}\left( \mathbf{\vec{r}}_{\ast }\right) }\leq 2R$ \ with a Lipschitz
constant $C^{\prime }\tau _{\ast }$ where $C^{\prime }$ which does not
depend on $\varrho ,\beta ,$\ and on $\mathbf{\vec{r}}_{\ast }$. We choose $%
\tau _{\ast \ast }$ \ so that $C^{\prime }\tau _{\ast \ast }\leq 1/2$ and
use Lemma \ref{Lemma contr}. According to this Lemma equation (\ref{eqav1F})
has a solution $\mathbf{\vec{v}}$ which satisfies (\ref{E1est}). This
solution is is a multi-particle wavepacket \ according to Proposition \ref%
{Proposition partnorm}.
\end{proof}

\begin{theorem}[particle wavepacket approximation]
\label{Theorem sumwavereg} Let the initial data $\mathbf{\hat{h}}$ in the
integral equation (\ref{ubaseq}) with solution $\mathbf{\hat{u}}\left( \tau
,\beta ;\mathbf{k}\right) $ be an multi-particle wavepacket $\mathbf{\hat{h}}%
\left( \beta ,\mathbf{k}\right) $ with $nk$-spectrum $S$ as in (\ref{P0}),
the regularity degree $s$ and with positions $\mathbf{r}_{\ast l}$ $\
l=1,\ldots ,N$, and components of $\mathbf{\hat{h}}\left( \beta ,\mathbf{k}%
\right) $ satisfy the inequality $\left\Vert \mathbf{\vec{h}}\right\Vert
_{\left( E^{1}\right) ^{2N}\left( \mathbf{\vec{r}}_{\ast }\right) }\leq R$.
Let $\tau _{\ast }\leq \tau _{\ast \ast }$. \ Assume that $S$ is universally
resonance invariant in the sense of Definition \ref{Definition omclos}. We
define $\mathbf{\hat{v}}\left( \tau ,\beta ;\mathbf{k}\right) $ by the
formula 
\begin{equation}
\mathbf{\hat{v}}\left( \tau ,\beta ;\mathbf{k}\right)
=\sum_{l=1}^{N}\dsum\nolimits_{\zeta =\pm }\mathbf{\hat{v}}_{l,\vartheta
}\left( \tau ,\beta ;\mathbf{k}\right) ,\ l=1,\ldots ,N,  \label{vsum}
\end{equation}%
where $\mathbf{\hat{v}}_{l,\vartheta }\left( \tau ,\beta ;\mathbf{k}\right) $
is a solution of (\ref{eqav}). Then every such $\mathbf{\hat{v}}_{l}\left( 
\mathbf{k};\tau ,\beta \right) $ is a particle-like wavepacket with the
position $\mathbf{r}_{\ast l}$ and%
\begin{equation}
\sup_{0\leq \tau \leq \tau _{\ast }}\left\Vert \mathbf{\hat{u}}\left( \tau
,\beta ;\mathbf{k}\right) -\mathbf{\hat{v}}\left( \tau ,\beta ;\mathbf{k}%
\right) \right\Vert _{L^{1}}\leq C_{1}\varrho +C_{2}\beta ^{s},
\label{uminv}
\end{equation}%
where the constant $C_{1}$ does not depend on $\varrho ,s$ and $\beta $, and
the constant $C_{2}$ does not depend on $\varrho ,\beta .$
\end{theorem}

\begin{proof}
Let $\mathbf{\vec{v}}\in \left( E^{1}\right) ^{2N}\left( \mathbf{\vec{r}}%
_{\ast }\right) $ be the solution of equation (\ref{eqav1F}) which exists by
Theorem \ref{Theorem E1sol}, it is a particle-like wavepacket. Note that \ 
\begin{equation*}
\Psi \left( \mathbf{\cdot },\mathbf{k}_{\ast l_{1}},\beta ^{1-\epsilon
}/2\right) \Psi \left( \mathbf{\cdot },\mathbf{k}_{\ast l_{1}},\beta
^{1-\epsilon }\right) =\Psi \left( \mathbf{\cdot },\mathbf{k}_{\ast
l_{1}},\beta ^{1-\epsilon }\right)
\end{equation*}%
and solution of (\ref{eqav1F}) has the form $\mathbf{\hat{v}}_{l,\vartheta
}\left( \tau ,\beta ;\mathbf{k}\right) =\Psi \left( \mathbf{\cdot },\mathbf{k%
}_{\ast l},\beta ^{1-\epsilon }\right) \left[ \ldots \right] $ and,
consequently, for such solutions $\Psi _{2}\mathbf{\vec{v}}_{\vec{\lambda}}=%
\mathbf{\vec{v}}_{\vec{\lambda}}$, the nonlinearity $\mathcal{F}%
_{n_{l},\vartheta ,\vec{\lambda},\Psi _{2}}^{\left( m\right) }\left( \mathbf{%
\vec{v}}\right) $ coincides with and equation $\Psi \left( \mathbf{\cdot },%
\mathbf{k}_{\ast l},\beta ^{1-\epsilon }\right) \mathcal{F}_{n_{l},\vartheta
,\vec{\xi}\left( \vec{\lambda}\right) }^{\left( m\right) }\left( \mathbf{%
\vec{v}}_{\vec{\lambda}}\right) $ and the equation (\ref{eqav1}) coincides
with (\ref{eqav}). Hence, $\mathbf{\vec{v}}$ is a solution of (\ref{eqav}).
Estimate (\ref{uminv}) follows from estimates (\ref{vminw}) and (\ref{uminw}%
).
\end{proof}

\begin{corollary}
\label{Corollary Theoremregwave}If conditions of Theorem \ref{Theorem
regwave} \ are satisfied, the statement of Theorem \ref{Theorem regwave}
holds.
\end{corollary}

\begin{proof}
Note that functions $\mathbf{\hat{w}}_{l,\vartheta }^{\prime }\left( \mathbf{%
k},\tau \right) =\Psi _{i_{l},\vartheta }\Pi _{n_{l},\vartheta }\mathbf{\hat{%
u}}\left( \mathbf{k},\tau \right) $, $\theta =\pm $, in Theorem \ref{Theorem
necessary} are the two components of $\mathbf{\hat{u}}_{l}\left( \tau ,\beta
;\mathbf{k}\right) \ $ in (\ref{ups1}). Hence, (\ref{wwp}) implies that%
\begin{equation}
\left\Vert \mathbf{\hat{u}}_{l}-\mathbf{\hat{w}}_{l,+}-\mathbf{\hat{w}}%
_{l,-}\right\Vert _{E}\leq C^{\prime }\beta ^{s},\ 0<\beta \leq \beta _{0},
\label{ulww}
\end{equation}%
where $\mathbf{\hat{w}}_{l,\vartheta }$ are solutions to \ (\ref{sysloc}).
According to (\ref{vminw}), if $\mathbf{\hat{v}}_{l,\vartheta }\left( 
\mathbf{k},\tau \right) \mathbf{\ }$is the solution of (\ref{eqav}) we have 
\begin{equation}
\left\Vert \mathbf{\hat{v}}_{l,\vartheta }-\mathbf{\hat{w}}_{l,\vartheta
}\right\Vert _{E}\leq C\varrho ,\ l=1,\ldots ,N;\ \vartheta =\pm .
\end{equation}%
Hence, 
\begin{equation}
\left\Vert \mathbf{\hat{u}}_{l}-\mathbf{\hat{v}}_{l,+}-\mathbf{\hat{v}}%
_{l,-}\right\Vert _{E}\leq C\varrho +C^{\prime }\beta ^{s},\ 0<\beta \leq
\beta _{0}.  \label{ulvv}
\end{equation}%
This inequality implies (\ref{uui}). We have proved that $\ \mathbf{\hat{v}}%
_{l,\vartheta }$ is a particle-like wavepacket as in Theorem \ref{Theorem
sumwavereg}. Estimate (\ref{ulvv}) implies that $\mathbf{\hat{u}}_{l}$\ is
equivalent to $\mathbf{\hat{v}}_{l}=\mathbf{\hat{v}}_{l,+}+\mathbf{\hat{v}}%
_{l,-}$ in the sense of (\ref{equiv}) of the degree $s_{1}=\min \left(
s,s_{0}\right) $.
\end{proof}

\section{Superposition principle and decoupling of the wavepacket
interaction system}

In this section we give the proof of the superposition principle of \cite%
{BF8} \ which is based on a study of the wavepacket interaction system \ (%
\ref{eqav}). We will show that when we omit cross-terms in the averaged
system wavepacket interaction system, the resulting error is estimated by $%
\frac{\varrho }{\beta ^{1+\epsilon }}\left\vert \ln \beta \right\vert $,
that is component wavepackets evolve essentially independently and the time
averaged wavepacket interaction system almost decouples.

Let $\ \mathcal{F}_{\limfunc{av},n_{l},\vartheta }$ be defined by (\ref{Fav}%
) and let a\ decoupled nonlinearity $\mathcal{F}_{\limfunc{av}%
,n_{l},\vartheta ,\mathrm{diag}}$ be defined by \ 
\begin{equation}
\mathcal{F}_{\limfunc{av},n_{l},\vartheta ,\mathrm{diag}}\left( \mathbf{\vec{%
w}}\right) =\sum\nolimits_{m\in \mathfrak{M}_{F}}\mathcal{F}%
_{n_{l},\vartheta }^{\left( m\right) },\ \mathcal{F}_{n_{l},\vartheta ,%
\mathrm{diag}}^{\left( m\right) }\left( \mathbf{\vec{w}}\right)
=\sum\nolimits_{\vec{\lambda}\in \Lambda _{n_{l},\vartheta }^{m,\mathrm{diag}%
}}\mathcal{F}_{n_{l},\vartheta ,\vec{\xi}\left( \vec{\lambda}\right)
}^{\left( m\right) }\left( \mathbf{\vec{w}}_{\vec{\lambda}}\right) .
\label{Favd}
\end{equation}%
where the set of indices $\Lambda _{n_{l},\vartheta }^{m,\mathrm{diag}}$
consists of 
\begin{equation}
\Lambda _{n_{l},\vartheta }^{m,\mathrm{diag}}=\left\{ \vec{\lambda}=\left( 
\vec{l},\vec{\zeta}\right) \in \Lambda _{n_{l},\vartheta }^{m,}:l_{j}=l,\
j=1,\ldots ,m\right\} .  \label{lamdiag}
\end{equation}%
Note that $\mathcal{F}_{n_{l},\vartheta ,\mathrm{diag}}^{\left( m\right) }%
\mathcal{\ }$in (\ref{Favd}) depends only on $\mathbf{w}_{l,+}$ and $\mathbf{%
w}_{l,-}$: 
\begin{equation}
\mathcal{F}_{n_{l},\vartheta ,\mathrm{diag}}^{\left( m\right) }\left( 
\mathbf{\vec{w}}\right) =\mathcal{F}_{\vartheta ,\mathrm{diag},l}^{\left(
m\right) }\left( \mathbf{w}_{l}\right) ,\ \mathbf{w}_{l}=\left( \mathbf{w}%
_{l,+},\mathbf{w}_{l,-}\right) .  \label{Fdiagl}
\end{equation}%
The coupling between different variables $\mathbf{v}_{l}$ \ in (\ref{eqav})
\ is caused by \ non-diagonal terms 
\begin{equation}
\mathcal{F}_{\limfunc{av},n_{l},\vartheta ,\mathrm{coup}}\left( \mathbf{\vec{%
w}}\right) =\mathcal{F}_{\limfunc{av},n_{l},\vartheta }\left( \mathbf{\vec{w}%
}\right) -\mathcal{F}_{\limfunc{av},n_{l},\vartheta ,\mathrm{diag}}\left( 
\mathbf{\vec{w}}\right) .  \label{Fcoup}
\end{equation}%
Obviously, equation (\ref{eqavF}) \ can be written in the form 
\begin{equation}
\mathbf{\vec{v}}=\mathcal{F}_{\text{$\limfunc{av}$},\Psi ,\mathrm{diag}%
}\left( \mathbf{\vec{v}}\right) +\mathcal{F}_{\text{$\limfunc{av}$},\Psi ,%
\mathrm{coup}}\left( \mathbf{\vec{v}}\right) +\mathbf{\vec{h}}_{_{\Psi }}.
\label{eqavF1}
\end{equation}%
The system of decoupled equations has the form 
\begin{equation}
\mathbf{\vec{v}}_{\mathrm{diag}}=\mathcal{F}_{\text{$\limfunc{av}$},\Psi ,%
\mathrm{diag}}\left( \mathbf{\vec{v}}_{\mathrm{diag}}\right) +\mathbf{\vec{h}%
}_{_{\Psi }}.  \label{eqavF1d}
\end{equation}%
or, when written in components, 
\begin{equation}
\mathbf{v}_{\mathrm{diag},l}=\mathcal{F}_{\text{$\limfunc{av}$},\Psi ,%
\mathrm{diag},l}^{\left( m\right) }\left( \mathbf{v}_{\mathrm{diag}%
,l}\right) +\mathbf{h}_{_{\Psi ,l}},l=1,\ldots ,N.  \label{eqdec}
\end{equation}%
We will prove that contribution of $\mathcal{F}_{\text{$\limfunc{av}$},\Psi ,%
\mathrm{coup}}$ in (\ref{eqavF1}) is small, the proof is based on the
following lemma.

\begin{lemma}[small coupling terms]
\label{Lemma GVM1}\textbf{\ }Let $\mathcal{F}_{n_{l},\vartheta ,\vec{\xi}%
\left( \vec{\lambda}\right) }^{\left( m\right) }\left( \mathbf{\vec{w}}_{%
\vec{\lambda}}\right) $\ be as in (\ref{Fav}), let all \ components $\mathbf{%
w}_{\lambda _{i}}$ of $\mathbf{\vec{w}}_{\vec{\lambda}}$ satisfy (\ref{hl0})
\ and be wavepackets in the sense of Definition \ref{dwavepack}, \ let also (%
\ref{rbb1}) hold. Assume also that: (i) the vector index $\vec{\lambda}$ has
at least two components $\lambda _{i}=\left( \zeta _{i},l_{i}\right) $ and $%
\lambda _{j}=\left( \zeta _{j},l_{j}\right) $ with $l_{i}\neq l_{j}$; (ii)
both $\mathbf{w}_{\lambda _{i}}$ and $\mathbf{w}_{\lambda _{j}}$ are
particle wavepackets in the sense of Definition \ref{Definition regwave};
(iii) either (\ref{NGVM}) or (\ref{NGVMe}) holds. \ Then for small $\beta $
and $\varrho $%
\begin{equation}
\left\Vert \mathcal{F}_{n_{l},\vartheta ,\vec{\xi}\left( \vec{\lambda}%
\right) }^{\left( m\right) }\left( \mathbf{\vec{w}}_{\vec{\lambda}}\right)
\right\Vert _{E^{N}}\leq C\frac{\varrho }{\beta ^{1+\epsilon }}\left\vert
\ln \beta \right\vert .  \label{NGVMest1}
\end{equation}
\end{lemma}

\begin{proof}
Since $\mathbf{k}_{\ast l}$ are not band-crossing points, according to
Definition \ref{Definition band-crossing point} \ and Condition \ref{clin}
the inequalities (\ref{Com2}) and (\ref{grchi}) hold. According to the
assumption of the theorem \ at least two $\mathbf{\hat{w}}_{l_{j}}$ are
different for different $j$. \ Let us assume that $l_{j_{1}}=l_{1}$, $%
l_{j_{2}}=l_{m}$ ,\ $l_{1}\neq l_{m}$ (the general case can be easily
reduced to this one by a relabeling of variables). Since $\mathbf{\hat{w}}%
_{l_{1}}$and $\mathbf{\hat{w}}_{l_{m}}$ are particle wavepackets, \ they
satisfy (\ref{gradker}) \ with $\mathbf{\mathbf{r}}$ replaced by $\mathbf{%
\mathbf{r}}_{l_{1}}$ \ and $\mathbf{\mathbf{r}}_{l_{m}}$ respectively. Let
us rewrite the integral with respect to $\tau _{1}$ in (\ref{Fm}) as 
\begin{gather}
\mathcal{F}_{n_{l},\vartheta ,\vec{\xi}\left( \vec{\lambda}\right) }^{\left(
m\right) }\left( \mathbf{\vec{w}}_{\vec{\lambda}}\right) \left( \mathbf{k}%
,\tau \right) =  \label{FSA} \\
\int_{0}^{\tau }\int_{\mathbb{D}_{m}}\exp \left\{ \mathrm{i}\phi _{\zeta ,%
\vec{\zeta}}\left( \mathbf{\mathbf{k}},\vec{k}\right) \frac{\tau _{1}}{%
\varrho }\right\} A_{\zeta ,\vec{\zeta}}^{\left( m\right) }\left( \mathbf{k},%
\vec{k}\right) \mathrm{\tilde{d}}^{\left( m-1\right) d}\vec{k}d\tau _{1} 
\notag
\end{gather}%
where 
\begin{equation}
A_{\zeta ,\vec{\zeta}}^{\left( m\right) }\left( \mathbf{k},\vec{k}\right)
=\chi _{\zeta ,\vec{\zeta}}^{\left( m\right) }\left( \mathbf{k},\vec{k}%
\right) \mathbf{w}_{l_{1}}\left( \mathbf{k}^{\prime }\right) \ldots \mathbf{w%
}_{l_{m}}\left( \mathbf{k}^{\left( m\right) }\right) ,  \label{Akk}
\end{equation}%
and then rewrite (\ref{FSA}) in the form 
\begin{gather}
\mathcal{F}_{n_{l},\vartheta ,\vec{\xi}\left( \vec{\lambda}\right) }^{\left(
m\right) }\left( \mathbf{\vec{w}}_{\vec{\lambda}}\right) \left( \mathbf{k}%
,\tau \right) =\mathcal{F}_{\zeta ,\vec{\zeta}}^{\left( m\right) }\left( 
\mathbf{w}_{l_{1}}\ldots \mathbf{w}_{l_{m}}\right) \left( \mathbf{k},\tau
\right) =  \label{FSA1} \\
\int_{0}^{\tau }\int_{\mathbb{D}_{m}}\exp _{\phi }\left( \mathbf{\mathbf{k}},%
\vec{k},\tau _{1},\varrho ,\mathbf{\mathbf{r}}_{l_{1}},\mathbf{\mathbf{r}}%
_{l_{m}}\right) A\left( \mathbf{k},\vec{k},\mathbf{\mathbf{r}}_{l_{1}},%
\mathbf{\mathbf{r}}_{l_{m}}\right) \mathrm{\tilde{d}}^{\left( m-1\right) d}%
\vec{k}d\tau _{1}.  \notag
\end{gather}%
where 
\begin{gather}
\exp _{\phi }\left( \mathbf{\mathbf{k}},\vec{k},\tau _{1},\varrho ,\mathbf{%
\mathbf{r}}_{l_{1}},\mathbf{\mathbf{r}}_{l_{m}}\right) =\exp \left\{ \mathrm{%
i}\phi _{\zeta ,\vec{\zeta}}\left( \mathbf{\mathbf{k}},\vec{k}\right) \frac{%
\tau _{1}}{\varrho }-\mathrm{i}\mathbf{\mathbf{r}}_{l_{1}}\mathbf{\mathbf{k}}%
^{\prime }-\mathrm{i}\mathbf{\mathbf{r}}_{l_{m}}\mathbf{\mathbf{k}}^{\left(
m\right) }\right\} ,  \label{expfi} \\
A\left( \mathbf{k},\vec{k},\mathbf{\mathbf{r}}_{l_{1}},\mathbf{\mathbf{r}}%
_{l_{m}}\right) =\mathrm{e}^{\mathrm{i}\mathbf{\mathbf{r}}_{l_{1}}\mathbf{%
\mathbf{k}}^{\prime }}\mathrm{e}^{\mathrm{i}\mathbf{\mathbf{r}}_{l_{m}}%
\mathbf{\mathbf{k}}^{\left( m\right) }}A_{\zeta ,\vec{\zeta}}^{\left(
m\right) }\left( \mathbf{k},\vec{k}\right) .  \notag
\end{gather}%
According to (\ref{conv}) $\mathbf{k}^{\left( m\right) }\left( \mathbf{k},%
\vec{k}\right) =\mathbf{k}-\mathbf{k}^{\prime }-\ldots -\mathbf{k}^{\left(
m-1\right) }$. Hence, \ picking a vector $\mathbf{p}$ with a unit length \
we obtain the formula 
\begin{equation}
\exp _{\phi }\left( \mathbf{\mathbf{k}},\vec{k},\tau _{1},\varrho ,\mathbf{%
\mathbf{r}}_{l_{1}},\mathbf{\mathbf{r}}_{l_{m}}\right) =\frac{\varrho 
\mathbf{p}\cdot \nabla _{\mathbf{k}^{\prime }}\exp _{\phi }\left( \mathbf{%
\mathbf{k}},\vec{k},\tau _{1},\varrho ,\mathbf{\mathbf{r}}_{l_{1}},\mathbf{%
\mathbf{r}}_{l_{m}}\right) }{\mathrm{i}\left[ \mathbf{p}\cdot \nabla _{%
\mathbf{k}^{\prime }}\phi _{\zeta ,\vec{\zeta}}\left( \mathbf{\mathbf{k}},%
\vec{k}\right) \tau _{1}-\varrho \mathbf{p}\cdot \left( \mathbf{\mathbf{r}}%
_{l_{1}}-\mathbf{\mathbf{r}}_{l_{m}}\right) \right] }.  \label{idexp1}
\end{equation}%
If we set 
\begin{gather}
\mathbf{\phi }^{\prime }=\nabla _{\mathbf{k}^{\prime }}\phi _{\zeta ,\vec{%
\zeta}}\left( \mathbf{\mathbf{k}}_{\ast l},\vec{k}_{\ast }\right) =\nabla _{%
\mathbf{k}^{\prime }}\omega \left( \zeta ^{\prime }\mathbf{k}_{\ast
}^{\prime }\right) -\nabla _{\mathbf{k}^{\left( m\right) }}\omega \left(
\zeta ^{\left( m\right) }\mathbf{k}_{\ast }^{\left( m\right) }\right) ,
\label{qcp} \\
c_{p}=\mathbf{p}\cdot \mathbf{\phi }^{\prime },q_{p}=\varrho \mathbf{p}\cdot
\left( \mathbf{\mathbf{r}}_{l_{1}}-\mathbf{\mathbf{r}}_{l_{m}}\right) , 
\notag
\end{gather}%
\begin{equation}
\theta _{0}\left( \mathbf{\mathbf{k}},\vec{k},\varrho ,\tau _{1}\right) =%
\frac{\left( c_{p}\tau _{1}-q_{p}\right) }{\left[ \mathbf{p}\cdot \nabla _{%
\mathbf{k}^{\prime }}\phi _{\zeta ,\vec{\zeta}}\left( \mathbf{\mathbf{k}},%
\vec{k}\right) \tau _{1}-\mathbf{p}\cdot \left( \mathbf{\mathbf{r}}_{l_{1}}-%
\mathbf{\mathbf{r}}_{l_{m}}\right) \right] },  \label{thet0}
\end{equation}%
then (\ref{idexp1}) can be recast as%
\begin{equation}
\exp _{\phi }\left( \mathbf{\mathbf{k}},\vec{k},\tau _{1},\varrho ,\mathbf{%
\mathbf{r}}_{l_{1}},\mathbf{\mathbf{r}}_{l_{m}}\right) =\frac{\varrho 
\mathbf{p}\cdot \nabla _{\mathbf{k}^{\prime }}\exp _{\phi }\left( \mathbf{%
\mathbf{k}},\vec{k},\tau _{1},\varrho ,\mathbf{\mathbf{r}}_{l_{1}},\mathbf{%
\mathbf{r}}_{l_{m}}\right) }{\mathrm{i}\left( c_{p}\tau _{1}-q_{p}\right) }%
\theta _{0}\left( \mathbf{\mathbf{k}},\vec{k},\varrho ,\tau _{1}\right) .
\label{idexp2}
\end{equation}%
If (\ref{NGVM}) holds \ the vector $\mathbf{\phi }^{\prime }\neq 0$ and to
get $\left\vert c_{p}\right\vert \neq 0$ we can take 
\begin{equation}
\mathbf{p}=\left\vert \mathbf{\phi }^{\prime }\right\vert ^{-1}\cdot \mathbf{%
\phi }^{\prime },\ \left\vert c_{p}\right\vert =p_{0}>0.  \label{pfi}
\end{equation}%
If (\ref{NGVMe}) holds we have $\mathbf{\phi }^{\prime }=0$ and we set 
\begin{equation}
\mathbf{p}=\left\vert \left( \mathbf{\mathbf{r}}_{l_{1}}-\mathbf{\mathbf{r}}%
_{l_{m}}\right) \right\vert ^{-1}\cdot \left( \mathbf{\mathbf{r}}_{l_{1}}-%
\mathbf{\mathbf{r}}_{l_{m}}\right) .  \label{pfie}
\end{equation}%
Let consider first the case when (\ref{NGVM}) holds. Notice that the
denominator in (\ref{idexp2}) vanishes for 
\begin{equation}
\tau _{10}=\frac{q_{p}}{c_{p}}.  \label{t10}
\end{equation}%
We split the integral with respect to $\tau _{1}$ in (\ref{FSA1}) \ into a
sum of two integrals, namely 
\begin{gather}
\mathcal{F}_{\zeta ,\vec{\zeta}}^{\left( m\right) }\left( \mathbf{w}%
_{l_{1}}\ldots \mathbf{w}_{l_{m}}\right) \left( \mathbf{k},\tau \right)
=F_{1}+F_{2},  \notag \\
F_{1}=\int_{\left\vert \tau _{10}-\tau _{1}\right\vert \geq c_{0}\beta
^{1-\epsilon }\left\vert \ln \beta \right\vert }\int_{\mathbb{D}_{m}}\exp
_{\phi }\left( \mathbf{\mathbf{k}},\vec{k},\tau _{1},\varrho ,\mathbf{%
\mathbf{r}}_{l_{1}},\mathbf{\mathbf{r}}_{l_{m}}\right) A\left( \mathbf{k},%
\vec{k},\mathbf{\mathbf{r}}_{l_{1}},\mathbf{\mathbf{r}}_{l_{m}}\right) 
\mathrm{\tilde{d}}^{\left( m-1\right) d}\vec{k}d\tau _{1},  \label{Fbet1} \\
F_{2}=\int_{\left\vert \tau _{10}-\tau _{1}\right\vert <c_{0}\beta
^{1-\epsilon }\left\vert \ln \beta \right\vert }\int_{\mathbb{D}_{m}}\exp
_{\phi }\left( \mathbf{\mathbf{k}},\vec{k},\tau _{1},\varrho ,\mathbf{%
\mathbf{r}}_{l_{1}},\mathbf{\mathbf{r}}_{l_{m}}\right) A\left( \mathbf{k},%
\vec{k},\mathbf{\mathbf{r}}_{l_{1}},\mathbf{\mathbf{r}}_{l_{m}}\right) 
\mathrm{\tilde{d}}^{\left( m-1\right) d}\vec{k}d\tau _{1},  \notag
\end{gather}%
where $c_{0}$ is a large enough constant which we estimate below \ in (\ref%
{ses}). Since $\mathbf{w}_{j}\ $are bounded in $E$ and (\ref{scale1}) holds,
we obtain similarly to (\ref{dtf}) the estimate 
\begin{equation}
\left\Vert F_{2}\right\Vert _{L^{1}}\leq Cc_{0}\beta ^{1-\epsilon
}\left\vert \ln \beta \right\vert \dprod\limits_{j=1}^{m}\left\Vert \mathbf{w%
}_{l_{j}}\right\Vert _{E}\leq C_{1}\left( R\right) \frac{\varrho \left\vert
\ln \beta \right\vert }{\beta ^{1+\epsilon }}.  \label{0beta2}
\end{equation}

To estimate the norm of $F_{1}$ we use (\ref{idexp1}) and integrate by parts
the integral in (\ref{Fbet1}) with respect to $\mathbf{k}^{\prime }$. We
obtain 
\begin{gather}
F_{1}=\int_{\left\vert \tau _{10}-\tau _{1}\right\vert \geq \beta
^{1-\epsilon }\left\vert \ln \beta \right\vert }I\left( \mathbf{\mathbf{k}}%
,\tau _{1}\right) d\tau _{1},  \label{Ik} \\
I\left( \mathbf{\mathbf{k}},\tau _{1}\right) =-\int_{\mathbb{D}_{m}}\frac{%
\varrho \exp _{\phi }\left( \mathbf{\mathbf{k}},\vec{k},\tau _{1},\varrho ,%
\mathbf{\mathbf{r}}_{l_{1}},\mathbf{\mathbf{r}}_{l_{m}}\right) }{\mathrm{i}%
\left( c_{p}\tau _{1}-\varrho q_{p}\right) }\mathbf{p}\cdot \nabla _{\mathbf{%
k}^{\prime }}\left[ \theta _{0}A\left( \mathbf{k},\vec{k},\mathbf{\mathbf{r}}%
_{l_{1}},\mathbf{\mathbf{r}}_{l_{m}}\right) \right] \,\mathrm{\tilde{d}}%
^{\left( m-1\right) d}\vec{k}.  \notag
\end{gather}%
According to (\ref{Akk}) and (\ref{conv}) the expansion of the gradient $\
\nabla _{\mathbf{k}^{\prime }}$ in the above formula involves derivatives of 
$\chi ,\theta _{0},$ $\mathrm{e}^{\mathrm{i}\mathbf{\mathbf{r}}_{l_{1}}%
\mathbf{\mathbf{k}}^{\prime }}\mathbf{w}_{l_{1}}$ and $\mathrm{e}^{\mathrm{i}%
\mathbf{\mathbf{r}}_{l_{m}}\mathbf{\mathbf{k}}^{\left( m\right) }}\mathbf{w}%
_{l_{m}}$. To estimate $\theta _{0}$ and $\nabla \theta _{0}$ \ we note that 
\begin{gather}
\theta _{0}\left( \mathbf{\mathbf{k}},\vec{k},\varrho ,\tau _{1}\right) =%
\frac{\left( \mathbf{p}\cdot \mathbf{\phi }^{\prime }\tau _{1}-q_{p}\right) 
}{\left( \mathbf{p}\cdot \mathbf{\phi }^{\prime }\tau _{1}-q_{p}\right)
+\tau _{1}\mathbf{p}\cdot \left[ \nabla _{\mathbf{k}^{\prime }}\phi _{\zeta ,%
\vec{\zeta}}\left( \mathbf{\mathbf{k}},\vec{k}\right) -\mathbf{\phi }%
^{\prime }\right] }  \label{frac} \\
=\frac{1}{1+\tau _{1}\mathbf{p}\cdot \left[ \nabla _{\mathbf{k}^{\prime
}}\phi _{\zeta ,\vec{\zeta}}\left( \mathbf{\mathbf{k}},\vec{k}\right) -%
\mathbf{\phi }^{\prime }\right] /\left( c_{p}\tau _{1}-q_{p}\right) }  \notag
\end{gather}%
Since $\left\vert \tau _{10}-\tau _{1}\right\vert \geq c_{0}\beta
^{1-\epsilon }\left\vert \ln \beta \right\vert $, from (\ref{t10}) \ we
infer 
\begin{equation}
\left\vert c_{p}\tau _{1}-q_{p}\right\vert \geq c_{p}c_{0}\beta ^{1-\epsilon
}\left\vert \ln \beta \right\vert .  \label{pfib}
\end{equation}%
From (\ref{hl0}), we see that in the integral (\ref{Ik}) integrands are
nonzero only if 
\begin{equation}
\left\vert \mathbf{\mathbf{k}}^{\left( j\right) }-\zeta ^{\left( j\right) }%
\mathbf{\mathbf{k}}_{\ast }^{\left( j\right) }\right\vert \leq \pi _{0}\beta
^{1-\epsilon },\left\vert \mathbf{\mathbf{k}}-\zeta \mathbf{\mathbf{k}}%
_{\ast }\right\vert \leq m\pi _{0}\beta ^{1-\epsilon },  \label{nzero}
\end{equation}%
where $\pi _{0}\leq 1$. Using the Taylor remainder estimate \ for $\nabla _{%
\mathbf{k}^{\prime }}\phi _{\zeta ,\vec{\zeta}}$ at $\vec{k}_{\ast }$ we
obtain the inequality 
\begin{equation}
\left\vert \nabla _{\mathbf{k}^{\prime }}\phi _{\zeta ,\vec{\zeta}}\left( 
\mathbf{\mathbf{k}},\vec{k}\right) -\mathbf{\phi }^{\prime }\right\vert \leq
2mC_{\omega ,2}\beta ^{1-\epsilon }.\text{ \ }  \label{Tayest}
\end{equation}%
\ Hence in (\ref{frac}) 
\begin{equation}
\left\vert \tau _{1}\mathbf{p}\cdot \left[ \nabla _{\mathbf{k}^{\prime
}}\phi _{\zeta ,\vec{\zeta}}\left( \mathbf{\mathbf{k}},\vec{k}\right) -%
\mathbf{\phi }^{\prime }\right] /\left( c_{p}\tau _{1}-q_{p}\right)
\right\vert \leq 2m\tau _{\ast }C_{\omega ,2}/\left( c_{p}c_{0}\left\vert
\ln \beta \right\vert \right) .  \label{pfib1}
\end{equation}%
Suppose that $\beta \leq 1/2$ is small and $c_{0}$ \ satisfies%
\begin{equation}
\frac{m\tau _{\ast }C_{\omega ,2}}{\left\vert \ln \beta \right\vert }\leq 
\frac{m\tau _{\ast }C_{\omega ,2}}{\ln 2}\leq \frac{1}{4}\left\vert
c_{p}\right\vert c_{0}.  \label{ses}
\end{equation}%
\ \ Then it follows from (\ref{frac}) with help of (\ref{ses}), (\ref{Com2}%
), (\ref{pfib}) and (\ref{pfib1}) that \ 
\begin{equation}
\left\vert \theta _{0}\left( \mathbf{\mathbf{k}},\vec{k},\varrho ,\tau
_{1}\right) \right\vert \leq 2.  \label{gradom}
\end{equation}%
Obviously, 
\begin{equation}
\nabla _{\mathbf{k}^{\prime }}\theta _{0}\left( \mathbf{\mathbf{k}},\vec{k}%
,\varrho ,\tau _{1}\right) =\frac{-\tau _{1}\nabla _{\mathbf{k}^{\prime }}%
\left[ \mathbf{p\cdot }\left( \nabla _{\mathbf{k}^{\prime }}\phi _{\zeta ,%
\vec{\zeta}}\left( \mathbf{\mathbf{k}},\vec{k}\right) -\mathbf{\phi }%
^{\prime }\right) \right] }{\left( c_{p}\tau _{1}-q_{p}\right) \left[ 1+\tau
_{1}\mathbf{p\cdot }\left[ \nabla _{\mathbf{k}^{\prime }}\phi _{\zeta ,\vec{%
\zeta}}\left( \mathbf{\mathbf{k}},\vec{k}\right) -\mathbf{\phi }^{\prime }%
\right] /\left( c_{p}\tau _{1}-q_{p}\right) \right] ^{2}}  \label{gkth}
\end{equation}%
Using (\ref{pfib1}), (\ref{ses}) \ and (\ref{Com2}) we obtain that 
\begin{equation}
\left\vert \nabla _{\mathbf{k}^{\prime }}\theta _{0}\left( \mathbf{\mathbf{k}%
},\vec{k},\varrho ,\tau _{1}\right) \right\vert \leq \frac{4\tau _{\ast }}{%
\left\vert c_{p}\tau _{1}-q_{p}\right\vert }\left\vert \nabla _{\mathbf{k}%
^{\prime }}\left[ \mathbf{p}\cdot \left( \nabla _{\mathbf{k}^{\prime }}\phi
_{\zeta ,\vec{\zeta}}\left( \mathbf{\mathbf{k}},\vec{k}\right) -\mathbf{\phi 
}^{\prime }\right) \right] \right\vert \leq \frac{4\tau _{\ast }C_{\omega ,2}%
}{\left\vert c_{p}\tau _{1}-q_{p}\right\vert }.  \label{gkest}
\end{equation}%
To estimate $\nabla _{\mathbf{k}^{\prime }}\chi $ we use (\ref{grchi}). We
conclude that the absolute value of the integral (\ref{Ik}) is not greater
than%
\begin{gather}
\left\vert I\left( \mathbf{\mathbf{k}},\tau _{1}\right) \right\vert \leq 
\frac{4\varrho \tau _{\ast }C_{\omega ,2}}{\left\vert \tau
_{1}c_{p}-q_{p}\right\vert ^{2}}\int_{\mathbb{D}_{m}}\left\vert A\left( 
\mathbf{k},\vec{k},\mathbf{\mathbf{r}}_{l_{1}},\mathbf{\mathbf{r}}%
_{l_{m}}\right) \right\vert \mathrm{\tilde{d}}^{\left( m-1\right) d}\vec{k}+
\label{Ikt1} \\
\frac{2\varrho \tau _{\ast }}{\left\vert \tau _{1}c_{p}-q_{p}\right\vert }%
\int_{\mathbb{D}_{m}}\left[ \left\vert \nabla _{\mathbf{k}^{\prime }}A\left( 
\mathbf{k},\vec{k},\mathbf{\mathbf{r}}_{l_{1}},\mathbf{\mathbf{r}}%
_{l_{m}}\right) \right\vert \right] \mathrm{\tilde{d}}^{\left( m-1\right) d}%
\vec{k} \\
\leq \left[ \frac{4C_{\omega ,2}\varrho \tau _{\ast }}{\left\vert \tau
_{1}c_{p}-q_{p}\right\vert ^{2}}\left\Vert \chi ^{\left( m\right) }\left( 
\mathbf{\mathbf{k}},\cdot \right) \right\Vert +\frac{2\varrho \tau _{\ast }}{%
\left\vert \tau _{1}c_{p}-q_{p}\right\vert }\left\Vert \left( \nabla
_{k^{\prime }}-\nabla _{k^{\left( m\right) }}\right) \chi ^{\left( m\right)
}\left( \mathbf{\mathbf{k}},\cdot \right) \right\Vert \right]
\dprod\limits_{j=1}^{m}\left\Vert \mathbf{\mathbf{w}}_{j}\right\Vert
_{L^{1}}+  \notag \\
\frac{2\varrho \tau _{\ast }\left\Vert \chi ^{\left( m\right) }\left( 
\mathbf{\mathbf{k}},\cdot \right) \right\Vert }{\left\vert \tau
_{1}c_{p}-q_{p}\right\vert }\left[ \dprod\limits_{j=2}^{m}\left\Vert \mathbf{%
\mathbf{w}}_{l_{j}}\right\Vert _{L^{1}}\left\Vert \nabla _{\mathbf{k}%
^{\prime }}\mathrm{e}^{\mathrm{i}\mathbf{\mathbf{r}}_{l_{1}}\mathbf{k}%
^{\prime }}\mathbf{\mathbf{w}}_{l_{1}}\right\Vert
_{L^{1}}+\dprod\limits_{j=1}^{m-1}\left\Vert \mathbf{\mathbf{w}}%
_{j}\right\Vert _{L^{1}}\left\Vert \nabla _{\mathbf{k}^{\left( m\right) }}%
\mathrm{e}^{\mathrm{i}\mathbf{\mathbf{r}}_{l_{m}}\mathbf{k}^{\left( m\right)
}}\mathbf{\mathbf{w}}_{m}\right\Vert _{L^{1}}\right] .  \notag
\end{gather}%
Note that $\left\Vert \mathbf{\mathbf{w}}_{j}\right\Vert _{L^{1}}$ are
bounded according to (\ref{L1b}) \ and $\ \nabla _{\mathbf{k}^{\left(
m\right) }}\mathrm{e}^{\mathrm{i}\mathbf{\mathbf{r}}_{l_{m}}\mathbf{k}%
^{\left( m\right) }}\mathbf{\mathbf{w}}_{l_{m}}$, $\nabla _{\mathbf{k}%
^{\prime }}\mathrm{e}^{\mathrm{i}\mathbf{\mathbf{r}}_{l_{1}}\mathbf{k}%
^{\prime }}\mathbf{\mathbf{w}}_{l_{1}}\mathbf{\ }$by (\ref{gradker}). Hence
we obtain 
\begin{equation}
\left\vert I\left( \mathbf{\mathbf{k}},\tau _{1}\right) \right\vert \leq 
\frac{C_{2}\varrho \beta ^{-1-\epsilon }}{\tau _{1}c_{p}-q_{p}}+\frac{%
\varrho C_{2}}{\left\vert \tau _{1}c_{p}-q_{p}\right\vert ^{2}}.  \label{II}
\end{equation}%
Obviously, 
\begin{gather*}
\int_{\left\vert \tau _{1}-q_{p}/c_{p}\right\vert \geq c_{0}\beta
^{1-\epsilon }\left\vert \ln \beta \right\vert }\left\vert \tau
_{1}c_{p}-q_{p}\right\vert ^{-1}d\tau _{1}=\frac{1}{c_{p}}\int_{c_{0}\beta
^{1-\epsilon }\left\vert \ln \beta \right\vert }^{\tau _{\ast }-q_{p}/c_{p}}%
\frac{d\tau _{1}}{\tau _{1}}=\frac{1}{c_{p}}\ln \frac{\tau _{\ast
}-q_{p}/c_{p}}{c_{0}\beta ^{1-\epsilon }\left\vert \ln \beta \right\vert } \\
\leq \frac{1}{c_{p}}\left( C+\left\vert \ln \left[ \beta ^{1-\epsilon
}\left\vert \ln \beta \right\vert \right] \right\vert \right) \leq \frac{1}{%
c_{p}}\left[ C+\left\vert \ln \beta \right\vert +\left\vert \ln \left\vert
\ln \beta \right\vert \right\vert \right] \leq \frac{1}{c_{p}}\left[
C+2\left\vert \ln \beta \right\vert \right] .
\end{gather*}%
Similarly, using (\ref{scale1}) we get 
\begin{gather*}
\int_{\left\vert \tau _{1}-q_{p}/c_{p}\right\vert \geq c_{0}\beta
^{1-\epsilon }\left\vert \ln \beta \right\vert }\left\vert \tau
_{1}c_{p}-q_{p}\right\vert ^{-2}d\tau _{1}=\frac{1}{c_{p}}\int_{c_{0}\beta
^{1-\epsilon }\left\vert \ln \beta \right\vert }^{\tau _{\ast }-q_{p}/c_{p}}%
\frac{d\tau _{1}}{\tau _{1}^{2}}=\frac{1}{c_{p}}\left[ \frac{1}{c_{0}\beta
^{1-\epsilon }\left\vert \ln \beta \right\vert }-\frac{1}{\tau _{\ast
}-q_{p}/c_{p}}\right] \\
\leq \frac{1}{c_{p}c_{0}\beta ^{1-\epsilon }\left\vert \ln \beta \right\vert 
}\leq \frac{C_{3}\varrho }{\beta ^{-1-\epsilon }\left\vert \ln \beta
\right\vert }
\end{gather*}%
Hence, we obtain \ for small $\beta $ 
\begin{equation}
\left\Vert \mathcal{F}_{\zeta ,\vec{\zeta}}^{\left( m\right) }\left( \mathbf{%
w}_{1}\ldots \mathbf{w}_{m}\right) \left( \mathbf{k},\tau \right)
\right\Vert _{E}\leq C_{4}\frac{\varrho }{\beta ^{1+\epsilon }}\left\vert
\ln \beta \right\vert .  \label{NGVM2}
\end{equation}%
Now let us consider the case when (\ref{NGVMe}) holds, $\ \mathbf{\phi }%
^{\prime }=0$ and $\mathbf{p}$ is defined by (\ref{pfie}). Turning to
expression (\ref{frac}) we notice that 
\begin{equation*}
c_{p}\tau _{1}-q_{p}=-\varrho \left\vert \mathbf{\mathbf{r}}_{l_{1}}-\mathbf{%
\mathbf{r}}_{l_{m}}\right\vert ,\ \tau _{\ast }\left\vert c_{p}\tau
_{1}-q_{p}\right\vert ^{-1}\leq \frac{1}{\beta ^{1+\epsilon }},
\end{equation*}%
and, according to (\ref{Tayest}), 
\begin{equation*}
\left\vert \nabla _{\mathbf{k}^{\prime }}\phi _{\zeta ,\vec{\zeta}}\left( 
\mathbf{\mathbf{k}},\vec{k}\right) -\mathbf{\phi }^{\prime }\right\vert \leq
C_{\omega ,2}\beta ^{1-\epsilon }.
\end{equation*}%
Then we estimate the denominator in (\ref{frac}) and (\ref{gkth}) using (\ref%
{NGVMe}):$\ $%
\begin{equation*}
\left\vert \tau _{1}\mathbf{p}\cdot \left[ \nabla _{\mathbf{k}^{\prime
}}\phi _{\zeta ,\vec{\zeta}}\left( \mathbf{\mathbf{k}},\vec{k}\right) -%
\mathbf{\phi }^{\prime }\right] /\left( c_{p}\tau _{1}-q_{p}\right)
\right\vert \leq \tau _{\ast }C_{\omega ,2}\beta ^{1-\epsilon }/\left(
\varrho \left\vert \mathbf{\mathbf{r}}_{l_{1}}-\mathbf{\mathbf{r}}%
_{l_{m}}\right\vert \right) \leq \frac{1}{2}.
\end{equation*}%
If $\beta $ is so small that (\ref{ses}) holds we again get (\ref{gradom})
and (\ref{gkest}). Hence, we obtain (\ref{NGVM2}) in this case as well (in
fact, in this case the logarithmic factor can be omitted). Finally, we
obtain (\ref{NGVMest0}) from (\ref{NGVM2}) after summing up over all $\vec{%
\lambda},\vec{\zeta}$.
\end{proof}

\begin{lemma}
\label{Lemma GVM0}\textbf{\ }Let the $nk$-spectrum $S$ be universally
resonance invariant. Let the operators $\ \mathcal{F}_{\limfunc{av}%
,n_{l},\vartheta }\left( \mathbf{\vec{w}}\right) $, $\mathcal{F}_{\limfunc{av%
},n_{l},\vartheta ,\mathrm{diag}}\left( \mathbf{\vec{v}}\right) $ and $%
\mathcal{F}_{\limfunc{av},n_{l},\vartheta ,\mathrm{coup}}$ be defined
respectively by (\ref{Fav}), (\ref{Fu}) and (\ref{Fcoup}). Let $\mathbf{\vec{%
v}},$ \ $\left\Vert \mathbf{\vec{v}}\right\Vert _{E^{N}}\leq 2R,$ be a
multi-wavepacket solution of (\ref{eqavF}) with the $nk$-spectrum $S$. Then
for small $\beta $ and $\varrho $ 
\begin{equation}
\left\Vert \mathcal{F}_{\limfunc{av},n_{l},\vartheta ,\mathrm{coup}}\left( 
\mathbf{\vec{v}}\right) \right\Vert _{E^{N}}\leq C\frac{\varrho }{\beta
^{1+\epsilon }}\left\vert \ln \beta \right\vert .  \label{NGVMest0}
\end{equation}
\end{lemma}

\begin{proof}
According to (\ref{Fav}), (\ref{Favd}) \ and (\ref{Fcoup}) $\mathcal{F}_{%
\limfunc{av},n_{l},\vartheta ,\mathrm{coup}}$ involves only terms with $\vec{%
\lambda}\in \Lambda _{n_{l},\vartheta }^{m,}\setminus \Lambda
_{n_{l},\vartheta }^{m,\mathrm{diag}}$ and it is sufficient to prove the
estimate (\ref{NGVMest1}) for indices $\vec{\lambda}\in \Lambda
_{n_{l},\vartheta }^{m,}\setminus \Lambda _{n_{l},\vartheta }^{m,\mathrm{diag%
}}$. Such indices involve at least two components $\lambda _{i}=\left( \zeta
_{i},l_{i}\right) $ and $\lambda _{j}=\left( \zeta _{j},l_{j}\right) $ with $%
l_{i}\neq l_{j}$ \ \ since the $nk$-spectrum is universally invariant, see (%
\ref{rearr0}). According to Theorem \ref{Theorem sumwavereg} \ the solution
\ $\mathbf{\vec{v}}$ is a particle-like wavepacket, therefore all components
of $\mathbf{\vec{v}}_{\vec{\lambda}}$ are particle-like; (\ref{hl0}) holds
according to (\ref{weq0}). Hence, all conditions of Lemma \ref{Lemma GVM1}
are fulfilled and (\ref{NGVMest0}) follows from (\ref{NGVMest1}).
\end{proof}

Note now that every equation (\ref{eqdec}) is an approximation of the
equation (\ref{varcu}) \ with single-wavepacket initial data $\mathbf{\hat{h}%
}_{l}$, namely%
\begin{equation}
\mathbf{\hat{u}}_{l}\left( \mathbf{k},\tau \right) =\mathcal{F}\left( 
\mathbf{\hat{u}}_{l}\right) \left( \mathbf{k},\tau \right) +\mathbf{\hat{h}}%
_{l}\left( \mathbf{k}\right) .  \label{uel}
\end{equation}
One can apply to this equation Theorems \ref{Theorem uminw} and \ref{Theorem
uminwav} \ formally restricted to the case $N=1$ of a single wavepacket.
Based on this observation and on the above Lemma we prove the following
theorem which implies previously formulated Theorems \ref{Theorem
Superposition} \ and \ref{Theorem Superposition1}.

\begin{theorem}
\label{Theorem superpositionsys} Assume that the multiwavepacket $\mathbf{%
\vec{h}}=\sum \mathbf{\hat{h}}_{l}$ is particle-like and its $nk$-spectrum
is universally resonance invariant. Assume also that either (\ref{NGVM}) or (%
\ref{NGVMe}) holds. Let $\mathbf{\hat{u}}$ be solution of equation (\ref%
{varcu}). Let $\ \mathbf{\hat{u}}_{l}$ be solutions of (\ref{uel}). Then the
superposition principle holds, namely%
\begin{equation}
\left\Vert \mathbf{\hat{u}}-\sum_{l=1}^{N}\mathbf{\hat{u}}_{l}\right\Vert
\leq C\frac{\varrho }{\beta ^{1+\epsilon }}\left\vert \ln \beta \right\vert
+C\beta ^{s}.  \label{usumvd}
\end{equation}
\end{theorem}

\begin{proof}
Let $\mathbf{v}_{\mathrm{diag},l}$ be a solution of the decoupled system (%
\ref{eqdec}). We compare systems (\ref{eqavF1}) and (\ref{eqavF1d}). The
difference between the systems is the term $\mathcal{F}_{\limfunc{av}%
,n_{l},\vartheta ,\mathrm{coup}}\left( \mathbf{\vec{v}}\right) .$ According
to Theorem \ref{Theorem E1sol} the solution $\mathbf{\vec{v}}$ is a
particle-like wavepacket and we can apply Lemma \ref{Lemma GVM0}. According
to this Lemma (\ref{NGVMest0}) holds. Applying Lemma \ref{Lemma contr} to
the equations (\ref{eqavF1}) and (\ref{eqavF1d}) and using (\ref{NGVMest0})\
we conclude that the difference of their solutions satisfies the inequality 
\begin{equation}
\left\Vert \mathbf{v}_{l}-\mathbf{v}_{\mathrm{diag},l}\right\Vert _{E}\leq
C^{\prime }\frac{\varrho }{\beta ^{1+\epsilon }}\left\vert \ln \beta
\right\vert +C^{\prime }\beta ^{s}.  \label{vvdi}
\end{equation}
According to Theorem \ref{Theorem sumwavereg} inequality (\ref{uminv}) holds
where $\mathbf{\vec{v}}$ \ is a solution of (\ref{eqavF}) which can be
rewritten in the form of (\ref{eqavF1}). From (\ref{uminv}) and (\ref{vvdi})
we infer 
\begin{equation}
\left\Vert \mathbf{\hat{u}}-\sum_{l=1}^{N}\mathbf{v}_{\mathrm{diag}%
,l}\right\Vert \leq C_{1}\frac{\varrho }{\beta ^{1+\epsilon }}\left\vert \ln
\beta \right\vert +C_{1}\beta ^{s}.  \label{usumvd1}
\end{equation}

Note that equation (\ref{eqdec}) for $\mathbf{v}_{\mathrm{diag},l}$\
coincides with the averaged equation (\ref{eqavF}) obtained for the wave
interaction system derived for (\ref{uel}). Therefore, applying Theorems \ref%
{Theorem uminw} and \ref{Theorem uminwav} to the case $N=1$ and $\ \mathbf{%
\hat{h}}=\mathbf{\hat{h}}_{l}$ we deduce from (\ref{uminw}) \ and (\ref%
{vminw}) the estimate 
\begin{equation}
\left\Vert \mathbf{\hat{u}}_{l}-\mathbf{v}_{\mathrm{diag},l}\right\Vert
_{E}\leq C_{2}\varrho +C_{2}^{\prime }\beta ^{s}.  \label{uvdi}
\end{equation}%
Finally, from (\ref{usumvd1}) and (\ref{uvdi}) we infer (\ref{usumvd}).
\end{proof}

\subsection{Generalizations}

In this section we show that the particle-like wavepacket invariance can be
extended to the case when $nk$-spectra $S$ are not universally resonance
invariant. So suppose that $nk$-spectrum $S$ is resonance invariant and
consider nonlinearities of the form similar to (\ref{Fav}) 
\begin{equation}
\mathcal{F}_{\text{res},n_{l},\vartheta }\left( \mathbf{\vec{w}}\right)
=\sum\nolimits_{m\in \mathfrak{M}_{F}}\mathcal{F}_{n_{l},\vartheta }^{\left(
m\right) },\ \mathcal{F}_{n_{l},\vartheta }^{\left( m\right)
}=\sum\nolimits_{\vec{\lambda}\in \Lambda _{n_{l},\vartheta }^{\prime }}%
\mathcal{F}_{n_{l},\vartheta ,\vec{\xi}\left( \vec{\lambda}\right) }^{\left(
m\right) }\left( \mathbf{\vec{w}}_{\vec{\lambda}}\right) ,  \label{restrnon}
\end{equation}%
where $\Lambda _{n_{l},\vartheta }^{\prime }\subseteq \Lambda
_{n_{l},\vartheta }^{m}$ is a given subset of $\Lambda ^{m}$. \ Obviously, $%
\mathcal{F}_{\text{av}}$ defined by (\ref{Fav}) has the form of (\ref%
{restrnon}) \ with $\Lambda _{n_{l},\vartheta }^{\prime }=\Lambda
_{n_{l},\vartheta }^{m}$.\ Let us introduce a multi-wavepacket 
\begin{equation}
\mathbf{\vec{w}}=\left( \mathbf{w}_{n_{1},+},\mathbf{w}_{n_{1},-},\ldots ,%
\mathbf{w}_{n_{N},+},\mathbf{w}_{n_{N},-}\right)  \label{mpw}
\end{equation}
with the $nk$-spectrum $S=\left\{ \left( n_{l},\theta \right) ,\ l=1,\ldots
,N;\theta =\pm \right\} $.

We call a subset $S^{\prime }\subset S$ \emph{sign-invariant} if with it has 
$\left( n_{l},\theta \right) $ as an element then $\left( n_{l},-\theta
\right) $ is also its element. Suppose that $S^{\prime }\subset S$ is
sign-invariant. It is easy to see that if a set $S^{\prime }\subset S$ is
sign-invariant then it is uniquely defined by a subset of indices $I^{\prime
}=I^{\prime }\left( S^{\prime }\right) \subset I=\left\{ 1,\ldots ,N\right\} 
$, namely 
\begin{equation*}
S^{\prime }=\left\{ \left( n_{l},\theta \right) :l\in I^{\prime }\left(
S^{\prime }\right) ,\ \theta =\pm \right\} .
\end{equation*}

\begin{definition}
\label{Definition GVM}We call index pair $\left( n_{l},\mathbf{k}_{\ast
l}\right) \ $ \emph{Group Velocity Matched (GVM)} with $\mathcal{F}_{\text{%
res},n_{l},\vartheta }$ if every nonzero term $\mathcal{F}_{n_{l},\vartheta ,%
\vec{\xi}\left( \vec{\lambda}\right) }^{\left( m\right) }$ in the sum (\ref%
{restrnon}) has the index $\vec{\lambda}$ such that \ for at least one \
component $\lambda _{j}=\left( \zeta ^{\left( j\right) },l_{j}\right) $ of
this index 
\begin{equation}
\nabla \omega _{n_{l}}\left( \mathbf{k}_{\ast l}\right) =\nabla \omega
_{n_{l_{j}}}\left( \mathbf{k}_{\ast l_{j}}\right) .  \label{GVMl}
\end{equation}%
We call $S^{\prime }$ a \emph{GVM set} with respect to the nonlinearity $%
\mathcal{F}_{\text{res}}$ defined by (\ref{restrnon}) if $S^{\prime }\subset
S$ is sign-invariant \ and every $\left( n_{l},\mathbf{k}_{\ast l}\right)
\in S^{\prime }$ is GVM.
\end{definition}

Obviously, if $S$ is universally resonance invariant and $\Lambda
_{n_{l},\vartheta }^{\prime }=\Lambda _{n_{l},\vartheta }^{m}$ as in (\ref%
{Fav}) then $S$ is a GVM set, and in this case $l_{j}=I_{0}$ as in
Definition \ref{Universal solution}. If $S^{\prime }\subset S$ is
sign-invariant we call a\ multi-wavepacket $\ \mathbf{\vec{w}}$ as in (\ref%
{mpw}) with the $nk$-spectrum $S=\left\{ \left( n_{l},\theta \right)
,l=1,\ldots ,N;\theta =\pm \right\} $ \emph{partially }$S^{\prime }$\emph{%
-localized}$\ $\emph{\ multi-wavepacket} if for every $\left( n_{l},\theta
\right) \in S^{\prime }$ the wavepacket $\mathbf{w}_{n_{1},\theta }\ $is a
spatially localized with a position $\mathbf{r}_{\ast l}$. Note that
according to Definition \ref{dmwavepack} if $S^{\prime }=S$ is a partially $%
S^{\prime }$-localized$\ $ multi-wavepacket then it is a multi-particle
wavepacket.

\ Theorem \ref{Theorem regwave} on the particle-like wavepacket preservation
can be generalized as follows.

\begin{theorem}[preservation of spatially localized wavepackets]
\label{Theorem regwaveg}Assume that conditions of Theorem \ref{Theorem
sumwave} hold, in particular the initial datum $\mathbf{\hat{h}}=\mathbf{%
\hat{h}}\left( \beta ,\mathbf{k}\right) $ is a multi-wavepacket with $nk$%
-spectrum $S.$ Assume also that $S^{\prime }\subset S$ is a GVM set, $%
\mathbf{\hat{h}}=\mathbf{\hat{h}}\left( \beta ,\mathbf{k}\right) $ is
partially $S^{\prime }$-localized wavepacket with positions $\mathbf{r}%
_{\ast l},$ $l\in I^{\prime }\left( S^{\prime }\right) $, and that (\ref%
{ros0}) holds. Then the solution $\mathbf{\hat{u}}\left( \tau ,\beta \right)
=\mathcal{G}\left( \mathcal{F}\left( \rho \left( \beta \right) \right) ,%
\mathbf{\hat{h}}\left( \beta \right) \right) \left( \tau \right) $ to (\ref%
{ubaseq}) for any $\tau \in \left[ 0,\tau _{\ast }\right] $ is a
multi-wavepacket with $nk$-spectrum $S$, and it is an $S^{\prime }$%
-localized wavepacket with positions $\mathbf{r}_{\ast l},$ $l\in I^{\prime
}\left( S^{\prime }\right) $. Namely, (\ref{uui}) holds where $\mathbf{\hat{u%
}}_{l}$ is wavepacket with $nk$-pair $\left( n_{l},\mathbf{k}_{\ast
l}\right) \in S^{\prime }$ defined by (\ref{ups1}), the constants \ $%
C,C_{1},C_{2}$ do not depend on $\mathbf{r}_{\ast l}$. and every $\mathbf{%
\hat{u}}_{l}$, $l\in I^{\prime }\left( S^{\prime }\right) $, is equivalent
in the sense of the equivalence (\ref{equiv}) of degree $s_{1}=\min \left(
s,s_{0}\right) $ to a spatially localized wavepacket with position $\mathbf{r%
}_{\ast l}$.
\end{theorem}

\begin{proof}
The proof of the Theorem is the same as the proof of Theorem \ref{Theorem
regwave} \ since it used only the fact that a universally resonance
invariant set is a GVM one, that allows to apply Lemma \ref{Lemma gradF}.
One also have to use the space $\left( E^{1}\right) ^{2N}\left( \mathbf{\vec{%
r}}_{\ast },S^{\prime }\right) $ with the norm defined \ by the formula
similar to (\ref{E12Nr}):%
\begin{equation}
\left\Vert \mathbf{\vec{w}}\right\Vert _{\left( E^{1}\right) ^{2N}\left( 
\mathbf{\vec{r}}_{\ast },S^{\prime }\right) }=\sum\nolimits_{l,\vartheta
}\left\Vert \mathbf{\hat{w}}_{l,\vartheta }\right\Vert _{E}+\beta
^{1+\epsilon }\sum\nolimits_{\vartheta =\pm }\sum\nolimits_{l\in I^{\prime
}\left( S^{\prime }\right) }\left\Vert \nabla _{\mathbf{k}}\left( \mathrm{e}%
^{-\mathrm{i}\mathbf{r}_{\ast l}\mathbf{k}}\mathbf{\hat{w}}_{l,\vartheta
}\right) \right\Vert _{E}.  \label{E1NS}
\end{equation}%
After replacing $\left( E^{1}\right) ^{2N}\left( \mathbf{\vec{r}}_{\ast
}\right) $\ with $\left( E^{1}\right) ^{2N}\left( \mathbf{\vec{r}}_{\ast
},S^{\prime }\right) $ we can literally repeat all the steps of the proof of
Theorem \ref{Theorem regwave} \ and obtain the statement of Theorem \ref%
{Theorem regwaveg}
\end{proof}

Below we prove that the superposition principle can hold not only for
universal resonance invariant multiwavepackets, but for other cases allowing
resonant processes such as second and third harmonic generations, three-wave
interaction etc. Here we prove a theorem applicable to such situations,
which is more general than Theorem \ref{Theorem Superposition}.

Let us consider a multi-wavepacket with \ resonance invariant $nk$-spectrum%
\begin{equation*}
S=\left\{ \left( n_{l},\mathbf{k}_{\ast l}\right) ,\ l=1,\ldots ,N\right\}
\end{equation*}%
as in (\ref{ssp1}), and assume that is a union of spectra $S_{p}$:%
\begin{equation}
S=S_{1}\cup \ldots \cup S_{K},\ S_{p}\cap S_{q}=\varnothing \text{ \ if \ }%
p\neq q.  \label{Sunion}
\end{equation}%
Recall that resonance interactions are defined in terms of vectors $\vec{%
\lambda}\in \Lambda ^{m}$ (see (\ref{setlam}), (\ref{lamprop})). We call a
vector $\vec{\lambda}=\left( \left( \zeta ^{\prime },l_{1}\right) ,\ldots
,\left( \zeta ^{\left( m\right) },l_{m}\right) \right) \in \Lambda ^{m}$ 
\emph{cross-interacting (CI)} if there exist at least two indices $\left(
\zeta ^{\left( i\right) },l_{i}\right) $ and $\left( \zeta ^{\left( j\right)
},l_{j}\right) $ such that $\left( \zeta ^{\left( i\right) },l_{i}\right)
\in S_{p_{i}},\left( \zeta ^{\left( j\right) },l_{j}\right) \in S_{p_{j}}\ \ 
$with $p_{i}\neq p_{j}$.

\begin{definition}[partially GVM decomposition]
\label{Catalytic decomposition}We call decomposition (\ref{Sunion}) \emph{%
partially GVM } with respect to $\mathcal{F}_{\text{res}}\ \ $defined by (%
\ref{restrnon}) if the following two conditions are satisfied: (i) every
spectrum $S_{j}$, $j=1,\ldots ,K,$ is resonance invariant; \ (ii) a solution 
$\left( m,\zeta ,n,\vec{\lambda}\right) \in P\left( S\right) $ of the
resonance equation (\ref{Omeq0}) with CI vector $\vec{\lambda}=\left( \left(
\zeta ^{\prime },l_{1}\right) ,\ldots ,\left( \zeta ^{\left( m\right)
},l_{m}\right) \right) \ $has at least two indices $\left( \zeta ^{\left(
i\right) },l_{i}\right) \in S_{p_{i}}$ and $\left( \zeta ^{\left( j\right)
},l_{j}\right) \in S_{p_{j}}$ with $p_{i}\neq p_{j}$ such that both $l_{i}$
and $l_{j}$ are GVM with respect to $\mathcal{F}_{\text{res}}\ $\ and 
\begin{equation}
\left\vert \nabla _{\mathbf{k}}\omega _{n_{l_{i}}}\left( \mathbf{k}_{\ast
l_{i}}\right) -\nabla _{\mathbf{k}}\omega _{n_{l_{j}}}\left( \mathbf{k}%
_{\ast l_{j}}\right) \right\vert \neq 0.\text{ }  \label{NGVMg}
\end{equation}
\end{definition}

Now we use Lemma \ref{Lemma GVM1} for small coupling. Being given a
partially GVM\ decomposition (\ref{Sunion}) we introduce the set of coupling
terms \ between $S_{p_{i}}$ \ and $\ S_{p_{j}}$as follows: 
\begin{equation}
\Lambda _{n_{l},\vartheta }^{m,\mathrm{coup}}=\left\{ \vec{\lambda}=\left( 
\vec{l},\vec{\zeta}\right) \in \Lambda _{n_{l},\vartheta }^{m,}:\exists \
i\neq j\text{ \ such that \ }l_{i}\in S_{p_{i}},l_{j}\in S_{p_{j}}\right\} ,
\label{weakcoup}
\end{equation}%
We also introduce a set of interactions reducible to \ every $S_{p}$
(block-diagonal) which is similar to (\ref{lamdiag}): 
\begin{equation}
\Lambda _{n_{l},\vartheta }^{m,\mathrm{red}}=\Lambda _{n_{l},\vartheta
}^{m}\setminus \Lambda _{n_{l},\vartheta }^{m,\mathrm{coup}},
\label{lamdiag1}
\end{equation}%
and the reduced operator \ 
\begin{equation}
\mathcal{F}_{\limfunc{av},n_{l},\vartheta ,\mathrm{red}}\left( \mathbf{\vec{w%
}}\right) =\sum\nolimits_{m\in \mathfrak{M}_{F}}\mathcal{F}_{n_{l},\vartheta
,\mathrm{red}}^{\left( m\right) }\left( \mathbf{\vec{w}}\right) ,\ \mathcal{F%
}_{n_{l},\vartheta ,\mathrm{red}}^{\left( m\right) }\left( \mathbf{\vec{w}}%
\right) =\sum\nolimits_{\vec{\lambda}\in \Lambda _{n_{l},\vartheta }^{m,%
\mathrm{red}}}\mathcal{F}_{n_{l},\vartheta ,\vec{\xi}\left( \vec{\lambda}%
\right) }^{\left( m\right) }\left( \mathbf{\vec{w}}_{\vec{\lambda}}\right) ,
\label{Favred}
\end{equation}%
where $\Lambda _{n_{l},\vartheta }^{m,\mathrm{red}}$ is defined by (\ref%
{lamdiag1}). Note that if the set $S$ is universal resonance invariant and \
every $S_{p_{i}}$ is a two-point set $\left\{ \left( +,l_{i}\right) ,\left(
+,l_{i}\right) \right\} $ then $\Lambda _{n_{l},\vartheta }^{m,\mathrm{red}%
}=\Lambda _{n_{l},\vartheta }^{m,\mathrm{diag}}$. We introduce also a
partially decoupled, reduced system similar to (\ref{eqavF1d})%
\begin{equation}
\mathbf{\vec{v}}_{\mathrm{red}}=\mathcal{F}_{\text{$\limfunc{av}$},\Psi ,%
\mathrm{red}}\left( \mathbf{\vec{v}}_{\mathrm{red}}\right) +\mathbf{\vec{h}}%
_{_{\Psi }},  \label{eqred1}
\end{equation}%
which can be rewritten in the decoupled form, similar to (\ref{eqdec}):%
\begin{equation}
\mathbf{v}_{\mathrm{red},p}=\mathcal{F}_{\text{$\limfunc{av}$},\Psi ,\mathrm{%
red},p}^{\left( m\right) }\left( \mathbf{v}_{\mathrm{red},p}\right) +\mathbf{%
h}_{_{\mathrm{red},\Psi ,p}},\ p=1,\ldots K.  \label{eqred}
\end{equation}%
Now $\mathbf{v}_{\mathrm{red},p}$ may include more than one wavepacket,
namely 
\begin{equation}
\mathbf{v}_{\mathrm{red},p}=\sum_{\left( n_{l},\theta \right) \in
S_{p}}\left( \mathbf{\vec{v}}_{\mathrm{red}}\right) _{n_{l},\theta },\;%
\mathbf{h}_{_{\mathrm{red},\Psi ,p}}=\sum_{\left( n_{l},\theta \right) \in
S_{p}}\left( \mathbf{\vec{h}}_{_{\Psi }}\right) _{n_{l},\theta },\
p=1,\ldots K.  \label{vhred}
\end{equation}%
The following theorem is a generalization of the Theorem \ref{Theorem
Superposition} on the superposition.

\begin{theorem}[general superposition principle ]
\label{Theorem Superposition gen}\ Suppose that the initial data $\mathbf{%
\hat{h}}$ of (\ref{ubaseq}) is a multi-wavepacket of the form%
\begin{equation}
\mathbf{\hat{h}}=\sum_{p=1}^{K}\mathbf{\hat{h}}_{\mathrm{red},p},\ 
\end{equation}%
where $\mathbf{\hat{h}}$ is a multi-wavepacket in the sense of Definition %
\ref{Definition omclos} with resonance invariant $nk$-spectrum $S$, $\mathbf{%
\hat{h}}_{\mathrm{red},p}$ is a multi-wavepacket with a resonance invariant $%
nk$-spectrum $S_{p}$\ and the decomposition (\ref{Sunion}) is partially GVM
in the sense of Definition \ref{Catalytic decomposition} with respect to the
nonlinearity $\mathcal{F}_{\text{av}}\ $defined by (\ref{Fav}). Suppose also
that (\ref{scale1}) holds. Then the solution $\mathbf{\hat{u}}$ $=\mathcal{G}%
\left( \mathbf{\hat{h}}\right) $ to the evolution equation (\ref{ubaseq})
satisfies the approximate superposition principle: 
\begin{equation}
\mathcal{G}\left( \sum_{p=1}^{K}\mathbf{\hat{h}}_{\mathrm{red},p}\right)
=\sum_{p=1}^{K}\mathcal{G}\left( \mathbf{\hat{h}}_{\mathrm{red},p}\right) +%
\mathbf{\tilde{D}},
\end{equation}%
with a small remainder $\mathbf{\tilde{D}}\left( \tau \right) $ satisfying
the following estimate 
\begin{equation}
\sup_{0\leq \tau \leq \tau _{\ast }}\left\Vert \mathbf{\tilde{D}}\left( \tau
\right) \right\Vert _{L^{1}}\leq C_{\epsilon }\frac{\varrho }{\beta
^{1+\epsilon }}\left\vert \ln \beta \right\vert ,  \label{Dgen}
\end{equation}%
where $\epsilon $ is the same as in Definition \ref{dwavepack} and can be
arbitrary small, $\tau _{\ast }$ does not depend on $\beta ,\varrho $ \ and $%
\epsilon $.
\end{theorem}

\begin{proof}
The proof of Theorem \ref{Theorem Superposition gen} \ is similar the proof
of Theorem \ref{Theorem superpositionsys}. Averaged system (\ref{eqavF}) can
be written similarly to (\ref{eqavF1}) \ in the form 
\begin{equation}
\mathbf{\vec{v}}=\mathcal{F}_{\text{$\limfunc{av}$},\Psi ,\mathrm{red}%
}\left( \mathbf{\vec{v}}\right) +\mathcal{F}_{\text{$\limfunc{av}$},\Psi ,%
\mathrm{coup}}\left( \mathbf{\vec{v}}\right) +\mathbf{\vec{h}}_{_{\Psi }}.
\label{redcoup}
\end{equation}%
Comparing now systems (\ref{redcoup}) and (\ref{eqred1}) we find that the
difference between them is the term $\mathcal{F}_{\limfunc{av}%
,n_{l},\vartheta ,\mathrm{coup}}\left( \mathbf{\vec{v}}\right) .$ According
to Theorem \ref{Theorem regwaveg} the solution $\mathbf{\vec{v}}$ is a
spatially localized wavepacket and, hence, we can apply Lemma \ref{Lemma
GVM0} getting the inequality (\ref{NGVMest0}). Applying Lemma \ref{Lemma
contr} to the equations (\ref{redcoup}) and (\ref{eqred1}) and using (\ref%
{NGVMest0})\ we conclude that \ the difference of their solutions satisfies
the inequality 
\begin{equation}
\left\Vert \mathbf{v}_{p}-\mathbf{v}_{\mathrm{red},p}\right\Vert _{E}\leq
C^{\prime }\frac{\varrho }{\beta ^{1+\epsilon }}\left\vert \ln \beta
\right\vert +C^{\prime }\beta ^{s},p=1,...,K  \label{vvred}
\end{equation}
According to Theorem \ref{Theorem sumwavereg} inequality (\ref{uminv}) holds
where $\mathbf{\vec{v}}$ \ is a solution of (\ref{redcoup}) and we infer
from (\ref{vvred}) 
\begin{equation}
\left\Vert \mathbf{\hat{u}}-\sum_{p=1}^{K}\mathbf{v}_{\mathrm{red}%
,p}\right\Vert \leq C_{1}\frac{\varrho }{\beta ^{1+\epsilon }}\left\vert \ln
\beta \right\vert +C_{1}\beta ^{s}.  \label{uminvg}
\end{equation}%
Similarly to (\ref{uel}) we introduce equation for $\mathbf{\hat{u}}_{%
\mathrm{red},p}=\mathcal{G}\left( \mathbf{\hat{h}}_{\mathrm{red},p}\right) $ 
\begin{equation}
\mathbf{\hat{u}}_{\mathrm{red},p}\left( \mathbf{k},\tau \right) =\mathcal{F}%
\left( \mathbf{\hat{u}}_{\mathrm{red},p}\right) \left( \mathbf{k},\tau
\right) +\mathbf{\hat{h}}_{\mathrm{red},p}\left( \mathbf{k}\right) .
\label{uelg}
\end{equation}%
Applying Theorems \ref{Theorem uminw} and \ref{Theorem uminwav} we infer
similarly to (\ref{uvdi}) the inequality 
\begin{equation}
\left\Vert \mathbf{\hat{u}}_{\mathrm{red},p}-\mathbf{v}_{\mathrm{red}%
,p}\right\Vert _{E}\leq C_{2}\varrho +C_{2}^{\prime }\beta ^{s}.
\label{uvdig}
\end{equation}%
Finally, from (\ref{uminvg}) and (\ref{uvdig}) we infer (\ref{Dgen}).
\end{proof}

\textbf{Acknowledgment:} Effort of A. Babin and A. Figotin is sponsored by
the Air Force Office of Scientific Research, Air Force Materials Command,
USAF, under grant number FA9550-04-1-0359.\newline

\end{document}